\documentclass[a4paper,12pt,leqno]{book}

\usepackage{hyperref}% ArXiv
\usepackage{makeidx}
\usepackage[T1]{fontenc}
\usepackage[english]{babel}
\usepackage{amsmath}
\input xy \xyoption{all}

\newlength\regolo

\newcounter{szn}
\newcounter{cnt}[szn]
\newcounter{itm}

\makeatletter

\newcommand\separazione{\bigskip\par}

\newcommand\sottosezione[1]{\@startsection{subparagraph}{2}{0pt}{-\baselineskip}{\baselineskip}{\centering\bfseries\large}{#1}}

\newcommand\sezione[2][]{\stepcounter{szn}\settowidth\regolo{#1}\ifdim\regolo>0pt\@startsection{section}{1}{0pt}{-\baselineskip}{\baselineskip}{\centering\bfseries\Large}[#1]{#2}\else\@startsection{section}{1}{0pt}{-\baselineskip}{\baselineskip}{\centering\bfseries\Large}[#2]{#2}\fi\addtocounter{szn}{-1}\refstepcounter{szn}}

\newcommand\capitolo\chapter

\newcommand\parte\part

\makeatother

\newenvironment{elenco}{\begin{list}{\roman{itm})}{\setlength\itemindent{0pt}\setlength\labelsep{0.5em}\setlength\labelwidth\leftmargin\addtolength\labelwidth{-\labelsep}\setlength\listparindent{.5\parindent}\setlength\parsep\parskip\setlength\itemsep\medskipamount\setlength\partopsep{0pt}\usecounter{itm}}}{\end{list}}

\newenvironment{equazione}{\addtocounter{cnt}{+1}\addtocounter{equation}{+1}\equation\tag{\thecnt}}{\endequation}

\newenvironment{equazione+}{\addtocounter{cnt}{+1}\addtocounter{equation}{+1}\flalign\tag{\thecnt}&\phantom!&&}{\endflalign}

\newenvironment{multiriga}{\addtocounter{cnt}{+1}\addtocounter{equation}{+1}\multline\tag{\thecnt}}{\endmultline}

\providecommand\qedsymbol{\textsl{q.e.d.}}
\newcommand\mathqed{\quad\hbox{\qedsymbol}}
\DeclareRobustCommand\qed{\ifmmode\mathqed\else\leavevmode\unskip\penalty9999\hbox{}\nobreak\hfill\quad\hbox{\qedsymbol}\fi}
\newenvironment{proof}[1][]{\begin{list}{}{\setlength\itemindent{0pt}\setlength\labelsep{0pt}\setlength\labelwidth{0pt}\setlength\leftmargin{0pt}\setlength\listparindent\parindent\setlength\parsep\parskip\setlength\partopsep{0pt}}\item\textsl{Proof}\settowidth\regolo{#1}\ifdim\regolo>0pt\textit{~~(#1)}\else\fi~~}{\qed\end{list}}

\newenvironment{definizione}[1][]{\refstepcounter{cnt}\begin{list}{}{\setlength\itemindent{0pt}\setlength\labelsep{0pt}\setlength\labelwidth{0pt}\setlength\leftmargin{0pt}\setlength\listparindent\parindent\setlength\parsep\parskip\setlength\partopsep{0pt}}\item\textbf{\thecnt~~Definition}~~}{\end{list}}

\newenvironment{enunciato}[2][]{\refstepcounter{cnt}\begin{list}{}{\setlength\itemindent{0pt}\setlength\labelsep{0pt}\setlength\labelwidth{0pt}\setlength\leftmargin\parindent\setlength\listparindent\parindent\setlength\parsep\parskip\setlength\partopsep{0pt}}\item\thecnt\settowidth\regolo{#2}\ifdim\regolo>0pt\settowidth\regolo{#1}\ifdim\regolo>0pt~\textsc{~#2~}~(#1)~\else~\textsc{~#2~}\fi\else\settowidth\regolo{#1}\ifdim\regolo>0pt~~(#1)~\else~\fi\fi\em~}{\end{list}}

\newenvironment{esempio}[1][]{\refstepcounter{cnt}\begin{list}{}{\setlength\itemindent\parindent\setlength\labelsep\parindent\setlength\labelwidth{0pt}\setlength\leftmargin{0pt}\setlength\listparindent\parindent\setlength\parsep\parskip\setlength\partopsep{0pt}}\item\thecnt\settowidth\regolo{#1}\ifdim\regolo>0pt\textit{~~Example:~#1}\else\textit{~~Example}\fi~~}{\end{list}}

\newenvironment{lemma}[1][]{\refstepcounter{cnt}\begin{list}{}{\setlength\itemindent{0pt}\setlength\labelsep{0pt}\setlength\labelwidth{0pt}\setlength\leftmargin\parindent\setlength\listparindent\parindent\setlength\parsep\parskip\setlength\partopsep{0pt}}\item\textbf{\thecnt~~Lemma}\settowidth\regolo{#1}\ifdim\regolo>0pt~~(#1)\slshape~~\else\slshape~~\fi}{\end{list}}

\newenvironment{nota}[1][]{\refstepcounter{cnt}\begin{list}{}{\setlength\itemindent\parindent\setlength\labelsep\parindent\setlength\labelwidth{0pt}\setlength\leftmargin{0pt}\setlength\listparindent\parindent\setlength\parsep\parskip\setlength\partopsep{0pt}}\item\thecnt\settowidth\regolo{#1}\ifdim\regolo>0pt\textit{~~Note:~#1}\else\textit{~~Note}\fi~~}{\end{list}}

\newenvironment{proposizione}[1][]{\refstepcounter{cnt}\begin{list}{}{\setlength\itemindent{0pt}\setlength\labelsep{0pt}\setlength\labelwidth{0pt}\setlength\leftmargin\parindent\setlength\listparindent\parindent\setlength\parsep\parskip\setlength\partopsep{0pt}}\item\textbf{\thecnt~~Proposition}\settowidth\regolo{#1}\ifdim\regolo>0pt~~(#1)\slshape~~\else\slshape~~\fi}{\end{list}}

\newenvironment{theorem}[1][]{\refstepcounter{cnt}\begin{list}{}{\setlength\itemindent{0pt}\setlength\labelsep{0pt}\setlength\labelwidth{0pt}\setlength\leftmargin\parindent\setlength\listparindent\parindent\setlength\parsep\parskip\setlength\partopsep{0pt}}\item\textbf{\thecnt~~Theorem}\settowidth\regolo{#1}\ifdim\regolo>0pt~~(#1)\slshape~~\else\slshape~~\fi}{\end{list}}

\newenvironment{theorem*}[1][]{\begin{list}{}{\setlength\itemindent{0pt}\setlength\labelsep{0pt}\setlength\labelwidth{0pt}\setlength\leftmargin\parindent\setlength\listparindent\parindent\setlength\parsep\parskip\setlength\partopsep{0pt}}\item\textbf{Theorem}\settowidth\regolo{#1}\ifdim\regolo>0pt~~(#1)\slshape~~\else\slshape~~\fi}{\end{list}}

\newenvironment{corollario}[1][]{\refstepcounter{cnt}\begin{list}{}{\setlength\itemindent{0pt}\setlength\labelsep{0pt}\setlength\labelwidth{0pt}\setlength\leftmargin\parindent\setlength\listparindent\parindent\setlength\parsep\parskip\setlength\partopsep{0pt}}\item\textbf{\thecnt~~Corollary}\settowidth\regolo{#1}\ifdim\regolo>0pt~~(#1)\slshape~~\else\slshape~~\fi}{\end{list}}

\usepackage{amssymb}
\usepackage{bbold}%*mbboard*
\usepackage{mathrsfs}
\DeclareMathAlphabet{\mathpzc}{T1}{pzc}{m}{it}
\newcommand\bydef{\stackrel{\text{\tiny\textrm{def}}}=}
\newcommand\gm[1]{\guillemotleft{\tiny~}{#1}{\tiny~}\guillemotright}
\newcommand\inciso[1]{\nobreakdash---\hspace{0pt}#1\nobreakdash---\hspace{0pt}}
\newcommand\refcnt[2][]{%
\settowidth\regolo{#1}\ifdim\regolo>0pt\ifmmode\text{\ref{#1}.\ref{#2}}\else{\ref{#1}.\ref{#2}}\fi\else\ifmmode\text{\ref{#2}}\else\ref{#2}\fi\fi}
\newcommand\refequ[2][]{%
\settowidth\regolo{#1}\ifdim\regolo>0pt\ifmmode\text{(\ref{#1}.\ref{#2})}\else(\ref{#1}.\ref{#2})\fi\else\ifmmode\text{(\ref{#2})}\else(\ref{#2})\fi\fi}
\newcommand\refsez[2][]{%
\settowidth\regolo{#1}\ifdim\regolo>0pt\S\S\ref{#1}--\ref{#2}\else\S\ref{#2}\fi}
\newcommand\refcpt[2][]{%
\settowidth\regolo{#1}\ifdim\regolo>0pt\textrm{\ref{#1}}--\textrm{\ref{#2}}\else\textrm{\ref{#2}}\fi}
\newcommand\Aut{\mathrm{Aut}}
\newcommand\dimens[1]{\ensuremath{\mathrm{dim}\,#1}}
\newcommand\End{\mathrm{End}}
\newcommand\GL{\mathit{GL}}
\newcommand\Hom{\ensuremath{\mathrm{Hom}}}
\newcommand\image[1]{\mathrm{Im}\,#1}
\newcommand\Iso{\mathrm{Iso}}
\newcommand\kernel[1]{\mathrm{Ker}\,#1}
\newcommand\Lis{\mathrm{Lis}}

\newcommand\C{\mathit C}
\newcommand\D{\mathit D}
\newcommand\Lebesgue[1][1]{\mathit L^{#1}}
\newcommand\support[2][]{\ensuremath{\mathrm{supp}_{#1}\,#2}}
\newcommand\T[2][]{\ensuremath{\mathit T_{#1}\,#2}}
\newcommand\test[1][\infty]{\mathit C^{#1}_{\mathit c}}
\newcommand\epito\twoheadrightarrow
\newcommand\from\leftarrow
\newcommand\infrom\hookleftarrow
\newcommand\into\hookrightarrow

\newcommand\isoto{\stackrel\thicksim\to}
\newcommand\longto\longrightarrow
\newcommand\xfrom[1]{\xleftarrow{#1}}
\newcommand\xto[1]{\xrightarrow{#1}}
\newcommand\Id{\mathit{Id}}
\newcommand\Ob{\mathrm{Ob}}

\newcommand\opposite[1]{\ensuremath{#1^{\mathrm{op}}}}
\newcommand\SpV{\ensuremath{\underline{\mathcal V\mspace{-2.8mu}\mathit{ec}}}}
\newcommand\VectorSpaces{\ensuremath{\{\mathrm{vector\:spaces}\}}}

\newcommand\ComplexVectorSpaces{\ensuremath{\{\mathrm{complex\:vector\:spaces}\}}}
\newcommand\SheavesOfModules[1]{\ensuremath{\left\{\mathrm{sheaves\:of\:#1\text-modules}\right\}}}
\newcommand\aeq\Leftrightarrow
\newcommand\seq\Rightarrow
\newcommand\pt{\ensuremath{\star}}
\newcommand\ev{\mathit{ev}}
\newcommand\id{\mathit{id}}
\newcommand\pr{\mathit{pr}}
\newcommand\nZ{\ensuremath{\mathbb Z}}

\newcommand\nR{\ensuremath{\mathbb R}}
\newcommand\nC{\ensuremath{\mathbb C}}
\newcommand\displaycap[3]{\ensuremath{\overset{#2}{\underset{#1}{\displaystyle\bigcap}}\,#3}}

\newcommand\displayprod[3]{\ensuremath{\overset{#2}{\underset{#1}{\displaystyle\prod}}\,#3}}

\newcommand\displaysum[3]{\ensuremath{\overset{#2}{\underset{#1}{\displaystyle\sum}}\,#3}}
\newcommand\dual[1]{\ensuremath{#1^\vee}}
\newcommand\bidual[1]{\ensuremath{#1^{\vee\vee}}}
\newcommand\inductivelim[2]{\ensuremath{\underset{#1}{\varinjlim}\,#2}}
\newcommand\modulo[1]{\ensuremath{\left|#1\right|}}
\newcommand\bigmod[1]{\ensuremath{\bigl|#1\bigr|}}

\newcommand\norma[1]{\ensuremath{\left\|#1\right\|}}

\newcommand\scalare[2]{\ensuremath{\left\langle#1,#2\right\rangle}}
\newcommand\bigsca[2]{\ensuremath{\bigl\langle#1,#2\bigr\rangle}}

\newcommand\parteRe[1]{\ensuremath{\Re\mathfrak e#1}}
\newcommand\trasposto[1]{\ensuremath{{^{\mathsf t}#1}}}
\newcommand\txtcap[3]{\ensuremath{\overset{#2}{\underset{#1}{\textstyle\cap}}\,#3}}
\newcommand\txtcup[3]{\ensuremath{\overset{#2}{\underset{#1}{\textstyle\cup}}\,#3}}
\newcommand\txtprod[3]{\ensuremath{\overset{#2}{\underset{#1}{\textstyle\prod}}\,#3}}
\newcommand\txtcoprod[3]{\ensuremath{\overset{#2}{\underset{#1}{\textstyle\coprod}}\,#3}}
\newcommand\txtsum[3]{\ensuremath{\overset{#2}{\underset{#1}{\textstyle\sum}}\,#3}}
\newcommand\can\cong
\newcommand\iso\approx
\newcommand\sheafconst[2][]{%
\settowidth{\regolo}{\ensuremath{#1}}\ifdim\regolo>0pt\ensuremath{\underline{#2}_{#1}}\else\ensuremath{\underline{#2}}\fi}
\newcommand\continuous[1][]{%
\settowidth{\regolo}{\ensuremath{#1}}\ifdim\regolo>0pt\ensuremath{\mathscr C^0_{#1}}\else\ensuremath{\mathscr C^0}\fi}
\newcommand\smooth[1][]{%
\settowidth{\regolo}{\ensuremath{#1}}\ifdim\regolo>0pt\ensuremath{\mathscr C^\infty_{#1}}\else\ensuremath{\mathscr C^\infty}\fi}
\newcommand\G{\ensuremath{\mathcal G}}
\renewcommand\H{\ensuremath{\mathcal H}}
\newcommand\K{\ensuremath{\mathcal K}}
\newcommand\s[1][]{%
\settowidth{\regolo}{\ensuremath{#1}}\ifdim\regolo>0pt\ensuremath{\mathit s\,#1}\else\ensuremath{\mathit s}\fi}
\renewcommand\t[1][]{%
\settowidth{\regolo}{\ensuremath{#1}}\ifdim\regolo>0pt\ensuremath{\mathit t\,#1}\else\ensuremath{\mathit t}\fi}
\renewcommand\c{\ensuremath{\mathit c}}
\newcommand\p{\ensuremath{\mathit p}}
\renewcommand\u{\ensuremath{\mathit u}}
\newcommand\imap{\ensuremath{\mathit i}}
\newcommand\mca[2][]{%
\settowidth{\regolo}{\ensuremath{#1}}\ifdim\regolo>0pt\ensuremath{#2^{\scriptscriptstyle(#1)}}\else\ensuremath{#2{_{\s}\times_{\t}}#2}\fi}
\newcommand\Kt[1][C]{\ensuremath{\mathcal{#1}}}
\newcommand\Tensor{\boldsymbol\otimes}
\newcommand\TU{\ensuremath{\mathsf1}}
\newcommand\Can[1]{\ensuremath{\underline{\mathcal C\mspace{-2.6mu}\mathit{an}}(#1)}}
\newcommand\stack[2][]{%
\settowidth{\regolo}{\ensuremath{#1}}\ifdim\regolo>0pt\ensuremath{\mathfrak{#2}(#1)}\else\ensuremath{\mathfrak{#2}}\fi}
\newcommand\xstack[2][]{%
\settowidth{\regolo}{\ensuremath{#1}}\ifdim\regolo>0pt\ensuremath{\mathfrak{#2}\left(#1\right)}\else\ensuremath{\mathfrak{#2}}\fi}
\newcommand\sheafhom[4][]{%
\ensuremath{\mathscr H\mspace{-5.4mu}\mathit{om}^{\stack{#1}}_{#2}(#3,#4)}}
\newcommand\sections[2][]{%
\settowidth{\regolo}{\ensuremath{#1}}\ifdim\regolo>0pt\ensuremath{\boldsymbol\Gamma_{#1}#2}\else\ensuremath{\boldsymbol\Gamma#2}\fi}
\newcommand\Des[2][]{%
\settowidth{\regolo}{\ensuremath{#1}}\ifdim\regolo>0pt\ensuremath{\underline{\mathcal D\mspace{-2.6mu}\mathit{es}}^{\mathfrak{#1}}(#2)}\else\ensuremath{\underline{\mathcal D\mspace{-2.6mu}\mathit{es}}(#2)}\fi}
\newcommand\V[2][]{%
\settowidth{\regolo}{\ensuremath{#1}}\ifdim\regolo>0pt\ensuremath{\mathit V^{\mathfrak{#1}}(#2)}\else\ensuremath{\mathit V(#2)}\fi}
\newcommand\R[2][]{%
\settowidth{\regolo}{\ensuremath{#1}}\ifdim\regolo>0pt\ensuremath{\mathit R^{\mathfrak{#1}}(#2)}\else\ensuremath{\mathit R(#2)}\fi}
\newcommand\xR[2][]{%
\settowidth{\regolo}{\ensuremath{#1}}\ifdim\regolo>0pt\ensuremath{\mathit R^{\mathfrak{#1}}\left(#2\right)}\else\ensuremath{\mathit R\left(#2\right)}\fi}
\newcommand\forget[2][]{%
\settowidth{\regolo}{\ensuremath{#1}}\ifdim\regolo=0pt\settowidth{\regolo}{\ensuremath{#2}}\ifdim\regolo>0pt\ensuremath{\boldsymbol{\omega}(#2)}\else\ensuremath{\boldsymbol{\omega}}\fi\else\settowidth{\regolo}{\ensuremath{#2}}\ifdim\regolo>0pt\ensuremath{\boldsymbol{\omega}^{\mathfrak{#1}}(#2)}\else\ensuremath{\boldsymbol{\omega}^{\mathfrak{#1}}}\fi\fi}
\newcommand\xforget[2][]{%
\settowidth{\regolo}{\ensuremath{#1}}\ifdim\regolo=0pt\ensuremath{\boldsymbol{\omega}\left(#2\right)}\else\ensuremath{\boldsymbol{\omega}^{\mathfrak{#1}}\left(#2\right)}\fi}
\newcommand\envelope[2][]{%
\settowidth{\regolo}{\ensuremath{#1}}\ifdim\regolo=0pt\settowidth{\regolo}{\ensuremath{#2}}\ifdim\regolo>0pt\ensuremath{\pi(#2)}\else\ensuremath{\pi}\fi\else\settowidth{\regolo}{\ensuremath{#2}}\ifdim\regolo>0pt\ensuremath{\pi^{\mathfrak{#1}}(#2)}\else\ensuremath{\pi^{\mathfrak{#1}}}\fi\fi}
\newcommand\xenvelope[2][]{%
\settowidth{\regolo}{\ensuremath{#1}}\ifdim\regolo=0pt\ensuremath{\pi\left(#2\right)}\else\ensuremath{\pi^{\mathfrak{#1}}\left(#2\right)}\fi}
\newcommand\Av[1]{%
\mathrm{Av}_{#1}}
\newcommand\f[2][]{\settowidth{\regolo}{\ensuremath{#1}}\ifdim\regolo=0pt\ensuremath{\sections{\mathscr{#2}}}\else\ensuremath{\mathscr{#2}}\fi}
\newcommand\E{\f[*]E}

\newcommand\fifu[1][\omega]{\ensuremath{\boldsymbol{#1}}}
\newcommand\tannakian[2][]{%
\settowidth{\regolo}{\ensuremath{#1}}\ifdim\regolo=0pt\settowidth{\regolo}{\ensuremath{#2}}\ifdim\regolo>0pt\ensuremath{\mathcal T(#2)}\else\ensuremath{\mathcal T}\fi\else\settowidth{\regolo}{\ensuremath{#2}}\ifdim\regolo>0pt\ensuremath{\mathcal T^{\mathfrak{#1}}(#2)}\else\ensuremath{\mathcal T^{\mathfrak{#1}}}\fi\fi}
\newcommand\xtannakian[2][]{%
\settowidth{\regolo}{\ensuremath{#1}}\ifdim\regolo=0pt\ensuremath{\mathcal T\left(#2\right)}\else\ensuremath{\mathcal T^{\mathfrak{#1}}\left(#2\right)}\fi}
\newcommand\Rf{\ensuremath{\mathscr R}}
\newcommand\Hfield[2][]{\settowidth{\regolo}{#1}\ifdim\regolo>0pt\ensuremath{\sections{\mathscr{#2}}}\else\ensuremath{\mathscr{#2}}\fi}%
\newcommand\Linear[1]{\ensuremath{\mathit L(#1)}}%
\newcommand\corners[1]{\ulcorner#1\urcorner}%
\newcommand\reprho{\rho}%
\newcommand\repsigma{\sigma}%
\newcommand\envehom{\ensuremath{\pi}}%
\newcommand\fibfunc[1][\omega]{\ensuremath{\boldsymbol{#1}}}%
\newcommand\tanngpd[2][]{\settowidth{\regolo}{#1}\ifdim\regolo>0pt{\settowidth{\regolo}{#2}\ifdim\regolo>0pt\ensuremath{\mathcal T(#2;#1)}\else\ensuremath{\mathcal T(#1)}\fi}\else{\settowidth{\regolo}{#2}\ifdim\regolo>0pt\ensuremath{\mathcal T(#2)}\else\ensuremath{\mathcal T}\fi}\fi}%
\newcommand\repfunc{\ensuremath{\mathscr{R}}}%
\renewcommand\i{\ensuremath{\mathit i}}

\makeindex

\begin{document}

\pagestyle{empty}\hbadness=1000\vbadness=1000

\begin{center}
{\LARGE Giorgio \ Trentinaglia}
\vskip 50pt\rule{\textwidth}{1pt}\vskip 10pt
{\LARGE\sffamily TANNAKA \ DUALITY \ FOR\\[10pt]PROPER \ LIE \ GROUPOIDS}\makebox[0pt][l]{\raisebox{1.55ex}\dag}
\vskip 10pt\rule{\textwidth}{1pt}\vskip 50pt
{\large\sffamily (PhD Thesis, \ Utrecht University, \ 2008)}
\end{center}
\vfill
\raisebox\depth\dag This work was financially supported by \mbox{Utrecht} University, the University of \mbox{Padua}, and a grant of the foundation ``Fon\-da\-zio\-ne \mbox{Ing.~Aldo~Gini}''

\newpage
\noindent\textbf{Abstract: } The main contribution of this thesis is a Tannaka duality theorem for proper Lie groupoids. This result is obtained by replacing the category of smooth vector bundles over the base manifold of a Lie groupoid with a larger category, the category of smooth Euclidean fields, and by considering smooth actions of Lie groupoids on smooth Euclidean fields. The notion of smooth Euclidean field that is introduced here is the smooth, finite dimensional analogue of the familiar notion of continuous Hilbert field. In the second part of the thesis, ordinary smooth representations of Lie groupoids on smooth vector bundles are systematically studied from the point of view of Tannaka duality, and various results are obtained in this direction.\vskip 10pt
\noindent\textbf{Keywords: } proper Lie groupoid, representation, tensor category, Tannaka duality, stack\vskip 10pt
\noindent\textbf{AMS Subject Classifications: } 58H05, 18D10
\vfill
\noindent\textbf{Acknowledgements: } I would like to thank my supervisor, \mbox{I.~Moerdijk}, for having suggested the research problem out of which the present work took shape and for several useful remarks, and also \mbox{M.~Crainic} and \mbox{N.~T.~Zung}, for their interest and for helpful conversations.

\pagestyle{headings}\hbadness=1000\vbadness=1000

\cleardoublepage
\setcounter{chapter}{-1}\refstepcounter{chapter}\addcontentsline{toc}{chapter}{Table of Contents}
\tableofcontents

\pagestyle{myheadings}\hbadness=1000\vbadness=10000

\cleardoublepage
\chapter*{Introduction}\markboth{INTRODUCTION}{INTRODUCTION}\addcontentsline{toc}{chapter}{Introduction}
Although a rigorous formulation of the problem with which this doctoral thesis is concerned will be possible only after the central ideas of Tannaka duality theory have been at least briefly discussed, I can nevertheless start with some comments about the general context where such a problem takes its appropriate place. Roughly speaking, my study aims at a better understanding of the relationship that exists between a given Lie groupoid and the corresponding category of representations. First of all, for the benefit of non-specialists, I want to explain the reasons of my interest in the theory of Lie groupoids (a precise definition of the notion of Lie groupoid can be found in \refsez{O.1.5} of this thesis) by drawing attention to the principal applications that justify the importance of this theory; in the second place, I intend to undertake a critical examination of the concept of representation in order to convince the reader of the naturalness of the notions I will introduce below.

\sottosezione{From Lie groups to Lie groupoids}\markright{From Lie groups to Lie groupoids}\addcontentsline{toc}{section}{From Lie groups to Lie groupoids}

Groupoids make their appearance in diverse mathematical contexts. As the name `groupoid' suggests, this notion generalizes that of group. In order to explain how and to make the definition more plausible, it is best to start with some examples.

The reader is certainly familiar with the notion of fundamental group of a topological space. The construction of this group presupposes the choice of a base point, and any two such choices give rise to the same group provided there exists a path connecting the base points (for this reason one usually assumes that the space is path connected). However, instead of considering only paths starting and ending at the same point, one might more generally allow paths with arbitrary endpoints; two such paths can still be composed as long as the one starts where the other ends. One obtains a well-defined associative partial operation on the set of homotopy classes of paths with fixed endpoints, for which the (classes of) constant paths are both left and right neutral elements. Observe that each path has a two-sided inverse, namely the path itself with reverse orientation.

In geometry, groups are usually groups of transformations\inciso{or symmetries}of some object or space. If $g$ is an element of a group $G$ acting on a space $X$ and $x$ is a point of $X$, one may think of the pair $(g,x)$ as an arrow going from $x$ to ${g\cdot x}$; again, two such arrows can be composed in an obvious way, by means of the group operation of $G$, provided one starts where the other ends. Composition of arrows is an associative partial operation on the set ${G\times X}$, which encodes both the multiplication law of the group $G$ and the $G$-action on $X$.

In the representation theory of groups, the linear group $\GL(V)$ associated with a finite dimensional vector space $V$ plays a fundamental role. If a vector bundle $E$ over a space $X$ is given instead of a single vector space $V$, one can consider the set $\GL(E)$ of all triples $(x,x',\lambda)$ consisting of two points of $X$ and a linear isomorphism $\lambda: E_x \isoto E_{x'}$ between the fibres over these points. As in the examples above, an element $(x,x',\lambda)$ of this set can be viewed as an arrow going from $x$ to $x'$; such an arrow can be composed with another one as long as the latter has the form $(x',x'',\lambda')$. Arrows of the form $(x,x,\id)$ are both left and right neutral elements for the resulting associative partial operation, and each arrow admits a two-sided inverse.

By abstraction from these and similar examples, one is led to consider small categories where every arrow is invertible. Such categories are referred to as groupoids. More explicitly, a groupoid consists of a space $X$ of ``base points'' (also called objects), a set \G\ of ``arrows'', endowed with source and target projections $\s, \t: \G \to X$, and an associative partial composition law $\mca\G \to \G$ (defined for all pairs of arrows $(g',g)$ with the property that the source of $g'$ equals the target of $g$), such that in correspondence with each point $x$ of $X$ there is a (necessarily unique) ``neutral'' or ``unit'' arrow, often itself denoted by $x$, and every arrow is invertible.

The notion of Lie groupoid generalizes that of Lie group. Much the same as a Lie group is a group endowed with a smooth manifold structure compatible with the multiplication law and with the operation of taking the inverse, a Lie groupoid is a groupoid where the sets $X$ and \G\ are endowed with a smooth manifold structure that makes the various maps which arise from the groupoid structure smooth. For instance, in each of the examples above one obtains a Lie groupoid when the space $X$ of base points is a smooth manifold, $G$ is a Lie group acting smoothly on $X$ and $E$ is a smooth vector bundle over $X$; these Lie groupoids are respectively called the {\em fundamental groupoid} of the manifold $X$, the {\em translation groupoid} associated with the smooth action of $G$ on $X$ and the {\em linear groupoid} associated with the smooth vector bundle $E$. There is also a more general notion of {\em $\C^\infty$\nobreakdash-structured groupoid,} about which we shall spend a few words later on in the course of this introduction, which we introduce in our thesis in order to describe certain groupoids that arise naturally in the study of Tannaka duality theory.

In the course of the second half of the twentieth century the notion of groupoid turned out to be very useful in many branches of mathematics, although this notion had in fact already been in the air since the earliest accomplishments of quantum mechanics\inciso{think, for example, of \mbox{Heisenberg's} formalism of matrices}or, more back in time, since the first investigations into classification problems in geometry. Nowadays, the theory of Lie groupoids constitutes the preferred language for the geometrical study of foliations \cite{Moerdijk&Mrcun'03}; the same theory has applications to noncommutative geometry \cite{Co94,CdSW} and quantization deformation theory \cite{Lm98}, as well as to symplectic and \mbox{Poisson} geometry \cite{We87,CDW87,DZ05}. Another source of examples comes from the study of orbifolds \cite{Mo02}; this subject is connected with the theory of stacks, which originated in algebraic geometry from \mbox{Grothendieck's} suggestion to use groupoids as the right notion to understand moduli spaces.

\separazione

\noindent When trying to extend representation theory from Lie groups to Lie groupoids, one is first of all confronted with the problem of defining a suitable notion of representation for the latter. As far as we are concerned, we would like to generalize the familiar notion of (finite dimensional) Lie group representation, by which one generally means a homomorphism $G \to \GL(V)$ of a Lie group $G$ into the group of automorphisms of some finite dimensional vector space $V$, so that as many constructions and results as possible can be adapted to Lie groupoids without essential changes; in particular, we would like to carry over Tannaka duality theory (see the next subsection) to the realm of Lie groupoids.

The notion of Lie group representation recalled above has an obvious naive extension to the groupoid setting. Namely, a representation of a Lie groupoid \G\ can be defined as a Lie groupoid homomorphism $\G \to \GL(E)$ (smooth functor) into the linear groupoid associated with some smooth vector bundle $E$ over the manifold of objects of \G. Any such representation assigns each arrow $x \to x'$ of \G\ a linear isomorphism $E_x \isoto E_{x'}$ in such a way that composition of arrows is respected. In our dissertation we will use the term `classical representation' to refer to this notion. Unfortunately, classical representations prove to be completely inadequate for the above-mentioned purpose of carrying forward Tannaka duality to Lie groupoids; we shall say something more about this matter later.

The preceding consideration leads us to introduce a different notion of representation for Lie groupoids. In doing this, however, we adhere to the point of view that the latter should be as close as possible to the notion of classical representation\inciso{in particular the new theory should extend the theory of classical representations}and that moreover in the case of groups one should recover the usual notion of representation recalled above.

\sottosezione{Historical perspective on Tannaka duality}\markright{Historical perspective on Tannaka duality}\addcontentsline{toc}{section}{Historical perspective on Tannaka duality}

It has been known for a long time, and precisely since the pioneer work of \textit{Pontryagin} and \textit{van~Kampen} in the 1930's, that a commutative locally compact group can be identified with its own bidual. Recall that if $G$ is such a group then its dual is the group formed by all the characters on $G$, that is to say the continuous homomorphisms of $G$ into the multiplicative group of complex numbers of absolute value one, the group operation being given by pointwise multiplication of complex functions; one may regard the latter group as a topological group\inciso{in fact, a locally compact one}by taking the topology of uniform convergence on compact subsets. There is a canonical pairing between $G$ and this dual, given by pointwise evaluation of characters at elements of $G$, which induces a continuous homomorphism of $G$ into its own bidual. Then one can prove that the latter correspondence is actually an isomorphism of topological groups; see for instance \textit{Dixmier~(1969)} \cite{Dix69}, \textit{Rudin~(1962)} \cite{Ru62}, or the book by \textit{Chevalley~(1946)} \cite{Che46}.

When one tries to generalize this duality result to non-\mbox{Abelian} locally compact groups, such as for instance Lie groups, it becomes evident that the whole ring of representations must be considered because characters are no longer sufficient to recapture the group. However, it is still an open problem to formulate and prove a general duality theorem for noncommutative Lie groups: even the case of simple algebraic groups is not well understood, despite the enormous accumulating knowledge on their irreducible representations. The situation is quite the opposite when the group is {\em compact,} because the dual object $\dual G$ of a compact group $G$ is discrete and so belongs to the realm of algebra: in this case, there is a good duality theory due to \mbox{\textit{H.~Peter,}} \mbox{\textit{H.~Weyl}} and \mbox{\textit{T.~Tannaka,}} which we now proceed to recall.

The early duality theorems of \mbox{\textit{Tannaka~(1939)}} \cite{Tan39} and \mbox{\textit{Krein~(1949)}} \cite{Kr49} concentrate on the problem of reconstructing a compact group from the ring of isomorphism classes of its representations. Owing to the ideas of \mbox{Grothendieck} \cite{Sa72}, these results can nowadays be formulated within an elegant categorical framework. Although we do not intend to enter into details now, these ideas are implicit in what we are about to say.

1. One starts by considering the category \R[\mathrm{0}]{G} of all continuous finite dimensional representations of the compact group $G$: the objects of \R[\mathrm{0}]{G} are the pairs $(V,\varrho)$ consisting of a finite dimensional real vector space $V$ and a continuous homomorphism $\varrho: G \to \GL(V)$; the morphisms are precisely the $G$\nobreakdash-equivariant linear maps.

2. There is an obvious functor \forget{} of the category \R[\mathrm{0}]{G} into that of finite dimensional real vector spaces, namely the forgetful functor $(V,\varrho) \mapsto V$. The natural endomorphisms of \forget{} form a topological algebra $\End(\forget{})$, when one endows $\End(\forget{})$ with the coarsest topology making each map $\lambda \mapsto \lambda(R)$ continuous as $R$ ranges over all objects of \R[\mathrm{0}]{G}.

3. The subset \tannakian{G} of this algebra, formed by the elements compatible with the tensor product operation on representations, in other words the natural endomorphisms $\lambda$ of \forget{} such that $\lambda({R\otimes R'}) = {\lambda(R) \otimes \lambda(R')}$ and $\lambda(\TU) = \id$, proves to be a compact group.

4. {\em (Tannaka)} The canonical map $\envelope{}: G \to \End(\forget{})$, defined by setting $\envelope{g}(R) = \varrho(g)$ for each object $R = (V,\varrho)$ of \R[\mathrm{0}]{G}, establishes an isomorphism of topological groups between $G$ and \tannakian{G}.

\sottosezione{What is new in this thesis}\markright{What is new in this thesis}\addcontentsline{toc}{section}{What is new in this thesis}

We are now ready to give a short summary of the original contributions of the present study.

Within the realm of Lie groupoids, proper groupoids play the same role as compact groups; for example, all isotropy groups of a proper Lie groupoid are compact (the isotropy group at a base point $x$ consists of all arrows $g$ with $\s(g) = \t(g) = x$). The main result of our research is a {\em Tannaka duality theorem for proper Lie groupoids,} which takes the following form.

To begin with, we construct, for each smooth manifold $X$, a category whose objects we call {\em smooth fields over $X$;} our notion of smooth field is the analogue, in the smooth and finite dimensional setting in which we are interested, of the familiar notion of continuous Hilbert field introduced by \textit{\mbox{Dixmier} and \mbox{Douady}} in the early 1960's \cite{DD63} (see also \mbox{\em Bos} \cite{Bos06} or \mbox{\em Kali\v snik} \cite{Kal07} for more recent work related to continuous Hilbert fields). The category of smooth fields is a proper enlargement of the category of smooth vector bundles. Like for vector bundles, one can define a notion of Lie groupoid representation on a smooth field in a completely standard way. Given a Lie groupoid \G, such representations and their obvious morphisms form a category that is related to the category of smooth fields over the base manifold $M$ of \G\ by means of a forgetful functor of the former into the latter category. To this functor one can assign, by generalizing the construction explained above in the case of groups, a groupoid over $M$, to which we shall refer as the {\em Tannakian groupoid associated with \G,} to be denoted by \tannakian\G, endowed with a natural candidate for a smooth structure on the space of arrows {\em ($\C^\infty$\nobreakdash-structured groupoid).} As for groups, there is a canonical homomorphism \envelope{} of \G\ into \tannakian{\G} that turns out to be compatible with this $\C^\infty$\nobreakdash-structure.

Our Tannaka duality theorem for proper Lie groupoids reads as follows:
\begin{theorem*}
Let \G\ be a proper Lie groupoid. The $\C^\infty$\nobreakdash-structure on the space of arrows of the Tannakian groupoid \tannakian{\G} is a genuine manifold structure so that \tannakian{\G} is a Lie groupoid. The canonical homomorphism \envelope{} is a Lie groupoid isomorphism $\G \can \tannakian\G$.
\end{theorem*}
The main point here is to prove the surjectivity of the homomorphism \envelope{}; the fact that \envelope{} is injective is a direct application of a theorem of \textit{N.T.~Zung.}

Actually, the reasonings leading to our duality theorem also hold, for the most part, for the representations of a proper Lie groupoid on vector bundles. Since from the very beginning of our research we were equally interested in studying such representations, we found it convenient to provide a general theoretical framework where the diverse approaches to the representation theory of Lie groupoids could take their appropriate place, so as to state our results in a uniform language. The outcome of such demand was the theory of {\em `smooth tensor stacks'.} Smooth vector bundles and smooth fields are two examples of smooth tensor stacks. Each smooth tensor stack gives rise to a corresponding notion of representation for Lie groupoids; then, for each Lie groupoid one obtains, by the same general procedure outlined above, a corresponding Tannakian groupoid, which will depend very much, in general, on the initial choice of a smooth tensor stack (for example, Tannaka duality fails in the context of representations on vector bundles).

Our remaining contributions are mainly concerned with the study of Tannakian groupoids arising from representations of \textsl{proper} Lie groupoids on \textsl{vector bundles.} Since in this case the reconstructed groupoid may not be isomorphic to the original one, the problem of whether the aforesaid standard $\C^\infty$\nobreakdash-structure on the space of arrows of the Tannakian groupoid turns the latter groupoid into a Lie groupoid becomes considerably more interesting and difficult than in the case of representations on smooth fields. Our principal result in this direction is that the answer to the indicated question is affirmative for all proper \textsl{regular} groupoids. In connection with this result we prove invariance of the solvability of the problem under Morita equivalence. Finally, we provide examples of {\em classically reflexive} proper Lie groupoids, i.e.\ proper Lie groupoids for which the groupoid reconstructed from the representations on vector bundles is isomorphic to the original one; however, our list is very short: failure of reflexivity is the rule rather than the exception when one deals with representations on vector bundles.

\sottosezione{Outline chapter by chapter}\markright{Outline chapter by chapter}\addcontentsline{toc}{section}{Outline chapter by chapter}

In order to help the reader find their own way through the dissertation, we give here a detailed account of how the material is organized.

\begin{center}
$* \quad * \quad *$
\end{center}

\noindent In Chapter~\refcpt1 we recall basic notions and facts concerning Lie groupoids.

The initial section is mainly about definitions, notation and conventions to be followed in the sequel.

The second section contains relatively more interesting material: after briefly recalling the familiar notion of a representation of a Lie groupoid on a vector bundle (classical representation), we supply a concrete example,\footnote{We discovered this counterexample independently, though it turned out later that the same had already been around for some time \cite{Lueck&Oliver'01}.} which motivates our introducing the notion of representation on a smooth field in Chapter~\refcpt{4}, showing that it is in general impossible to distinguish two Lie groupoids from one another just on the basis of knowledge of the corresponding categories of representations on vector bundles; more precisely, we shall explicitly construct a principal $\mathit T^2$\nobreakdash-bundle over the circle (where $\mathit T^k$ denotes the $k$\nobreakdash-torus), together with a homomorphism onto the trivial $\mathit T^1$\nobreakdash-bundle over the circle, such that the obvious pull-back of representations along this homomorphism yields an isomorphism between the categories of classical representations of these two bundles of Lie groups.

In Section~\ref{O.3.2} we review the notion of a {\em (normalized) Haar system} on a Lie groupoid; this is the analogue, for Lie groupoids, of the notion of (probability) Haar measure on a group. Like probability Haar measures, normalized Haar systems can be used to obtain invariant functions, metrics etc.\ by means of the usual averaging technique. The possibility of constructing equivariant maps lies at the heart of our proof that the homomorphism \envelope{} mentioned above is surjective for every proper Lie groupoid.

Section~\ref{N.4} introduces the reader to a relatively recent result obtained by \mbox{\em N.T.~Zung} about the local structure of proper Lie groupoids; this general result was first conjectured by \mbox{\em A.~Weinstein} in his famous paper about the local linearizability of proper regular groupoids \cite{Weinstein'02} (where the result is proved precisely under the additional assumption of regularity). Zung's local linearizability theorem states that each proper Lie groupoid \G\ is, locally in the vicinity of any given \G\nobreakdash-invariant point of its base manifold, isomorphic to the translation groupoid associated with the induced linear action of the isotropy group of \G\ at the point itself on the respective tangent space. As a consequence of this, every proper Lie groupoid is locally \mbox{Morita} equivalent to the translation groupoid associated with some compact Lie group action. The local linearizability of proper Lie groupoids accounts for the injectivity of the homomorphism \envelope{}.

Finally, in Section~\ref{N.5}, we prove a statement relating the global structure up to Morita equivalence of a proper Lie groupoid and the existence of globally faithful representations: precisely, we show that a proper Lie groupoid admits a globally faithful representation on a smooth vector bundle if and only if it is Morita equivalent to the translation groupoid of a compact Lie group action. Although this result is not elsewhere used in our work, we present a proof of it here because we believe that the same technique, applied to representations on smooth fields, may be used to obtain nontrivial information about the global structure of arbitrary proper Lie groupoids (since every such groupoid trivially admits globally faithful representations on smooth fields).

\begin{center}
$* \quad * \quad *$
\end{center}

\noindent Chapter~\refcpt{2} is mainly concerned with the background notions needed in order to formulate precisely the reconstruction problem in full generality. The formal categorical framework within which this problem is most conveniently stated in the language of tensor categories and tensor functors.

Section~\ref{O.1.1+N.6} introduces the pivotal notion of a {\em tensor category:} this will be, for us, an additive $k$\nobreakdash-linear category \Kt\ ($k$ = real or complex numbers) endowed with a bilinear bifunctor $(A,B) \mapsto {A\otimes B}: {\Kt\times\Kt} \to \Kt$ called a {\em tensor product,} a distinguished object \TU\ called the {\em tensor unit} and various natural isomorphisms called \textit{$\mathit{ACU}$~constraints} which, roughly speaking, make the product $\otimes$ associative and commutative with neutral element \TU. The notion of {\em rigid} tensor category is also briefly recalled: this is a tensor category with the property that each object $R$ admits a dual, that is an object $R'$ for which there exist morphisms ${R'\otimes R} \to \TU$ and $\TU \to {R\otimes R'}$ compatible with one another in an obvious sense; the category of finite dimensional vector spaces\inciso{or, more generally, smooth vector bundles over a manifold}is an example.

In Section~\ref{O.1.2} we review the notions of a {\em tensor functor} (morphism of tensor categories) and a {\em tensor preserving natural transformation} (morphism of tensor functors): one obtains a tensor functor by attaching, to an ordinary functor $F$, (natural) isomorphisms ${F(A) \otimes F(B)} \can F({A\otimes B})$ and $\TU \can F(\TU)$, called {\em tensor functor constraints,} compatible with the \mbox{$\mathit{ACU}$}~constraints of the two tensor categories involved; a tensor preserving natural transformation of tensor functors is simply an ordinary natural transformation $\lambda$ such that $\lambda({A\otimes B}) = {\lambda(A) \otimes \lambda(B)}$ and $\lambda(\TU) = \id$ up to the obvious identifications provided by the tensor functor constraints. If an object $R$ admits a dual $R'$ in the above sense, then $\lambda(R)$ is an isomorphism for any tensor preserving $\lambda$ (a tensor preserving functor will preserve duals whenever they exist). A fundamental example of tensor functor is the pull-back of smooth vector bundles along a smooth mapping of manifolds.

Section~\ref{N.8} hints at the relationship between real and complex theory: to mention one example, in the case of groups one can either consider linear representations on real vector spaces and then take the group of all tensor preserving natural automorphisms of the standard forgetful functor or, alternatively, consider linear representations on complex vector spaces and then take the group of all self-conjugate tensor preserving natural automorphisms; these two groups, of course, will turn out to be the same. We indicate how these comments may be generalized to the abstract categorical setting we have just outlined to the reader.

Section~\ref{O.1.3} is devoted to a concise exposition, without any ambition to completeness, of the algebraic geometer's point of view on Tannaka duality. In fact, many fundamental aspects of the algebraic theory are omitted here; we refer more demanding readers to \textit{\mbox{Saavedra} (1972)} \cite{Sa72}, \textit{\mbox{Deligne} and \mbox{Milne} (1982)} \cite{DeMi82} and \textit{\mbox{Deligne} (1990)} \cite{De90}. We thought it necessary to include this exposition with the intent of providing adequate grounds for understanding certain questions reaised in Chapter~\refcpt{5}.

Contrary to the rest of the chapter, Section~\ref{O.1.4+.3.3} is entirely based on our own work. In this section we prove a key technical lemma which we exploit later on, in Section~\ref{N.20}, to establish the surjectivity of the envelope homomorphism \envelope{} (see above) for all proper Lie groupoids; this lemma reduces the latter problem to that of checking that a certain extendability condition for morphisms of representations is satisfied. The proof of our result makes use of the classical Tannaka duality theorem for compact (Lie) groups, though for the rest it is purely algebraic and it does not reproduce any known argument.

\begin{center}
$* \quad * \quad *$
\end{center}

\noindent In Chapter~\refcpt{3}, we introduce our abstract systematization of representation theory. Our ideas took shape gradually, during the attempt to make the treatment of various inequivalent approaches to the representation theory of Lie groupoids uniform. A collateral benefit of this abstraction effort was a gain in simplicity and formal elegance, along with a general better understanding of the mathematical features of the theory itself.

We begin with the description of a certain categorical structure, that we shall call {\em fibred tensor category,} which permits to make sense of the notion of `Lie groupoid action' in a natural way. Smooth vector bundles and smooth fields provide examples of such a structure. A fibred tensor category \stack{C} may be defined as a correspondence that assigns a tensor category \stack[X]{C} to each smooth manifold $X$ and a tensor functor $f^*: \stack[X]{C} \to \stack[Y]{C}$ to each smooth mapping $f: Y \to X$, along with a coherent system of tensor preserving natural isomorphisms $({g\circ f})^* \can {f^*\circ g^*}$ and $\id^* \can \Id$. Most notions needed in representation theory can be defined purely in terms of the fibred tensor category structure, provided this enjoys some additional properties which we now proceed to summarize.

In Section \ref{N.11}, we make from the outset the assumption that \stack{C} is a {\em prestack,} in other words that the obvious presheaf $U \mapsto \Hom_{\stack[U]{C}}(E|_U,F|_U)$ is a sheaf on $X$ for all objects $E$, $F$ of the category \stack[X]{C}. We also require \stack{C} to be {\em smooth,} that is to say, roughly speaking, that for each $X$ there is an isomorphism of complex algebras $\End(\TU_X) \simeq \C^\infty(X)$, where $\TU_X$ denotes the tensor unit in \stack[X]{C}.

Let \smooth[X] denote the sheaf of smooth functions on $X$. For each smooth prestack \stack{C} one can associate to every object $E$ of the category \stack[X]{C} a sheaf of \smooth[X]\nobreakdash-modules, \sections E, to be called the {\em sheaf of smooth sections} of $E$. The latter operation yields a functor of \stack[X]{C} into the category of sheaves of \smooth[X]\nobreakdash-modules. One has a natural transformation ${\sections{E} \otimes_{\smooth[X]} \sections{E'}} \to \sections{({E\otimes E'})}$, which need not be an isomorphism, and an isomorphism $\smooth[X] \simeq \sections{(\TU_X)}$ of \smooth[X]\nobreakdash-modules, that behave much as usual tensor functor constraints do. The compatibility of the operation $E \mapsto \sections{E}$ with the pullback along a smooth map $f: Y \to X$ is measured by a canonical natural morphism of sheaves of \smooth[Y]\nobreakdash-modules ${f^*(\sections E)} \to \sections{({f^*E})}$. For each point $x$ of $X$, there is a functor which assigns, to every object $E$ of the category \stack[X]{C}, a complex vector space $E_x$ to be referred to as the {\em fibre of $E$ at $x$;} a local smooth section $\zeta \in \sections E(U)$, defined over an open neighbourhood $U$ of $x$, will determine a vector $\zeta(x) \in E_x$ to be referred to as the {\em value of $\zeta$ at $x$.}

In order to show that Morita equivalences have the usual property of inducing a categorical equivalence between the categories of representations, we further need to impose the condition that \stack C is a {\em stack.} This condition, examined in Section \ref{N.12}, means that when one is given an open cover $\{U_i\}$ of a (paracompact) manifold $M$, along with a family of objects $E_i \in {\Ob\,\stack[U_i]{C}}$ and a cocycle of isomorphisms $\theta_{ij}: E_i|_{U_i\cap U_j} \isoto E_j|_{U_i\cap U_j}$, there must be some object $E$ in \stack[M]{C} which admits a family of isomorphisms $E|_{U_i} \isoto E_i \in \stack[U_i]{C}$ compatible with $\{\theta_{ij}\}$. Naively speaking, one can glue objects in \stack C together. When \stack C is a smooth stack, the category \stack[M]{C} will essentially contain the category of all smooth vector bundles over $M$ as a full subcategory.

In Section \ref{N.13}, we lay down the foundations of the representation theory of Lie groupoids relative to a {\em type \stack T,} for an arbitrary smooth stack of tensor categories \stack T. A {\em representation of type \stack T} of a Lie groupoid \G\ is a pair $(E,\varrho)$ consisting of an object $E$ of the category \stack[M]{T} (where $M$ is the base of \G) and an arrow $\varrho: {\s^*E} \to {\t^*E}$ in the category \stack[\G]{T} (where $\s, \t: \G \to M$ are the source resp.\ target map of \G) such that ${\u^*\varrho} = \id_E$ (where $\u: M \to \G$ denotes the unit section) and ${\mathit m^*\varrho} = {{{\p_1}^*\varrho} \circ {{\p_2}^*\varrho}}$ (where $\mathit m, \p_1, \p_2: \mca\G \to \G$ respectively denote multiplication, first and second projection). With the obvious notion of morphism, representations of type \stack T of a Lie groupoid \G\ form a category \R[T]\G. This category inherits an additive linear tensor structure from the base category \stack[M]{T}, making the forgetful functor $(E,\varrho) \mapsto E$ a strict linear tensor functor of \R[T]{\G} into \stack[M]{T}. The latter functor will be denoted by \forget[T]{\G} and will be called the {\em standard fibre functor of type \stack T} associated with \G.

Each homomorphism of Lie groupoids $\phi: \G \to \H$ induces a linear tensor functor $\phi^*: \R[T]{\H} \to \R[T]{\G}$ that we call the {\em pullback along $\phi$.} One has tensor preserving natural isomorphisms $({\psi\circ\phi})^* \can {\phi^* \circ \psi^*}$. In Section \ref{N.14} we show that for every Morita equivalence $\phi: \G \to \H$ the pullback functor $\phi^*$ is an equivalence of tensor categories.

\begin{center}
$* \quad * \quad *$
\end{center}

\noindent Chapter~\refcpt{4} is the core of our dissertation. This is the place where we describe the general duality theory for Lie groupoids in the abstract framework of Chapters~\refcpt[2]{3} and where we prove our most important results, culminating in the above-mentioned reconstruction theorem for proper Lie groupoids.

Section~\ref{N.15} contains a detailed description of in what type of Lie groupoid representations one should be interested, from our point of view, when dealing with duality theory of Lie groupoids. Namely, we say that a type \stack T is a {\em stack of smooth fields} if it meets a number of extra requirements, called `axioms', which we now proceed to summarize.

\sloppy
Our first axiom says that the canonical morphisms ${\sections E \otimes_{\smooth[X]} \sections{E'}} \to \sections{({E\otimes E'})}$ and ${f^*(\sections E)} \to \sections{({f^*E})}$ (cfr.\ the summary of Ch.~\refcpt{3}, \refsez{N.11}) are surjective; this axiom conveys information about the smooth sections of ${E\otimes E'}$ and ${f^*E}$ and it implies that the fibre at $x$ of an object $E$ is spanned, as a vector space, by the values $\zeta(x)$ as $\zeta$ ranges over all germs of local smooth sections of $E$ at $x$.

\fussy
Next, recall that any arrow $a: E \to E'$ in \stack[X]T induces a morphism of sheaves of \smooth[X]\nobreakdash-modules $\sections a: \sections E \to \sections{E'}$ and a bundle of linear maps $\{a_x: E_x \to {E'}_x\}$; these are mutually compatible, in an obvious sense. Our second and third axioms completely characterize the arrows in \stack[X]T in terms of their effect on smooth sections and the bundles of linear maps they induce; namely, an arrow $a: E \to E'$ vanishes if and only if $a_x$ vanishes for all $x$, and every pair formed by a morphism of \smooth[X]\nobreakdash-modules $\alpha: \sections E \to \sections{E'}$ and a compatible bundle of linear maps $\{\lambda_x: E_x \to {E'}_x\}$ gives rise to a (unique) arrow $a: E \to E'$ such that $\alpha = \sections a$ or, equivalently, $\lambda_x = a_x$ for all $x$.

Then there is an axiom requiring the existence of local Hermitian metrics on the objects of \stack[X]T. A {\em Hermitian metric on $E$} is an arrow ${E\otimes E^*} \to \TU$ inducing a positive definite Hermitian sesquilinear form on each fibre $E_x$; the axiom says that for any paracompact $M$, each object of \stack[M]T admits Hermitian metrics. This assumption has many useful consequences: for example, it implies various continuity principles for smooth sections and a fundamental extension property for arrows.

The remaining two axioms impose various finiteness conditions on \stack T: roughly speaking, finite dimensionality of the fibres of an arbitrary object $E$ and local finiteness of the sheaf of modules \sections E. More precisely, one axiom canonically identifies \stack[\pt]T, as a tensor category, with the category of finite dimensional vector spaces\inciso{where \pt\ denotes the one-point manifold}so that, for instance, the functor $E \mapsto E_x$ becomes a tensor functor of \stack[X]T into the category of such spaces; the other axiom requires the existence, for each point $x$, of an open neighbourhood $U$ such that $\sections E(U)$ is spanned, as a $\C^\infty(U)$\nobreakdash-module, by a finite set of sections of $E$ over $U$.

In Section~\ref{N.16}, we introduce our fundamental example of a stack of smooth fields (which is to play a role in our reconstruction theorem for proper Lie groupoids in \refsez{N.20}), to which we refer as the type \stack{E^\infty} of {\em smooth Euclidean fields.} The notion of smooth Euclidean field over a manifold $X$ generalizes that of smooth vector bundle over $X$ in that the dimension of the fibres is allowed to vary discontinuously over $X$ or, in other words, the sheaf of smooth sections is no longer a locally free \smooth[X]\nobreakdash-module. Our theory of smooth Euclidean fields may be regarded as the counterpart, in the smooth setting, of the well-established theory of continuous Hilbert fields \cite{DD63}.

In Section~\ref{N.17} we prove various results about the equivariant extension of morphisms of Lie groupoid representations whose type is a stack of smooth fields; in combination with the technical lemma of \refsez{O.1.4+.3.3}, these extension results allow one to establish the surjectivity of the envelope homomorphism \envelope{} associated with representations on an arbitrary stack of smooth fields. The proofs are based on the usual averaging technique\inciso{which makes sense for any proper Lie groupoid because of the existence of normalized Haar systems}and, of course, on the axioms for stacks of smooth fields.

In Sections~\ref{N.18}--\ref{N.19}, we delve into the formalism of fibre functors with values in an arbitrary stack of smooth fields. A {\em fibre functor,} with values in a stack of smooth fields \stack F, is a faithful linear tensor functor \fifu\ of some additive tensor category \Kt\ into \stack[M]F, for some fixed paracompact manifold $M$ to be referred to as the {\em base} of \fifu. This notion is obtained by abstracting the fundamental features, which allow one to make sense of the construction of the Tannakian groupoid, from the concrete example provided by the standard forgetful functor associated with the representations of type \stack F of a Lie groupoid over $M$. To any fibre functor \fifu\ with base $M$, one can assign a groupoid \tannakian{\fifu} over $M$ to which we refer as the {\em Tannakian groupoid} associated with \fifu\ constructed, like in the case of groups, by taking all tensor preserving natural automorphisms of \fifu. The set of arrows of \tannakian{\fifu} comes naturally equipped with a topology and a {\em smooth functional structure} that is a sheaf $\Rf^\infty$ of algebras of continuous real valued functions on \tannakian{\fifu} closed under composition with arbitrary smooth functions $\nR^d \to \nR$; the notion of smooth functional structure is analogous to that of \textit{$\C^\infty$-ring,} cfr \cite{MoeRey86,MoeRey91}.

In Section~\ref{N.20}, we reap the fruits of all our previous work and prove several statements of fundamental importance about the Tannakian groupoid \tannakian{\G} associated with the standard forgetful functor \forget{\G} on the category of representations of an arbitrary proper Lie groupoid \G. (We are still dealing with a situation where the type is an arbitrary stack of smooth fields.) Recall that there is a canonical homomorphism $\envelope{}: \G \to \tannakian{\G}$ defined by setting $\envelope{g}(E,\varrho) = \varrho(g)$, which, as previously mentioned, turns out to be surjective for proper \G; the proof of this theorem is based on the results of Sections~\ref{O.1.4+.3.3} and~\ref{N.17}. Moreover, when \G\ is proper, the Tannakian groupoid \tannakian{\G} becomes a topological groupoid and \envelope{} a homomorphism of topological groupoids: then we show that injectivity of \envelope{} implies that \envelope{} is an isomorphism of topological groupoids and that this in turn implies that the above-mentioned functional structure on \tannakian{\G} is actually a Lie groupoid structure for which \envelope{} becomes an isomorphism of Lie groupoids. Accordingly, we say that a Lie groupoid \G\ is {\em reflexive}\inciso{relative to a certain type}if \envelope{} induces a homeomorphism between the spaces of arrows of \G\ and \tannakian{\G}. Our main theorem, which concludes the section, states that every proper Lie groupoid is reflexive relative to the type \stack{E^\infty} of smooth Euclidean fields. The injectivity of \envelope{} for this particular type of representations is an easy consequence of Zung's local linearizabilty result for proper Lie groupoids.

\begin{center}
$* \quad * \quad *$
\end{center}

\noindent Besides establishing a Tannaka duality theory for proper Lie groupoids, the work described above also leads to results concerning the classical theory of representations of Lie groupoids on vector bundles. Chapter~\refcpt{5} concentrates on what can be said about the latter case exclusively from the abstract standpoint of the theory of fibre functors outlined in \refsez[N.18]{N.19}. The main objects of study here are certain fibre functors, which will be referred to as {\em classical fibre functors,} enjoying formal properties analogous to those possessed by the standard forgetful functor associated with the category of classical representations of a Lie groupoid.

The distinctive features of classical fibre functors are the rigidity of the domain tensor category \Kt\ and the type being equal to the stack of smooth vector bundles. Section~\ref{N.21} collects some general remarks about such fibre functors and some basic definitions. For any classical fibre functor \fifu, the Tannakian groupoid \tannakian{\fifu} proves to be a $\C^\infty$\nobreakdash-structured groupoid over the base $M$ of \fifu; this means that all structure maps of \tannakian{\fifu} are morphisms of functionally structured spaces with respect to the $\C^\infty$\nobreakdash-functional structure $\Rf^\infty$ on \tannakian{\fifu} introduced in \refsez{N.18}. One can define, for every $\C^\infty$\nobreakdash-structured groupoid \tannakian{}, an obvious notion of $\C^\infty$\nobreakdash-representation on a smooth vector bundle; such representations form a tensor category \R[\infty]{\tannakian{}}. Every object $R$ of the domain category \Kt\ of a classical fibre functor \fifu\ determines a $\C^\infty$\nobreakdash-representation $\ev_R$, which we call {\em evaluation at $R$,} of the Tannakian groupoid \tannakian{\fifu} on the vector bundle $\fifu(R)$. The operation $R \mapsto \ev_R$ provides a tensor functor of \Kt\ into the category of $\C^\infty$\nobreakdash-representations of \tannakian{\fifu}, the {\em evaluation functor} associated with \fifu.

Section~\ref{O.2.1} is preliminary to Section~\ref{O.2.2+N.23}. It is devoted to a discussion of the technical notion of a {\em tame submanifold} which we introduce in order to define representative charts in the subsequent section. All the reader needs to know about tame submanifolds is that these are particular submanifolds of Lie groupoids with the property that whenever a Lie groupoid homomorphism establishes a bijective correspondence between two of them, the induced bijection is actually a diffeomorphism and that Morita equivalences preserve tame submanifolds.

The fact that \tannakian{\fifu} is a $\C^\infty$\nobreakdash-structured groupoid for every classical \fifu\ poses the question of whether \tannakian{\fifu} is actually a Lie groupoid. In Section~\ref{O.2.2+N.23} we start tackling this issue by providing a necessary and sufficient criterion, which proves to be convenient enough to use in practice, for the answer to the latter question being positive for a given \fifu. This criterion is expressed in terms of the notion of a {\em representative chart,} that is a pair $(\Omega,R)$ consisting of an open subset $\Omega$ of \tannakian{\fifu} and an object $R$ of the domain category \Kt\ of \fifu\ such that the evaluation representation at $R$ induces a homeomorphism between $\Omega$ and a tame submanifold of the linear groupoid $\GL({\fifu R})$; then \tannakian{\fifu} is a Lie groupoid if, and only if, representative charts cover \tannakian{\fifu} and $(\Omega,{R\oplus S})$ is a representative chart for every representative chart $(\Omega,R)$ and for every object $S$ of \Kt.

Section~\ref{O.2.3} introduces a notion of morphism for (classical) fibre functors. Roughly speaking, a morphism of \fifu[\omega] into \fifu[\omega'], over a smooth mapping $f: M \to M'$ of the base manifolds, is a tensor functor of \Kt[C'] into \Kt[C] compatible with the pullback of vector bundles along $f$; every morphism $\fifu[\omega] \to \fifu[\omega']$ over $f$ induces a homomorphism of $\C^\infty$\nobreakdash-structured groupoids $\tannakian{\fifu[\omega]} \to \tannakian{\fifu[\omega']}$ over $f$.

Section~\ref{O.2.4} is devoted to the study of {\em weak equivalences} of (classical) fibre functors: we define them as those morphisms over a surjective submersion which have the property of being a categorical equivalence. As an application of the criterion of \refsez{O.2.2+N.23}, we show that if \fifu[\omega] is weakly equivalent to \fifu[\omega'], then \tannakian{\fifu[\omega]} is a Lie groupoid if and only if \tannakian{\fifu[\omega']} is; when this is the case, the Lie groupoids \tannakian{\fifu[\omega]} and \tannakian{\fifu[\omega']} turn out to be Morita equivalent.

\begin{center}
$* \quad * \quad *$
\end{center}

\noindent In Chapter~\refcpt{6}, we apply the general abstract theory of the preceding chapter to the motivating example provided by the standard forgetful functor on the category of classical representations of a proper Lie groupoid \G. The Tannakian groupoid associated with the latter classical fibre functor will be denoted by \tannakian[\infty]{\G}; in fact, this construction can be extended to a functor $\text- \mapsto \tannakian[\infty]{\text-}$ of the category of Lie groupoids into the category of $\C^\infty$\nobreakdash-structured groupoids so that the envelope homomorphism \envelope{\text-} becomes a natural transformation $(\text-) \to \tannakian[\infty]{\text-}$. We will focus our attention on the following two problems: in the first place, we want to understand whether the Tannakian groupoid \tannakian[\infty]{\G} is a Lie groupoid, let us say for \G\ proper; secondly, we are interested in examples of {\em classically reflexive} Lie groupoids, that is to say Lie groupoids \G\ for which the envelope homomorphism \envelope{} is an isomorphism of topological groupoids between \G\ and \tannakian[\infty]{\G} (recall that, under the assumption of properness, it is sufficient that \envelope{} is injective).

In Section~\ref{O.3.5}, we collect what we know about the first of the two above-mentioned problems in the general case of an arbitrary proper Lie groupoid. Namely, we show that the condition, in the criterion for smoothness of \refsez{O.2.2+N.23}, that $(\Omega,{R\oplus S})$ should be a representative chart for every representative chart $(\Omega,R)$ and object $S$, is always satisfied by the standard forgetful functor on the category of classical representations of a proper Lie groupoid \G\ so that \tannakian[\infty]{\G} is a (proper) Lie groupoid if and only if one can find enough representative charts; if this is the case, then the envelope map \envelope{} is a full submersion of Lie groupoids whose associated pullback functor $\envelope{}^*$ establishes an isomorphism of the corresponding categories of classical representations inverse to the evaluation functor of \refsez{N.21}.

Section~\ref{O.3.6} prosecutes the study initiated in the previous section by providing a proof of the fact that \tannakian[\infty]{\G} is a Lie groupoid for every proper regular groupoid \G. We conjecture that the same statement holds true for every proper \G, that is even without the regularity assumption.

Section~\ref{O.4.3} contains a list of examples of classically reflexive (proper) Lie groupoids; since, as \refsez{O.3.1+.4.1} exemplifies, most Lie groupoids fail to be classically reflexive, this list cannot be very long. To begin with, translation groupoids associated with compact Lie group actions are evidently classically reflexive. Next, we observe that any \'etale Lie groupoid whose source map is proper is necessarily classically reflexive because, for such groupoids, one can make sense of the regular representation. Finally, orbifold groupoids\inciso{by which we mean proper effective groupoids}are classically reflexive because the standard action on the tangent bundle of the base manifold yields a globally faithful classical representation.

\sottosezione{Some possible applications}\markright{Some possible applications}\addcontentsline{toc}{section}{Some possible applications}

The study of classical fibre functors in Chapter~\refcpt{5} was originally motivated by the example treated in Chapter~\refcpt{6}, namely the standard forgetful functor associated with the category of classical representations of a Lie groupoid. However, examples of classical fibre functors can also be found by looking into different directions.

To begin with, one could consider representations of {\em Lie algebroids} \cite{Moerdijk&Mrcun'03,Crainic'03,ELW98}. Recall that a representation of a Lie algebroid $\mathfrak g$ over a manifold $M$ is a pair $(E,\nabla)$ consisting of a vector bundle $E$ over $M$ and a flat $\mathfrak g$\nobreakdash-connection $\nabla$ on $E$, that is, a bilinear map ${\Gamma(\mathfrak g) \times \Gamma(E)} \to \Gamma(E)$ (global sections), $\C^\infty(M)$\nobreakdash-linear in the first argument, \mbox{Leibnitz} in the second and with vanishing curvature. Such representations naturally form a tensor category.

Another example of the same sort is provided by the {\em singular foliations} introduced by \mbox{I.~Androulidakis} and \mbox{G.~Skandalis} \cite{AS06}. Here one is given a locally finite sheaf $\mathscr F$ of modules of vector fields over a manifold $M$, closed under the Lie bracket; this is to be thought of as inducing a `singular' foliation of $M$, in that $\mathscr F$ is no longer necessarily locally free and so the dimension of the leaves may jump. Again, one can consider pairs $(E,\nabla)$ formed by a vector bundle $E$ over $M$ and a morphism of sheaves $\nabla: {\mathscr F\otimes \sections E} \to \sections E$ enjoying formal properties analogous to those defining a flat connection.

In his paper about the local linearizability of proper Lie groupoids \cite{Zu06}, N.T.~Zung poses the question of whether a space, which is locally isomorphic to the orbit space of a compact Lie group action, is necessarily the orbit space $M/\G$ associated with a proper Lie groupoid \G\ over a manifold $M$. Of course, this question is not stated very precisely; its rigorous formulation, as far as we can see, should be given in the following terms. Let us call a $\C^\infty$\nobreakdash-structured space $(X,\mathscr F^\infty)$ a {\em generalized orbifold} if the space $X$ is Hausdorff, paracompact and locally isomorphic, as a functionally structured space, to an orbit space associated with some linear compact Lie group action\nobreakdash---in other words, locally isomorphic to a space of the form $(V/G,\smooth[V/G])$ for some representation $G \to \GL(V)$ of a compact Lie group $G$ on a finite dimensional vector space $V$. The theory of functionally structured spaces suggests the right notion of smooth map of generalized orbifolds and hence the right notion of isomorphism. Zung's theorem implies that the orbit space $(M/\G,\smooth[M/\G])$ of a proper Lie groupoid \G\ over a manifold $M$ is a generalized orbifold: then the question is whether an arbitrary generalized orbifold is actually of this precise form.

\sloppy
Classical fibre functors make their natural appearance in connection with any given generalized orbifold $X$. (Conventionally, we will refer to the $\C^\infty$\nobreakdash-structure of $X$, when necessary, by means of the notation \smooth[X].) Let \V[\infty]X denote the category of locally free sheaves of \smooth[X]\nobreakdash-modules (of locally finite rank), endowed with the standard linear tensor structure; one may refer to the objects of this category as {\em vector bundles over $X$.} Choose a locally finite cover $\{U_i\}$ of $X$ by open subsets $U_i$ such that for each $i$ there is an isomorphism $V_i/G_i \iso U_i$; we regard the maps $\phi_i: V_i \to U_i$ as fixed once and for all, and we assume, for simplicity, that the $V_i$ all have the same dimension. Letting $M$ be the disjoint union ${\coprod V_i}$, one has an obvious classical fibre functor $\fifu^X_M = \fifu^X_{\{V_i,\phi_i\}}$ over $M$ sending each object $\mathscr E$ of the category \V[\infty]X to the smooth vector bundle ${\oplus_i {\phi_i}^*\mathscr E}$ over $M$.

\fussy
The Tannakian groupoid $\tannakian[\infty]{X} = \tannakian{\fifu^X_M}$ is a $\C^\infty$\nobreakdash-structured groupoid with the property that the obvious map $\phi: M \to X$ induces an isomorphism of functionally structured spaces between $M/\tannakian[\infty]{X}$ and $X$; thus, the study of this groupoid might be relevant to the above-mentioned problem. Similarly, the study of the Tannakian groupoids associated with the other examples might lead to interesting information about the underlying geometrical objects, at least when the situation involves some kind of properness. In this connection, it is natural to hope for a general result relating the domain category of a classical fibre functor with the category of $\C^\infty$\nobreakdash-representations of the corresponding Tannakian groupoid, for example via the standard evaluation functor described in \refsez{N.21}.

\separazione

\noindent A well-known conjecture, which has been raising some interest recently \cite{HeMe04,Kal07}, states that every proper \'etale Lie groupoid is Morita equivalent to the translation groupoid associated with some compact Lie group action or, equivalently, that every such groupoid admits a globally faithful classical representation (cfr.\ Ch.~\refcpt{1}, \refsez{N.5}). This conjecture is related to the question of whether proper \'etale Lie groupoids are classically reflexive (we have already observed that the answer is affirmative in the effective case, see Ch.~\refcpt{6}, \refsez{O.4.3}). It is known that for each groupoid \G\ of this kind, there exist a proper effective Lie groupoid $\tilde\G$ and a submersive epimorphism $\G \to \tilde\G$; the kernel of this homomorphism is necessarily a bundle of finite groups $\mathcal B$ embedded into \G, hence, one gets an exact sequence of Lie groupoids $1 \to \mathcal B \into \G \to \tilde\G \to 1$ where $\mathcal B$ and $\tilde\G$ are both classically reflexive. These considerations strongly suggest that one should investigate how the property of reflexivity behaves with respect to Lie groupoid extensions.

\pagestyle{headings}\hbadness=1000\vbadness=1000

\capitolo[Lie Groupoids, Classical Representations]{Lie Groupoids and their Classical Representations}\label{1}

The present chapter is essentially introductory: we regard all the material thereof as well-known. Our purpose is, first of all, to fix some notational conventions and some standard terminology concerning Lie groupoids; this is done in \refsez{O.1.5}. Next, in \refsez{O.3.1+.4.1}, we provide a detailed discussion of a concrete example which is to serve as motivation for the approach we will adopt in Chapters \refcpt[3]{4}. In \refsez[O.3.2]{N.4} we treat the two fundamental pillars on to which our main result holds: Haar systems and \mbox{Zung's} linearizability theorem; we decided to include a presentation of these topics here because we found it difficult to provide adequate references for them. The chapter ends with a digression on the problem of representing a proper Lie groupoid as a global quotient arising from a smooth compact Lie group action.

\sezione{Generalities about Lie Groupoids}\label{O.1.5}
The term \index{groupoid|emph}{groupoid} refers to a small category where every arrow is invertible. A \index{Lie groupoid|see{groupoid}}\textit{Lie groupoid} can be approximately described as an internal groupoid in the category of smooth manifolds. To construct a Lie groupoid \G\ one has to give a pair of manifolds of class $\C^\infty$ \mca[0]\G\ and \mca[1]\G, respectively called \index{manifold of arrows, objects}\index{objects, manifold of -}\textit{manifold of objects} and \index{arrows, manifold of -}\textit{manifold of arrows,} and a list of smooth maps called \index{structure maps}\textit{structure maps.} The basic items in this list are the \index{source}\textit{source} map $\s: \mca[1]\G \to \mca[0]\G$ and the \index{target}\textit{target} map $\t: \mca[1]\G \to \mca[0]\G$; these have to meet the requirement that the fibred product $\mca[2]\G = \mca{\mca[1]\G}$ exists in the category of $\C^\infty$-manifolds. Then one has to give a \index{composition}\textit{composition} map $\c: \mca[2]\G \to \mca[1]\G$, a \index{unit map, section}\textit{unit} map $\u: \mca[0]\G \to \mca[1]\G$ and an \index{inverse}\textit{inverse} map $\imap: \mca[1]\G \to \mca[1]\G$, for which the familiar algebraic laws must be satisfied.

Terminology and Notation: The points $x = \s(g)$ and $x' = \t(g)$ are resp.\ called the \textit{source} and the \textit{target} of the \textit{arrow} $g$. We let $\G(x,x')$ denote the set of all the arrows whose source is $x$ and whose target is $x'$; we shall use the abbreviation $\G|_x$ for the \index{isotropy group}\textit{isotropy} or \index{vertex group|see{isotropy group}}\textit{vertex} group $\G(x,x)$. Notationally, we will often identify a point $x \in \mca[0]\G$ and the corresponding unit arrow $\u(x) \in \mca[1]\G$. It is costumary to write ${g'\cdot g}$ or $g'g$ for the composition $\c(g',g)$ and $g^{-1}$ for the inverse $\imap(g)$.

Our description of the notion of Lie groupoid is still incomplete. It turns out that a couple of additional requirements are needed in order to get a reasonable definition.

Recall that a manifold $M$ is said to be \index{paracompact}\textit{paracompact} if it is Hausdorff and there exists an ascending sequence of open subsets with compact closure $\cdots \subset U_i \subset \overline U_i \subset U_{i+1} \subset \cdots$ such that $M = \txtcup{i=0}\infty{U_i}$. A Hausdorff manifold is paracompact if and only if it possesses a countable basis of open subsets. Any open cover of a paracompact manifold admits a locally finite refinement. Any paracompact manifold admits partitions of unity of class $\C^\infty$ (subordinated to an open cover; cf.\ for instance \mbox{Lang} \cite{La01}).

In order to make the fibred product \mca{\mca[1]\G} meaningful as a manifold and for other purposes related to our studies, we shall include the following additional conditions in the definition of Lie groupoid:
\begin{elenco}
\item[\textsl{1.}]The source map $\s: \mca[1]\G \to \mca[0]\G$ is a submersion with Hausdorff fibres;
\item[\textsl{2.}]The manifold \mca[0]\G\ is paracompact.
\end{elenco}
Note that we do not require that the manifold of arrows \mca[1]\G\ is Hausdorff or paracompact; actually, this manifold is neither Hausdorff nor second countable in many examples of interest. The definition here differs from that in \mbox{Moerdijk} and \mbox{Mr\v cun} \cite{Moerdijk&Mrcun'03} in that we additionally require that the manifold \mca[0]\G\ is paracompact. The first condition implies at once that the domain of the composition map is a submanifold of the Cartesian product ${\mca[1]\G \times \mca[1]\G}$ and that the target map is a submersion with Hausdorff fibres; thus, the source fibres $\G(x,\text-) = \s^{-1}(x)$ and the target fibres $\G(\text-,x') = \t^{-1}(x')$ are closed Hausdorff submanifolds of \mca[1]\G. A Lie groupoid \G\ is said to be \index{Hausdorff groupoid}\index{groupoid!Hausdorff}\textit{Hausdorff} if the manifold \mca[1]\G\ is Hausdorff.

Some more Terminology: The manifold \mca[0]\G\ is usually called the \index{base}\textit{base} of the groupoid \G; one also says that \G\ is a groupoid \index{groupoid!over a manifold|see{base}}over the manifold \mca[0]\G. We shall often use the notation $\G^x = \G(x,\text-) = \s^{-1}(x)$ for the fibre of the source map over a point $x \in \mca[0]\G$. More generally, we shall write
\begin{equazione}\label{i1}
\G(S,S') = \bigl\{g \in \mca[1]\G: \s(g) \in S\text{~\&~} \t(g) \in S'\bigr\}\text, \quad \G|_S = \G(S,S)
\end{equazione}
and $\G^S = \G(S,\text-) = \G(S,\mca[0]\G) = \s^{-1}(S)$ for all subsets $S, S' \subset \mca[0]\G$.

A \index{homomorphism of groupoids|emph}\textit{homomorphism of Lie groupoids} is a smooth functor. More precisely, a homomorphism $\varphi: \G \to \H$ consists of two smooth maps $\mca[0]\varphi: \mca[0]\G \to \mca[0]\H$ and $\mca[1]\varphi: \mca[1]\G \to \mca[1]\H$, compatible with the groupoid structure in the sense that ${\s\circ \mca[1]\varphi} = {\mca[0]\varphi \circ \s}$, ${\t\circ \mca[1]\varphi} = {\mca[0]\varphi \circ \t}$ and $\mca[1]\varphi({g'\cdot g}) = {\mca[1]\varphi(g') \cdot \mca[1]\varphi(g)}$. Lie groupoids and their homomorphisms form a category.

There is also a notion of \index{topological groupoid}\textit{topological groupoid:} this is just an internal groupoid in the category of topological spaces and continuous mappings. In the continuous case the definition is much simpler and one need not worry about the domain of definition of the composition map. With the obvious notion of homomorphism, topological groupoids constitute a category.

\begin{esempio}\label{i2}
Every smooth manifold $M$ can be regarded as a Lie groupoid by taking $M$ itself as the manifold of arrows and the identity map $\id: M \to M$ as the unit section. Alternatively, one can form the \textit{pair groupoid over $M$}; this is the Lie groupoid whose manifold of arrows is ${M\times M}$ and whose source and target map are the two projections.
\end{esempio}

\begin{esempio}\label{i3}
Any Lie group $G$ can be regarded as a Lie groupoid over the one-point manifold by taking $G$ itself as the manifold of arrows.
\end{esempio}

\begin{esempio}[linear groupoids]\label{i4}
If $E$ is a real or complex smooth vector bundle (of locally finite rank) over a manifold $M$, one can form the \index{linear groupoid@linear groupoid \ensuremath{\GL(E)}|emph}\index{GL(E)@\ensuremath{\GL(E)} (linear groupoid)|emph}\textit{linear groupoid $\GL(E)$ associated with $E$}. This is defined as the groupoid over $M$ whose arrows $x \to x'$ are the linear isomorphisms $E_x \isoto E_{x'}$ between the fibres of $E$ over the points $x$ and $x'$. There is an obvious smooth structure turning $\GL(E)$ into a Lie groupoid.
\end{esempio}

\begin{esempio}[action groupoids]\label{i5}
Let $G$ be a Lie group acting smoothly (from the left) on a manifold $M$. Then one can define the \index{action groupoid@action groupoid \ensuremath{G\ltimes M}|emph}\index{G ltimes M@\ensuremath{G\ltimes M} (action groupoid)|emph}\textit{action} (or \index{translation groupoid|see{action groupoid}}\textit{translation}) \textit{groupoid ${G\ltimes M}$} as the Lie groupoid over $M$ whose manifold of arrows is the Cartesian product ${G\times M}$, whose source and target map are respectively the projection onto the second factor $(g,x) \mapsto x$ and the action $(g,x) \mapsto {g x}$ and whose composition law is the operation
\begin{equazione}\label{i6}
{(g',x')(g,x)} = (g'g,x)\text.
\end{equazione}
There is a similar construction ${M\rtimes G}$ associated with right actions.
\end{esempio}

Let \G\ be a Lie groupoid and let $x$ be a point of its base manifold \mca[0]\G. The \index{orbit}\textit{orbit of \G} (or \textit{\G-orbit}) \textit{through $x$} is the subset
\begin{equazione}\label{i7}
{\G x} \bydef {\G\cdot x} \bydef \t\bigl(\G^x\bigr) = \{x' \in \mca[0]\G| \exists g: x \to x'\}\text.
\end{equazione}
Note that the isotropy group $\G|_x$ acts from the the right on the manifold $\G^x$; this action is clearly free and transitive along the fibres of the restriction of the target map \t\ to $\G^x$. The following result holds (see \cite{Moerdijk&Mrcun'03} p.~115):
\begin{theorem}\label{i8}
Let \G\ be a Lie groupoid and let $x, x' \in \mca[0]\G$. Then
\begin{elenco}
\item[1.]$\G(x,x')$ is a closed submanifold of \mca[1]\G;
\item[2.]$\G|_x$ is a Lie group;
\item[3.]the \G-orbit through $x$ is an immersed submanifold of \mca[0]\G;
\item[4.]the target map $\t: \G^x \to {\G x}$ proves to be a principal $\G|_x$-bundle.
\end{elenco}
\end{theorem}
It is worthwhile spending a couple of words about the manifold structure that is asserted to exist on the \G-orbit through $x$. The set ${\G x}$ can obviously be identified with the homogeneous space ${\G^x/(\G|_x)}$. Now, it can be proved that there exists a (possibly non-Hausdorff) manifold structure on this quotient space, such that the quotient map turns out to be a principal bundle.

\separazione

\noindent We say that a Lie (or topological) groupoid \G\ is \index{proper groupoid}\index{groupoid!proper}\textit{proper} if \G\ is Hausdorff and the combined source--target map $(\s,\t): \mca[1]\G \to {\mca[0]\G \times \mca[0]\G}$ is proper (in the familiar sense: the inverse image of a compact subset is compact).

The manifold of arrows \mca[1]{\G} of a proper Lie groupoid \G\ is always paracompact. Indeed, by the definition of Lie groupoid, the base $M$ of \G\ is a paracompact manifold and therefore there exists an invading sequence $\cdots \subset U_i \subset \overline U_i \subset U_{i+1} \subset \cdots \subset M$ of pre-compact open subsets; the inverse images $\Gamma_i = \G|_{U_i} = (\s,\t)^{-1}({U_i\times U_i})$ form an analogous sequence inside the (Hausdorff) manifold \mca[1]\G.

Let $x_0$ be a point of $M$. We know the orbit $S = {\G x_0}$ is an immersed submanifold of $M$ (precisely, there exists a unique manifold structure on $S$ such that $\t: \G^{x_0} \to S$ is a principal right $\G|_{x_0}$-bundle and the inclusion $S \into M$ an immersion). Now, \textsl{it follows from the properness of \G\ that $S$ is actually a submanifold of $M$}. To see this, fix a point $s_0 \in S$. Since there exists a local equivariant chart $\G(x_0,W) \iso {W\times \G|_{x_0}}$ where $W$ is both an open neighborhood of $s_0$ in $S$ and a submanifold of $M$, it will be enough to prove the existence of an open ball $B \subset M$ at $s_0$ such that ${S\cap B} \subset W$. To do this, take a sequence of open balls $B_i$ shrinking to $s_0$: the decreasing sequence $\Sigma_i = {\G(x_0,\overline B_i) - \G(x_0,W)}$ of closed subsets of the manifold $\G(x_0,\text-)$ is contained in the compact subset $\G(x_0,\overline B_1)$ and therefore, since ${\bigcap \Sigma_i} = \varnothing$, there exists some $i$ such that $\G(x_0,B_i) \subset \G(x_0,W)$.

\sezione{Classical Representations}\label{O.3.1+.4.1}
In this section we introduce the costumary notion of representation of a Lie groupoid on a smooth vector bundle and we explain, by means of a counterexample, why this notion is inadequate for the purpose of building a possible Tannaka duality theory for proper Lie groupoids.

Let \G\ be a Lie groupoid and let $M$ be its base. We let \index{R(T;k)@\R[\infty]{\mathcal T;k} (category of smooth representations on vector bundles)}\R[\infty]{\G;\nC} denote the category of all \nC-linear \index{C infinity representation@\ensuremath{C^\infty}-representation}\index{representation!classical|see{\ensuremath{C^\infty}- or smooth}}\index{representation!C infinity or smooth@\ensuremath{C^\infty}- or smooth}\index{classical representation|see{\ensuremath{C^\infty}-representation}}\textit{classical representations} of \G. The objects of this category are the pairs $(E,\varrho)$ consisting of a smooth complex vector bundle $E$ (of locally finite rank) over $M$ and a Lie groupoid homomorphism
\begin{equazione}\label{O.esp12}
\begin{split}
\xymatrix@C=30pt@R=25pt{\G\ar[d]_{(\s,\t)}\ar[r]^-\varrho & \GL(E)\ar[d]^{(\s,\t)} \\ {M\times M}\ar[r]^{\id\times\id} & {M\times M}\text;\!\!}
\end{split}
\end{equazione}
the arrows, let us say those $a: (E,\varrho) \to (F,\varsigma)$, are the morphisms of vector bundles $a: E \to F$ such that the square
\begin{equazione}\label{O.esp13}
\begin{split}
\xymatrix@C=35pt{E_x\ar[d]_{a_x}\ar[r]^-{\varrho(g)} & E_{x'}\ar[d]^{a_{x'}} \\ F_x\ar[r]^-{\varsigma(g)} & F_{x'}}
\end{split}
\end{equazione}
commutes for all $x, x' \in M$ and $g \in \G(x,x')$. There is an entirely analogous notion of \nR-linear classical representation of \G, where real vector bundles are used instead of complex ones. One obtains a corresponding category \R[\infty]{\G;\nR}. Insofar as a particular choice of coefficients is not relevant to the subject matter of a discussion, we shall write simply \R[\infty]\G\ and suppress any further reference to coefficients.

\separazione

\noindent Lie groupoids cannot always be distinguished from one another just on the basis of knowledge of the respective categories of classical representations; this consideration motivates our approach to Tannaka duality as described in Chapter \refcpt{4}. We are going to substantiate our assertion by means of a counterexample which we discovered independently in 2005: only recently \mbox{A.~Henriques} pointed out to us that the same counterexample was already known in the context of orbispace theory, see \mbox{L\"uck} and \mbox{Oliver} (2001) \cite{Lueck&Oliver'01}.

Recall that a \index{Lie bundle}\textit{Lie bundle} (also known as \index{bundle of Lie groups|see{Lie bundle}}\textit{bundle of Lie groups}) is a Lie groupoid whose source and target map coincide.

Fix a Lie group $H$ and choose an automorphism $\chi \in \Aut(H)$. There is a general procedure\inciso{completely analogous to the construction of \mbox{M\"obius} bands, \mbox{Klein} bottles {\em et similia}}by means of which one can obtain a locally trivial Lie bundle $\G = \G_{H;\chi} \to \mathrm S^1$ with fibre $H$ over the unit circle. Put $\mca[1]\G = {({\nR\times H})/\thicksim}$ where $\thicksim$ is the equivalence relation
\begin{equazione}\label{O.esp32}
(t,h) \thicksim (t',h') \quad \aeq \quad t'-t = \ell \in \nZ \quad \text{and} \quad h' = \chi^\ell(h)\text.
\end{equazione}
The bundle fibration $\mca[1]\G \to \mathrm S^1$ (= source map of \G\ = target map of \G) is defined as the unique map that makes the square
\begin{equazione}\label{O.esp33}
\begin{split}
\xymatrix{{\nR\times H}\ar[r]\ar[d]_{\text{quot.~proj.}} & \nR\ar[d]^{t\mapsto\mathit e^{2\pi it}} \\ \mca[1]\G\ar@{-->}[r] & \mathrm S^1}
\end{split}
\end{equazione}
commute. In terms of representatives of elements of \mca[1]\G, the composition law $\c: {\mca[1]\G \times_{\mathrm S^1} \mca[1]\G} \to \mca[1]\G$ can be defined by setting
\begin{equazione}\label{O.esp34}
{[t',h'] \cdot [t,h]} = [t',{h'\cdot \chi^k(h)}]\text,
\end{equazione}
where $k = t'-t \in \nZ$ and the square bracket notation indicates that we are taking equivalence classes. This operation turns $\G \to \mathrm S^1$ into a bundle of groups over the circle, with fibre $H$.

\sloppy
Consider the open cover of $\mathrm S^1$ determined by the local exponential parametrizations $(0,1) \isoto U$ and $(-\tfrac12,\tfrac12) \isoto V$. One has two corresponding mutually compatible trivializing charts for \mca[1]\G\ over $\mathrm S^1$, namely
\begin{equazione}\label{O.esp35}
\tau_U: \mca[1]\G|_U \isoto {U\times H} \quad \text{and} \quad \tau_V: \mca[1]\G|_V \isoto {V\times H}\text:
\end{equazione}
the former sends $g \in \mca[1]\G|_U$ to the pair $(\mathit e^{2\pi it},h)$ with $[t,h] = g$ and $0<t<1$, the latter sends $g \in \mca[1]\G|_V$ to the pair $(\mathit e^{2\pi it},h)$ with $[t,h] = g$ and $-\tfrac12<t<\tfrac12$. These charts determine the differentiable structure. Notice, by the way, that the transition map between them, namely
\begin{equazione}\label{O.esp36}
{\tau_U \circ {\tau_V}^{-1}}: {({U\cap V}) \times H} \isoto {({U\cap V}) \times H}\text,
\end{equazione}
is given by the identity over ${W\times H}$ and by $(w',h) \mapsto (w',\chi(h))$ over ${W'\times H}$, if one lets $(0,\tfrac12) \isoto W$ and $(\tfrac12,1) \isoto W'$ denote the two connected components of the intersection ${U\cap V}$.

\fussy
We start by studying the complex classical representations of the Lie bundle $\G_{H;\chi}$, which are technically easier to handle. The analogous result for real representations will be deduced as a corollary.

% La revisione finisce qui
Fix a classical representation $(E,\varrho) \in {\Ob\,\R[\infty]{\G;\nC}}$ on a smooth complex vector bundle $E$ of rank $\ell$ over $\mathrm S^1$. Since $U$ and $V$ are contractible open subsets of $\mathrm S^1$, the vector bundle $E$ will be trivial over each of them i.e.\ there will exist smooth vector bundle isomorphisms%
\begin{equazione}\label{O.esp37}%
E|_U \isoto {U\times \nC^\ell} \quad\text{and} \quad E|_V \isoto {V\times \nC^\ell}\text.%
\end{equazione}%
These will form a trivializing atlas for $E$ over $\mathrm S^1$, whose unique transition mapping will be given by, let us say,%
\begin{equazione}\label{O.esp38}%
Q: W \to \GL(\ell;\nC) \quad\text{and} \quad Q': W' \to \GL(\ell;\nC)\text.%
\end{equazione}%
Accordingly, the Lie bundle $\GL(E)$ over $\mathrm S^1$ (that is, by abuse of notation, the restriction of the linear groupoid $\GL(E)$ to the diagonal $\mathrm S^1 \into {\mathrm S^1\times \mathrm S^1}$) will be described by trivializing charts of the following form%
\begin{equazione}\label{O.esp39}%
\GL(E)|_U \isoto {U\times \GL(\ell;\nC)} \quad\text{and} \quad \GL(E)|_V \isoto {V\times \GL(\ell;\nC)}\text,%
\end{equazione}%
whose transition map ${({U\cap V})\times \GL(\ell;\nC)} \isoto {({U\cap V})\times \GL(\ell;\nC)}$ will send $w \in W$ to $A \mapsto {Q(w) A Q(w)^{-1}}$ and $w' \in W'$ to $A \mapsto {Q'(w') A Q'(w')^{-1}}$.%

In this situation one can write down corresponding local expressions for $\varrho$, namely $\varrho_U(u,h) = \bigl(u,A_U(u,h)\bigr)$ over $U$ and $\varrho_V(v,h) = \bigl(v,A_V(v,h)\bigr)$ over $V$ with $A_U: {U\times H} \to \GL(\ell;\nC)$ a smooth family of representations of $H$ etc., which make the following squares%
\begin{equazione}\label{O.esp40}%
\begin{split}%
\xymatrix{\mca[1]\G|_U\ar[r]^-{\varrho|_U}\ar[d]^{\tau_U}_\iso & \GL(E)|_U\ar[d]^{\iso_U} & \mca[1]\G|_V \ar[r]^-{\varrho|_V}\ar[d]^-{\tau_V}_\iso & \GL(E)|_V\ar[d]^{\iso_V} \\ {U\times H}\ar@{-->}[r]^-{\varrho_U} & {U\times \GL(\ell;\nC)} & {V\times H}\ar@{-->}[r]^-{\varrho_V} & {V\times \GL(\ell;\nC)}}%
\end{split}%
\end{equazione}%
commute. If we take their restrictions to $W$, $W'$ respectively, we obtain%
\begin{equazione}\label{O.esp41}%
\begin{split}%
\xymatrix{{W\times H}\ar@{-->}[r]^-{\varrho_V} & {W\times \GL(\ell;\nC)} & {W'\times H}\ar@{-->}[r]^-{\varrho_V} & {W'\times \GL(\ell;\nC)} \\ \mca[1]\G|_W\ar[r]^-{\varrho|_W}\ar[d]^{\tau_U}_\iso\ar[u]_{\tau_V}^\iso & \GL(E)|_W\ar[d]^{\iso_U}\ar[u]_{\iso_V} & \mca[1]\G|_{W'}\ar[r]^-{\varrho|_{W'}}\ar[d]^{\tau_U}_\iso\ar[u]_{\tau_V}^\iso & \GL(E)|_{W'}\ar[d]^{\iso_U}\ar[u]_{\iso_V} \\ {W\times H}\ar@{-->}[r]^-{\varrho_U} & {W\times \GL(\ell;\nC)} & {W'\times H}\ar@{-->}[r]^-{\varrho_U} & {W'\times \GL(\ell;\nC)}}%
\end{split}%
\end{equazione}%
and hence, making use of the explicit expression \refequ{O.esp36} for the transition map ${\tau_U\circ {\tau_V}^{-1}}$, we are led to the following relations: for all $h \in H$%
\begin{equazione}\label{O.esp42}%
\left\{\begin{aligned}%
A_U(w,h) &= {Q(w) A_V(w,h) Q(w)^{-1}} & &\phantom= \text{~for all~} w \in W%
\\%
A_U\bigl(w',\chi(h)\bigr) &= {Q'(w') A_V(w',h) Q'(w')^{-1}} & &\phantom= \text{~for all~} w' \in W'\text.%
\end{aligned}\right.%
\end{equazione}%

From now on, we assume that $H$ is compact. We also fix two points $w_0 \in W$ and $w_0' \in W'$. There is a continuous path $\gamma_U: [0,1] \to U$ from $w_0$ to $w_0'$. This gives a continuous map%
\begin{equazione}\label{O.esp43}%
{[0,1]\times H} \xto{\gamma_U\times\id} {U\times H} \xto{A_U} \GL(\ell;\nC)%
\end{equazione}%
which is clearly a homotopy of representations of $H$ connecting $A_U(w_0,\text-)$ to $A_U(w_0',\text-)$. Then, as remarked in Note~\refcnt{ii1}, there will be an invertible matrix $R \in \GL(\ell;\nC)$ such that%
\begin{equazione}\label{O.esp44}%
A_U(w_0,\text-) = {R A_U(w_0',\text-) R^{-1}}\text.%
\end{equazione}%
A second path $\gamma_V: [0,1] \to V$ connecting $w_0$ to $w_0'$ will analogously yield a matrix $S \in \GL(\ell;\nC)$ such that%
\begin{equazione}\label{O.esp45}%
A_V(w_0,\text-) = {S A_V(w_0',\text-) S^{-1}}\text.%
\end{equazione}%
Making the appropriate substitutions in \refequ{O.esp42}, we finally find an invertible matrix $P \in \GL(\ell;\nC)$ such that%
\begin{equazione}\label{O.esp46}%
A_U\bigl(w_0,\chi(h)\bigr) = {P A_U(w_0,h) P^{-1}} \quad\text{for all~} h \in H\text.%
\end{equazione}%

Next, we further specialize down to the case where $H$ is abelian and connected. Motivated by Eq.~\refequ{O.esp46}, we focus our attention on those matrix representations $A: H \to \GL(\ell;\nC)$ such that%
\begin{equazione}\label{O.esp47}%
\exists P \in \GL(\ell;\nC) \quad\text{for which} \quad A(\chi(h)) = {P A(h) P^{-1}}\text.%
\end{equazione}%
By \mbox{Schur's} Lemma, every irreducible matrix representation of an Abelian Lie group must be one-dimensional (cf.\ for instance \mbox{Br\"ocker} and \mbox{tom Dieck} p.~69) and therefore, because of the compactness of $H$, necessarily a {\em character} i.e.\ a Lie group homomorphism of $H$ into the 1\nobreakdash-torus $\mathrm T^1$. Since every representation of a compact Lie group is a direct sum of irreducible ones ({\em ibid.\ }p.~68), it is no loss of generality to assume Eq.~\refequ{O.esp47} to be of the following form%
\begin{equazione}\label{O.esp48}%
\left(\begin{matrix}%
({\alpha_1\circ \chi})(h) & \cdots & 0 \\ \vdots & \ddots & \vdots \\ 0 & \cdots & ({\alpha_\ell\circ \chi})(h)%
\end{matrix}\right) = {P%
\left(\begin{matrix}%
\alpha_1(h) & \cdots & 0 \\ \vdots & \ddots & \vdots \\ 0 & \cdots & \alpha_\ell(h)%
\end{matrix}\right) P^{-1}}\text,%
\end{equazione}%
where $\alpha_1, \ldots, \alpha_\ell: H \to \mathrm T^1$ are characters of $H$.%

The two complex diagonal matrices occurring in Eq.~\refequ{O.esp48} must have the same characteristic polynomial $p(h,X) \in \nC[X]$. Thus, if we put%
\begin{equazione}\label{O.esp50}%
\beta_j = {\alpha_j\circ \chi}: H \to \mathrm T^1 \quad\text{and} \quad F_{ij} = \bigl\{h \in H: \alpha_i(h) = \beta_j(h)\bigr\}\text,%
\end{equazione}%
we can in particular express $H$ as a finite union $F_{11}\cup \cdots\cup F_{1\ell}$ of closed subsets. Now, it follows by a standard inductive argument that one of them, let us say $F_{11}$, must have nonempty interior; therefore, the two characters $\alpha_1$ and $\beta_1$ coincide on all of $H$, because a homomorphism of connected Lie groups is determined by its differential at the neutral element ({\em ibid.\ }p.~24). Cancelling the two corresponding linear factors in $p(h,X)$ we obtain%
\begin{equazione}\label{O.esp51}%
{\bigl(X-\beta_2(h)\bigr) \cdots \bigl(X-\beta_\ell(h)\bigr)} = {\bigl(X-\alpha_2(h)\bigr) \cdots \bigl(X-\alpha_\ell(h)\bigr)}\text.%
\end{equazione}%
Then, arguing by induction on the degree of the polynomial, we conclude that there is a permutation $\sigma$ on $\ell$ letters such that  $\alpha_i = \beta_{\sigma(i)} = {\alpha_{\sigma(i)} \circ \chi}$ for all $i = 1, \ldots, \ell$.%

Now, consider for instance $\alpha_1$. Write $\sigma$ as a product of disjoint cycles and consider the cycle $\bigl(1,\sigma(1),\ldots,\sigma^r(1)\bigr)$ where $r \geqq 0$ and $\sigma^{r+1}(1) = 1$. Then we have $\alpha_1 = {\alpha_{\sigma(1)} \circ \chi} = {\bigl({\alpha_{\sigma(\sigma(1))} \circ \chi}\bigr) \circ \chi} = {\alpha_{\sigma^2(1)} \circ \chi^2} = \cdots = {\alpha_{\sigma^r(1)} \circ \chi^r} = {\bigl({\alpha_{\sigma(\sigma^r(1))} \circ \chi}\bigr) \circ \chi^r} = {\alpha_{\sigma^{r+1}(1)} \circ \chi^{r+1}} = {\alpha_1\circ \chi^{r+1}}$. Therefore $\alpha_1$ is an example of a character $\alpha: H \to \mathrm T^1$ with the special property%
\begin{equazione}\label{O.esp52}%
\exists r \geqq 0 \quad\text{such that} \quad \alpha = {\alpha\circ \chi^{r+1}}\text.%
\end{equazione}%

Finally, let us take $H = \mathrm T^2 = {\mathrm T^1\times \mathrm T^1}$ to be the 2\nobreakdash-torus. Fix an arbitrary $\ell \in \nZ$, and consider the map%
\begin{equazione}\label{O.esp53}%
\chi_\ell: \mathrm T^2 \to \mathrm T^2 \quad \text{defined by the rule} \quad (s,t) \mapsto (s,{s^\ell t})\text.%
\end{equazione}%
This is an automorphism of the Lie group $\mathrm T^2$, with inverse $\chi_{-\ell}$.%

Any 2\nobreakdash-character $\alpha: \mathrm T^2 \to \mathrm T^1$ can be written as the product $\alpha(s,t) = {\mu(s)\nu(t)}$ of the two 1\nobreakdash-characters $\mu(s) = \alpha(s,1)$ and $\nu(t) = \alpha(1,t)$. If we assume that $\alpha$ enjoys the property \refequ{O.esp52} then we get ${\mu(s)\nu(t)} = \alpha(s,t) = \alpha\bigl(s,s^{\ell(r+1)}t\bigr) = {\mu(s) \nu(s)^{\ell(r+1)} \nu(t)}$ and therefore $\nu(s)^{\ell(r+1)} = 1$ for all $s \in \mathrm T^1$. Now, if $\ell\neq0$ then $\nu$ must be trivial, because $r+1>0$. It follows that%
\begin{equazione}\label{O.esp55}%
\alpha(s,t) = \mu(s)%
\end{equazione}%
does not depend on $t$.%

\begin{proposizione}\label{O.cor1}
Fix any integer $0 \neq \ell \in \nZ$ and let $\G_{\mathrm T^2;\chi_\ell} \to \mathrm S^1$ be the locally trivial Lie bundle with fibre $\mathrm T^2$ over the circle, constructed as explained above by using $\chi_\ell \in \Aut(\mathrm T^2)$ as twisting automorphism.

Then there exists an embedding of Lie bundles over the circle
\begin{equazione}\label{O.esp56}
\begin{split}
\xymatrix@C=35pt{{\mathrm S^1 \times \mathrm T^1}\ar@{^{(}->}[r]^-\varphi\ar[d] & \G_{\mathrm T^2;\chi_\ell}\ar[d] \\ {\mathrm S^1 \times \mathrm S^1}\ar[r]^{\id\times\id} & {\mathrm S^1 \times \mathrm S^1}}
\end{split}
\end{equazione}
with the property that every classical representation $(E,\varrho)$ in $\R[\infty]{\G_{\mathrm T^2;\chi_\ell}}$ pulls back to a trivial representation $(E,{\varrho\circ\varphi})$ of ${\mathrm S^1 \times \mathrm T^1}$.
\end{proposizione}
\begin{proof}
Define the embedding $\varphi$ as follows. Given $(x,z) \in {\mathrm S^1 \times \mathrm T^1}$, send it to the equivalence class $[t,(1,z)]$, no matter what $t$ you choose as long as $\mathit e^{2\pi it} = x$. With respect to either of the two charts $\tau_U$ and $\tau_V$ of Eq.~\refequ{O.esp35}, the local expression of this embedding is simply $(x,z) \mapsto (x;1,z)$.

Now, let $(E,\varrho)$ be a \nC\nobreakdash-linear representation of \G\ and let $w_0 \in W$ be the point we selected in the course of the discussion above. In the chart $\tau_U$ the isotropy group $\G|_{w_0}$ and the torus $\mathrm T^2$ are identified by the induced Lie group isomorphism $\G|_{w_0} \iso \mathrm T^2$. The subgroups $\varphi(\{w_0\}\times \mathrm T^1) \subset \G|_{w_0}$ and ${\{1\}\times \mathrm T^1} \subset \mathrm T^2$ correspond to one another under this isomorphism; moreover, the homomorphism $\varrho_{w_0}: \G|_{w_0} \to \GL(E_{w_0})$ is given the matrix representation $A = A_U(w_0,\text-): \mathrm T^2 \to \GL(\ell;\nC)$ of Eq.~\refequ{O.esp47}. Therefore, since from Eq.~\refequ{O.esp55} we know that ${\{1\}\times \mathrm T^1}$ is contained in $\kernel A$, we conclude that the image $\varphi(\{w_0\}\times \mathrm T^1)$ is contained in $\kernel{\varrho_{w_0}}$. By the standard homotopy argument of Note~\refcnt{ii1} we finally get $\varphi(\{x\}\times \mathrm T^1) \subset \kernel{\varrho_x}$ for all $x \in \mathrm S^1$. This completes the proof in the \nC\nobreakdash-linear case.%

Finally, let $R = (E,\varrho)$ be any \nR-linear classical representation of \G. It will be enough to take the complexification ${R\otimes\nC} = ({E\otimes\nC},{\varrho\otimes\nC})$ and observe that $\kernel{\varrho_x} = \kernel{\varrho_x\otimes\nC} = \kernel{({\varrho\otimes\nC})_x}$ for all $x$.
\end{proof}

Consider the map ${\nR \times \mathrm T^2} \to {\mathrm S^1 \times \mathrm T^1}$ given by $(t;z,z') \mapsto (\mathit e^{2\pi it},z)$. This induces an epimorphism of Lie bundles over $\mathrm S^1$
\begin{equazione+}\label{O.esp57}
\psi: \G_{\mathrm T^2;\chi_\ell} \longto \underline{\mathrm T}^1 & \bigl(\underline{\mathrm T}^1 \bydef {\mathrm S^1 \times \mathrm T^1}\bigr)
\end{equazione+}
whose kernel is precisely the image of the embedding $\varphi$ of the preceding proposition. The latter map yields an identification of forgetful functors
\begin{equazione}\label{O.equ10}
\begin{split}
\xymatrix{\R[\infty]{\underline{\mathrm T}^{\mathrm1}}\ar[d]_{\text{forg.~func.}}\ar[r]^-{\psi^*}_-\simeq & \R[\infty]{\G_{\mathrm T^2;\chi_\ell}}\ar[d]^{\text{forg.~func.}} \\ \V[\infty]{\mathrm S^1}\ar@{=}[r] & \V[\infty]{\mathrm S^1}}
\end{split}
\end{equazione}
defined as $\psi^*(E,\varrho) \bydef (E,{\psi\circ\varrho})$. One easily recognizes that the functor $\psi^*$ is an isomorphism of categories. Indeed, its inverse $\psi_*$ can be constructed explicitly by means of the familiar universal property of the quotient (which in the present case follows immediately from Proposition \refcnt{O.cor1}), namely
\begin{equazione}\label{O.esp58}
\begin{split}
\xymatrix@C=30pt@R=25pt{\G_{\mathrm T^2;\chi_\ell}\ar[d]_\psi\ar[r]^-\varrho & \GL(E) \\ \underline{\mathrm T}^1\ar@{-->}[ur]_{\psi_*\varrho} &}
\end{split}
\end{equazione}
for every $(E,\varrho) \in {\Ob\,\R[\infty]{\G_{\mathrm T^2;\chi_\ell}}}$, so that $(E,{\psi_*\varrho})$ is an object of \R[\infty]{\underline{\mathrm T}^{\mathrm1}} (one obviously sets $\psi_*(a) = a$ for all morphisms $a$).

The existence of the identification of categories \refequ{O.equ10} shows in a very convincing way that, in general, a category of classical representations does not provide enough information to recover the Lie groupoid from which it originates; this is true independently of the recipe one might invent for a possible reconstruction theory. Note also that this failure already occurs under circumstances where the Lie groupoid is a very reasonable one (compact, abelian and connected). Of course, what we are saying does not exclude the possibility that in some special cases the reconstruction may be feasible; we shall give a few elementary examples of this sort later on in \refsez{O.4.3}.%

\begin{nota}\label{ii1}
(Compare also \mbox{Br\"ocker} and \mbox{tom Dieck} \cite{Broecker&tomDieck'85} p.~84) Let $G$ be a Lie group and let $\varrho_t: G \to \GL(V)$ be a family of representations $\varrho_t$ depending continuously on $g \in G$ and $t \in [0,1]$, in other words, a homotopy of representations. We claim that \textsl{if $G$ is compact, the representations $\varrho_0$ and $\varrho_1$ are isomorphic}\nobreakdash---i.e.\ there exists some $A \in \GL(V)$ which conjugates $\varrho_0$ into $\varrho_1$. To begin with, let $\dual G$ denote the set of isomorphism classes of irreducible $G$-modules. For each $\gamma \in \dual G$, select a representative $V_\gamma$. Then for every $t \in [0,1]$ one can decompose the $G$-module $V_t = (V,\varrho_t)$ into a direct sum $V_t \iso {\underset{\gamma\in\dual G}\oplus {n_\gamma^t V_\gamma}}$ in which the integer $n_\gamma^t$ = multiplicity of $V_\gamma$ in $V_t$ = ${\int {\chi_t \overline{\chi_\gamma}}}$, where $\chi_\gamma$ is the character of $V_\gamma$ and $\chi_t = \txtsum{\gamma\in\dual G}{}{n_\gamma^t \chi_\gamma}$ is the character of $V_t$, depends continuously on $t$ and is therefore constant. This proves the claim.

More generally, one has that if $f_t: G \to H$ is any homotopy of homomorphisms of a {\em compact} Lie group $G$ into a Lie group $H$ then $f_0$ and $f_1$ are conjugate: see \mbox{Conner} and \mbox{Floyd} (1964) \cite{Conner&Floyd'64} Lemma~38.1. Their result is a consequence of the following theorem of \mbox{Montgomery} and \mbox{Zippin} (1955) (which can be found in \cite{MoZi55} p.~216):
\begin{theorem*}
Let $G$ be a Lie group and $F$ a compact subgroup of $G$. Then there exists an open set $O$ in $G$, $F \subset O$, with the property that if $H$ is a compact subgroup of $G$ and $H \subset O$, then there is a $g$ in $G$ such that
$$
g^{-1}Hg \subset F\text.
$$
Moreover given any neighborhood $W$ of $e$, $O$ can be so chosen that for every $H \subset O$ the desired $g$ can be selected in $W$.
\end{theorem*}
Compare \mbox{Bredon} (1972) \cite{Bredon'72} II.5.6.
\end{nota}

\sezione{Normalized Haar Systems}\label{O.3.2}
Normalized Haar systems on proper Lie groupoids are the analogue of Haar probability measures on compact Lie groups. In the present section we review the basic definitions and give some details about the construction of Haar systems on proper Lie groupoids; an exposition of this material can also be found in \mbox{Crainic} \cite{Crainic'03}. Let \G\ be a Lie groupoid over a manifold $M$.

\begin{definizione}\label{iii1}
A \index{positive Haar system}\index{Haar system!positive}\textit{positive Haar system} on \G\ is a family of positive Radon measures $\{\mu^x\}$ ($x \in M$), each one with support concentrated in the respective source fibre $\G^x = \G(x,\text-) = \s^{-1}(x)$, satisfying the following conditions:
\begin{elenco}%
\item${\int \varphi\,\mathit d\mu^x} > 0$ for all nonnegative $\varphi \in \test(\G^x)$, $\varphi \neq 0$;%
\item for every $\varphi \in \test({\mca[1]\G;\nC})$, the function $\Phi: M \to \nC$ defined by%
\begin{equazione}\label{iii3}%
\Phi(x) \bydef {\int_{\G^x} \varphi|_{\G^x}\:\mathnormal d\mu^x}%
\end{equazione}%
is of class $\C^\infty$;%
\item{\em right invariance:} for arbitrary $g \in \G(x,x')$ and $\varphi \in \test(\G^x)$,%
\begin{equazione}\label{iii4}%
{\int_{\G^{\smash{x'}}} {\varphi\circ\tau^g}\:\mathit d\mu^{x'}} = {\int_{\G^x} \varphi\,\mathit d\mu^x}%
\end{equazione}%
where $\tau^g: \G(x',\text-) \to \G(x,\text-)$ denotes right translation $h \mapsto hg$.%
\end{elenco}
In this definition the term `positive' refers to the first condition whereas the term `smooth' is occasionally used to emphasize the second condition.
\end{definizione}

The existence of positive (smooth) Haar systems on a Lie groupoid \G\ can be established if \G\ is {\em proper.} (Recall that \G\ is proper if it is Hausdorff and the map $(\s,\t): \G \to {M\times M}$ is proper in the usual sense.) One way to do this is the following. One starts by fixing a Riemann metric on the vector bundle $\mathfrak g \to M$, where $\mathfrak g$ is the Lie \index{algebroid}\index{Lie algebroid}algebroid of \G\ (cfr.\ \mbox{Crainic} \cite{Crainic'03} or \mbox{Moerdijk} and \mbox{Mr\v cun} \cite{Moerdijk&Mrcun'03}, Chapter 6; note the use of paracompactness). Right translations determine isomorphisms $\T{\G(x,\text-)} \iso {\t^*\mathfrak g}|_{\G(x,\text-)}$ for all $x \in M$. These can be used to induce, on the source fibres $\G(x,\text-)$, Riemann metrics whose associated volume forms provide the desired system of measures.

Positive Haar systems are not entirely adequate for our purposes. We will find the following notion more useful:
\begin{definizione}\label{iii2}
A \index{Haar system!normalized}\index{normalized Haar system}\textit{normalized Haar system} on \G\ is a family of positive Radon measures $\{\mu^x\}$ ($x \in M$), each with support concentrated in the respective source fibre $\G(x,\text-)$, with the following properties: \textsl{(a)}~all smooth functions on $\G(x,\text-)$ are integrable with respect to $\mu^x$, that is to say
\begin{equazione}\label{iii.5}
\C^\infty\bigl({\G(x,\text-);\nC}\bigr) \subset \Lebesgue(\mu^x;\nC)\text;
\end{equazione}
\textsl{(b)}~Conditions \textsl{ii)} and \textsl{iii)} of the preceding definition hold for an arbitrary smooth function $\varphi$ on \mca[1]\G, respectively $\G(x,\text-)$; \textsl{(c)}~the following normality condition is satisfied:
\begin{elenco}
\item[\textsl{i$\,'$)}]${\int \mathit d\mu^x} = 1$, for all $x \in M$.
\end{elenco}
\end{definizione}

Every proper Lie groupoid admits normalized (smooth) Haar systems. For such a groupoid \G, one can prove this by using a \index{cut-off function}cut-off function, namely a positive, smooth function $c$ on the base $M$, such that the source map \s\ restricts to a proper map on $\support{({c\circ\t})}$ and ${\int ({c\circ\t})\, \mathit d\nu^x} = 1$ for all $x \in M$, where $\{\nu^x\}$ is a fixed positive (smooth) Haar system on \G. The system of positive measures $\mu^x = {({c\circ\t}) \nu^x}$ has the desired properties.%

Observe that if $E \in {\Ob\,\V[\infty]M}$ is a smooth vector bundle of locally finite rank over the base of \G\ and $\psi: \G \to E$ is a smooth mapping such that for each $x \in M$ the fibre $\G(x,\text-)$ is mapped into the vector space $E_x$, then the integral%
\begin{equazione}\label{iii.6}%
\Psi(x) \bydef {\int \psi_x\,\mathit d\mu^x}%
\end{equazione}%
makes sense and defines a smooth section of $E$. This follows easily from the properties listed in Definition~\refcnt{iii2}, by working in local coordinates.%

\sezione{The Local Linearizability Theorem}\label{N.4}
Let \G\ be a Lie groupoid and let $M$ be its base manifold. We say that a submanifold $N$ of $M$ is a \index{slice}\textit{slice} at the point $z \in N$ if the orbit immersion ${\G z} \into M$ is transversal to $N$ at $z$. A submanifold $S$ of $M$ will be called a slice if it is a slice at all of its points. The following remark will be used very often: \textsl{Let $N$ be a submanifold of $M$ and let $g \in \G^N \equiv \s^{-1}(N)$; then $N$ is a slice at $z = \s(g)$ if and only if the intersection ${\G^N\cap \t^{-1}(z')}$, $z' = \t(g)$ is transversal at $g$.} Indeed, from the equalities%
$$%
\T[g]{\G^N} = {\T[z]N \oplus \T[g]{\G^z}} \quad\text{and} \quad\T[g]{\t^{-1}(z')} = {\T[z]{\G z'} \oplus W} = {\T[z]{\G z}\oplus W}\text,%
$$%
where $W$ is a linear subspace of \T[g]{\G^z}, it follows immediately that%
\begin{equazione}\label{iv.2}%
{\T[g]{\G^N} + \T[g]{\t^{-1}(z')}} = {\bigl({\T[z]N + \T[z]{\G z}}\bigr) \oplus \T[g]{\G^z}}\text.%
\end{equazione}%
By virtue of this fact, one obtains that \textsl{for each submanifold $N$ of $M$, the subset of all points at which $N$ is a slice is an open subset of $N$.} In order to ascertain it, fix a point $z$ belonging to this subset. Since the intersection of $\G^N$ with the fibre $\t^{-1}(z)$ must be transversal at $\u(z) \in \G(z,z)$, there will be a neighbourhood $\Gamma^N$ of $\u(z)$ in $\G^N$ such that for all $g \in \Gamma^N$ the intersection ${\G^N\cap \t^{-1}(\t[g])}$ is transversal at $g$. Now, if $S$ is an open neighbourhood of $z$ in $N$ such that $\u(S) \subset \Gamma^N$, one has that $S$ is a slice.%

Let $R$, $S$ be mutually transversal submanifolds of a manifold $N$: then ${R\cap S}$ is a submanifold of $N$, of dimension ${r+s-n}$.%

Next, let $p: Y \to X$ be a submersion, let $S$ be any submanifold of $Y$ and fix $s_0 \in S$. Put $x_0 = p(s_0)$. Then \textsl{$S$ intersects the fibre $p^{-1}(x_0)$ transversally at $s_0$ if and only if the restriction $p|_S: S \to X$ is submersive at that point;} from this, it immediately follows that when the intersection ${S\cap p^{-1}(x_0)}$ is transversal at $s_0$, there exists a neighbourhood $A$ of $s_0$ in $S$ such that at all points $a \in A$ the intersection ${S\cap p^{-1}(x)}$, $x = p(a)$ is also transversal. In order to check the previous claim, it is not restrictive to assume that $Y = {X\times Z}$ is a Cartesian product and that $p = \pr$ is the projection onto the first factor. Setting $s_0 = (x_0,z_0)$, one obtains for the tangent spaces the picture%
\begin{equazione}\label{iv.3}%
{\T[s_0]S + \T[z_0]Z}\, \subset\, \T[s_0]({X\times Z})\, =\, {\T[x_0]X \oplus \T[z_0]Z}\, \xto{\pr_*}\, \T[x_0]X\text,%
\end{equazione}%
from which it is evident that \T[s_0]S contains a linear subspace $W$ such that $\pr_*(W) = \T[x_0]X$ if and only if the inclusion~\refequ{iv.3} is an equality.%

\sloppy%
\begin{nota}\label{iv.4}%
\textsl{If a submanifold $S$ of $M$ is a slice then the intersection ${\s^{-1}(S) \cap \t^{-1}(S)}$ is transversal and the restriction $\G|_S$ is a Lie groupoid over $S$.} Indeed, let us fix $g \in \G(z,z')$ with $z, z' \in S$. Since%
\begin{equazione}\label{iv.5}%
{\T[g]{\s^{-1}(S)} + \T[g]{\t^{-1}(z')}} \subset {\T[g]{\s^{-1}(S)} + \T[g]{\t^{-1}(S)}}\text,%
\end{equazione}%
one immediately obtains the transversality at $g$ of the intersection written above. The target map \t\ will induce a submersion of $\s^{-1}(S)$ onto an open subset of $M$ and this submersion will in turn induce a submersion of ${\s^{-1}(S) \cap \t^{-1}(S)}$ onto $S$.%
\end{nota}%

\fussy%
\begin{nota}\label{iv.6}%
\textsl{Let $S$ be a slice; then ${\G\cdot S}$ is an open subset of $M$.} To verify this it will be enough to show that given any point $z \in S$ there exists a neighbourhood $U$ of $z$ in $M$ such that ${\s^{-1}(S)\cap \t^{-1}(u)} \neq \emptyset$ for all $u \in U$. This is true because the intersection ${\s^{-1}(S)\cap \t^{-1}(z)}$ is nonempty and transversal. Then $U \subset {\G\cdot S}$, from which the inclusion ${\G\cdot z} \subset {\G\cdot U} \subset {\G\cdot S}$ finally follows.%
\end{nota}%

\begin{theorem*}[N.T.~Zung]\label{N.viii1}%
Let \G\ be a proper Lie groupoid and let $X$ be its base manifold. Let $x_0 \in X$ be a point which is not moved by the tautological action of $\G$ on its own base.%

Then there exists a continuous linear representation $G \to \GL(V)$ of the isotropy group $G \equiv \G|_{x_0}$ on a finite dimensional vector space $V$, such that for some open neighbourhood $U$ of $x_0$ one can find an isomorphism of Lie groupoids $\G|_U \iso {G \ltimes V}$ which makes $x_0$ correspond to $0$.%
\end{theorem*}%
\begin{proof}%
See \mbox{Zung's} paper \cite{Zu06}.%
\end{proof}%

We want to give a careful proof of the statement that any proper Lie groupoid is locally \index{Morita equivalent|see{Morita equivalence}}Morita equivalent to the translation groupoid associated with a (linear) compact Lie group action; this will of course follow from \mbox{Zung's} theorem. The latter statement is a key ingredient in the proof of our \gm{main reconstruction theorem}, Theorem \refcnt[N.20]{xx.18}. Let us begin with a technical observation about slices.%

Let $S$, $T$ be two slices in $M$. Let $g_0 \in \G(S,T)$; put $s_0 \equiv \s(g_0) \in S$ and $t_0 \equiv \t(g_0) \in T$. To fix ideas, suppose $\dimens S \leqq \dimens T$. Then we claim that \textsl{there exists a smooth section $\tau: B \to \mca[1]\G$ to the target map of \G, defined over some open neighbourhood $B$ of $t_0$ in $T$, such that $\tau(t_0) = g_0$ and the composite map ${\s\circ\tau}$ induces a submersion of $B$ onto an open neighbourhood of $s_0$ in $S$.} To begin with, let us notice\inciso{in general}that if one is given a couple of smooth submersions $Y \xfrom p X \xto q Z$ with $\dimens Y \geqq \dimens Z$ then for each point $x \in X$ there exists a smooth $p$\nobreakdash-section $\pi: U \to X$, defined over some open neighbourhood $U$ of $p(x)$, such that $\pi(p(x)) = x$ and the composite ${q\circ\pi}: U \to N$ is a submersion onto an open neighbourhood of $q(x)$ in $Z$; this is seen by means of an obvious argument based on elementary linear algebra: there exists a complementary subspace $F$ to $\kernel{\T[x]p}$ in \T[x]X such that ${F + \kernel{\T[x]q}} = \T[x]X$. Now, the intersection%
\begin{equazione}\label{iv.7}%
X \equiv {\s^{-1}(S)\cap \t^{-1}(T)} \subset \mca[1]\G%
\end{equazione}%
is transversal, because for all $g \in \G(s,t)$ with $s \in S$ and $t \in T$, $\s^{-1}(S)$ will intersect $\t^{-1}(t)$ and hence a fortiori $\t^{-1}(T)$ transversally at $g$ (since $S$ is a slice). Moreover, the source map $\s: \G \to M$ restricts to a submersion of $X$ onto $S$, for\inciso{since $T$ is a slice}the submanifold $\t^{-1}(T)$ is transversal to every \s\nobreakdash-fibre it intersects and therefore the restriction $\s: \t^{-1}(T) \to M$ is a submersion. Symmetrically, the induced mapping $\t|_X: X \to T$ will be submersive. Thus we can apply the foregoing general remark about submersions to get a smooth target section $\tau$ with the desired properties.%
\begin{corollario}\label{iv.8}%
Let \G\ be a proper Lie groupoid over a manifold $M$.%

Then for each point $x_0 \in M$ there exist a finite dimensional linear representation $G \to \GL(V)$ of a compact Lie group $G$, and a \G\nobreakdash-invariant open neighbourhood $U$ of $x_0$ in $M$ along with a \index{Morita equivalence}Morita equivalence $\iota: {G\ltimes V} \into \G|_U$, such that $\mca[0]\iota: V \into U$ is an embedding of manifolds mapping $0$ to $x_0$.%
\end{corollario}%
\begin{proof}%
By properness, we can find a slice $S \subset M$ such that ${S\cap {\G\cdot x_0}} = \{x_0\}$. Then $\G|_S$ is a proper Lie groupoid for which the point $x_0$ is invariant. By \mbox{Zung's} theorem, we can assume that there exists an isomorphism of Lie groupoids ${G\ltimes V} \iso \G|_S$, $0 \mapsto x_0$, for some linear compact Lie group action $G \to \GL(V)$. We contend that ${G\ltimes V} \iso \G|_S \into \G|_U$, where $U$ is the open subset ${\G\cdot S} \subset M$, is the Morita equivalence $\iota$ we are looking for. It will be sufficient to prove that the surjective mapping ${V \times_U \G|_U} \to U$, $(v,g) \mapsto \t(g)$ is a submersion. This clearly follows from the preceding observation about slices when we take $T \equiv U$.%
\end{proof}%

We conclude this section with some remarks relating the groupoids $\G|_S$ and $\G|_T$ induced on two different slices $S$ and $T$. Suppose $\dimens S \leqq \dimens T$. Let $s_0 \in S$ and $t_0 \in T$ be two points lying on the same \G\nobreakdash-orbit. Then%
\begin{elenco}%
\item for some open neighbourhoods $B \subset T$ of $t_0$ and $A \subset S$ of $s_0$ there exists a Morita equivalence $\G|_B \epito \G|_A$ mapping $t_0$ to $s_0$ and inducing a submersion of $B$ onto $A$;%
\item for some open neighbourhood $A \subset S$ of $s_0$ there exists an embedding of Lie groupoids $\G|_A \into \G|_T$ mapping $s_0$ to $t_0$ and inducing a slice embedding $A \into T$ (ie an embedding whose image is a slice);%
\item if in particular $\dimens S = \dimens T$ then the Lie groupoids $\G|_S$ and $\G|_T$ are locally isomorphic about the points $s_0$ and $t_0$.%
\end{elenco}%
Let us verify the assertion \textsl{i)}. Choose any $g_0 \in \G(s_0,t_0)$. By the technical observations preceding Corollary~\refcnt{iv.8}, we can find a smooth target section $\tau: B \to \mca[1]\G$ so that ${\s\circ\tau}$ is a submersion onto an open neighbourhood $A \subset S$ of $s_0$. The latter map can be lifted to%
\begin{equazione}\label{iv.9}%
\G|_B \to \G_A\text, \quad h \:\mapsto \:{\tau(\t[h])^{-1} \cdot h \cdot \tau(\s[h])}\text;%
\end{equazione}%
this formula sets up the required Morita equivalence. In an entirely analogous manner assertion \textsl{ii)} can be proved by considering a suitable smooth source section $\sigma: A \to \mca[1]\G$ such that ${\t\circ\sigma}$ is a slice embedding of $A$ into $T$ and then by lifting this embedding to one of Lie groupoids%
\begin{equazione}\label{iv.10}%
\G|_A \into \G|_T\text, \quad g \:\mapsto \:{\sigma(\t[g]) \cdot g \cdot \sigma(\s[g])^{-1}}\text.%
\end{equazione}%
\begin{nota}\label{iv.11}%
\textsl{Let $\sigma: U \to \mca[1]\G$ be a local bisection. Suppose $S \subset U$ is a slice. Then $T \equiv \t\bigl(\sigma(S)\bigr)$ is also a slice; moreover, there exists a Lie groupoid isomorphism $\G|_S \xto\iso \G|_T$ which lifts the map ${\t\circ\sigma}$.}%

Let us prove that $T$ is a slice. Put $V = \t\bigl(\sigma(U)\bigr)$. Fix a point $s_0 \in S$ and let $t_0 \equiv \t(\sigma(s_0))$. Then%
\begin{equazione}\label{iv.12}%
\t\bigl(\sigma({{\G\cdot s_0}\cap U})\bigr) = {{\G\cdot t_0}\cap V}%
\end{equazione}%
and therefore, since ${\t\circ\sigma}$ is a diffeomorphism of $U$ onto $V$, the orbit ${\G\cdot s_0}$ intersects the submanifold $S$ transversally at $s_0$ if and only if ${\G\cdot t_0}$ intersects $T$ transversally at $t_0$; our claim follows. Next, observe that ${\t\circ\sigma}$ is certainly a diffeomorphism of $S$ onto $T$, which can be lifted\inciso{via $\sigma$, as in \refequ{iv.10}}to a Lie groupoid isomorphism with the expected properties.%
\end{nota}%

\sezione{Global Quotients}\label{N.5}
The material presented in this section is not directly relevant to the problem discussed in the thesis; if the reader wishes to do so, he may go directly to the next chapter. As before, we lay no claim to originality.%

\begin{lemma}\label{v.1}%
Let \H\ be a proper Lie groupoid, acting without isotropy on its own base $F$ (i.e.\ all isotropy groups of \H\ are assumed to be trivial).%

Then the orbit space ${F/\H}$ has a unique manifold structure such that the quotient map $q: F \to {F/\H}$ is a submersion.%
\end{lemma}%
\begin{proof}%
The mapping $(\s,\t): \H \to {F\times F}$ is an injective immersion. Indeed, for a fixed $h \in \H(f,f')$, $f, f' \in F$, the tangent map%
\begin{equazione}\label{v.2}%
\T[h]\H \xto{\:\:\T[h]{(\s,\t)}\:\:} \T[(f,f')]{({F\times F})} \can {\T[f]F\oplus \T[f']F}%
\end{equazione}%
equals the linear map ${\T[h]{\s}\oplus \T[h]{\t}}$; therefore%
\begin{equazione}\label{v.3}%
\kernel{\T[h]{(\s,\t)}} = {\kernel{\T[h]{\s}}\cap \kernel{\T[h]{\t}}} = \T[h]{\H(f,f')} = 0%
\end{equazione}%
(cfr.\ for example \cite{Moerdijk&Mrcun'03}, \textit{proof of Thm.~5.4,} p.~117; by the triviality of the isotropy groups of $\H$, the latter tangent space must be zero).%

Moreover, because of properness, $(\s,\t): \H \to {F\times F}$ is also a closed map, hence in fact an embedding of smooth manifolds.%

It follows that the equivalence relation $R = \image{(\s,\t)} = \{(f,f')|\exists h: f \to f'\text{~in~}\H\}$ is a submanifold of ${F\times F}$; the projection onto the second factor clearly restricts to a submersion of $R$ onto $F$. Therefore, by \mbox{Godement's} Theorem (\textit{see}~\cite{Serre'64}, p.~92), there is a manifold structure on the quotient space ${F/R} = {F/\H}$, making the quotient map $q: F \to {F/\H}$ a submersion.%
\end{proof}%

This lemma applies when a proper Lie groupoid \G\ with base $M$ {\em acts freely} from the left on a manifold $F$ along some smooth mapping $p: F \to M$. By definition, this means that the corresponding action groupoid $\H \equiv {\G\ltimes F}$ has trivial isotropy groups. In order to conclude that there exists a smooth manifold structure on the quotient space ${F/\G}$, for which the projection $F \to {F/\G}$ is submersive, one needs to check that the groupoid ${\G\ltimes F}$ is also proper. So, let $C \subset {F\times F}$ be any compact subset and put $C_1 = \pr_1(C) \subset F$; since $F$ is a Hausdorff manifold, the inverse image $(\s_{\H},\t_{\H})^{-1}(C)$ will be a closed subset of the manifold ${\G\times F}$ and hence of the compact set%
\begin{equazione}\label{v.4}%
{{(\s_\G,\t_\G)^{-1}}\bigl(({p\times p})(C)\bigr) \times C_1} \subset {\G\times F}\text,%
\end{equazione}%
where ${p\times p}$ denotes the smooth map $(f,f') \mapsto (p(f),p(f'))$.%

Now, suppose that a Lie group $K$ acts smoothly on $F$ from the right, in such a way that $p: F \to M$ turns out to be a principal $K$\nobreakdash-bundle. Assume that this action commutes with the given left action of \G. Then there is a well-defined induced right action of $K$ on the quotient manifold ${F/\G}$. This is in fact a smooth action because of an elementary property of submersions (see e.g.\ p.~\pageref{O.esp101} below): the action map ${{F/\G}\times K} \to {F/\G}$ has to be smooth because upon composing it with the submersion ${F\times K} \to {{F/\G}\times K}$ one obtains a smooth map, namely ${F\times K} \to F \to {F/\G}$.%

\separazione%

The next result should probably be regarded as folklore. Its statement, along with the key idea for the proof presented here, was suggested to me by \mbox{I.~Moerdijk} as early as the beginning of 2006.%

\begin{theorem}\label{v.5}%
Suppose a proper Lie groupoid \G\ admits a global faithful representation on a smooth vector bundle.%

Then \G\ is \index{Morita equivalence}Morita equivalent to the translation groupoid associated with a compact, connected Lie group action.%
\end{theorem}%
\begin{proof}%
Let $\varrho: \G \into \GL(E)$ be a faithful representation on a\inciso{let us say, real}smooth vector bundle $E$ over the base $M$ of \G. By properness of \G, we can find a $\varrho$\nobreakdash-invariant metric\footnote{This can be proved in a standard way, very much like in the case of groups, by using Haar systems as a substitute for Haar measures. Cfr.\ Proposition~\ref{N.17}.\ref{N.iv6}.} on $E$, which we fix once and for all. Then let $F = \mathrm{Fr}(E) \xto p M$ be the {\em orthonormal frame bundle} of $E$ (relative to the chosen invariant metric): recall that the fibre $F_x$ above $x$ is the space of all linear isometries $f: \nR^d \isoto E_x$, where $d$ is the rank of $E_x$. The total space $F$ of this fibre bundle is a paracompact Hausdorff manifold; moreover, the fibration $p$ is a principal bundle for the canonical right action of the orthogonal group $K = \mathit{O}(d)$ on $F$ (defined by $fk = {f\circ k}$). Since $\varrho$ acts on $E$ by isometries, the Lie groupoid \G\ will act on the manifold $F$ from the left\inciso{via the representation $\varrho$, that is by the rule $gf = {\varrho(g)\circ f}$}along the bundle map $p$. Clearly, the two actions commute.%

Next, let the ``double action groupoid'' ${\G\ltimes F \rtimes K}$ be the Lie groupoid over the manifold $F$ that is obtained as follows. Its manifold of arrows is ${({\G\ltimes F})\times K}$, viz.\ the submanifold of the Cartesian product ${({\G\times F})\times K}$ consisting of all triples $(g,f,k)$ with $\s(g) = p(f)$. The source map sends the arrow $(g,f,k)$ to $f$ and the target map to $gfk$. The composition of arrows is defined to be ${(g',f',k')\cdot (g,f,k)} = (g'g,f,kk')$. Then the identity arrow at $f$ is $(p(f),f,\id)$ and the inverse must be given by $(g,f,k)^{-1} = (g^{-1},gfk,k^{-1})$. All these structure maps are obviously smooth.%

Now, we claim that there are Morita equivalences%
\begin{equazione}\label{v.6}%
\G \xfrom{\quad\tilde p\quad} {\G\ltimes F \rtimes K} \xto{\quad\tilde q\quad} {{F/\G}\rtimes K}%
\end{equazione}%
from the double action groupoid. This will show that \G\ is Morita equivalent to the action groupoid ${{F/\G}\rtimes K}$, as contended. Perhaps it is good to spend a couple of words to state the formulas for {\em right} action groupoids; these are obtained by regarding a given right action of a Lie group $H$ on a manifold $X$ as a left action of the opposite group. Thus $(x,h) \mapsto x$, resp.\ $\mapsto {x\cdot h}$ is the source, resp.\ target map, and ${(x',h')\cdot (x,h)} = (x,hh')$ is the composition.%

We start with the construction of the equivalence to the left%
\begin{equazione}\label{v.7}%
\tilde p: {\G\ltimes F \rtimes K} \longto \G\text.%
\end{equazione}%
As the notation $\tilde p$ suggests, this equivalence is to be given by the surjective submersion $p: F \to M$ on base manifolds; as to arrows, we put $\tilde p(g,f,k) = g$. It is immediate to check that $\tilde p$ defines a Lie groupoid homomorphism of ${\G\ltimes F \rtimes K}$ onto \G. All one needs to do now in order to show that $\tilde p$ is a Morita equivalence is to solve, within the category of smooth manifolds, the universal problem stated in the left-hand diagram below:%
\begin{equazione}\label{v.8}%
\begin{split}%
\xymatrix@C=13pt{X\ar@/_1.5pc/[ddr]_{(f,f')}\ar@{-->}[dr]\ar@/^1pc/[drr]^(.6)g & & & X\ar@/_1.5pc/[ddr]_{(f,f')}\ar@{-->}[dr]\ar@/^1pc/[drr]^(.5){({q\circ f},k)} & & \\ & {\G\ltimes F \rtimes K}\ar[r]^(.64){\tilde p}\ar[d] & \G\ar[d] & & {\G\ltimes F \rtimes K}\ar[r]^(.5){\tilde q}\ar[d] & {{F/\G}\rtimes K}\ar[d] \\ & {F\times F}\ar[r]^-{p\times p} & {M\times M} & & {F\times F}\ar[r]^-{q\times q} & {{F/\G} \times {F/\G}}\text.\!\!}%
\end{split}%
\end{equazione}%
It will be enough to notice that the map $X \to K\text,\: x \mapsto \kappa(x)$, which assigns the linear isometry $\kappa(x) = {f'(x)^{-1} \circ \varrho(g(x)) \circ f(x)}$ to each $x$, is of class $\C^\infty$. Then we can define the dashed arrow in the aforesaid diagram to be $x \mapsto (g(x),f(x),\kappa(x))$. This is clearly the unique possible solution.%

We turn our attention now to the other equivalence%
\begin{equazione}\label{v.9}%
\tilde q: {\G\ltimes F \rtimes K} \longto {{F/\G}\rtimes K}\text.%
\end{equazione}%
This is given by $q$ on objects and by $\tilde q(g,f,k) = (q(f),k)$ on arrows. Clearly, the map $\tilde q$ so defined is a homomorphism of Lie groupoids. Since the base mapping $q: F \to {F/\G}$ is known to be a surjective submersion by Lemma~\refcnt{v.1}, in order to show that $\tilde q$ is a Morita equivalence it will be enough to solve the right-hand universal problem of \refequ{v.8}. We observe that from the properness of \G\ and the faithfulness of $\varrho$ it follows\inciso{see for instance Corollary~\refcnt[O.2.2+N.23]{O.cor4} below}that the image $\varrho(\G) \subset \GL(E)$ is a submanifold; moreover, it can be shown\inciso{cfr.\ Lemma~\refcnt[O.3.5]{O.lem7}, for example}that $\varrho: \G \xto\iso \varrho(\G)$ is actually a diffeomorphism. Now, the map $X \to \GL(E)\text,\: x \mapsto \gamma(x)$, that sends $x$ to the isometry $\gamma(x) = {f'(x)\circ k(x)\circ f(x)^{-1}}$, is clearly smooth and factors through the submanifold $\varrho(\G)$. Then we may use the fact that $\varrho$ is a diffeomorphism of \G\ onto $\varrho(\G)$ and define the dashed arrow as $x \mapsto \bigl(\varrho^{-1}(\gamma(x)),f(x),k(x)\bigr)$; this is of course a smooth correspondence.%
\end{proof}%

\capitolo{The Language of Tensor Categories}\label{2}

With the exception of \refsez{O.1.4+.3.3}, the present chapter offers an introduction to the categorical setting of the modern theory of Tannaka duality originating from the ideas of \mbox{A.~Grothendieck} and \mbox{N.~Saavedra Rivano;} cfr \cite{Sa72,DeMi82,De90,JoSt91}.

In Section \ref{O.1.4+.3.3} we prove a key technical lemma which will be used in the proof of our reconstruction theorem in \refsez{N.20}; since this lemma deals with a fairly abstract categorical situation, we thought it was more appropriate to include it in this chapter.

\sezione{Tensor Categories}\label{O.1.1+N.6}
A \index{tensor structure|see{tensor category}}\textit{tensor structure} on a category \Kt\ consists of the following data:
\index{tensor unit@tensor unit \ensuremath{\TU}}\begin{equazione}\label{vi1}
\text{a bifunctor~} \boldsymbol\otimes: {\Kt\times\Kt} \longto \Kt\text, \quad \text{a distinguished object~} \TU \in \Ob(\Kt)
\end{equazione}
and a list of natural isomorphisms, called \index{ACU constraints@\textit{ACU} constraints|emph}\index{constraint}\textit{ACU constraints:}
\begin{equazione}\label{vi2}
\begin{split}
\begin{array}{c}
\alpha_{R,S,T}: {R\otimes ({S\otimes T})} \isoto {({R\otimes S})\otimes T}\text,
\\[\medskipamount]
\gamma_{R,S}: {R\otimes S} \isoto {S\otimes R}\text,
\\[\medskipamount]
\lambda_R: R \isoto {\TU\otimes R} \qquad \text{and} \qquad \rho_R: R \isoto {R\otimes\TU}
\end{array}
\end{split}
\end{equazione}
satisfying MacLane's \index{coherence conditions}\textit{coherence conditions} (cfr for example MacLane (1971), pp.~157~ff.\ and especially p.~180 for a detailed exposition). A \index{tensor category|emph}\textit{tensor category} is a category endowed with a tensor structure. In the terminology of \cite{MacLa71}, the present notion corresponds to that of ``{symmetric monoidal category}''. The natural isomorphism $\alpha$ resp.\ $\gamma$ is called the \index{associativity constraint}\textit{associativity} resp.\ \index{commutativity constraint}\textit{commutativity} constraint; $\lambda$ and $\rho$ are the \index{tensor unit constraints}\textit{tensor unit} constraints.

In order to state MacLane's Coherence Theorem for tensor categories, it will be convenient to introduce the concepts of \gm{canonical multi-functor} and \gm{canonical transformation}. These will constitute respectively the objects and the morphisms of a category \index{Can(C)@\Can{\Kt} (category of canonical multi-functors)}\Can\Kt.

A \index{multi-functor}\textit{multi-functor} on \Kt\ is a functor $\Phi: \Kt^I \to \Kt$ for some finite set $I$. The cardinality $\modulo I = \mathrm{Card}\,I$ will be called the \gm{-ariety} of $\Phi$.

The \index{canonical multi-functor}\textit{canonical} multi-functors are, roughly speaking, those obtained as {\em formal} iterates of $\boldsymbol\otimes$, possibly involving \TU. The adjective `formal' here means that a `{canonical multi-functor}' is not just a certain type of multi-functor, in that one should regard the particular inductive construction, by which a canonical multi-functor is obtained, as part of the defining data; we do not want to go into details here: the interested reader may consult \cite{MacLa71}. The recursive rules for generating canonical multi-functors are listed below:
\begin{elenco}
\item the unique $0$-ary canonical multi-functor is $\TU: \Kt^\varnothing = \{\pt\} \to \Kt$, $\pt \mapsto \TU$;
\item the ``identity'': $\Kt^{\{\pt\}} \to \Kt$ is canonical;%
\item if $\Phi: \Kt^I \to \Kt$ and $\Psi: \Kt^J \to \Kt$ are canonical then so is ${\Phi\otimes\Psi}: \Kt^{I\sqcup J} \to \Kt$ where ${I\sqcup J}$ indicates disjoint union;%
\item if $I \xto\sigma J$ is a bijection of finite sets and $\Phi: \Kt^I \to \Kt$ is canonical then $\Phi^\sigma: \Kt^J \to \Kt^I \to \Kt$ is also canonical.
\end{elenco}
Canonical multi-functors are the objects of \Can\Kt. As to \index{canonical transformation}canonical natural transformations, they are recursively generated as follows:
\begin{elenco}
\item[\textsl{a)}]the identity $\id: \Phi \to \Phi$ is canonical; if $\eta: \Phi \to \Phi'$, with $\Phi, \Phi': \Kt^I \to \Kt$, and $\theta: \Psi \to \Psi'$, with $\Psi, \Psi': \Kt^J \to \Kt$, are canonical transformations of canonical multi-functors, then so is ${\eta\otimes\theta}: {\Phi\otimes\Psi} \to {\Phi'\otimes\Psi'}$ (natural transformations of multi-functors $\Kt^{I\sqcup J} \to \Kt$); if $I \xto\sigma J$ is a bijection of sets then $\theta^\sigma: \Phi^\sigma \to \Psi^\sigma$ is also canonical;
\item[\textsl{b)}]$\alpha_{\Phi,\Psi,\mathrm X}: \smash{\bigl[{\Phi \otimes ({\Psi \otimes \mathrm X})}\bigr]}^\sigma \isoto \smash{\bigl[{({\Phi \otimes \Psi}) \otimes \mathrm X}\bigr]}^\tau$ and its inverse ${\alpha_{\Phi,\Psi,\mathrm X}}^{-1}$ are canonical transformations, where $\sigma$, $\tau$ are the bijections ${I\sqcup ({J\sqcup K})} \to {I\sqcup J\sqcup K} \from {({I\sqcup J})\sqcup K}$;%
\item[\textsl{c)}]$\gamma_{\Phi,\Psi}: {\Phi\otimes\Psi} \isoto [{\Psi\otimes\Phi}]^\sigma$ (along with its inverse) is canonical, where ${I\sqcup J} \xfrom\sigma {J\sqcup I}$ is the obvious bijection;%
\item[\textsl{d)}]$\lambda_\Phi: \Phi \isoto ({\TU\otimes\Phi})^\sigma$ and $\rho_\Phi: \Phi \isoto ({\Phi\otimes\TU})^\tau$ (along with their inverses) are canonical, where ${\varnothing\sqcup I} \xto\sigma I \xfrom\tau {I\sqcup\varnothing}$ are the obvious bijections.%
\end{elenco}%
It is clear that all canonical transformations are isomorphisms.%

\index{coherence theorem}\textit{MacLane's Coherence Theorem} for \gm{symmetric monoidal categories} (\gm{tensor categories} in our terminology) can now be stated as follows:
\begin{theorem*}
The category \Can\Kt\ is a preorder. That is to say, for any canonical multi-functors $\Phi$ and $\Psi$ there is at most one canonical natural transformation $\Phi \to \Psi$.%
\end{theorem*}%
\begin{proof}%
See {[MacLane]}, \textsc{xi.1} p.~253.%
\end{proof}%

This theorem says that any diagram of canonical multi-functors and canonical natural transformations one can possibly construct will commute. When one is given such a diagram, let us say of multi-functors $\Kt^I \to \Kt$, one may choose an identification $\{1,\ldots,i\} \isoto I$ and denote a generic object of $\Kt^I$ by $(R_1,\ldots,R_i)$, $R_1, \ldots, R_i \in \Ob(\Kt)$. Evaluating the given diagram at this $i$\nobreakdash-tuple of objects\inciso{so that $\Phi \xto\theta \Psi$ becomes $\Phi(R_1,\ldots,R_i) \xto{\theta(R_1,\ldots,R_i)} \Psi(R_1,\ldots,R_i)$, for instance}one obtains a commutative diagram in \Kt.%

\begin{nota}\label{vi.6}%
(See also \textit{Saavedra,} 1.3.3.1) \textsl{Let $(\Kt,\boldsymbol\otimes,\TU)$ be a tensor category. Then $\End_{\Kt}(\TU)$ is a commutative ring.} To see this, observe that the tensor unit constraint $\TU \can {\TU\otimes\TU}$ establishes a canonical isomorphism of rings between $\End(\TU)$ and $\End({\TU\otimes\TU})$. Now, if $e, e' \in \End(\TU)$ then ${ee'} \can {({1\otimes e})({e'\otimes 1})} = {e'\otimes e} = {({e'\otimes 1})({1\otimes e})} \can {e'e}$ in this isomorphism and hence ${ee'} = {e'e}$. Note that this proof only uses the coherence identity $\lambda_\TU = \rho_\TU$ for the tensor unit constraints; the commutativity constraint plays no role.%
\end{nota}%

\sottosezione{Rigid tensor categories}%

A tensor category $(\Kt,\boldsymbol\otimes)$ is said to be \index{closed tensor category}\index{tensor category!closed}\textit{closed,} whenever one can exhibit a bifunctor $\underline{\mathrm{hom}}: {\opposite\Kt \times \Kt} \longto \Kt$, called `\index{internal hom (bifunctor)}internal hom' and denoted by
$$
(X,Y) \mapsto Y^X \equiv \underline{\mathrm{hom}}(X,Y)\text,
$$
along with natural transformations (in the variable $Y$)
$$
\eta^X_Y: Y \to ({Y\otimes X})^X\qquad \text{and}\qquad \varepsilon^X_Y: {Y^X\otimes X} \to Y\text,
$$
satisfying the triangular identities for an adjunction
$$
\Kt\bigl({X\otimes T},Y\bigr) \isoto \Kt\bigl(X,\underline{\mathrm{hom}}(T,Y)\bigr)\quad \text(\text{in~the~variables~} (X,Y) \in {\opposite\Kt\times\Kt}\text)%
$$%
between the functors \gm{${\text-\otimes T}$} and \gm{$\underline{\mathrm{hom}}(T,\text-)$} and making%
\begin{equazione}\label{O.equ200}%
\begin{split}
\xymatrix@C=40pt{{Y^{X'}\otimes X}\ar[d]^{\id\otimes a}\ar[r]^-{Y^a\otimes \id} & {Y^X\otimes X}\ar[d]^\varepsilon & ({Y \otimes X})^X\ar[r]^-{({\id\otimes a})^\id} & ({Y\otimes X'})^X \\ {Y^{X'}\otimes X'}\ar[r]^-\varepsilon & Y & Y\ar[u]_\eta\ar[r]^-\eta & ({Y\otimes X'})^{X'}\ar[u]_{\id^a}}%
\end{split}
\end{equazione}%
commute for every arrow $a: X \to X'$.%

Suppose now that an `internal hom' bifunctor and natural transformations $\eta$, $\varepsilon$ with these properties have been fixed. Then there is an obvious arrow%
\begin{equazione}\label{O.equ201}%
\delta^{S,T}_{X,Y}: {X^S\otimes Y^T} \to ({X \otimes Y})^{S\otimes T}\text,%
\end{equazione}%
namely the unique solution $d$ to the equation%
$$%
{\varepsilon\circ({d\otimes \id})} = {({\varepsilon\otimes\varepsilon})\circ\can},%
$$%
where $\can$ is the unique canonical isomorphism. Because of \refequ{O.equ200}, the arrow $\delta$ must be natural in all variables. By the same reason, the solution%
\begin{equazione}\label{O.equ202}%
\iota_X: X \to \bidual{X}%
\end{equazione}%
(where we put $\dual X \equiv \underline{\mathrm{hom}}(X,\TU)$, to be called the \index{dual}\textit{dual} of $X$) to the equation
$$%
{\varepsilon\circ({\iota_X\otimes\id})} = {\varepsilon\circ\can}%
$$%
is natural in $X$.%

A different choice of internal hom bifunctor and natural transformations $\eta$ and $\varepsilon$ will yield the same natural arrows $\delta$ and $\iota$ up to isomorphism: thus it makes sense to call a closed tensor category \index{rigid tensor category}\index{tensor category!rigid}\textit{rigid} when these natural arrows are {\em isomorphisms.}%

One can also formulate this notion in terms of duals, since for a rigid tensor category one has the identification%
\begin{equazione}\label{vi.7}%
\underline{\mathrm{hom}}(X,Y) \iso {\dual X\otimes Y}\text,%
\end{equazione}%
cf.\ {\em Deligne~(1990),} \cite{De90}~2.1.2.%

Let $(\Kt,\otimes)$ be a rigid tensor category. The contravariant functor%
$$%
X \mapsto {\dual X}, \quad f \mapsto \trasposto{f}%
$$%
is an equivalence between $\Kt$ and its opposite category \opposite\Kt\ (because it is {\em involutive}, ie its composite with itself is naturally isomorphic to the identity, since rigidity implies that (\ref{O.equ202}) is a natural isomorphism).%

This gives in particular a bijection between the hom-sets%
$$%
f \mapsto \trasposto{f}: \Hom_{\Kt}(X,Y) \isoto \Hom_{\Kt}({\dual Y},{\dual X})\text{.}%
$$%
One also has an ``internal'' isomorphism%
$$%
Y^X \isoto {\dual X}^{\dual Y},%
$$%
namely the composite%
$$%
Y^X \xleftarrow{\, \iso \,} {{\dual X} \otimes Y} \xto{\, {\id \otimes \iota_Y} \,} {{\dual X} \otimes \bidual{Y}} \xto{\, \iso \,} {\bidual{Y} \otimes {\dual X}} \xto{\, \iso \,} {\dual X}^{\dual Y}.%
$$%

For every object of $\Kt$ there is an arrow $X^X \isoto {{\dual X} \otimes X} \xto{\varepsilon} \TU$. If we apply the functor $\Hom_{\Kt}(\TU,\cdot)$ to this, we obtain the \index{trace}\textit{trace} map%
\begin{equazione}\label{O.esp109}%
\mathrm{Tr}_X: \End_{\Kt}(X) \to \End_{\Kt}(\TU).%
\end{equazione}%
The \index{rank}\textit{rank} of $X$ is defined as $\mathrm{Tr}_X(1_X)$. There are the formulas%
\begin{equazione}\label{O.esp108}%
\begin{split}
\begin{array}{l}%
\mathrm{Tr}_{X \otimes X'}(f \otimes f') = \mathrm{Tr}_X(f) \mathrm{Tr}_{X'}(f'),%
\\[\medskipamount]%
\mathrm{Tr}_\TU(f) = f.%
\end{array}%
\end{split}
\end{equazione}%

A tensor category $(\Kt,\otimes)$ is said to be \index{additive tensor category}\index{tensor category!additive}\textit{additive} if the category $\Kt$ is endowed with an additive structure such that the bifunctor $\otimes$ is biadditive, that is additive in each variable separately. Moreover, if the hom-sets $\Kt(A,B)$ are endowed with a real (or complex) vector space structure in such a way that composition of arrows and the bifunctor $\otimes$ are bilinear, then we say that $(\Kt,\otimes)$ is a \index{linear tensor category}\index{tensor category!linear}\textit{linear} tensor category.

\begin{esempio}\label{vi3}
Let $\SpV_\nC$ be the category of vector spaces over $\nC$ of finite dimension. Then this is an abelian rigid tensor category, and all the preceding definitions have their usual meaning.%
\end{esempio}%

\begin{esempio}\label{vi4}
Let $M$ be a smooth manifold. Let \index{V(X;k)@\V[\infty]{X;k}, \stack[X]{V^\infty} (category of vector bundles)|emph}\Kt[C] = \V[\infty]{M;\nC} be the category of smooth complex vector bundles of locally finite rank over $M$. The direct sum operation $(E,F) \mapsto {E\oplus F}$ makes it into an additive \nC-linear category, although in general not an abelian one, since a map of vector bundles over $M$ need not have a kernel, for instance. We shall identify the category of finite dimensional vector spaces over $\nC$ with \V[\infty]{\pt;\nC} where $\pt$ is the one-point manifold.

The category $\V[\infty]{M;\nC}$ is endowed with a canonical rigid tensor structure, obtained from the rigid tensor structure of $\SpV_\nC$ by means of the general procedure described in Lang 2001 \cite{La01} p.~58, as follows. Recall that a \index{multi-functor}multi-functor
$$%
\Phi: \underset{n \text{ times}}{\SpV_\nC \times \cdots \times \SpV_\nC} \longto \SpV_\nC%
$$%
(where case $n = 0$ corresponds to the choice of an object $\Phi(\cdot) \in \Ob(\SpV_\nC)$, and we allow $\Phi$ to be contravariant in some variables), such that the induced mappings%
$$%
{\Linear{V_1,W_1} \times \cdots \times \Linear{V_n,W_n}} \to \Linear{\Phi(V_1,\ldots,V_n),\Phi(W_1,\ldots,W_n)}%
$$%
are of class $\mathit{C}^\infty$, determines a corresponding multi-functor%
$$%
\overline{\Phi}: {\V[\infty]{M;\nC} \times \cdots \times \V[\infty]{M;\nC}} \longto \V[\infty]{M;\nC}%
$$%
with the same variance and satisfying the following properties:%
\begin{elenco}%
\item for every $x \in M$, the fiber above $x$ is%
\begin{equazione}\label{O.esp111}%
\begin{array}{ll}%
{\overline{\Phi}(E_1,\ldots,E_n)}_x &= {\{x\} \times \Phi({E_1}_x,\ldots,{E_n}_x)} \\ &\iso \Phi({E_1}_x,\ldots,{E_n}_x);%
\end{array}%
\end{equazione}%
\item for arbitrary morphisms of vector bundles $a_i: E_i \to F_i, i = 1, \ldots, n$, ${\overline{\Phi}(a_1,\ldots,a_n)}_x$ corresponds to $\Phi({a_1}_x,\ldots,{a_n}_x)$ up to the canonical identifications (\ref{O.esp111});
\item If the vector bundles $E_i \iso {M \times \boldsymbol{E}_i}$ are trivial, then these trivializations $\iso_i$ determine a trivialization%
$$%
\overline{\Phi}(E_1,\ldots,E_n) \iso {M \times \Phi(\boldsymbol{E}_1,\ldots,\boldsymbol{E}_n)}%
$$%
in a canonical way; in particular, in the case $n = 0$, $\overline{\Phi}(\text-) \iso {M \times \Phi(\text-)}$ (the standard notation is then $\underline{\Phi(\text-)} = \overline{\Phi}(\text-)$).%
\end{elenco}%
A natural transformation $\lambda: \Phi \to \Psi$ of multi-functors with the same variance induces a natural transformation $\overline{\lambda}: \overline{\Phi} \to \overline{\Psi}$, such that ${\overline{\lambda}(E_1,\ldots,E_n)}_x$ corresponds to $\lambda({E_1}_x,\ldots,{E_n}_x)$ up to the identifications (\ref{O.esp111}). Observe that $\overline{\lambda \circ \mu} = {\overline{\lambda} \circ \overline{\mu}}$ and $\overline{\id} = \id$.

We can apply these constructions to the multifunctors and natural transformations which define the rigid tensor structure of $\SpV_\nC$, in order to obtain a similar structure on $\V[\infty]{M;\nC}$.
\end{esempio}%

\sezione{Tensor Functors}\label{O.1.2}
Let \Kt[C], \Kt[D] be tensor categories. A \index{tensor functor|emph}\textit{tensor functor} $:\Kt[C] \longto \Kt[D]$ consists of the data $(F,\tau,\upsilon)$, where
$$
F: \Kt[C] \longto \Kt[D]
$$
is a functor, $\tau$ is a \index{constraint|emph}\index{tensor functor constraints|emph}natural isomorphism of bifunctors
$$%
\tau_{R,S}^{}: {F(R)\otimes F(S)} \isoto F({R\otimes S})%
$$%
such that the \index{coherence conditions|emph}diagrams
$$%
\xymatrix@C=35pt{{FR \otimes (FS \otimes FT)}\ar[d]^\alpha\ar[r]^-{\id \otimes \tau} & {FR \otimes F(S \otimes T)}\ar[r]^-\tau & F(R \otimes (S \otimes T))\ar[d]^{F(\alpha)} \\ {(FR \otimes FS) \otimes FT}\ar[r]^-{\tau \otimes \id} & {F(R \otimes S) \otimes FT}\ar[r]^-\tau & F((R \otimes S)\otimes T)}%
$$%
and%
$$%
\xymatrix{{F(R)\otimes F(S)}\ar[d]^\tau\ar[r]^\gamma & {F(S)\otimes F(R)}\ar[d]^\tau \\ F(R \otimes S)\ar[r]^{F(\gamma)} & F(S \otimes R)}%
$$%
commute, and%
$$%
\upsilon: \TU \isoto F(\TU)%
$$%
is an isomorphism in \Kt[D] such that%
$$%
\xymatrix{F(R)\ar[r]^-{F(\lambda)}\ar[d]_\lambda & F(\TU \otimes R) && F(R)\ar[r]^-{F(\rho)}\ar[d]_\rho & F(R \otimes \TU) \\ {\TU \otimes F(R)}\ar[r]^-{\upsilon \otimes \id} & {F\TU \otimes F(R)}\ar[u]^\tau && {F(R) \otimes \TU}\ar[r]^-{\id \otimes \upsilon} & {F(R) \otimes F\TU}\ar[u]^\tau}%
$$%
commute. (Commutativity of one square implies commutativity of the other, because of the symmetry of the monoidal structure.)%

Now suppose that $\Kt[C]$ and $\Kt[D]$ are closed tensor categories. Let $F: \Kt[C] \longto \Kt[D]$ be a tensor functor. (We shall usually omit writing down the full triple of data.) Then there is a canonical arrow%
$$%
p^R_S: F(S^R) \to {F S}^{F R},%
$$%
namely the unique solution $p$ to the problem%
$$%
\xymatrix{{F(S^R) \otimes FR}\ar[d]^\tau\ar[r]^-{p \otimes \id} & {{F S}^{F R} \otimes FR}\ar[d]^\varepsilon \\ F(S^R \otimes R)\ar[r]^-{F(\varepsilon)} & FS.}%
$$%
This arrow is natural in the variables $R, S$. A {\em rigid functor} is a tensor functor between closed tensor categories such that this natural arrow is an isomorphism. If $\Kt[C]$ and $\Kt[D]$ are both rigid, then a tensor functor $F: \Kt[C] \longto \Kt[D]$ is automatically rigid.%
\begin{esempio}\label{vii1}%
Let $f: M \to N$ be a $\C^\infty$-mapping of smooth manifolds. This map induces the \textit{base change} or \index{pullback!along a smooth map}\textit{pullback} functor
$$%
f^*: \V[\infty]N \longto \V[\infty]M.%
$$%
Recall that for $x \in M$ the fiber ${(f^* F)}_x$ coincides with ${\{x\} \times F_{f(x)}}$, since $f^* F$ is by construction a subset of ${M \times F}$. For every functor of several variables $\Phi$ as in the last example of Section \ref{O.1.1+N.6}, we have a canonical natural isomorphism%
\begin{equazione}\label{O.esp113}%
{f^* \overline{\Phi}(E_1,\ldots,E_n)} \iso {\overline{\Phi}(f^*E_1,\ldots,f^*E_n)}.%
\end{equazione}%
It follows at once from the existence of these canonical natural isomorphisms that $f^*$ can be regarded as a tensor functor (with respect to the standard tensor structure described in the last example of the preceding section). It is also clear from (\ref{O.esp113}) that this tensor functor is rigid. (Of course, rigidity of the pullback functor follows also indirectly from rigidity of the categories $\V[\infty]{M}, \V[\infty]{N}$.)%
\end{esempio}%
\begin{definizione}\label{vii2}%
Let $\lambda: F \to G$ be a natural transformation of tensor functors. $\lambda$ is said to be \index{tensor preserving}\textit{tensor-preserving,} or a {morphism of tensor functors}, whenever the diagrams
$$%
\label{O.esp110}\xymatrix{{FR \otimes FS}\ar[d]^\tau\ar[rr]^{\lambda(R) \otimes \lambda(S)} & & {GR \otimes GS}\ar[d]^\tau & & \TU\ar[d]^\upsilon\ar[r]^{\id} & \TU\ar[d]^\upsilon \\ F(R \otimes S)\ar[rr]^{\lambda(R \otimes S)} & & G(R \otimes S) & & {F \TU}\ar[r]^{\lambda(\TU)} & {G \TU}}%
$$%
commute. The collection of all such $\lambda$'s will be denoted by $\Hom^\otimes(F,G)$.%
\end{definizione}%

\sezione{Complex Tensor Categories}\label{N.8}
An \index{anti-involution}\textit{anti-involution} on a \nC-linear tensor category $\Kt = (\Kt,\otimes)$ is an anti-linear tensor functor
\begin{equazione}\label{N.ix1}
\boldsymbol*: \Kt \to \Kt\text, \quad R \mapsto R^*%
\end{equazione}%
for which there exists a tensor preserving natural \index{constraint}isomorphism
\begin{equazione}\label{N.ix2}
\iota_R: R^{**} \isoto R \qquad \text{with} \qquad \iota(R^*) = \iota(R)^*\text.
\end{equazione}
By fixing one such isomorphism, one obtains a mathematical structure which we call \index{complex tensor category|emph}\textit{complex tensor category.} A {morphism} of complex tensor categories, or \index{complex tensor functor|emph}\textit{complex tensor functor,} is obtained by attaching, to an ordinary \nC-linear tensor functor $F$, a tensor preserving natural \index{constraint}isomorphism
\begin{equazione}\label{N.ix3}
\xi_R: F(R)^* \isoto F(R^*)
\end{equazione}
such that the following \index{coherence conditions|emph}diagram commutes
\begin{equazione}\label{N.ix4}
\begin{split}
\xymatrix@R=13pt@C=15pt{F(R)^{**}\ar[r]^-{\can^*}\ar[dr]_\can & F(R^*)^*\ar[r]^-\can & F(R^{**})\ar[dl]^{F\,\can} \\ & FR\text.\!\! &}%
\end{split}%
\end{equazione}%

\begin{esempio}[the category of vector spaces]\label{N.ix5}%
If $V$ is a complex vector space, we let $V^*$ denote the space obtained by retaining the additive structure of $V$ but changing the scalar multiplication into $zv^* = (\overline zv)^*$; the star here indicates that a vector of $V$ is to be regarded as one of $V^*$. Since any linear map $f: V \to W$ will map $V^*$ linearly into $W^*$, we can also regard $f$ as a linear map $f^*: V^* \to W^*$. Moreover, the unique linear map of ${V^*\otimes W^*}$ into $({V\otimes W})^*$ sending ${v^*\otimes w^*} \mapsto ({v\otimes w})^*$ is an isomorphism, and complex conjugation sets up a linear bijection between $\nC$ and $\nC^*$. This turns vector spaces into a complex tensor category with $V^{**} = V$.%
\end{esempio}%

\begin{esempio}[the category of vector bundles over a manifold]\label{N.ix6}%
By using the procedure described in Example~\refcnt[O.1.1+N.6]{vi4} one can transport the complex tensor structure of the preceding example to the category \V[\infty]{M;\nC} of smooth complex vector bundles (of locally finite rank) over a manifold $M$.
\end{esempio}%

\separazione%

Consider a complex tensor category $(\Kt,\otimes,*)$. By a \index{sesquilinear form}\textit{sesquilinear form} on an object $R \in \Ob(\Kt)$ we mean any arrow $b: {R\otimes R^*} \to \TU$. A sesquilinear form $b$ on the object $R$ will be said to be \index{Hermitian form}\textit{Hermitian} when the sesquilinear form $\tilde b$ on $R$, defined as the composite
\begin{equazione}\label{viii.8}%
{R\otimes R^*} \can {R^{**}\otimes R^*} \can ({R\otimes R^*})^* \xto{\:\:b^*\:\:} \TU^* \can \TU\text,%
\end{equazione}%
coincides with $b$ itself, i.e.\ $\tilde b = b$. Note that one always has the equality $\tilde{\tilde b} = b$. Clearly, in the examples above one recovers the familiar notions.%

Suppose now that our complex tensor category is \index{rigid tensor category|emph}\index{tensor category!rigid|emph}\textit{rigid.} Then for each object $R$ we can find another \index{dual|emph}object $R'$, along with arrows $e^R: {R'\otimes R} \to \TU$ and $d^R: \TU \to {R\otimes R'}$, such that the following compositions are identities:
\begin{equazione}\label{viii.9}
\begin{array}{c}%
R \can {\TU\otimes R} \xto{\:\:d^R\otimes R\:\:} {R\otimes R'\otimes R} \xto{\:\:R\otimes e^R\:\:} {R\otimes\TU} \can R%
\\[\medskipamount]%
R' \can {R'\otimes\TU} \xto{\:\:R'\otimes d^R\:\:} {R'\otimes R\otimes R'} \xto{\:\:e^R\otimes R'\:\:} {\TU\otimes R'} \can R'\text.%
\end{array}%
\end{equazione}%
We make the assumption that for each object $R$ we have selected one such triple $(\dual R,e^R,d^R)$. Then for each $R$ we obtain a well-defined isomorphism $q^R: \smash{\dual R}^* \isoto \dual{\smash{(R^*)}}$, namely the unique arrow $q$ such that%
\begin{multiriga}\label{viii.10}%
{\smash{\dual R}^* \otimes R^*} \xto{\:q\otimes R^*\:} {\dual{\smash{R^*}} \otimes R^*} \xto{\:e^{R^*}\:} \TU\quad \text{equals}\\%
{\smash{\dual R}^* \otimes R^*} \can ({\dual R\otimes R})^* \xto{\:(e^R)^*\:} \TU^* \can \TU\text.%
\end{multiriga}%
We say that a sesquilinear form $b$ on $R$ is \index{nondegenerate form}\textit{nondegenerate,} when the arrows $b_{\text-}: R \to \dual{\smash{R^*}}$ and $b^{\text-}: R^* \to \dual R$, defined as the unique solutions to
\begin{multiriga}\label{viii.12}%
{R\otimes R^*} \xto{\:b_{\text-}\otimes R^*\:} {\dual{\smash{R^*}} \otimes R^*} \xto{\:e^{R^*}\:} \TU\quad \text{equals}\quad b\qquad \text{and}\\%
b\quad \text{equals}\quad {R\otimes R^*} \xto{\:R\otimes b^{\text-}\:} {R\otimes\dual R} \can {\dual R\otimes R} \xto{\:e^R\:} \TU\text,%
\end{multiriga}%
are {\em isomorphisms.} \textsl{If $b$ is Hermitian then $b_{\text-}$ is an isomorphism if and only if so is $b^{\text-}$.} Indeed, the diagrams%
\begin{equazione}\label{viii.13}%
\begin{split}%
\xymatrix@C=35pt{R^*\ar[d]^{b^{\text-}}\ar[r]^-{\left({\tilde b_{\text-}}\right)^*} & \smash{\dual{\smash{R^*}}}^* & & R^{**}\ar[d]^\can\ar[r]^-{\left(\tilde b\right)^{\text-*}} & \smash{\dual R}^*\ar[d]^{q^R}_\iso \\ \dual R\ar[r]^-\can & \smash{\dual R}^{**}\ar[u]^\iso_{\left(q^R\right)^*} & & R\ar[r]^{b_{\text-}} & \dual{\smash{R^*}}}%
\end{split}%
\end{equazione}%
commute for an arbitrary sesquilinear form $b$ on $R$.%

\separazione%

Let $(\Kt,\otimes,*)$ be a complex tensor category. By a \index{descent datum}\textit{descent datum} on an object $R \in \Ob(\Kt)$ we mean an isomorphism $\mu: R \isoto R^*$ such that the composition $R \stackrel\mu\iso R^* \stackrel{\mu^*}\iso R^{**} \can R$ equals $\id_R$. We let ${\nR\,\Kt}$ denote the category whose objects are the pairs $(R,\mu)$ consisting of an object $R$ of \Kt\ and a descent datum $\mu$ on $R$ and whose morphisms $a: (R,\mu) \to (R',\mu')$ are the morphisms $a: R \to R'$ such that ${\mu'\cdot a} = {a^*\cdot\mu}$. Note that ${\nR\,\Kt}$ is naturally an \nR\nobreakdash-linear category; moreover, there is an obvious induced tensor structure, which turns ${\nR\,\Kt}$ into an \nR\nobreakdash-linear tensor category.

As an example of this construction, observe that one has an obvious equivalence of real tensor categories between $\SpV_{\nR}$ and $\nR(\SpV_{\nC})$: in one direction, to any real vector space $V$ one can assign the pair $({\nC\otimes V},{{z\otimes v} \mapsto {\overline z\otimes v}})$; conversely, any descent datum $\mu: U \isoto U^*$ on a complex vector space $U$ determines the real subspace $U^\mu \subset U$ of $\mu$\nobreakdash-invariant vectors. More generally, one has analogous equivalences of real tensor categories between \V[\infty]{M;\nR} and $\nR\bigl(\V[\infty]{M;\nC}\bigr)$, \R[\infty]{M;\nR} and $\nR\bigl(\R[\infty]{M;\nC}\bigr)$ and so on.%

Notice that any complex tensor functor $F: \Kt[C] \to \Kt[D]$ will induce a linear tensor functor ${\nR\,F}: {\nR\,\Kt[C]} \to {\nR\,\Kt[D]}$. By using the fact that the isomorphism ${R\oplus R^*} \iso ({R\oplus R^*})^*$ is a descent datum on ${R\oplus R^*}$ for each $R$, one can easily show that setting $\hat\lambda(R,\mu) = \lambda(R)$ defines a \textsl{bijection}%
\begin{equazione}\label{viii.14}%
\Hom^{\otimes,\boldsymbol*}(F,G) \isoto \Hom^\otimes({\nR\,F},{\nR\,G})\text, \quad \lambda \mapsto \hat\lambda%
\end{equazione}%
between the \index{self-conjugate}\textit{self-conjugate} tensor preserving transformations $F \to G$ and the tensor preserving transformations ${\nR\,F} \to {\nR\,G}$, for any complex tensor functors $F, G: \Kt[C] \to \Kt[D]$.

\sezione{Review of Groups and Tannaka Duality}\label{O.1.3}
Throughout the present section, $k$ is a fixed field. We let $\SpV_k$ denote the category of finite dimensional vector spaces over $k$; this is a rigid abelian linear tensor category with $\End(\TU) = k$. All $k$-algebras are understood to be commutative.%

Let $G = {\mathrm{Spec}\,A}$ be an affine group scheme over $k$, ie a group object in the category ${\mathfrak{Sch}(k)}$ of (affine) schemes over $k$ (schemes endowed with a morphism $G \to {\mathrm{Spec}\,k}$, in other words with $A$ a $k$-algebra). This means that we have morphisms of schemes: ``multiplication'' ${G \times_k G} \to G$, ``unit element'' ${\mathrm{Spec}\,k} \to G$, ``inverse'' $G \to G$ (over $k$), satisfying the usual group laws; equivalently, one is given morphisms of $k$-algebras $\Delta: A \to {A \otimes_k A}$, $\varepsilon: A \to k$ and $\sigma: A \to A$ (the comultiplication, counit and coinverse maps) such that the following axioms hold: coassociativity, coidentity%
$$%
\xymatrix{A\ar[d]^\Delta\ar[r]^-\Delta & {A \otimes A}\ar[d]^{\id \otimes \Delta} && A\ar[dr]_\iso\ar[r]^-\Delta & {A \otimes A}\ar[d]^{\varepsilon \otimes \id} \\ {A \otimes A}\ar[r]^-{\Delta \otimes \id} & {A \otimes A \otimes A} &&& {k \otimes A}}%
$$%
and coinverse%
$$%
\xymatrix{A\ar[d]^\varepsilon\ar[r]^-\Delta & {A \otimes A}\ar[d]^{(\sigma,\id)} \\ k\,\ar@{^(->}[r] & A.}%
$$%

If $A$ is a finitely generated $k$-algebra, we say that $G$ is {\em algebraic} or that it is an {\em algebraic group}. One defines a {\em coalgebra} over $k$ to be a vector space $C$ over $k$ endowed with linear maps $\Delta: C \to {C \otimes_k C}$ and $\varepsilon: C \to k$ satisfying the coassociativity and coidentity axioms. A {\em (right) comodule} over a coalgebra $C$ is a vector space $V$ over $k$ together with a linear map $\rho: V \to {V \otimes C}$ such that the following diagrams commute%
$$%
\xymatrix{V\ar[d]^\rho\ar[r]^-\rho & {V \otimes C}\ar[d]^{\rho \otimes \Delta} && V\ar[dr]_\iso\ar[r]^-\rho & {V \otimes C}\ar[d]^{\id \otimes \varepsilon} \\ {V \otimes C}\ar[r]^-{\rho \otimes \id} & {V \otimes C \otimes C} &&& {V \otimes k}}$$%
For example, $\Delta$ defines a $C$-comodule structure on $C$ itself.%

An affine group scheme $G = {\mathrm{Spec}\,A}$ can be regarded as a functor $G: {k\text{-}\mathfrak{alg}} \longto \mathfrak{groups}$ of $k$-algebras with values into groups (cf.\ also Waterhouse 1979 \cite{Wa79}):%
$$%
G(R) = \Hom_{k\text{-}\mathfrak{alg}}^{}(A,R), \qquad\qquad \text{for every $k$-algebra } R,%
$$%
so in particular, when $R = k$,%
\begin{align*}%
G(k) &= \Hom_{k\text{-}\mathfrak{alg}}^{}(A,k) \\ &= \Hom_{\mathfrak{Sch}(k)}^{}({\mathrm{Spec}\,k},G)%
\end{align*}%
is the set of closed $k$-rational points of $G$. The group structure on $G(R)$ is obtained as follows: for $s,t \in G(R)$, the product ${s \cdot t}$, the neutral element and the inverse $s^{-1}$ are respectively defined as%
$$%
A \xto{\Delta} {A \otimes_k A} \xto{s \otimes_k t} {R \otimes_k R} \xto{\text{mult.}} R,%
$$%
$$%
A \xto{\varepsilon} k \xto{\text{unit}} R,%
$$%
$$%
A \xto{\sigma} A \xto{s} R.%
$$%

Let $\Kt$ be a rigid abelian $k$-linear tensor category, and let $\omega: \Kt \longto \SpV_k$ be an exact faithful $k$-linear tensor functor. Then one can define a functor%
$$%
\underline{\Aut}^\otimes(\omega): {k\text{-}\mathfrak{alg}} \longto \mathfrak{groups},%
$$%
as follows. For $R$ a $k$-algebra, there is a canonical tensor functor $\phi_R: \SpV_k \longto \mathfrak{Mod}_R, V \mapsto {V \otimes_k R}$ into the category of $R$-modules (this is an abelian tensor category with $\End(\TU) = R$, but in general it will not be rigid because not all $R$-modules will be reflexive). If $F, G: \Kt \longto \SpV_k$ are tensor functors, then we can define $\underline{\Hom}^\otimes(F,G)$ to be the functor of $k$-algebras%
$$%
\underline{\Hom}^\otimes(F,G)(R) = \Hom^\otimes(\phi_R \circ F, \phi_R \circ G).%
$$%
Thus $\underline{\Aut}^\otimes(\omega)(R)$ consists of families $(\lambda_X), X \in \Ob(\Kt)$ where $\lambda_X$ is an $R$-linear automorphism of ${\omega(X) \otimes_k R}$ such that $\lambda_{X_1 \otimes X_2} = \lambda_{X_1} \otimes \lambda_{X_2}$, $\lambda_\TU$ is the identity mapping of $R$, and%
$$%
\xymatrix{{\omega(X) \otimes R}\ar[d]^{\omega(a) \otimes \id}\ar[r]^{\lambda_X} & {\omega(X) \otimes R}\ar[d]^{\omega(a) \otimes \id} \\ {\omega(Y) \otimes R}\ar[r]^{\lambda_Y} & {\omega(Y) \otimes R}}%
$$%
commutes for every arrow $a: X \to Y$ in $\Kt$. In the special case where $\Kt = \R{\mathnormal G;\mathnormal k}$ for some affine group scheme $G$ over $k$, and $\omega$ is the forgetful functor $\R{\mathnormal G; \mathnormal k} \longto \SpV_k$, it is clear that every element of $G(R)$ defines an element of $\underline{\Aut}^\otimes(\omega)(R)$. One has the following result%
\begin{proposizione}\label{ix1}%
The natural transformation $G \to \underline{\Aut}^\otimes(\omega)$ (of functors of $k$-algebras with values into groups) is an isomorphism.%
\end{proposizione}%
\begin{theorem}\label{O.mainthm}%
Let $\Kt$ be a rigid abelian tensor category such that $\End(\TU) = k$, and let $\omega: \Kt \longto \SpV_k$ be an exact faithful $k$-linear tensor functor. Then%
\begin{elenco}%
\item the functor $\underline\Aut^\otimes(\omega)$ of $k$-algebras is representable by an affine group scheme $G$;
\item $\omega$ defines an equivalence of tensor categories%
$$%
\Kt \longto \R{\mathnormal G; \mathnormal k}.%
$$%
\end{elenco}%
\end{theorem}%
\begin{definizione}
A \textit{neutral Tannakian category} over $k$ is a rigid abelian $k$-linear tensor category $\Kt$ for which there exists an exact faithful $k$-linear tensor functor $\omega: \Kt \longto \SpV_k$. Any such functor is said to be a {\em fibre functor} for $\Kt$.
\end{definizione}

\sezione{A Technical Lemma on Compact Groups}\label{O.1.4+.3.3}
Throughout the present section, let \SpV\ denote the complex tensor category of complex vector spaces of finite dimension (see Note~\refcnt[N.8]{N.ix5}).%

\separazione%

Let \Kt\ be an arbitrary additive complex tensor category. Let $F: \Kt \to \SpV$ be a complex tensor functor. Moreover, let $H$ be a topological group. Suppose we are given a homomorphism of monoids%
\begin{equazione}\label{O.esp79}%
\pi: H \to \End^{\otimes,\boldsymbol*}(F)\text.%
\end{equazione}%
We shall say that $\pi$ is {\em continuous} if for every object $R \in \Ob(\Kt)$ the induced representation%
\begin{equazione}\label{O.esp80}%
\pi_R: H \to \End\bigl(F(R)\bigr)%
\end{equazione}%
defined by $h \mapsto \pi_R(h) \equiv {\pi(h)}(R)$ is continuous.%

\begin{proposizione}[Technical Lemma.]\label{O.prp5}%
Let \Kt, $F$ and $H$ be as above. Suppose in addition that $H$ is a compact Lie group. Finally, let $\pi: H \to \End^{\otimes,\boldsymbol*}(F)$ be a continuous homomorphism.%

Assume the following condition holds:%
\begin{elenco}%
\item[\textsl{(*)}]for any couple of objects $R, S \in \Ob(\Kt)$ and for each homomorphism $A: F(R) \to F(S)$ of the corresponding $H$\nobreakdash-modules\inciso{in other words, for each \nC\nobreakdash-linear map $A$ such that the diagram%
\begin{equazione}\label{O.esp82}%
\begin{split}%
\xymatrix@C=40pt@R=23pt{F(R)\ar[d]^A\ar[r]^-{\pi_R(h)} & F(R)\ar[d]^A \\ F(S)\ar[r]^-{\pi_S(h)} & F(S)}%
\end{split}%
\end{equazione}%
commutes $\forall h \in H$}there is an arrow $R \xto a S$ such that $A = F(a)$.%
\end{elenco}%
Then $\pi$ is surjective; in particular, $\End^{\otimes,\boldsymbol*}(F) = \Aut^{\otimes,\boldsymbol*}(F)$ is necessarily a group.%
\end{proposizione}%
\begin{proof}%
Put $K \bydef \kernel\pi \subset H$. This is a closed normal subgroup, because it coincides with the intersection ${\bigcap \kernel{\pi_X}}$ over all objects $X$ of \Kt. On the quotient $G \bydef H/K$ there is a unique (compact) Lie group structure such that the quotient homomorphism $H \epito G$ is a Lie group homomorphism. Every $\pi_X$ can indifferently be thought of as a continuos representation of $H$ or a continuous representation of $G$, and every linear map $A: F(X) \to F(Y)$ is a morphism of $G$\nobreakdash-modules if and only if it is a morphism of $H$\nobreakdash-modules. Being continuous, every $\pi_X$ is also smooth.%

We claim there exists an object $R$ of \Kt\ such that the corresponding $\pi_R$ is faithful as a representation of $G$. This can be seen in a completely standard way, cf.\ for instance {\em Br\"ocker and tom Dieck (1985),} pp.~136\nobreakdash--137; nonetheless, in the present more abstract situation it will be useful to have a look at the argument in detail anyway. The `Noetherian' property of the compact Lie group $G$ allows us to find $X_1, \ldots, X_\ell \in \Ob(\Kt)$ with the property that%
\begin{equazione}\label{O.esp83}%
{\kernel{\pi_{X_1}} \cap \cdots \cap \kernel{\pi_{X_\ell}}} = \{e\}%
\end{equazione}%
as representations of $G$, where $e$ denotes the neutral element. Then, setting $R \bydef {X_1 \oplus \cdots \oplus X_\ell}$, the representation $\pi_R$ will be faithful because of the existence of an isomorphism of $G$\nobreakdash-modules%
\begin{equazione}\label{O.equ12}%
F({X_1 \oplus \cdots \oplus X_\ell}) \iso {F(X_1) \oplus \cdots \oplus F(X_\ell)}\text.%
\end{equazione}%
(The existence of such isomorphisms follows from the remark that a natural transformation of additive functors is additive: for instance, when $\ell=2$,%
$$%
\xymatrix@C=28pt{FX\ar[d]^{F\,i_X}\ar[rr]^-{\pi(h)(X)} & & FX\ar[d]^{F\,i_X} & {FX\oplus FY}\ar[d]^\iso\ar@{}[dl]|*{\seq}\ar[rr]^-{\pi_X(h)\,\oplus\,\pi_Y(h)} & & {FX\oplus FY}\ar[d]^\iso \\ F({X\oplus Y})\ar[rr]^-{\pi(h)({X\oplus Y})} & & F({X\oplus Y}) & F({X\oplus Y})\ar[rr]^-{\pi_{X\oplus Y}(h)} & & F({X\oplus Y}) \\ FY\ar[u]_{F\,i_Y}\ar[rr]^-{\pi(h)(Y)} & & FY\ar[u]_{F\,i_Y} &&&}%
$$%
shows that the canonical isomorphism ${F(X)\oplus F(Y)} \iso F({X\oplus Y})$ is also an isomorphism of $H$\nobreakdash-modules or, equivalently, $G$\nobreakdash-modules.)%

It follows that the $G$\nobreakdash-module $F(R)$ is a tensor generator for the complex tensor category \R{G;\nC} of continuous finite dimensional complex $G$\nobreakdash-modules. Indeed, every irreducible such $G$\nobreakdash-module $V$ embeds as a submodule of some tensor power ${F(R)^{\otimes k} \otimes \smash{\left(F(R)^*\right)}^{\otimes\ell}}$ (see for instance {\em Br\"ocker and tom Dieck, 1985}); since by assumption each $\pi(h)$ is a self-conjugate tensor preserving natural transformation, this tensor power will be naturally isomorphic to $F\left({R^{\otimes k} \otimes (R^*)^{\otimes\ell}}\right)$ as a $G$\nobreakdash-module and hence, as a consequence of the existence of the $G$\nobreakdash-module isomorphisms \refequ{O.equ12}, for each object $V$ of \R{G;\nC} there will be some object $X$ of \Kt\ such that $V$ embeds into $F(X)$ as a submodule.%

Next, consider an arbitrary natural transformation $\lambda \in \End(F)$. Let $X$ be an object of \Kt\ and let $V \subset FX$ be a submodule. The choice of a complement to $V$ in $FX$ determines a module endomorphism $P: FX \to V \into FX$ which, by condition~\textsl{(*),} comes from some endomorphism $X \xto p X \in \Kt$. Therefore%
\begin{equazione}\label{x.6}%
\begin{split}%
\xymatrix@C=35pt{FX\ar[d]^P\ar[r]^-{\lambda(X)} & FX\ar[d]^P \\ FX\ar[r]^-{\lambda(X)} & FX}%
\end{split}%
\end{equazione}%
commutes and, consequently, $\lambda(X)$ maps $V$ into itself. I will usually omit $X$ from the notation and simply write $\lambda_V: V \to V$ for the linear map that $\lambda(X)$ induces on $V$ by restriction.%

Given any other submodule $W \subset FY$ and any module homomorphism $B: V \to W$, the diagram%
\begin{equazione}\label{O.equ11}%
\begin{split}%
\xymatrix@C=35pt@R=23pt{V\ar[d]^B\ar[r]^{\lambda_V} & V\ar[d]^B \\ W\ar[r]^{\lambda_W} & W}%
\end{split}%
\end{equazione}%
is necessarily commutative. To prove this, extend the given homomorphism $B: V \to W$ to a homomorphism $B': FX \to FY$ (for instance, by choosing a complement to $V$ in $FX$ and then by taking the composite map $FX \to V \xto B W \into FY$) and then argue as before, by invoking the assumption~\textsl{(*).}%

Next, we define an isomorphism of complex algebras%
\begin{equazione}\label{x.8}%
\theta: \End(F) \isoto \End(\omega_G)%
\end{equazione}%
so that the following diagram commutes%
\begin{equazione}\label{O.equ13}%
\begin{split}%
\xymatrix{H\ar[d]_\pr\ar[r]^-\pi & \End(F)\ar[d]^\theta \\ G\ar[r]^-{\pi_G} & \End(\omega_G)\text,\!\!}%
\end{split}%
\end{equazione}%
where $\omega_G: \R{G;\nC} \to \SpV$ is the standard forgetful functor (which to any $G$\nobreakdash-module associates the underlying complex vector space) and $\pi_G(g)$ is the natural transformation $\varrho \mapsto {\pi_G(g)}(\varrho) \equiv \varrho(g)$. Given a module $V$, there exists an object $X$ of \Kt\ together with an embedding $V \into FX$, so we may define ${\theta(\lambda)}(V)$ to be the restriction of $\lambda(X)$ to $V$ (this makes sense in view of the above remarks). Of course, it is necessary to check that $\theta$ is well-defined.%

Suppose we are given two objects $X, Y \in \Ob(\Kt)$, along with $G$\nobreakdash-module embeddings of $V$ into $FX$, $FY$ respectively. Since it is always possible to embed everything equivariantly into $F({X\oplus Y})$ and since doing this does not affect the induced $\lambda_V$'s, it will be no loss of generality to assume that $X=Y$. Let $W, W' \subset FX$ be the submodules corresponding to the two different embeddings of $V$ into $FX$. Then by our remark~\refequ{O.equ11} there is a commutative diagram%
\begin{equazione}\label{x.10}%
\begin{split}%
\xymatrix@C=30pt{V\ar@{=}[d]\ar[r]^-\iso & W\ar[r]^-{\lambda_W}\ar[d]^\iso & W\ar[d]^\iso\ar[r]^-{\iso^{-1}} & V\ar@{=}[d] \\ V\ar[r]^-\iso & W'\ar[r]^-{\lambda_{W'}} & W'\ar[r]^-{\iso^{-1}} & V\text,\!\!}%
\end{split}%
\end{equazione}%
which shows that the two different embeddings precisely determine the same linear endomorphism of $V$.%

Clearly, \refequ{O.equ11} implies that $\theta(\lambda) \in \End(\omega_G)$. For $\mu \in \End(\omega_G)$ and $X \in \Kt$, put $\mu^F(X) = \mu(FX)$; then $\mu^F \in \End(F)$ and $\theta(\mu^F) = \mu$, because of the existence of embeddings $V \into FX$ and because of naturality of $\mu$: hence $\theta$ is surjective. The latter map is also injective since $\lambda(X) = \theta(\lambda)(FX)$. It is straightforward to check that the diagram \refequ{O.equ13} commutes.%

Now, to conclude the proof, it will be enough to show that $\theta$ induces a bijection between $\End^{\otimes,\boldsymbol*}(F)$ and $\End^{\otimes,\boldsymbol*}(\omega_G) = \tannakian G$, because then from \refequ{O.equ13} we get at once the following commutative square%
\begin{equazione}\label{x.11}%
\begin{split}%
\xymatrix@C=30pt@R=30pt{H\ar[d]_\pr\ar[r]^-\pi & \End^{\otimes,\boldsymbol*}(F)\ar[d]^\theta_\iso \\ G\ar[r]^-{\pi_G} & \tannakian G\text,\!\!}%
\end{split}%
\end{equazione}%
where the map on the bottom is a bijection (by the classical Tannaka duality theorem for compact groups), whence surjectivity of $\pi$ is evident.%

For instance, suppose $\lambda \in \End^\otimes(F)$ and let $V$ and $W$ be $G$\nobreakdash-modules that admit equivariant embeddings $V \into FX$ and $W \into FY$ for some $X, Y \in \Ob(\Kt)$. Since we are dealing with finite dimensional spaces, ${V\otimes W} \into {FX\otimes FY} \can F({X\otimes Y})$ will be also an embedding of $G$\nobreakdash-modules. Then, by the definition of $\theta$ and the assumption that $\lambda$ is tensor preserving, we see that the diagram%
\begin{equazione}\label{x.12}%
\begin{split}%
\xymatrix@C=40pt@R=23pt{F({X\otimes Y})\ar[r]^{\lambda({X\otimes Y})} & F({X\otimes Y}) \\ {V\otimes W}\ar@{^{(}->}[u]\ar[r]^-{\lambda_V\otimes\lambda_W} & {V\otimes W}\ar@{^{(}->}[u]}%
\end{split}%
\end{equazione}%
must commute. This shows that ${\theta(\lambda)}({V\otimes W}) = {{\theta(\lambda)}(V) \otimes {\theta(\lambda)}(W)}$. The reverse direction is straightforward.%
\end{proof}%

\separazione%

The argument that we used above in order to find the tensor generator $R$ admits the following generalization to the non-compact case. Let \Kt\ and $F$ be as in the statement of the preceding proposition.%
\begin{proposizione}\label{O.lem8}%
Let $G$ be a Lie group. Suppose that%
\begin{equazione}\label{x.13}%
\pi: G \longto \Aut(F)%
\end{equazione}%
is a faithful continuous homomorphism\nobreakdash---in other words, a continuous homomorphism such that for each $g \neq e \in G$ there exists an object $X$ in \Kt\ with $\pi_X(g) \neq \id_{FX}$.%

Then there exists an object $R \in \Ob(\Kt)$ for which $\kernel{\pi_R}$ is a discrete subgroup of $G$ or, equivalently, for which the continuous representation%
\begin{equazione}\label{x.14}%
\pi_R: G \to \GL(FR)%
\end{equazione}%
is faithful\inciso{i.e.\ injective}on some open neighbourhood of $e$.%
\end{proposizione}%
\begin{proof}%
Let $X$ be an arbitrary object of \Kt. Then $K \bydef \kernel{\pi_X}$ is a closed Lie subgroup of $G$. The connected component $K_e$ of $e$ in $K$ is also a closed Lie subgroup of $G$; in particular, the inclusion map $K_e \into G$ is an embedding of Lie groups (that is, a Lie subgroup and an embedding of manifolds). So, if $Y$ is another object, the continuous representation $\pi_Y: G \to \GL(FY)$ induces by restriction a continuous representation of $K_e$.%

The kernel $D \bydef {K_e \cap \kernel{\pi_Y}}$ is a closed Lie subgroup\inciso{in particular, a closed submanifold}of $K_e$ again. Thus, either $\dimens D < \dimens{K_e}$ or $D = K_e$, because $K_e$ is connected. Since $\pi$ is faithful, when $\dimens{K_e} > 0$ we can always find some object $Y$ such that $D \subsetneqq K_e$.%

Then it follows that for each $X \in \Ob(\Kt)$ one can always find another object $Y$ such that the submanifold $\kernel{\pi_{X\oplus Y}}$ has dimension strictly smaller than the dimension of $\kernel{\pi_X}$, unless $\dimens{\kernel{\pi_X}} = 0$. Hence an inductive argument using additivity of the category \Kt\ will yield an object $R$ such that $\dimens{\kernel{\pi_R}} = 0$ i.e.\ $\kernel{\pi_R}$ is discrete, as contended.%
\end{proof}%

\capitolo{Representation Theory Revisited}\label{3}

In the present chapter we introduce our language of {\em smooth stacks of (additive, real or complex) tensor categories,} or briefly {\em smooth (real or complex) tensor stacks.} We propose this language as the general foundational framework for the theory of representations of Lie groupoids.%

\textit{Some general conventions.} We use the expressions `{smooth}' and `{of class $\C^\infty$}' as synonyms. The capital letters $X, Y$ and $Z$ stand for manifolds of class $\C^\infty$, the corresponding lower-case letters $x, x', \ldots, y,$ etc.\ denote points on these manifolds. \smooth[X] indicates the sheaf of smooth functions on $X$ (we usually omit the subscript). Sheaves of \smooth[X]\nobreakdash-modules will also be referred to as {\em sheaves of modules over $X$.} For practical purposes, we need to consider manifolds which are possibly neither Hausdorff nor paracompact.%

\sezione{The Language of Fibred Tensor Categories}\label{N.11}
\textit{Fibred tensor categories.} \index{fibred tensor category}Fibred tensor categories will be denoted by means of capital Gothic type variables. Of course, as in \refsez{N.8}, we have to distinguish between the notions of real and complex fibred tensor category. We do the complex version; the real case is entirely analogous.

A fibred complex tensor category \stack T assigns, to each smooth manifold $X$, an additive complex tensor category%
\begin{equazione}\label{N.i1}%
\stack[X]T = \bigl(\stack[X]T,\boldsymbol\otimes_X,\TU_X,\boldsymbol*_X\bigr)%
\end{equazione}%
or $\bigl(\stack[X]T,\boldsymbol\otimes,\TU,\boldsymbol*\bigr)$ for short\inciso{omitting subscripts when they are clear from the context}and, to each smooth mapping $X \xto f Y$, a complex tensor functor%
\begin{equazione}\label{N.i2}%
f^*: \stack[Y]T \longto \stack[X]T%
\end{equazione}
called \gm{\index{pullback!along a smooth map}pull-back along $f$}. Moreover, for each pair of composable smooth maps $X \xto f Y \xto g Z$ and for each manifold $X$, any fibred complex tensor category provides self-conjugate tensor preserving natural \index{constraint}\index{fibred tensor category constraints}isomorphisms
\begin{equazione}\label{N.i3}%
\left\{\begin{aligned}%
\delta&: {f^*\circ g^*} \isoto ({g\circ f})^*%
\\%
\varepsilon&: \Id \isoto {\id_X}^*\text.%
\end{aligned}\right.%
\end{equazione}%
These are altogether required to make the following \index{coherence conditions|emph}diagrams commute
\begin{equazione}\label{N.i4}%
\begin{split}%
\xymatrix@C=17pt{f^*g^*h^*\ar[d]^(.47){\delta\cdot h^*}\ar[r]^-{f^*\delta} & f^*(hg)^*\ar[d]^(.47)\delta &&& {\id_X}^*f^*\ar[d]^(.47)\delta & f^*\ar[d]^(.47){f^*\varepsilon}\ar@{=}[dl]\ar[l]_(.34){\varepsilon\cdot f^*} \\ (gf)^*h^*\ar[r]^-\delta & (hgf)^* &&& f^* & f^*{\id_Y}^*\text{ .}\mspace{-10mu}\ar[l]_(.5)\delta}%
\end{split}%
\end{equazione}%
This is all of the mathematical data we need to introduce in order to speak about smooth tensor stacks and, later on, representations of Lie groupoids. All the required concepts can\inciso{and will}be defined in terms of the given categorical structure \stack T, i.e.\ \index{canonical isomorphism|see{constraint}}\textit{canonically.} We now explain how.

\sottosezione{Smooth tensor prestacks}%

Throughout the present subsection we let \stack P denote a fibred complex tensor category, fixed once and for all.%

\textit{Notation.}~For $i_U: U \into X$ the inclusion of an open subset, we shall put $E|_U = {i_U}^*E$ and $a|_U = {i_U}^*a$ for any object $E$ and morphism $a$ of the category \stack[X]P. (More generally, we shall adopt this abbreviation for the inclusion $i_S: S \into X$ of any {\em submanifold.})%

For any pair of objects $E, F \in {\Ob\,\stack[X]P}$, we let \index{Hom(E,F)@\sheafhom[C]XEF (sheaf hom)}\index{sheaf hom@sheaf hom \sheafhom[C]XEF}\sheafhom[P]XEF denote the presheaf of complex vector spaces over $X$ defined by
\begin{equazione}\label{N.i5}%
U \mapsto \Hom_{\stack[U]P}^{}(E|_U,F|_U)\text,%
\end{equazione}%
with the obvious restriction maps $a \mapsto {j^*a}$ corresponding to the inclusions $j: V \into U$ of open subsets. (To be precise, restriction along $j$ sends $a$ to the unique morphism $E|_V \to F|_V$ which corresponds to $j^*a$ up to the canonical isomorphisms $j^*(E|_U) \can E|_V$ and $j^*(F|_U) \can F|_V$ of \refequ{N.i3}.) Now, the requirement that \stack P be a \index{tensor prestack}\index{prestack}\textit{prestack} means exactly that any such presheaf is in fact a {\em sheaf;} in particular, it entails that one can glue any family of compatible local morphisms over $X$. Two special cases will be of particular interest to us: the sheaf $\sections E = \sheafhom[P]X\TU E$, to be referred to as the \index{Gamma(E)@\sections E, \f H (sheaf of sections)|emph}\index{sheaf sections@sheaf of sections \sections E, \f H|emph}\textit{sheaf of smooth sections} of $E \in {\Ob\,\stack[X]P}$, and the sheaf $\dual E = \sheafhom[P]XE\TU$, to be referred to as the {\em sheaf dual} of $E$. For any open subset $U$, the elements of $\sections E(U)$ will be of course referred to as the \index{section|emph}\index{smooth section|see{section}}\textit{smooth sections} of $E$ over $U$; it is perhaps useful to point out that it makes sense, for smooth sections over $U$, to take linear combinations with complex coefficients, because $\sections E(U)$ has a canonical vector space structure.

Since a morphism $a: E \to F$ in \stack[X]P yields a morphism $\sections a: \sections E \to \sections F$ of sheaves of complex vector spaces over $X$ (by composing $\TU|_U \to E|_U \xto{a|_U} F|_U$), we obtain a canonical functor%
\begin{equazione}\label{N.i6}%
\sections = \sections[X]: \stack[X]P \longto \SheavesOfModules{\sheafconst[X]\nC}\text,%
\end{equazione}%
where \sheafconst[X]\nC\ denotes the constant sheaf over $X$ of value \nC. (Note that a sheaf of complex vector spaces over a topological space $X$ is exactly the same thing as a sheaf of \sheafconst[X]\nC\nobreakdash-modules.)%

This functor is certainly linear. Moreover, there is an evident way to make it a {\em pseudo-tensor} functor of the tensor category $\bigl(\stack[X]P,\boldsymbol\otimes_X,\TU_X\bigr)$ into the category of sheaves of \sheafconst[X]\nC\nobreakdash-modules (with the standard tensor structure). In detail, a natural transformation $\tau_{E,F}: {\sections[X]E \otimes_{\sheafconst[X]\nC} \sections[X]F} \to \sections[X]{(E\otimes F)}$ arises, in the obvious manner, from the local pairings%
\begin{equazione}\label{N.i7}%
\begin{aligned}%
{\sections E(U)\mspace{5mu}\times\mspace{7.5mu}\sections F(U)} \quad&\longto \quad\sections{(E\otimes F)}(U)%
\\%
({\scriptstyle \TU|_U \xto a E|_U}\,,\;{\scriptstyle \TU|_U \xto b F|_U}) \quad&\mapsto \quad{\scriptstyle \TU|_U \mspace{8mu}\can \mspace{11mu}\TU|_U\otimes\TU|_U \mspace{8mu}\xto{a\otimes b} \mspace{11mu}{E|_U\otimes F|_U} \mspace{8mu}\can \mspace{11mu}(E\otimes F)|_U}%
\end{aligned}%
\end{equazione}%
(which are bilinear with respect to locally constant coefficients), and a morphism $\upsilon: \sheafconst[X]\nC \to \sections[X]\TU$ can be easily defined as follows%
\begin{equazione}\label{N.i8}%
\begin{aligned}%
{\scriptscriptstyle\left\{\begin{array}{c}\mathrm{locally\:constant\:complex}\\\mathrm{valued\;functions\;on}\;U\end{array}\right\}} \quad&\longto \quad\sections\TU(U)%
\\%
{t: U \to \nC} \qquad\qquad&\mapsto \quad{{t\cdot 1_U}: \TU|_U \to \TU|_U}%
\end{aligned}%
\end{equazione}%
(where $1_U = \id: \TU|_U \to \TU|_U$ is the ``unity constant section''); the operation of multiplication by $t$ in \refequ{N.i7} and \refequ{N.i8} is well-defined because $t$ is a complex constant, at least locally. It is easy to check that these morphisms of sheaves make all the diagrams in the definition of a tensor functor commute.%

Note that for $X = \pt$, where \pt\ is the one-point manifold, one has the standard identification \SheavesOfModules{\sheafconst[\pt]\nC} = \ComplexVectorSpaces\ of complex tensor categories. One may therefore regard, for $X = \pt$, the functor \refequ{N.i6} as a linear pseudo-tensor functor%
\begin{equazione}\label{N.i9}%
\stack[\pt]P \;\longto \;\ComplexVectorSpaces\text.%
\end{equazione}%
It will be convenient to have a short notation for this; making the above identification of categories explicit, we put, for all objects $E \in {\Ob\,\stack[\pt]P}$,%
\begin{equazione}\label{N.i10}%
E_* \,= \,(\sections[\pt]E)(\pt)%
\end{equazione}%
(so this is a complex vector space), and do the same for morphisms. Now, as a part of the definition of the general notion of smooth tensor stack, {\em we ask that the following condition be satisfied:} the morphism of sheaves \refequ{N.i8} is an isomorphisms for $X = \pt$. Let us record an immediate consequence of this requirement: there is a {\em canonical} isomorphism of complex vector spaces%
\begin{equazione}\label{N.i11}%
\nC \isoto \TU_*\text.%
\end{equazione}%

\begin{nota}\label{xi34}%
When dealing with the case of fibred complex tensor categories, one also has a natural morphism of sheaves of modules over $X$%
\begin{equazione}\label{xi35}%
(\sections[X]E)^* \longto \sections[X]{(E^*)}%
\end{equazione}%
defined by means of the anti-involution and the obvious related canonical isomorphisms. Since $\zeta^{**} = \zeta$ (up to canonical isomorphism), it follows at once that \refequ{xi35} is a natural isomorphism for an arbitrary complex tensor prestack; in fact, \refequ{xi35} is an isomorphism of pseudotensor functors viz.\ it is compatible\inciso{in the sense of \S\ref{O.1.2}}with the natural transformations \refequ{N.i7} and \refequ{N.i8}. Because of these considerations, we will not need to worry about complex structure in our subsequent discussion of ``axioms'' in \S\ref{N.15}.%
\end{nota}%

\textit{Notation.}~(Fibres of an object) Besides the fundamental notion of \gm{sheaf of smooth sections} we are now able to introduce a second one, that of \gm{fibre at a point}. Namely, given an object $E \in {\Ob\,\stack[X]P}$, we define the \index{fibre}\textit{fibre of $E$ at $x$} to be the finite dimensional complex vector space $E_x = (x^*E)_*$; we use the same name for the point $x$ and for the (smooth) mapping $\pt \to X, \pt \mapsto x$, so that $x^*$ is just the ordinary notation \refequ{N.i2} for the pull-back, $x^*E$ belongs to \stack[\pt]P and we can apply our notation \refequ{N.i10}. Similarly, whenever $a: E \to F$ is a morphism in \stack[X]P, we let $a_x: E_x \to F_x$ denote the linear map $(x^*a)_*$. Since $\text{-} \mapsto (\text{-})_x$ is by construction the composite of two complex pseudo-tensor functors, it may itself be regarded as a complex pseudo-tensor functor. If in particular we apply this to a local smooth section $\zeta \in \sections E(U)$ and make use of the canonical identification \refequ{N.i11}, we get, for $u$ in $U$, a linear map%
\begin{equazione}\label{N.i12}%
\nC \,\isoto \,(\TU_{\pt})_* \,\can \,({u^*\,\TU|_U})_* \,\xto{(u^*\zeta)_*} \,({u^*\,E|_U})_* \,\can \,(u^*E)_* \,= \,E_u\text,%
\end{equazione}%
which corresponds to a vector $\zeta(u) \in E_u$ (the image of the unity $1 \in \nC$) to be called the \index{value}\textit{value} of $\zeta$ at $u$. One has the intuitive formula
\begin{equazione}\label{N.i13}%
{a_u\cdot\zeta(u)} = [\sections a(U)\zeta](u)\text.%
\end{equazione}%
Notice also that the vectors ${\zeta(u) \otimes \eta(u)}$ and $({\zeta\otimes\eta})(u)$ correspond to one another in the canonical linear map ${E_u\otimes F_u} \to ({E\otimes F})_u$ (we may state this loosely by saying they are equal).%

We have not explained yet what we mean when we say that a tensor prestack is \gm{smooth}. This was not necessary before because all we have said so far does not depend on that specific property. However, from this precise moment we begin to develop systematically concepts which, even in order to be defined, presuppose the smoothness of the tensor prestack, so it becomes necessary to fill the gap.%

Consider the tensor unit $\TU \in {\Ob\,\stack[X]P}$ and let $x$ be any point. There is a canonical isomorphism $\nC \can \TU_x$ analogous to \refequ{N.i11}, namely the composite $\nC \can (\TU_{\pt})_* \can ({x^*\TU})_* = \TU_x$. This identification allows us to define a {\em canonical} homomorphism of complex algebras%
\begin{equazione}\label{N.i14}%
\End_{\stack[X]P}(\TU) \;\longto \;\left\{\mathrm{functions\:}X\to\nC\right\}\text, \quad e \mapsto \tilde e%
\end{equazione}%
by putting $\tilde e(x)$ = the complex constant such that the linear map ``scalar multiplication by $\tilde e(x)$'' (of \nC\ into itself) corresponds to $e_x: \TU_x \to \TU_x$ under the linear isomorphism $\nC \can \TU_x$. We shall say that the tensor prestack \stack P is \index{smooth tensor prestack}\index{tensor prestack!smooth|emph}\textit{smooth} if the homomorphism \refequ{N.i14} determines a {\em one-to-one} correspondence onto the subalgebra of {\em smooth} functions on $X$
\begin{equazione}\label{N.i15}%
\End_{\stack[X]P}(\TU) \can \C^\infty(X)\text.%
\end{equazione}%

A first consequence of the smoothness of \stack P is the possibility to endow each space $\Hom_{\stack[X]P}(E,F)$ with a $\C^\infty(X)$\nobreakdash-module structure, canonical and compatible with the already defined operation of multiplication by locally constant functions. Indeed, the natural action%
\begin{equazione}\label{N.i16}%
\begin{aligned}%
{\End_{\stack[X]P}(\TU) \mspace{5mu}\times \mspace{5mu}\Hom_{\stack[X]P}(E,F)} \mspace{+5mu}\longto \mspace{+5mu}\Hom_{\stack[X]P}(E,F)\text,%
\\%
(e,a) \mspace{+30mu}\mapsto \mspace{+30mu}{\scriptstyle E \mspace{10mu}\can \mspace{10mu}{\TU\otimes E} \mspace{10mu}\xto{e\otimes a} \mspace{10mu}{\TU\otimes F} \mspace{10mu}\can \mspace{10mu}F} \mspace{+20mu}%
\end{aligned}%
\end{equazione}%
turns $\Hom_{\stack[X]P}(E,F)$ into a left $\End_{\stack[X]P}(\TU)$\nobreakdash-module, hence we can use the identification of \nC\nobreakdash-algebras \refequ{N.i15} to make $\Hom_{\stack[X]P}(E,F)$ a $\C^\infty(X)$\nobreakdash-module; in short, the module multiplication can be written as $(\tilde e,a) \mapsto {e\otimes a}$.%

Accordingly, $\sheafhom[P]XEF(U) = \Hom_{\stack[U]P}(E|_U,F|_U)$ inherits a canonical structure of $\C^\infty(U)$\nobreakdash-module, for each open subset, and one verifies at once that this makes \sheafhom[P]XEF a sheaf of \smooth[X]\nobreakdash-modules. Of course, the remark applies in particular to any sheaf of {\em `smooth'} sections \sections[X]E, partly justifying the terminology; moreover, one readily sees that any morphism $a: E \to F$ in the category \stack[X]P induces a morphism $\sections[X]a: \sections[X]E \to \sections[X]F$ of sheaves of \smooth[X]\nobreakdash-modules. So we get a $\C^\infty(X)$\nobreakdash-linear functor%
\begin{equazione}\label{N.i17}%
\stack[X]P \;\longto \;\SheavesOfModules{\smooth[X]}\text,%
\end{equazione}%
still denoted by \sections[X]{}. (Notice that both categories have \Hom\nobreakdash-sets enriched with a $\C^\infty(X)$\nobreakdash-module structure\footnote{Such that the composition of morphisms is $\C^\infty(X)$\nobreakdash-bilinear.}. The $\C^\infty(X)$\nobreakdash-linearity of the functor amounts by definition to the $\C^\infty(X)$\nobreakdash-linearity of all the maps%
\begin{displaymath}%
\Hom_{\stack[X]P}(E,F) \;\to \;\Hom_{\smooth[X]}(\sections[X]E,\sections[X]F)\text, \quad a \mapsto \sections[X]a\text{. )}%
\end{displaymath}%

If one also takes into account the tensor structure then the process of \gm{upgrading} the functor \refequ{N.i6} can be pursued further by observing that the operations described in \refequ{N.i7}, \refequ{N.i8} may now be used to define morphisms of sheaves of \smooth[X]\nobreakdash-modules%
\begin{equazione}\label{N.i18}%
\left\{\begin{aligned}%
\tau&: {\sections[X]E \otimes_{\smooth[X]} \sections[X]F} \to \sections[X]{(E\otimes F)}\text,%
\\%
\upsilon&: \smooth[X] \to \sections[X]\TU\text;%
\end{aligned}\right.%
\end{equazione}%
the morphism $\tau = \tau_{E,F}$ is natural in the variables $E,F$ and, along with $\upsilon$, makes \refequ{N.i17} a {\em pseudo-tensor} functor of the tensor category \stack[X]P into the tensor category of sheaves of \smooth[X]\nobreakdash-modules. This is closer than \refequ{N.i6} to being a {\em tensor} functor, in that the morphism $\upsilon$ is evidently an {\em isomorphism} of sheaves of \smooth[X]\nobreakdash-modules.%

Consider next a smooth mapping of manifolds $f: X \to Y$. Suppose that $U \subset X$ and $V \subset Y$ are open subsets with $f(U) \subset V$, and let $f_U$ denote the induced mapping of $U$ into $V$. For any object $F$ of the category \stack[Y]P, we obtain a correspondence of local smooth sections%
\begin{equazione}\label{N.i19}%
(\sections[Y]F)(V) \;\longto \;\sections[X]{(f^*F)}(U)\text, \quad \eta \mapsto {\eta\circ f}%
\end{equazione}%
by putting ${\eta\circ f}$ equal by definition to the composite%
\begin{equazione}\label{N.i20}%
\TU|_U \,\can \,{(f^*\TU)}|_U \,\can \,f_U^*(\TU|_V) \,\xto{f_U^*(\eta)} \,f_U^*(F|_V) \,\can \,{(f^*F)}|_U\text.%
\end{equazione}%
One easily verifies that for $U$ fixed and $V$ variable, the maps \refequ{N.i19} form an inductive system indexed over the inclusions of neighbourhoods $V \supset V' \supset f(U)$, and that eventually they induce a morphism of sheaves of \smooth[X]\nobreakdash-modules
\begin{equazione}\label{N.i21}%
f^*(\sections[Y]F) \longto \sections[X]{(f^*F)}\text,%
\end{equazione}%
where $f^*(\sections[Y]F)$ is the ordinary pull-back in the sense of sheaves of modules over smooth manifolds. It is also clear that the morphism \refequ{N.i21} is natural in $F$, and also a morphism of pseudo-tensor functors (in other words, it is tensor preserving). To conclude, let us give some motivation for the notation \gm{${\eta \circ f}$}. There is an obvious canonical isomorphism of vector spaces%
\begin{equazione}\label{N.i22}%
(f^*F)_x = (x^*f^*F)_* \can ({f(x)}^*F)_* = F_{f(x)}\text.%
\end{equazione}%
Now, we have the two vectors $\eta(f(x)) \in F_{f(x)}$ and $({\eta\circ f})(x) \in (f^*F)_x$, and you can easily check that they correspond to one another in the above isomorphism. We can state this loosely as%
\begin{equazione}\label{N.i23}%
({\eta\circ f})(x) = \eta(f(x))\text.%
\end{equazione}%
The last expression evidently justifies our notation.%

\sezione{Smooth Tensor Stacks}\label{N.12}
It will be convenient to regard the open coverings of a manifold $X$ as smooth mappings onto $X$. This can be made precise as follows. Borrowing some standard terminology from algebraic geometers, we shall say that a smooth mapping $p: X' \to X$ is \index{flat map}\textit{flat,} if it is surjective and it restricts to an open embedding $p_{U'}: U' \into X$ on each connected component $U'$ of $X'$; we may think of $p$ as codifying a certain open covering of $X$, indexed by the set of connected components of $X'$. A \index{refinement}\textit{refinement} of $X' \xto p X$ will be obtained by composing $p$ backwards with another flat mapping $X'' \xto{p'} X'$. The fundamental property of flat mappings is that they can be pulled back, preserving flatness, along any smooth map: precisely, for any $Y \xto f X$ the pull-back%
\begin{equazione}\label{N.i24}%
{Y\times_X^{}X'} \:= \:\bigl\{(y,x'): f(y) = p(x')\bigr\}%
\end{equazione}%
will make sense in the category of $\C^\infty$\nobreakdash-manifolds and the first projection $\pr_1: {Y\times_X^{}X'} \to Y$ will be a flat mapping. Particularly relevant is the case where $f$ is also a flat mapping, leading to the ``standard'' common refinement for $f$ and $p$.%

Some standard abbreviations. For any flat mapping $p: X' \to X$, let%
\begin{equazione}\label{N.i25}%
X'' \:= \:{X'\times_{\mspace{-3mu}X}X'} \:= \:\bigl\{(x'_1,x'_2): p(x'_1) = p(x'_2)\bigr\}\text,%
\end{equazione}%
with the two projections $p_1, p_2: X'' \to X'$; and the triple fibred product%
\begin{equazione}\label{N.i26}%
X''' \:= \:{X'\times_{\mspace{-3mu}X}X'\times_{\mspace{-3mu}X}X'} \:= \:\bigl\{(x'_1,x'_2,x'_3): p(x'_1) = p(x'_2) = p(x'_3)\bigr\}%
\end{equazione}%
with its projections $p_{12}, p_{23}, p_{13}: X''' \to X''$ resp.\ given by $(x'_1,x'_2,x'_3) \mapsto (x'_1,x'_2)$ and so forth.%

A \index{descent datum}\textit{descent datum} for a smooth complex tensor prestack \stack P, {\em relative to} the flat mapping $p: X' \to X$, will be a pair $(E',\theta)$ consisting of an object $E' \in \stack[X']P$ and an isomorphism $\theta: {{p_1}^*E'} \isoto {{p_2}^*E'}$ in \stack[X'']P, such that ${p_{13}}^*(\theta) = {{p_{12}}^*(\theta) \circ {p_{23}}^*(\theta)}$ up to the canonical isos. A morphism of descent data, let us say of $(E',\theta)$ into $(F',\xi)$, will be a morphism $a': E' \to F'$ in \stack[X']P compatible with $\theta$ and $\xi$ in the sense that ${{p_2}^*(a') \circ \theta} = {\xi\circ {p_1}^*(a')}$. Descent data of type \stack P and relative to $X' \xto p X$ (and their morphisms) form a category \index{Des(X'/X)@\Des[C]{X'/X} (category of descent data)}\Des[P]{X'/X}. There is an obvious functor
\begin{equazione}\label{N.i31}%
\stack[X]P \longto \Des[P]{X'/X}\text, \quad E \mapsto ({p^*E},\phi_E)\text{, } a \mapsto {p^*a}%
\end{equazione}%
defined by letting $\phi_E$ be the canonical isomorphism ${p_1}^*({p^*E}) \can {({p\circ p_1})^*E} = {({p\circ p_2})^*E} \can {p_2}^*({p^*E})$. Whenever the functor \refequ{N.i31} is an equivalence of categories for every flat mapping of manifolds $p: X' \to X$, one says that the prestack \stack P is a \index{smooth tensor stack}\index{tensor stack}\index{stack}\textit{stack.}

%\begin{esempio}[smooth complex vector bundles]\label{N.i32}%
%\end{esempio}%

\begin{nota}\label{N.i33}
Depending on one's purposes, the condition that the functors \refequ{N.i31} be equivalences of categories for all flat mappings $X' \to X$ can be weakened to some extent. For example, one could ask it to be satisfied just for all flat $X' \to X$ over a Hausdorff, paracompact $X$. In fact, the latter condition will prove to be sufficient for all our purposes: no relevant aspect of the theory seems to depend on the stronger requirement. We propose to use the term \gm{\index{smooth tensor parastack}\index{tensor parastack}\index{parastack}parastack} for the weaker notion; we will often be sloppy and use `stack' as a synonym to `parastack'.
\end{nota}

\sottosezione{Locally trivial objects}%

Let \stack S be any smooth tensor prestack. An object $E \in {\Ob\,\stack[X]S}$ will be called \index{trivial object}\textit{trivial} if there exists some $V \in {\Ob\,\stack[\pt]S}$ for which one can find an isomorphism $E \stackrel\alpha\iso {{c_X}^*V}$ in \stack[X]S, where $c_X: X \to \pt$ denotes the collapse map. Any such pair $(V,\alpha)$ will be said to constitute a \index{trivialization}\textit{trivialization} of $E$.

For an arbitrary manifold $X$, let \index{V(X)@\V[C]X (subcategory of locally trivial objects)}\V[S]X denote the full subcategory of \stack[X]S formed by the \index{locally trivial object}locally trivial objects of locally finite rank; more explicitly, $E \in {\Ob\,\stack[X]S}$ will be an object of \V[S]X provided one can cover $X$ with open subsets $U$ such that $E|_U$ trivializes in \stack[U]S by means of a trivialization of the form $({\TU\oplus\cdots\oplus\TU},\alpha)$ or, equivalently, such that in \stack[U]S there exists an isomorphism $E|_U \iso {\TU_U\oplus\cdots\oplus\TU_U}$. It follows at once from the bilinearity of $\Tensor$, the triviality of \TU\ and the linearity of $f^*$ that the operation $X \mapsto \V[S]X$ determines a fibred (additive, complex) tensor subcategory of \stack S. Hence $X \mapsto \V[S]X$ inherits a fibred tensor structure from \stack S. It is easy to see that one gets in fact a smooth tensor prestack $\mathit V^{\stack S}$; moreover, it is obvious that $\mathit V^{\stack S}$ is a parastack resp.\ a stack if such is \stack S.

The complex tensor category \V[S]X very closely relates to that of smooth complex vector bundles over $X$. Let us make this precise. Clearly, every object $E \in \V[S]X$ yields a smooth complex vector bundle over $X$: just put $\tilde E = \{(x,e): x\in X, e \in E_x\}$, with the local trivializing charts obtained from local trivializations $E|_U \stackrel\alpha\iso {\TU_U\oplus\cdots\oplus\TU_U}$, $\alpha = (\alpha_1,\ldots,\alpha_d)$ by setting $(u,e) = \bigl(u;\alpha_{1,u}(e),\ldots,\alpha_{d,u}(e)\bigr) \in {U\times\nC^d}$. Since any morphism $a: E \to E'$ in \V[S]X can be locally described in terms of ``matrix expressions'' with smooth coefficients, setting ${\tilde a\cdot(x,e)} = (x,{a_x\cdot e})$ defines a morphism of smooth vector bundles $\tilde a: \tilde E \to \tilde{E'}$. It is an exercise to show that the assignment $E \mapsto \tilde E$ defines a faithful complex tensor functor of \V[S]X into smooth complex vector bundles. Under extremely mild hypotheses, this functor will actually prove to be an equivalence of complex tensor categories; this will happen, for example, when \stack S is a parastack and $X$ is paracompact, or when \stack S is stack.%

In conclusion, we see that for \stack S a smooth tensor (para-)stack (and $X$ a reasonable manifold), the category \stack[X]S will essentially include\inciso{as a full tensor subcategory}all smooth vector bundles over $X$. One arrives at the same results, alternatively, by considering the functor \sections[X]\ and the category of locally free sheaves of \smooth[X]\nobreakdash-modules of locally finite rank. This last remark can be summarized in the diagram%
\begin{equazione}\label{xii.33}%
\begin{split}%
\xymatrix@M=5pt@R=11pt{\V[S]X\ar@{^{(}->}[dr]!UL_(.7){\sections[X]{}}\ar[rr]^-*{\text-\mapsto\widetilde{(\text-)}}_-*{\simeq} & & \stack[X]{V^\infty}\ar@{_{(}->}[dl]!UR^(.7){\sections[X]{}} \\ & \SheavesOfModules{\smooth[X]} & }%
\end{split}%
\end{equazione}%
(commutative up to canonical natural isomorphism). The smooth tensor stack \stack{V^\infty} is therefore, in a very precise sense, the ``smallest'' possible.%

\sezione{Foundations of Representation Theory}\label{N.13}
We develop our theory of representations relative to a \gm{\index{type}type}. This can be any smooth complex tensor parastack \stack S, in the sense of Note~\refcnt[N.12]{N.i33}. Once a type \stack S has been fixed, one can associate to any Lie groupoid a mathematical object called \gm{fibre functor}.

\sloppy
This is done as follows. Let \G\ be a Lie groupoid, let us say, with base $M$. We are going to construct a category \index{R(G)@\R[T]{\G} (category of type \stack T representations)}\R[S]\G, along with a functor \forget[S]\G\ of \R[S]\G\ into \stack[M]S that we shall call the \gm{standard fibre functor} of \G\ {(of type \stack S).} An object of the category \R\G\ = \R[S]\G\ (every time we like we can omit writing the type \stack S, as this is fixed) is defined to be a \index{representation!type T@of type \stack T}pair $(E,\varrho)$ with $E$ an object of \stack[M]S and $\varrho$ a morphism in \stack[\G]S
\begin{equazione}\label{N.ii1}
\varrho: \s^*E \to \t^*E%
\end{equazione}%
(where $\s, \t: \G \to M$ denote the source, resp.\ target map of \G), such that the appropriate conditions for $\varrho$ to be an action\inciso{in other words, for it to be compatible with the groupoid structure}are satisfied, namely:%
\begin{elenco}%
\item$\corners{\u^*\varrho} = \id_E$, where $\u: M \to \G$ denotes the unit section. (Here and in the sequel we adopt the device of putting corners around a morphism to indicate the morphism\inciso{which one, will always be clear from the context}that corresponds to it up to some canonical identifications; for instance, the last equality, spelled out explicitly, means that the diagram%
\begin{equazione}\label{N.equR2}%
\begin{split}%
\xymatrix@C=13pt@R=19pt{{\u^*\s^*E}\ar[dr]_-\can^-{\text{can.}}\ar[rr]^{\u^*\varrho} && {\u^*\t^*E}\ar[dl]_-\can^-{\text{can.}} \\ & E &}%
\end{split}%
\end{equazione}%
commutes, where we use the identifications ${\u^*\s^*E} \can {({\s\circ\u})^*E} = {{\id_M}^*E} \can E$ etc.\ provided by the fibred tensor structure constraints associated with \stack S);%
\item if we let $\mca[2]\G = \mca\G$ denote the manifold of composable arrows of \G, $\c: \mca[2]\G \to \G\text,\, (g',g) \mapsto g'g$ the composition of arrows and $\p_0, \p_1: \mca[2]\G \to \G$ the two projections $(g',g) \mapsto g'\text,\, \mapsto g$ onto the first and second factor respectively, we have the identity $\corners{\c^*\varrho} = {\corners{\p_0^*\varrho} \cdot \corners{\p_1^*\varrho}}$; that is to say, according to our convention, we have the commutativity of the following diagram in the category \stack[{\mca[2]\G}]S:%
\begin{equazione}\label{N.equR3}%
\begin{split}%
\xymatrix@C=37pt@R=17pt{ & \save[]+<-12pt,0pt>*+<8pt>{\c^*\s^*E}!<+24pt,0pt>\ar@{=}[dl]!<+16pt,+10pt>\ar[rrr]!<-12pt,0pt>^-{\c^*\varrho}\restore & & & \save[]+<+12pt,0pt>*+<8pt>{\c^*\t^*E}!<-24pt,0pt>\ar@{=}[dr]!<-16pt,+10pt>\restore & \\ \save[]*+<8pt>{{\p_1}^*\s^*E}!<-4pt,0pt>\ar[ddrr]!<-16pt,+8pt>_-{{\p_1}^*\varrho}\restore & & & & & \save[]*+<8pt>{{\p_0}^*\t^*E}!<+4pt,0pt>\ar@{<-}[ddll]!<+16pt,+8pt>^-{{\p_0}^*\varrho}\restore \\ & & & & & \\ & & \save[]+<-6pt,0pt>*+<8pt>{{\p_1}^*\t^*E}\ar@{=}[r]!<-16pt,0pt>\restore & \save[]+<+6pt,0pt>*+<8pt>{{\p_0}^*\s^*E}\restore & & }%
\end{split}%
\end{equazione}%
(which involves the canonical identifications ${\c^*\s^*E} \can {({\s\circ\c})^*E} = {({\s\circ\p_1})^*E} \can {{\p_1}^*\s^*E}$ etc.\ provided by the structure constraints of \stack S). We shall also write \gm{${\c^*\varrho} = {{\p_0^*\varrho} \cdot {\p_1^*\varrho}}$~(mod~$\can$)}.%
\end{elenco}%
This concludes the description of the objects of \R[S]\G; we shall call them {\em representations} of \G, or \index{action|see{representation}}\textit{\G-actions} (of type \stack S). As morphisms of \G\nobreakdash-actions $a: (E,\varrho) \to (E',\varrho')$ we take all those morphisms $a: E \to E'$ in \stack[M]S which make the following square commutative
\begin{equazione}\label{N.ii4}%
\begin{split}%
\xymatrix@C=23pt@R=23pt{{\s^*E}\ar[d]^(.49){\s^*a}\ar[r]^\varrho & {\t^*E}\ar[d]^(.49){\t^*a} \\ {\s^*E'}\ar[r]^{\varrho'} & {\t^*E'}\text{~.}\!\!\!\!}%
\end{split}%
\end{equazione}
We endow the category \R[S]\G\ with the linear structure of \stack[M]S. Then the forgetful functor
\begin{equazione}\label{N.ii5}%
\forget[S]\G: \R[S]\G \longto \stack[M]S\text, \quad (E,\varrho) \mapsto E%
\end{equazione}%
is linear and faithful. We call it the \index{omega(G)@\forget[T]\G, \forget[\infty]{\G} (forgetful functor)|emph}\index{standard fibre functor!type T@(of type \stack T) \forget[T]\G|emph}\textit{standard fibre functor} of \G\ \textit{(of type \stack S).}

\fussy
Observe that the linear category \R[S]\G\ is {\em additive.} Indeed, fix any objects $R_0, R_1 \in \R\G$, let us say $R_i = (E_i,\varrho_i)$, and choose a representative $E_0 \stackrel{i_0}\into {E_0\oplus E_1} \stackrel{i_1}\infrom E_1$ for the direct sum in \stack[M]S. Then, since the linear functors $\s^*$, $\t^*$ have to preserve direct sums (cf.\ \textit{MacLane (1998),} p.~197), there will be a unique `universal' isomorphism in \stack[\G]S%
$$%
\s^*({E_0\oplus E_1}) = {{\s^*E_0} \oplus {\s^*E_1}} \xto{\:\:\varrho_0\,\oplus\,\varrho_1\:\:} {{\t^*E_0} \oplus {\t^*E_1}} = \t^*({E_0\oplus E_1})\text.%
$$%
One checks that the pair ${R_0\oplus R_1} = ({E_0\oplus E_1},{\varrho_0\oplus\varrho_1})$ is a \G\nobreakdash-action, that $R_0 \stackrel{i_0}\into {R_0\oplus R_1} \stackrel{i_1}\infrom R_1$ are morphisms of \G\nobreakdash-actions, and that they yield a direct sum in \R\G. The process to obtain a null representation is entirely analogous, starting from a null object in \stack[M]S.%

\begin{lemma}\label{N.ii7}%
For an arbitrary \G\nobreakdash-action $(E,\varrho) \in \R[S]\G$, the morphism $\varrho: \s^*E \to \t^*E$ is necessarily an isomorphism in \stack[\G]S.%
\end{lemma}%
\begin{proof}%
Let \Kt\ be any category. Define two arrows $a, a'$ to be `equivalent', and write $a \thicksim a'$, if they are isomorphic as objects of the arrow category $\mathit{Ar}(\Kt)$ (in other words, if there exist isomorphisms between their domains and codomains which transform the one arrow into the other). Then the following assertions hold: \textsl{a)}~for any functor $F: \Kt[C] \to \Kt[D]$, $a \thicksim a'$ implies $Fa \thicksim Fa'$; \textsl{b)}~the existence of a natural iso $F \isoto G$ implies $Fa \thicksim Ga$ for every $a$; \textsl{c)}~if $a \thicksim a'$ and $a$ is left (resp.\ right) invertible, then the same is true of $a'$; \textsl{d)}~$ba \thicksim \id$ implies that $a$ is left invertible and $b$ right invertible.%

Let $\i: \G \to \G\text,\, g \mapsto g^{-1}$ be the inverse, and consider the two maps $(\i,\id), (\id,\i): \G \to \mca[2]\G$ given by $g \mapsto (g^{-1},g)\text,\, \mapsto (g,g^{-1})$ respectively. Then one has the following equivalences of arrows in the category \stack[\G]S%
\begin{align*}%
\id_{\s^*E} &= {\s^*\id_E} \stackrel{\text{\tiny\textsl{a)}}}\thicksim {\s^*\u^*\varrho} \stackrel{\text{\tiny\textsl{b)}}}\thicksim {({\u\circ\s})^*\varrho} = {[{\c \circ (\i,\id)}]^*\varrho}%
\\%
&\stackrel{\text{\tiny\textsl{b)}}}\thicksim {(\i,\id)^*\c^*\varrho} \stackrel{\text{\tiny\textsl{a)}}}\thicksim {(\i,\id)^*\corners{\c^*\varrho}} \stackrel{\text{\tiny(\ref{N.equR3})}} = {(\i,\id)^*({\corners{\p_0^*\varrho} \cdot \corners{\p_1^*\varrho}})}%
\\%
&= {{(\i,\id)^*\corners{\p_0^*\varrho}} \cdot {(\i,\id)^*\corners{\p_1^*\varrho}}}\text,%
\end{align*}%
hence ${(\i,\id)^*\corners{\p_1^*\varrho}}$ is left invertible in \stack[\G]S, by \textsl{d)}. Since this is in turn equivalent to ${(\i,\id)^*{\p_1}^*\varrho} \thicksim [{\p_1 \circ (\i,\id)}]^*\varrho = {{\id_{\G}}^*\varrho} \thicksim \varrho$, $\varrho$ itself will be left invertible in \stack[\G]S, by \textsl{c)}. An analogous reasoning will establish the right invertibility of $\varrho$. It follows that $\varrho$ is invertible.%
\end{proof}%

Next, we discuss the standard \index{tensor category}\index{complex tensor category}tensor structure on the category \R\G. This structure makes \R\G\ an additive linear tensor category. The standard fibre functor $\fifu = \forget\G$ turns out to be a {\em strict} \index{tensor functor}\index{complex tensor functor}tensor functor of \R\G\ into \stack[M]S, in the sense that the identities $\fifu({R\otimes S}) = {\fifu(R) \otimes \fifu(S)}$ and $\fifu(\TU) = \TU$ hold, so that they can be taken respectively as the natural constraints $\tau$ and $\upsilon$ in the definition of tensor functor.

We start with the construction of the bifunctor $\otimes: {\R\G \times \R\G} \to \R\G$. For two arbitrary representations $R, S \in \R\G$, let us say $R = (E,\varrho)$ and $S = (F,\repsigma)$, we put ${R\otimes S} = ({E\otimes F},\corners{\varrho\otimes\repsigma})$, where\inciso{following the usual convention}$\corners{\varrho\otimes\repsigma}$ stands for the composite morphism%
\begin{equazione}\label{N.ii10}%
\s^*({E\otimes F}) \can {{\s^*E} \otimes {\s^*F}} \xto{\:\:\varrho\,\otimes\,\repsigma\:\:} {{\t^*E} \otimes {\t^*F}} \can \t^*({E\otimes F})\text.%
\end{equazione}%
It is easy to recognize that the pair ${R\otimes S}$ is itself a \G\nobreakdash-action, i.e.\ an object of the category \R\G; moreover, if $(E,\varrho) \xto a (E',\varrho')$ and $(F,\repsigma) \xto b (F',\repsigma')$ are morphisms in \R\G\ then so is ${a\otimes b}: {R\otimes S} \to {R'\otimes S'}$.%

We define the tensor unit of \R\G\ to be the pair $(\TU_M,\corners\id)$, where $\TU_M$ the tensor unit of \stack[M]S and $\corners\id$ is the composite canonical isomorphism%
\begin{equazione}\label{N.ii11}%
{\s^*\TU_M} \can \TU_{\G} \can {\t^*\TU_M}\text.%
\end{equazione}%

The \textit{ACU} natural \index{ACU constraints@\textit{ACU} constraints}\index{constraint}constraints $\alpha\text,\, \gamma\text,\, \lambda\text,\, \rho$ for the tensor structure of the base category \stack[M]S will provide analogous constraints for the tensor product we just introduced on \R\G. (For example, consider three representations $R, S, T \in \R\G$ and let $E, F, G \in \stack[M]S$ be the respective supports; then the isomorphism $\alpha_{E,F,G}: {E \otimes ({F\otimes G})} \isoto {({E\otimes F}) \otimes G}$ is also an isomorphism $\alpha_{R,S,T}: {R \otimes ({S\otimes T})} \isoto {({R\otimes S}) \otimes T}$ in \R\G.) A~fortiori, the coherence diagrams for such `inherited' constraints will commute.

\sezione{Homomorphisms and Morita Invariance}\label{N.14}
We now proceed to study the operation of taking the inverse image of a representation along a \index{homomorphism of groupoids}homomorphism of Lie groupoids. Then we concentrate on the special case of Morita equivalences; in order to give a satisfactory treatment of these, it will be necessary to analyze natural transformations of Lie groupoid homomorphisms first.

Let $\varphi: \G \to \H$ be a homomorphism of Lie groupoids and let $M \xto f N$ be the smooth map induced by $\varphi$ on the base manifolds.%

Suppose $(F,\sigma) \in \R[S]\H$. Consider the morphism\inciso{which we also denote by ${\varphi^*\sigma}$, slightly abusing notation}defined as follows:%
\begin{equazione}\label{N.ii15}%
{\s_{\G}}^*({f^*F}) \can {\varphi^*{\s_{\H}}^*F} \xto{\:\varphi^*\sigma\:} {\varphi^*{\t_{\H}}^*F} \can {\t_{\G}}^*({f^*F})\text;%
\end{equazione}%
the equalities ${f\circ\s_{\G}} = {\s_{\H}\circ\varphi}$ etc.\ account, of course, for the existence of the canonical isomorphisms occurring in \refequ{N.ii15}. It is straightforward to check that the pair $({f^*F},{\varphi^*\sigma})$ constitutes an object of the category \R[S]\G\ and that if $(F,\sigma) \xto b (F',\sigma')$ is a morphism of \H\nobreakdash-actions then ${f^*b}$ is a morphism of $({f^*F},{\varphi^*\sigma})$ into $({f^*F'},{\varphi^*\sigma'})$ in \R[S]\G. Hence we get a functor%
\begin{equazione}\label{N.ii16}%
\varphi^*: \R[S]\H \longto \R[S]\G\text,%
\end{equazione}%
which we agree to call the \index{inverse image}\textit{inverse image} or \index{pullback!of representations}\textit{pull-back} (of representations) \textit{along $\varphi$.}

It is fairly easy to check that the \index{tensor functor constraints}\index{constraint}constraints
\begin{equazione}\label{N.ii17}%
\left\{\begin{aligned}%
\upsilon&: \TU_M \isoto {f^*\TU_N}%
\\%
\tau_{F,F'}&: {{f^*F}\otimes{f^*F'}} \isoto f^*({F\otimes F'})\text,%
\end{aligned}\right.%
\end{equazione}%
associated with the tensor functor $f^*$, can also function as isomorphisms of \G\nobreakdash-actions, $\upsilon: \TU \isoto \varphi^*(\TU)$ and $\tau_{S,S'}: {\varphi^*(S) \otimes \varphi^*(S')} \isoto \varphi^*({S\otimes S'})$, for all $S, S' \in \R\H$ with, let us say, $S = (F,\sigma)$ and $S' = (F',\sigma')$. A fortiori, these isomorphisms are natural and they provide appropriate tensor functor constraints for $\varphi^*$, thus making $\varphi^*$ a \index{tensor functor}tensor functor of the tensor category \R\H\ into the tensor category \R\G.%

Let $\G \xto\varphi \H \xto\psi \K$ be two composable homomorphisms of Lie groupoids and let $X \xto{\varphi_0} Y \xto{\psi_0} Z$ denote the respective maps on bases. Note that for an arbitrary action $T = (G,\tau) \in \R\K$ the canonical isomorphism ${{\varphi_0}^*\smash{\psi_0}^*G} \can {({\psi_0\circ\varphi_0})^*G} = {\smash{({\psi\circ\varphi})_0}^*G}$ is actually a morphism $\varphi^*({\psi^*T}) \isoto {({\psi\circ\varphi})^*T}$ in the category \R\G. Hence we get an isomorphism of tensor functors%
\begin{equazione}\label{xiv.4}%
{\varphi^* \circ \psi^*} \xto\can ({\psi\circ\varphi})^*\text.%
\end{equazione}%

It is worthwhile remarking that $\varphi^*$ fits in the following diagram%
\begin{equazione}\label{N.ii18}%
\begin{split}%
\xymatrix@C=30pt{\R[S]\H\ar[d]_{\xforget[S]\H}\ar[r]^-{\varphi^*} & \R[S]\G\ar[d]^{\xforget[S]\G} \\ \stack[N]S\ar[r]^-{f^*} & \stack[M]S\text,\!\!}%
\end{split}%
\end{equazione}%
whose commutativity is to be interpreted as an equality of composite {\em tensor functors}\nobreakdash---thus, involving also the constraints.%

The notion from Lie groupoid theory we want to dualize next is that of natural transformation; this comes about especially when one considers Morita equivalences, as we shall see soon. Recall that a \index{transformation}\textit{transformation} $\tau: \varphi_0 \to \varphi_1$ (between two Lie groupoid homomorphisms $\varphi_0, \varphi_1: \G \to \H$) is a smooth mapping $\tau$ of the base manifold $M$ of \G\ into the manifold of arrows of \H, such that $\tau(x): f_0(x) \to f_1(x)$ $\forall x \in M$ and the familiar diagram
\begin{equazione}\label{N.ii19}%
\begin{split}%
\xymatrix{f_0(x)\ar[d]_{\varphi_0(g)}\ar[r]^-{\tau(x)} & f_1(x)\ar[d]^{\varphi_1(g)} \\ f_0(x')\ar[r]^-{\tau(x')} & f_1(x')}%
\end{split}%
\end{equazione}%
is commutative for all $g \in \mca[1]\G$, $g: x \to x'$. Suppose an action $S = (F,\sigma) \in \R[S]\H$ is given. Then one can apply $\tau^*$ to the isomorphism $\sigma: {\s^*F} \xto\iso {\t^*F}$ to obtain an isomorphism ${f_0^*F} \xto\iso {f_1^*F}$ in the category \stack[M]S%
\begin{equazione}\label{N.ii20}%
{f_0^*F} \can {\tau^*\s^*F} \xto{\:\tau^*\sigma\:} {\tau^*\t^*F} \can {f_1^*F}\text,%
\end{equazione}%
which may be denoted by the symbol ${\sigma\circ\tau}$. (Here one uses the identities $f_0 = {\s_{\H}\circ\tau}$ and $f_1 = {\t_{\H}\circ\tau}$.) By expressing \refequ{N.ii19} as an identity between suitable smooth maps, one can check that ${\sigma\circ\tau}$ is actually an isomorphism of \G\nobreakdash-actions between ${\varphi_0^*S}$ and ${\varphi_1^*S}$: in detail, consider the maps $({\tau\circ\t},\varphi_0)$ and $(\varphi_1,{\tau\circ\s})$, of \mca[1]\G\ (manifold of arrows) into $\mca[2]\H \equiv \mca\H$ (manifold of composable arrows), respectively given by $g \mapsto (\tau(\t[g]),\varphi_0(g))$ and $g \mapsto (\varphi_1(g),\tau(\s[g]))$; the commutativity of \refequ{N.ii19} implies that upon composing these maps with multiplication $\c: \mca[2]\H \to \H$ one gets the same result, ${\c \circ ({\tau\circ\t},\varphi_0)} = {\c \circ (\varphi_1,{\tau\circ\s})}$; from the latter identity it is easy to see that \refequ{N.ii20} is a morphism in \R[S]\G. Then the rule $(F,\sigma) \mapsto {\sigma\circ\tau}$ defines a natural isomorphism\inciso{in fact, a tensor preserving one}between the functors $\varphi_0^*, \varphi_1^*: \R[S]\H \to \R[S]\G$; we will use the notation%
\begin{equazione}\label{N.ii21}%
\tau^*: \varphi_0^* \xto\thicksim \varphi_1^*\text, \quad \tau^* \in \Iso^\otimes(\varphi_0^*,\varphi_1^*)\text.%
\end{equazione}%

We are now ready to discuss Morita equivalences. Recall that a homomorphism $\varphi: \G \to \H$ is said to be a \index{Morita equivalence|emph}\textit{Morita equivalence} in case
\begin{equazione}\label{xiv.8}%
\begin{split}%
\xymatrix@C=29pt@R=23pt{\G\ar[d]_{(\s,\t)}\ar[r]^-\varphi & \H\ar[d]^{(\s,\t)} \\ {M\times M}\ar[r]^-{f\times f} & {N\times N}}%
\end{split}%
\end{equazione}%
is a pullback diagram in the category of $\C^\infty$ manifolds and the mapping%
\begin{equazione}\label{xiv.9}%
{\t\circ\pr_2}: {M {_f\times_{\s}} \H} \to N\text,%
\end{equazione}%
which, loosely speaking, sends $f(x) \xto h y$ to $y$, is a surjective submersion. Our main goal in this section is to show that \textsl{the pull-back functor $\varphi^*: \R\H \to \R\G$ associated with a Morita equivalence $\varphi$ is an \index{equivalence of tensor categories}equivalence of tensor categories.}\footnote{Recall that a tensor functor $\Phi: \Kt[C] \to \Kt[D]$ is said to be a \index{tensor equivalence|emph}\textit{tensor equivalence} in case there exists a tensor functor $\Psi: \Kt[D] \to \Kt[C]$ along with tensor preserving natural isomorphisms ${\Psi\circ\Phi} \simeq \Id_{\Kt[C]}$ and ${\Phi\circ\Psi} \simeq \Id_{\Kt[D]}$.} Clearly, it will be enough to show that $\varphi^*$ is a categorical equivalence (in the familiar sense): this means that we have to look for a functor $\varphi_!: \R\G \to \R\H$ such that natural isomorphisms ${\varphi_! \circ \varphi^*} \simeq \Id_{\R\H}$ and ${\varphi^* \circ \varphi_!} \simeq \Id_{\R\G}$ exist.

Notice that the condition that the map~\refequ{xiv.9} should be a surjective submersion will of course be satisfied when $f$ itself is a surjective submersion. As a first step, we show how the task of constructing a quasi-inverse for the pullback functor $\varphi^*$ associated with an arbitrary Morita equivalence $\varphi$ may be reduced to the special case where $f$ is precisely a surjective submersion. To this end, consider the \index{weak pullback}\textit{weak pullback} (see~\cite{Moerdijk&Mrcun'03}, pp.~123--132)
\begin{equazione}\label{xiv.10}%
\begin{split}%
\xymatrix@C=40pt{\mathcal P\ar[d]^\chi\ar[r]^\psi & \G\ar[d]^\varphi_{}="b" \\ \H\ar[r]^\Id^(.26)\quad="a"\ar@{=>}@/^.8pc/"a";"b"^\tau & \H\text.\!\!}%
\end{split}%
\end{equazione}%
Let $P$ be the base manifold of the Lie groupoid $\mathcal P$. It is well-known (ibid.\ p.~130) that the Lie groupoid homomorphisms $\psi$ and $\chi$ are Morita equivalences with the property that the respective base maps $\mca[0]\psi: P \to M$ and $\mca[0]\chi: P \to N$ are surjective submersions. Now, if we succeed in proving that $\psi^*$ and $\chi^*$ are categorical equivalences then, since by \refequ{xiv.4} and \refequ{N.ii21} above we have a natural isomorphism (actually, a tensor preserving one)%
\begin{equazione}\label{xiv.11}%
\chi^* \xto{\:\iso\:} ({\varphi\circ\psi})^* \xfrom{\:\can\:} {\psi^* \circ \varphi^*}\text,%
\end{equazione}%
the same will be true of $\varphi^*$.%

From now on we will work under the hypothesis that the given Morita equivalence $\varphi$~\refequ{xiv.8} determines a surjective submersion $f: M \epito N$ on base manifolds. This being the case, there exists an open cover $N = \txtcup{i\in I}{}{V_i}$ of the manifold $N$ by open subsets $V_i$ such that for each of them one can find a smooth section $s_i: V_i \into M$ to $f$. We fix such a cover and such sections once and for all.%

Let an arbitrary object $R = (E,\varrho) \in \R[S]\G$ be given. For each $i \in I$ one can take the pull-back $E_i \equiv {{s_i}^*E} \in \stack[V_i]S$. Fix a couple of indices $i, j \in I$. Then, since \refequ{xiv.8} is a pull-back diagram, for each $y \in {V_i\cap V_j}$ there is exactly one arrow $g(y): s_i(y) \to s_j(y)$ such that $\varphi(g(y)) = y$. More precisely, let $y \mapsto g(y) = g_{ij}(y)$ be the smooth mapping defined as the unique solution to the following universal problem (in the $\C^\infty$ category)%
\begin{equazione}\label{xiv.12}%
\begin{split}%
\xymatrix@C=30pt{V_{ij}\ar@/_1.5pc/[ddr]_(.4){(s_i,s_j)}\ar@{-->}[dr]^(.55){g_{ij}}\ar@/^1.3pc/[drr]^(.6){\u|V_{ij}} \\ & \G\ar[d]^{(\s,\t)}\ar[r]^\varphi & \H\ar[d]^{(\s,\t)} \\ & {M\times M}\ar[r]^{f\times f} & {N\times N}\text,\!\!}%
\end{split}%
\end{equazione}%
where $\u: N \to \H$ denotes the unit section and $V_{ij} \equiv {V_i\cap V_j}$. Then, putting $E_i|_j = E_i|_{V_i\cap V_j}$ and $E_j|_i = E_j|_{V_i\cap V_j}$, one may pull the action $\varrho$ back along the map $g_{ij}$ so as to get an isomorphism $\theta_{ij}: E_i|_j \isoto E_j|_i$ in the category \stack[V_{ij}]S:%
\begin{equazione}\label{xiv.13}%
E_i|_j \can {({\s\circ g_{ij}})^*E} \can {{g_{ij}}^*\s^*E} \xto{\:\:{g_{ij}}^*\varrho\:\:} {{g_{ij}}^*\t^*E} \can {({\t\circ g_{ij}})^*E} \can E_j|_i%
\end{equazione}%
or, as an identity up to canonical isomorphisms,%
\begin{equazione+}\label{xiv.13'}%
\theta_{ij} = {{g_{ij}}^*\varrho}\text. & \bigl(\text{mod $\can$}\bigr)%
\end{equazione+}%
(Note that the fact that $\varrho$ is an isomorphism in the category \stack[\G]S, that is to say Lemma~\refcnt[N.13]{N.ii7}, is used in an essential way.) Next, from the obvious remark that for an arbitrary third index $k \in I$ one has $g_{ik}(y) = {g_{jk}(y) g_{ij}(y)}$ $\forall y \in V_{ijk} \equiv {V_i\cap V_j\cap V_k}$ (or better $g_{ik}|_j = {\c\circ (g_{jk}|_i,g_{ij}|_k)}$, where $g_{ik}|_j$ denotes the restriction of $g_{ik}$ to $V_{ijk}$ etc.), and from the multiplicative axiom \refequ[N.13]{N.equR3} for $\varrho$, it follows that the system of isomorphisms $\{\theta_{ij}\}$ constitutes a ``cocycle'' or ``descent datum'' for the family $\{E_i\}_{i\in I} \in \xstack[{\smash[b]{\txtcoprod{i\in I}{}{V_i}}}]S$, relative to the flat mapping $\txtcoprod{i\in I}{}{V_i} \to N$. Since $N$ is a paracompact manifold and \stack S is a smooth parastack, there exists some object ${\varphi_!E}$ of \stack[N]S along with isomorphisms $\theta_i: ({\varphi_!E})|_i \equiv ({\varphi_!E})|_{V_i} \xto\iso E_i$ in \stack[V_i]S, compatible with $\{\theta_{ij}\}$ in the sense that, modulo the identification $({\varphi_!E})_i|_{V_{ij}} \can ({\varphi_!E})_j|_{V_{ij}}$, one has the identity%
\begin{equazione+}\label{xiv.14}%
\theta_j|_i = \theta_j|_{V_{ij}} = {\theta_{ij} \cdot \theta_i|_{V_{ij}}} = {\theta_{ij} \cdot \theta_i|_j}\text. & \bigl(\text{mod $\can$}\bigr)%
\end{equazione+}%

For simplicity, let us put $F \equiv {\varphi_!E}$. Our next step will be to define a morphism $\sigma = {\varphi_!\varrho}: {{\s_{\H}}^*F} \to {{\t_{\H}}^*F}$, which is to provide the \H\nobreakdash-action on $F$. For each pair $V_i,V_{i'}$ we introduce the abbreviation $\H_{i,i'} \equiv \H(V_i,V_{i'})$; we also write $\H_{ij,i'j'} \equiv \H(V_{ij},V_{i'j'})$. Then the subsets $\H_{i,i'} \subset \mca[1]\H$ form an open cover of the manifold \mca[1]\H. Now, let $g_{i,i'}: \H_{i,i'} \to \G$ be the smooth map obtained by solving the following universal problem%
\begin{equazione}\label{xiv.15}%
\begin{split}%
\xymatrix@C=30pt{\H_{i,i'}\ar[d]_{(\s,\t)}\ar@{-->}[dr]^(.55){g_{i,i'}}\ar@/^1.5pc/[drr]^(.55){\text{~~~~inclusion}} \\ {V_i\times V_{i'}}\ar@/_1.3pc/[dr]_(.3){s_i\times s_{i'}} & \G\ar[d]^{(\s,\t)}\ar[r]^\varphi & \H\ar[d]^{(\s,\t)} \\ & {M\times M}\ar[r]^{f\times f} & {N\times N}\text.\!\!}%
\end{split}%
\end{equazione}%
We can use this map to define a morphism $\sigma_{i,i'}: ({{\s_{\H}}^*F})|_{i,i'} \to ({{\t_{\H}}^*F})|_{i,i'}$ in the category \stack[\H_{i,i'}]S, as follows:%
\begin{multiriga}\label{xiv.16}%
({{\s_{\H}}^*F})|_{i,i'} \can (\s_{\H}|_{i,i'})^*(F|_i) \xto{\:\:(\s_{\H}|_{i,i'})^*\theta_i\:\:} {(\s_{\H}|_{i,i'})^*E_i}%
\\%
\can {{g_{i,i'}}^*{\s_{\G}}^*E} \xto{\:\:{g_{i,i'}}^*\varrho\:\:} {{g_{i,i'}}^*{\t_{\G}}^*E}%
\\%
\can {(\t_{\H}|_{i,i'})^*E_{i'}} \xto{\:\:(\t_{\H}|_{i,i'})^*{\theta_i}^{-1}\:\:} (\t_{\H}|_{i,i'})^*(F|_{i'}) \can ({{\t_{\H}}^*F})|_{i,i'}%
\end{multiriga}%
or, in the form of an identity modulo canonical identifications,%
\begin{equazione+}\label{xiv.16'}%
\sigma_{i,i'} = {{(\t_{\H}|_{i,i'})^*{\theta_i}^{-1}} \cdot {{g_{i,i'}}^*\varrho} \cdot {(\s_{\H}|_{i,i'})^*\theta_i}}\text. & \bigl(\text{mod $\can$}\bigr)%
\end{equazione+}%
Starting from the equality of mappings%
\begin{equazione}\label{xiv.17}%
g_{i,i'}|_{j,j'} = {\left({g_{j'i'} \circ \t_{\H}|_{ij,i'j'}}\right) g_{j,j'}|_{i,i'} \left({g_{ji} \circ \s_{\H}|_{ij,i'j'}}\right)}%
\end{equazione}%
(note that $g_{j'i'} = {\i_{\G} \circ g_{i'j'}}$ where $\i_{\G}$ is the inverse map of \G) and the ``mod $\can$'' identities \refequ{xiv.13'}, \refequ{xiv.14} and \refequ{xiv.16'}, one can check that $\sigma_{i,i'}|_{j,j'} = \sigma_{j,j'}|_{i,i'}$ in \stack[\H_{ij,i'j'}]S; hence the morphisms $\sigma_{i,i'}$ glue together into a unique morphism $\sigma = {\varphi_!\varrho}$ of \stack[{\mca[1]\H}]S, with the property that $\sigma|_{i,i'} = \sigma_{i,i'}$.%

Next, suppose we are given a morphism $a: R \to R'$ in \R[S]\G, where $R' = (E',\varrho')$, let us say. Then we can obtain a morphism ${\varphi_!a}: {\varphi_!R} \to {\varphi_!R'}$, where ${\varphi_!R} = ({\varphi_!E},{\varphi_!\varrho})$ etc., by first letting $b_i = {{s_i}^*a}$ and the observing that%
\begin{equazione}\label{xiv.18}%
{\theta_{ij}' \cdot b_i|_j} = {b_j|_i \cdot \theta_{ij}} \quad \text{in~}\stack[V_{ij}]S%
\end{equazione}%
(because of the definition of $\theta_{ij} = \theta_{ij}^R$ and $\theta_{ij}' = \theta_{ij}^{R'}$ and because $a$ is a \G\nobreakdash-equivariant morphism). In this way we get a functor of \R[S]\G\ into \R[S]\H.%

The construction of the isomorphisms ${\varphi^* \circ \varphi_!} \simeq \Id_{\R\G}$ and ${\varphi_! \circ \varphi^*} \simeq \Id_{\R\H}$ is left as an exercise, to be done along the same lines.%

\capitolo{General Tannaka Theory}\label{4}

In the preceding chapter we laid down the foundations of Representation Theory in the abstract setting of smooth tensor stacks. The assumptions on the type \stack S were quite mild there, nothing more than just smoothness and the property of being a stack. However, in order to get our reconstruction theory to work effectively, we need to impose further restrictions on the type \stack S. We will call a smooth tensor stack a {\em stack of smooth fields} when it meets such additional requirements.%

The additional properties which characterize stacks of smooth fields are introduced in \refsez{N.15}. The stack of smooth vector bundles is an example. In the subsequent section we provide another fundamental example, the stack of {\em smooth (Euclidean) fields,} which will play a major role in the achievement of our Tannaka duality theorem for proper Lie groupoids in \refsez{N.20}. This stack is a nontrivial extension of the stack of smooth vector bundles, but its definition is as simple.%

\sezione{Stacks of Smooth Fields}\label{N.15}
The expression \gm{\index{stack of smooth fields}stack of smooth fields} will be employed to indicate a smooth (real or complex) tensor stack\footnote{In accordance with the philosophy of Note~\ref{N.12}.\ref{N.i33}, we use the word `stack' but we really mean `parastack'.} for which the axiomatic conditions listed below are satisfied. When dealing specifically with stacks of smooth fields we shall prefer them to be represented by the letter \stack F, which is more suggestive than the usual \stack S.

\sottosezione{The axioms}%

\sloppy%
Our first axiom is about the tensor product and pull-back operations. Roughly speaking, it states that the sections of a tensor product or a pull-back are exactly what one would expect them to be on the basis of the standard definition of tensor product and pull-back of sheaves of \smooth\nobreakdash-modules; however, for such sections the relation of equality may be coarser, in the sense that more sections may be regarded as being identical.%
\begin{enunciato}[tensor product \& pull-back]{Axiom I}\label{N.iiiA1}%
The canonical natural morphisms \refequ[N.11]{N.i18} and \refequ[N.11]{N.i21}%
$$%
\left\{\begin{aligned}%
&{\sections E\otimes_{\smooth[X]}^{}\sections{E'}} \to \sections{(E\otimes E')}%
\\%
&f^*(\sections[Y]F) \to \sections[X]{(f^*F)}%
\end{aligned}\right.%
$$%
are {\em surjective} (= epimorphisms of sheaves).%
\end{enunciato}%
Thus, every local smooth section of ${E\otimes E'}$ will possess, in the vicinity of each point, an expression as a finite linear combination, with smooth coefficients, of sections of the form ${\zeta\otimes\zeta'}$. Similarly, given any partial smooth section of $f^*F$, it will be possible to express it locally as a finite linear combination, with coefficients in \smooth[X], of sections of the form ${\eta\circ f}$.%

\fussy%
Suppose $E \in \stack[X]F$. Let us go back for a moment to the map $\sections E(U) \to E_x, \zeta \mapsto \zeta(x)$ we defined in \S\ref{N.11} (for each open neighbourhood $U$ of the point $x$). These maps are evidently compatible with the restriction to a smaller open neighbourhood of $x$, hence on passing to the inductive limit they will determine a linear map%
\begin{equazione}\label{N.iii1}%
(\sections E)_x \to E_x\text, \quad \zeta \mapsto \zeta(x)%
\end{equazione}%
of the stalk of \sections E at $x$ into the fibre of $E$ at the same point. We call this map the \index{evaluation of germs}\textit{evaluation (of germs)} at $x$. Notice, by the way, that the identity
\begin{equazione}\label{N.iii2}%
({\alpha\zeta})(x) = {\alpha(x)\zeta(x)}%
\end{equazione}%
holds for all germs of smooth sections $\zeta \in (\sections E)_x$ and of smooth functions $\alpha \in \smooth[X,x]$. It follows from Axiom~\textsc i (pull-back) that {\em for any stack of smooth fields, the evaluation of germs at a point is a surjective linear map.} Indeed, the stalk $(\sections E)_x$ coincides, as a vector space, with the space of global sections of $x^*(\sections E)$ (recall that $(\sections E)_x = \inductivelim{U\ni x}{\sections E(U)} = {x^{-1}(\sections E)}(\pt)$, actually as a \smooth[X,x]\nobreakdash-module), and the fibre $E_x$ is defined as the space of global sections of \sections{(x^*E)}; it is immediate to recognize that the evaluation of germs is just the map of global sections induced by \refequ[N.11]{N.i21}.%

\separazione%

The second axiom says that a difference between any two morphisms can be detected by looking at the linear maps they induce on the fibres.%
\begin{enunciato}[vanishing]{Axiom II}\label{N.iiiA2}%
Let $a: E \to E'$ be a morphism in \stack[X]F. Suppose that $a_x: E_x \to E'_x$ is zero $\forall x \in X$. Then $a = 0$.%
\end{enunciato}%

As a first, immediate consequence, an arbitrary section $\zeta \in \sections E(U)$ will vanish if and only if all its values $\zeta(u)$ will be zero as $u$ ranges over $U$: thus, one sees that {\em smooth sections are characterized by their values;} intuitively, one can think of the elements of $\sections E(U)$ as sections\inciso{in the usual sense}of the `bundle' of fibres $\{E_u\}$.%

Furthermore, by combining Axioms \textsc{ii} and \textsc i, it follows that {\em the functor $\sections[X]: \stack[X]F \to \SheavesOfModules{\smooth[X]}$ is faithful.} This is an easy consequence of the surjectivity of the evaluation of germs at a point; the argument we propose now will also be preparatory to the next axiom.%

For each morphism $a: E \to F$ in \stack[X]F, consider the `bundle' of linear maps $\{a_x: E_x \to F_x\}$ and the morphism $\alpha = \sections a: \sections E \to \sections F$ of sheaves of \smooth[X]\nobreakdash-modules. We start by asking what relation there is between these data. The link between the two is obviously provided by the above canonical evaluation maps of the stalks onto the fibres $(\sections E)_x \epito E_x$: it is clear that the stalk homomorphism $\alpha_x$ and the linear map $a_x$ have to be compatible, in the sense that the following square should commute%
\begin{equazione}\label{N.iii5}%
\begin{split}%
\xymatrix{(\sections E)_x\ar@{->>}[d]^{\text{eval.}}\ar[r]^{\alpha_x} & (\sections F)_x\ar@{->>}[d]^{\text{eval.}} \\ E_x\ar[r]^{a_x} & F_x\text.\!\!}%
\end{split}%
\end{equazione}%
In general, we shall say that a morphism of sheaves of modules $\alpha: \sections E \to \sections F$ and a `bundle' of linear maps $\{a_x: E_x \to F_x\}$ are {\em compatible,} whenever the diagram \refequ{N.iii5} commutes for all $x \in X$. Notice that, in view of the preceding axioms, compatibility implies that the morphism of sheaves and the bundle of linear maps determine each other unambiguously. (Indeed, in one direction, the morphism $\alpha$ clearly determines the maps $a_x$ through the commutativity of \refequ{N.iii5}. Conversely, the commutativity of \refequ{N.iii5} for all $x$ entails that for any smooth section $\zeta \in \sections E(U)$ one has the formula $[\alpha(U)\mspace{.6mu}\zeta](x) = a_x\bigl(\zeta(x)\bigr)$, and therefore, if $\alpha$ and $\beta$ are both compatible with $\{a_x\}$, it follows by Axiom \textsc{ii} that ${\alpha(U)\mspace{.6mu}\zeta} = {\beta(U)\mspace{.6mu}\zeta}$ for all $\zeta$ and hence that $\alpha = \beta$.) In particular, from $\sections a = \sections b$ it will follow that $a_x = b_x$ for all $x$ and therefore that $a = b$.%

\separazione%

Let us call a morphism of sheaves of modules $\alpha: \sections E \to \sections F$ {\em representable,} if it admits a compatible bundle of linear maps $\{a_x: E_x \to F_x\}$. Our next axiom, which complements the preceding one by providing a general criterion for the existence of morphisms in \stack[X]F, states that the collection of such morphisms is ``as big as possible'':%
\begin{enunciato}[morphisms]{Axiom III}\label{N.iiiA3}%
For every representable $\alpha: \sections E \to \sections F$, there exists a morphism $a: E \to F$ in \stack[X]F such that $\sections a = \alpha$.%
\end{enunciato}%

This axiom will not be used anyhere in the present section. It will play a role only in \S\ref{N.17}, where it is needed in order to construct morphisms of representations by means of fibrewise integration.%

\separazione%

We cannot yet deduce, from the axioms we have introduced so far, certain very intuitive properties that are surely reasonable for a ``smooth section''; for instance, if a section\inciso{or, more generally, a morphism}vanishes over a dense open subset of its domain of definition, it would be natural to expect it to be zero everywhere. Analogously, if the value of a section is non zero at a point then it should be non zero at all nearby points. The next axiom yields such properties, among many other consequences.%

We shall say that a Hermitian\inciso{or, in the real case, symmetric}form $\phi: {E\otimes E^*} \to \TU$ in \stack[X]F is a \index{metric|emph}\textit{Hilbert metric} on $E$, when for every point $x$ the induced form $\phi_x$ on the fibre $E_x$
\begin{equazione}\label{N.iii6}%
{E_x\otimes {E_x}^{\!*}} \;\xto{\text{can.}} \;({E\otimes E^*})_x \;\xto{\;\phi_x\,} \;\TU_x \;\can \;\nC%
\end{equazione}%
is a Hilbert metric (in the familiar sense, viz.\ positive definite).%

\begin{enunciato}[metrics]{Axiom IV}\label{N.iiiA4}%
Any object $E \in {\Ob\,\stack[X]F}$ supports \index{local metric}\textit{local} metrics; that is to say, the open subsets $U$ such that one can find a Hilbert metric on $E|_U$ cover $X$.
\end{enunciato}%
In general, one can only assume {\em local} metrics to exist, think e.g.\ of smooth vector bundles; however, as for vector bundles, global metrics can be constructed from local ones as soon as smooth partitions of unity are available on the manifold $X$ (e.g.\ when $X$ is paracompact).%

Let $E \in {\Ob\,\stack[X]F}$ and let $\phi$ be a Hilbert metric on $E$. By a \index{orthonormal frame|emph}\textit{$\phi$\nobreakdash-orthonormal frame for $E$} about a point $x$ of $X$ we mean a list of sections $\zeta_1, \ldots, \zeta_d \in \sections E(U)$, defined over a neighbourhood of $x$, such that for all $u$ in $U$ the vectors $\zeta_1(u), \ldots, \zeta_d(u)$ are orthonormal in $E_u$ (with respect to $\phi_u$) and
\begin{equazione}\label{N.iii7}%
\mathrm{Span}\,\bigl\{\zeta_1(x), \ldots, \zeta_d(x)\bigr\} \,= \,E_x\text.%
\end{equazione}%
{\em Orthonormal frames for $E$ exist about each point $x$ for which the fibre $E_x$ is finite dimensional.} Indeed, over some neighbourhood $N$ of $x$ we can first of all find local smooth sections $\zeta_1, \ldots, \zeta_d$ with the property that the vectors $\zeta_1(x), \ldots, \zeta_d(x)$ form a basis of the space $E_x$ (Axiom \textsc{i}). Since for all $n \in N$ the vectors $\zeta_1(n), \ldots, \zeta_d(n)$ are linearly dependent if and only if there is a $d$\nobreakdash-tuple of complex numbers $(z_1,\ldots,z_d)$ with ${\modulo{z_1}^2 + \cdots + \modulo{z_d}^2} = 1$ and $\txtsum{i=1}d{z_i\mspace{.6mu}\zeta_i(n)} = 0$, the continuous function%
$$%
{N\times \mathrm S^{2d-1}} \to \nR\text, \qquad (n;s_1,t_1,\ldots,s_d,t_d) \mapsto \modulo{\txtsum{\ell=1}d{{(s_\ell + {i t_\ell})} \zeta_\ell(n)}}%
$$%
must have a minimum $c>0$ at $n=x$, hence a lower bound $\tfrac{c}{2}$ on a suitable neighbourhood $U$ of $x$ so that the $\zeta_i(u)$ must be linearly independent for all $u \in U$. At this point it is enough to apply the Gram\nobreakdash--Schmidt process in order to obtain an orthonormal frame about $x$. This elementary observation (existence of orthonormal frames) will prove to be very useful. Let us start to illustrate its importance with some basic applications.%

Consider an \index{embedding}\textit{embedding $e: E' \into E$} in the category \stack[X]F, that is to say, a morphism such that the linear map $e_x: \smash{E'}_x \into E_x$ is injective for all $x$.\footnote{It follows immediately from Axiom \textsc{ii} that an embedding is a monomorphism. The converse need not be true because the functor $E \mapsto E_x$ doesn't have any exactness properties. For example, let $a$ be a smooth function on \nR\ such that $a(t) = 0$ if and only if $t = 0$. Then $a$, regarded as an element of $\End(\TU)$, is both mono and epi in $\stack F(\nR)$ while $a_0 = 0: \nC \to \nC$ is neither injective nor surjective.} Suppose there exists a global metric $\phi$ on the object $E$; also assume that $E' \in {\Ob\,\V[F]X}$ is locally trivial of (locally) finite rank. Then {\em $e$ admits a co-section,} i.e.\ {\em there exists a morphism $p: E \to E'$ with ${p\circ e} = \id$} (so $e$ is a section in the categorical sense). To prove this, note first of all that the metric $\phi$ will induce a metric $\phi'$ on $E' \into E$. Fix any point $x \in X$. Since $E'_x$ is finite dimensional, there exists a $\phi'$\nobreakdash-orthonormal frame for $E'$ about $x$, let us say $\zeta'_1, \ldots, \zeta'_d \in \sections{E'}(U)$. Put $\zeta_i^{} = {\sections e(U) \zeta'_i} \in \sections E(U)$, let $\phi_U$ be the metric induced on $E|_U$, and consider%
\begin{equazione}\label{N.iii10}%
\zeta^i: \,E|_U \,\can \,{E|_U\otimes\TU_U} \,\can \,{E|_U\otimes{\TU|_U}^*} \,\xto{E|_U^{}\otimes\zeta_i^*} \,{E|_U\otimes{E|_U}^*} \,\xto{\phi_U^{}} \,\TU_U\text.%
\end{equazione}%
Define $p_U: E|_U \to E'|_U$ as the composite of $E|_U \xto{\zeta^1\oplus\,\cdots\,\oplus\,\zeta^d} {\TU\oplus\cdots\oplus\TU}$ and ${\TU\oplus\cdots\oplus\TU} \xto{\zeta'_1\oplus\,\cdots\,\oplus\,\zeta'_d} E'|_U$. Note that ${(p_U^{})}_u: E_u^{} \to E'_u$ is the orthogonal projection, with respect to $\phi_u$, onto $\smash{E'}_u \into E_u$: it follows by Axiom \textsc{ii} that $p_U$ does not actually depend on $U$ or the other choices involved, so that we get a well-defined morphism $p: E \to E'$, by the prestack property; moreover, we have ${p\circ e} = \id$ for similar reasons.%

Another application: let $E \in {\Ob\,\stack[X]F}$, and suppose that the dimension of the fibres is (finite and) locally constant over $X$; {\em then $E \in {\Ob\,\V[F]X}$} i.e.\ {\em $E$ is locally trivial, of locally finite rank.} Indeed, fix an arbitrary point $x$. By Axiom \textsc{iv}, there exists an open neighbourhood $U$ of $x$ such that $E|_U$ supports a metric $\phi_U$. Since $E_x$ is finite dimensional, it is no loss of generality to assume that a $\phi_U$\nobreakdash-orthonormal system $\zeta_1, \ldots, \zeta_d \in \sections E(U)$ can be found; one can also assume ${\mathrm{dim}\,E_u} = d$ constant over $U$. Take $e \bydef {\zeta_1\oplus\cdots\oplus\zeta_d}: E' \bydef {\TU\oplus\cdots\oplus\TU} \into E|_U$ and $p: E|_U \to E'$ as above. It is immediate to see that $e$ and $p$ are fibrewise inverse to one another.%

\begin{lemma}\label{N.iii4}%
Let $X$ be a paracompact manifold and let $S \stackrel{i_S}{\into} X$ be a {\em closed} submanifold. Let \stack F be a stack of smooth fields.%

Let $E, F \in {\Ob\,\stack[X]F}$, and suppose that $E' = E|_S$ belongs to $\V[F]S$, i.e.\ is locally free, of locally finite rank.%

Then every morphism $a': E' \to F'$ in $\stack[S]F$ can be extended to a morphism $a: E \to F$ in $\stack[X]F$, i.e.\ $a' = a|_S$ for such an $a$.%
\end{lemma}%
\begin{proof}%
Fix a point $s \in S$. Then there exists an open neighbourhood $A$ of $s$ in $S$ such that over $A$ we can find a trivialization ($d$ summands)%
\begin{equazione}\label{N.iii3}%
E'|_A \iso {\TU_A\oplus\cdots\oplus\TU_A}\text.%
\end{equazione}%
Let $\zeta'_1, \ldots, \zeta'_d \in \sections{E'}(A)$ be the sections corresponding to this trivialization (so for instance $\zeta'_1$ is the composite ${\TU_S}|_A \can \TU_A \stackrel{\text{\tiny 1st}}\into {\TU_A \oplus \cdots \oplus \TU_A} \iso E'|_A$). Also, let $U$ be any open subset of $X$ such that ${U\cap S} = A$.%

Now, by Axiom \textsc i (pull-back case), taking smaller $U$ and $A$ about $s$ if necessary, it is no loss of generality to assume that there exist local sections $\zeta_1, \ldots, \zeta_d \in \sections E(U)$ with $\zeta'_k = {\zeta_k^{}\circ i_S}$, $\,k = 1, \ldots, d$. To see this, observe that locally about $s$ each $\zeta'_k$ is a finite linear combination $\txtsum{j}{}{\alpha_{j,k} ({\zeta_{j,k} \circ i_S})}$ with $\zeta_{j,k} \in \sections E(U)$ and $\alpha_{j,k} \in \C^\infty(A)$, by the cited axiom; hence if $U$ is chosen conveniently, let us say say so that there exists a diffeomorphism of $U$ onto a product ${A \times \nR^n}$, the coefficients $\alpha_{j,k}$ will extend to some smooth functions $\tilde\alpha_{j,k} \in \C^\infty(U)$ and $\zeta_k = \txtsum{j}{}{\tilde\alpha_{j,k}\mspace{.6mu}\zeta_{j,k}}$ will meet our requirements.%

We have already observed \refequ[N.11]{N.i22} that there is a canonical isomorphism of vector spaces ${(i_S^*E)}_s \can E_{i(s)}^{}$ which makes $({\zeta_k\circ i_S})(s)$ correspond to $\zeta_k(x)$, where we put $x = i_S(s)$. Hence the values $\zeta_k(x)$, $\,k = 1, \ldots, d\,$ are linearly independent in the fibre $E_x$, because the same is true of the values $\zeta'_k(s)$, $\,k = 1, \ldots, d\,$ in $\smash{E'}_s$ (the trivializing isomorphism \refequ{N.iii3} above yields a linear isomorphism $(E')_s \iso \nC^d$ which, as one can easily check, makes $\zeta'_k(s)$ correspond to the $k$\nobreakdash-th standard basis vector of $\nC^d$). This implies that if $U$ is small enough then the morphism $\zeta = {\zeta_1 \oplus \cdots \oplus \zeta_d}: {\TU_U \oplus \cdots \oplus \TU_U} \to E|_U$ is an embedding and admits a cosection $p: E|_U \to {\TU_U \oplus \cdots \oplus \TU_U}$, by Axiom~\textsc{iv} (existence of local metrics).%

Next, set $\eta'_k = {\sections{a'}(A)\mspace{.6mu}\zeta'_k} \in \sections{F'}(A)$. As remarked earlier in the proof, it is no loss of generality to assume that there exist partial sections $\eta_1, \ldots, \eta_d$ in $\sections F(U)$ with $\eta'_k = {\eta_k^{}\circ i_S}$. Again, these sections can be combined into a morphism $\eta: {\TU_U \oplus \cdots \oplus \TU_U} \to F|_U$ ($d$\nobreakdash-fold direct sum).%

Finally, we can take the composite%
$$%
E|_U \;\xto{\;p\;} \;\underset{\text{$d$ summands}}{\underbrace{\TU_U \oplus \cdots \oplus \TU_U}} \;\xto{\;\eta\;} \;F|_U\text.%
$$%
It is immediate to check that the restriction of this morphism to the submanifold $A \into U$ coincides with $a'|_A$, up to the canonical identifications $(E|_U)|_A \can E'|_A$ and $(F|_U)|_A \can F'|_A$. Let us summarize briefly what we have done so far: starting from an arbitray point $s \in S$, we have found an open neighbourhood $U = U^s$ of $x = i_S(s)$ in $X$, along with a morphism $a^s: E|_U \to F|_U$ whose restriction to $A = {U\cap S}$ agrees with $a'|_A$. This means that we have solved our problem locally.%

To conclude the proof, consider the open cover of $X$ formed by the open subsets $\{U^s: s \in S\}$ and the complement $U = {\complement_XS}$. (Here we use, of course, the closedness of $S$.) Since $X$ is a paracompact manifold, we can find a smooth partition of unity $\{\theta_i: i \in I\} \cup \{\theta\}$ subordinated to this open cover. Then\inciso{by the prestack property}the sum $a \bydef \txtsum{i\in I}{}{\theta_i\mspace{.6mu}a^{s_i}}$ corresponds to a well-defined morphism $E \to F$ in \stack[X]F, clearly extending $a'$.%
\end{proof}%

\separazione%

The last two axioms impose various finiteness requirements, both on the fibres and on the sheaf of smooth sections of an object.%

To begin with, there is a stock of conditions we shall impose on \stack F in order that the category \stack[\pt]F may be equivalent, as a tensor category, to the category of vector spaces of finite dimension. We gather these conditions into what we call the ``{\index{dimension axiom}dimension axiom}'':

\begin{enunciato}[dimension]{Axiom V}\label{N.iiiA5}%
It is required of the canonical pseudo-tensor functor $\refequ[N.11]{N.i9}: \stack[\pt]F \to \VectorSpaces$ that%
\begin{elenco}%
\item[\textsl{a)}]it is {\em fully faithful;}%
\item[\textsl{b)}]it factors through the subcategory whose objects are the finite dimensional vector spaces, in other words $E_*$ \refequ[N.11]{N.i10} is finite dimensional for all $E \in \stack[\pt]F$;%
\item[\textsl{c)}]it is a genuine {\em tensor} functor, i.e.\ \refequ[N.11]{N.i7} and \refequ[N.11]{N.i8} become isomorphisms of sheaves for $X = \pt$.%
\end{elenco}%
\end{enunciato}%

In particular, for each object $V \in \stack[\pt]F$ there exists a {\em trivialization} of $V$, i.e.\ an isomorphism $V \iso {\TU\oplus\cdots\oplus\TU}$ ({\em finite} direct sum). The number of copies of $\TU$ in any such decomposition determines the {\em dimension} of an object.%

Moreover, it follows from this axiom, and precisely from \textsl{c)}, that the functor `fibre at $x$', $E \mapsto E_x$ is a {\em complex tensor functor.} (In general, it is only a {\em complex pseudo-tensor functor,} see \S\ref{N.11}.)%

\separazione%

An object $E$ of \stack[X]F is \index{locally finite object, sheaf}\textit{locally finite,} if \sections E is a locally finitely generated \smooth[X]\nobreakdash-module. In other words, $E$ is locally finite if the manifold $X$ admits a cover by open subsets $U$ such that there exist local sections $\zeta_1, \ldots, \zeta_d \in \sections E(U)$ with the property
\begin{equazione}\label{N.iii8}%
\sections E|_U \;= \;\smooth[U]\,\{\zeta_1, \ldots, \zeta_d\}\text.%
\end{equazione}%
(The expression on the right-hand side has a clear meaning as a presheaf of sections over $U$; since it is always possible to assume $U$ paracompact, this presheaf is in fact a sheaf, as one can easily see by means of partitions of unity.) The condition on $U$ amounts to the existence of an epimorphism of sheaves of modules%
\begin{equazione}\label{N.iii9}%
\underset{\text{$d$ summands}}{\underbrace{\smooth[U]\oplus\cdots\oplus\smooth[U]}} \,\epito \,\sections E|_U\text.%
\end{equazione}%
\begin{enunciato}[local finiteness]{Axiom VI}\label{N.iiiA6}%
Let $X$ be a smooth manifold. Every object $E \in {\Ob\,\stack[X]F}$ is locally finite.%
\end{enunciato}%

The present axiom, like Axiom \textsc{iii} above, will play a role in the proof of the `Averaging Lemma' only, in \S\ref{N.17}.%

\sezione{Smooth Euclidean Fields}\label{N.16}
Our next goal in this section is to elaborate a concrete model for the axioms we just proposed. Of course, in order to be useful, such a model ought to contain much more than just vector bundles: in fact, we intend to exploit it later on, in \S\ref{N.20}, to prove a general reconstruction theorem for proper Lie groupoids. We first introduce a somewhat weaker notion which, however, is of some interest on its own.%

\begin{definizione}\label{xvi.1}
By a \index{smooth Hilbert field}\textit{smooth Hilbert field} we mean an object \Hfield H consisting of \textsl{(a)}~a family $\{H_x\}$ of Hilbert spaces, indexed over the set of points of a manifold $X$, and \textsl{(b)}~a sheaf $\,\Hfield[s]H$ of $\,\smooth[X]$\nobreakdash-modules of local sections of $\{H_x\}$, subject to the following conditions:
\begin{elenco}%
\item $\bigl\{\zeta(x): \zeta \in {(\Hfield[s]H)}_x\bigr\}$, where ${(\Hfield[s]H)}_x$ indicates the stalk at $x$, is a {\em dense} linear subspace of $H_x$;%
\item for each open subset $U$, and for all sections $\zeta, \zeta' \in \Hfield[s]H(U)$, the function \scalare{\zeta}{\zeta'} on $U$ defined by $u \mapsto \bigsca{\zeta(u)}{\zeta'(u)}$ turns out to be smooth.%
\end{elenco}%
We refer to the manifold $X$ as the {\em base} of \Hfield H; we can also say that \Hfield H is a smooth Hilbert field {\em over $X$.}%
\end{definizione}%

Some explanations are perhaps in order. By a \gm{local section of $\{H_x\}$} we mean here an element of the product \txtprod{x\in U}{}{H_x} of all the spaces over some open subset $U$ of $X$. The definition establishes in particular that for each open subset $U$ the set of sections $\Hfield[s]H(U)$ is a submodule of the $\C^\infty(U)$\nobreakdash-module of all the sections of $\{H_x\}$ over $U$. \Hfield[s]H will be called the \index{Gamma(E)@\sections E, \f H (sheaf of sections)}\index{sheaf sections@sheaf of sections \sections E, \f H}\textit{sheaf of smooth sections of \Hfield H} and the elements of $\,\Hfield[s]H(U)$ will be accordingly referred to as the \index{section}\textit{smooth sections of \Hfield H over $U$.} This terminology, overlapping with that of \S\ref{N.11}, has been introduced intentionally and will be justified soon.

Next, we need a suitable notion of morphism. There are various possibilities here. We choose the notion which seems to fit our purposes better: a bundle of bounded linear maps inducing a morphism of sheaves of modules. Precisely, let \Hfield E and \Hfield F be smooth Hilbert fields over $X$. A morphism of \Hfield E into \Hfield F is a family of bounded linear maps $\{a_x: E_x \to F_x\}$, indexed over the set of points of $X$, such that for each open subset $U \subset X$ and for all $\zeta \in \Hfield[s]E(U)$ the section over $U$ given by $u \mapsto {a_u\cdot\zeta(u)}$ belongs to $\Hfield[s]F(U)$.%

Smooth Hilbert fields over $X$ and their morphisms form a category which will be denoted by \stack[X]{H^\infty}. We want to turn the operation $X \mapsto \stack[X]{H^\infty}$ into a fibred (complex) tensor category \stack{H^\infty}, in the sense of \S\ref{N.11}. This fibred tensor category will prove to be a smooth tensor parastack (but not a stack: this is the reason why we work with the weaker notion of parastack) satisfying some of the axioms, although\inciso{of course}not all of them: for this reason, \stack{H^\infty} constitutes a source of interesting examples.%

Let us start with the definition of the tensor structure on the category \stack[X]{H^\infty} of smooth Hilbert fields.%

We shall concern ourselves with the tensor product of Hilbert fields in a moment; before doing that however we review the tensor product of Hilbert spaces. Let $V$ be a complex vector space. We denote by $V^*$ the space obtained by retaining the additive structure of $V$ while changing the scalar multiplication into ${zv^*} = ({\overline zv})^*$; the star here indicates that a vector of $V$ is to be regarded as one of $V^*$. If $\phi: {E\otimes E^*} \to \nC$ and $\psi: {F\otimes F^*} \to \nC$ are sesquilinear forms then we can combine them into a sesquilinear form on the tensor product ${E\otimes F}$%
\begin{equazione}\label{N.iiiH1}%
{({E\otimes F})\otimes {({E\otimes F})}^*} \,\can \,{({E\otimes E^*})\otimes ({F\otimes F^*})} \,\xto{\phi\,\otimes\,\psi} \,{\nC\otimes\nC} \,\can \,\nC\text.%
\end{equazione}%
If we compute this form on the generators of ${E\otimes F}$ we get%
\begin{equazione}\label{N.iiiH2}%
\scalare{e\otimes f}{e'\otimes f'} = {\scalare{e}{e'}\scalare{f}{f'}}\text.%
\end{equazione}%
Suppose now that both $\phi$ and $\psi$ are Hilbert space inner products. Then this formula shows that the form \refequ{N.iiiH1} is Hermitian. Moreover, if we express an arbitrary element $w$ of ${E\otimes F}$ as a linear combination $\txtsum{i=1}k{\txtsum{j=1}\ell{a_{i,j}\mspace{.6mu}{e_i \otimes f_j}}}$ with $e_1, \ldots, e_k$, resp.\ $f_1, \ldots, f_\ell$ orthonormal in $E$, resp.\ $F$, we see from \refequ{N.iiiH2} that $a_{i,j} = \scalare{w}{e_i\otimes f_j} = 0$ for all $i, j$ implies $w = 0$. Hence the form is non degenerate. The same expression can be used to show that the form is positive definite:%
$$%
\scalare{w}{w} \,= \,\txtsum{i,i'}{}{\,\txtsum{j,j'}{}{\,a_{i,j} \:\overline{a_{i',j'}} \:\delta^{j,j'}_{i,i'}}} \,= \,\txtsum{i,j}{}{\modulo{a_{i,j}}^2} \,\geqq \,0\text.%
$$%
The space ${E\otimes F}$ can be completed with respect to the pre-Hilbert inner product \refequ{N.iiiH1} to a Hilbert space called the \gm{\index{tensor product!of Hilbert spaces}Hilbert tensor product} of $E$ and $F$. We agree that from now on, when $E$ and $F$ are Hilbert spaces, the symbol ${E\otimes F}$ will denote the {\em Hilbert} tensor product of $E$ and $F$. It is equally easy to see that if $a: E \to E'$ and $b: F \to F'$ are bounded linear maps of Hilbert spaces then their tensor product extends by continuity to a bounded linear map of ${E\otimes F}$ into ${E'\otimes F'}$ that we still denote by ${a\otimes b}$. Moreover, the canonical isomorphisms of vector spaces ${u\otimes ({v\otimes w})} \mapsto {({u\otimes v})\otimes w}$ etc.\ extend by continuity to unitary isomorphisms ${E\otimes ({F\otimes G})} \isoto {({E\otimes F})\otimes G}$ etc.\ of Hilbert spaces.%

Suppose now that \Hfield E and \Hfield F are Hilbert fields over $X$. Consider the bundle of tensor products $\{{E_x\otimes F_x}\}$. For arbitrary local sections $\zeta \in \Hfield[s]E(U)$ and $\eta \in \Hfield[s]F(U)$, we let ${\zeta\otimes\eta}$ denote the section of $\{{E_x\otimes F_x}\}$ given by $u \mapsto {\zeta(u)\otimes\eta(u)}$. The law%
\begin{equazione}\label{N.iiiH3}%
U \;\mapsto \;{\smooth(U)\,\bigl\{{\zeta\otimes\eta}: \,\zeta \in \Hfield[s]E(U), \,\eta \in \Hfield[s]F(U)\bigr\}}%
\end{equazione}%
defines a sub-presheaf of the sheaf of local sections of $\{{E_x\otimes F_x}\}$. (We use expressions of the form ${\smooth(U)\{\cdots\}}$ to indicate the $\smooth(U)$\nobreakdash-module spanned by a collection of sections over $U$.) Let ${\Hfield E\otimes\Hfield F}$ denote the Hilbert field over $X$ consisting of the bundle $\{{E_x\otimes F_x}\}$ and the sheaf (of sections of this bundle) generated by the presheaf \refequ{N.iiiH3}, in other words, the smallest subsheaf of the sheaf of local sections of $\{{E_x\otimes F_x}\}$ containing \refequ{N.iiiH3}. We call ${\Hfield E\otimes\Hfield F}$ the \index{tensor product!of smooth Hilbert fields}\textit{tensor product} of \Hfield E and \Hfield F. Observe that for all morphisms $\Hfield E \xto\alpha \Hfield{E'}$ and $\Hfield F \xto\beta \Hfield{F'}$ of Hilbert fields over $X$, the bundle of bounded linear maps $\{{a_x\otimes b_x}\}$ yields a morphism ${\alpha\otimes\beta}$ of ${\Hfield E\otimes\Hfield F}$ into ${\Hfield{E'}\otimes\Hfield{F'}}$.%

Another operation which applies to Hilbert spaces is conjugation. This operation sends a Hilbert space $E$ to the conjugate vector space $E^*$ endowed with the Hermitian product $\scalare{v^*}{w^*} = \scalare{w}{v}$. We now carry conjugation of Hilbert spaces over to a functorial construction on Hilbert fields. Let \Hfield E be a Hilbert field over $X$. We get the conjugate field ${\Hfield E}^*$ by taking the bundle $\{{E_x}^*\}$ of conjugate spaces, along with the local smooth sections of \Hfield E regarded as local sections of $\{{E_x}^*\}$. If $\alpha = \{a_x\}: \Hfield E \to \Hfield F$ is a morphism of Hilbert fields over $X$ then, since a linear map $a_x: E_x \to F_x$ also maps ${E_x}^*$ linearly into ${F_x}^*$, we get a morphism $\alpha^* = \{{a_x}^*\}: {\Hfield E}^* \to {\Hfield F}^*$. Observe that the correspondence $\alpha \mapsto \alpha^*$ is {\em anti-}linear. Note also that ${\Hfield E}^{**} = \Hfield E$.%

The rest of the construction (tensor unit, the various constraints \ldots) is completely obvious. One obtains a complex tensor category, that is easily recognized to be additive as a \nC\nobreakdash-linear category. It remains to construct the complex tensor functor $f^*: \stack{H^\infty}(Y) \to \stack{H^\infty}(X)$ associated with a smooth map $f: X \to Y$, and to define the constraints \refequ[N.11]{N.i3}.%

Let \Hfield H be a Hilbert field over $Y$. The \index{pullback!along a smooth map}\index{pullback!of smooth Hilbert fields}\textit{pull-back} of \Hfield H along $f$, denoted by ${f^*\Hfield H}$, is the Hilbert field over $X$ whose description is as follows: the underlying bundle of Hilbert spaces, indexed by the points of $X$, is $\bigl\{H_{f(x)}\bigr\}$; the sheaf of smooth sections is generated\inciso{as a subsheaf of the sheaf of all local sections of the bundle $\bigl\{H_{f(x)}\bigr\}$}by the presheaf
\begin{equazione}\label{N.iiiH4}%
U \;\mapsto \;\smooth[X](U)\bigl\{{\eta\circ f}: \,\eta \in \Hfield[s]H(V), \,V \supset f(U)\bigr\}\text.%
\end{equazione}%
Since this is a presheaf of \smooth[X]\nobreakdash-modules (of sections), it follows that $\Hfield[s]{(\mathnormal f^*\Hfield H)}$ is a sheaf of \smooth[X]\nobreakdash-modules (of sections). Moreover, it is clear that for any morphism $\beta: \Hfield H \to \Hfield{H'}$ of Hilbert fields over $Y$, the family of bounded linear maps $\{b_{f(x)}\}$ defines a morphism $f^*\beta: {f^*\Hfield H} \to {f^*\Hfield{H'}}$ of Hilbert fields over $X$.%

\sloppy%
Observe that ${{f^*\Hfield H} \otimes {f^*\Hfield{H'}}}$ and $f^*({\Hfield H \otimes \Hfield{H'}})$ are exactly the same smooth Hilbert field over $X$, essentially because ${{(\eta\otimes\eta')} \circ f} = {(\eta\circ f)} \otimes {(\eta'\circ f)}$; also $\smooth[X] = {f^*\smooth[Y]}$. These identities can function as tensor functor constraints. Similarly $f^*(\Hfield H^*) = ({f^*\Hfield H})^*$ can be taken as a constraint, so we get a complex tensor functor $f^*: \stack{H^\infty}(Y) \to \stack{H^\infty}(X)$.%

\fussy%
Since the identities $f^*(g^*\Hfield H) = {({g\circ f})^*\Hfield H}$ and ${\smash{\id_X}^*\Hfield H} = \Hfield H$ hold, the operation $X \mapsto \stack[X]{H^\infty}$ is a ``strict'' fibred complex tensor category.%

Note that the `sheaf of sections'\inciso{defined abstractly only in terms of the prestack structure of \stack{H^\infty}, as explained in \S\ref{N.11}}turns out to be precisely the `sheaf of smooth sections' which we introduced in the above definition as one of the two constituent data of a smooth Hilbert field. However, note that the fibre $\Hfield H_x$ (in the sense of \S\ref{N.11}) will be in general only a dense subspace of the Hilbert space $H_x$ (this is the reason why we use two distinct notations); of course, $\Hfield H_x = H_x$ whenever $H_x$ is finite dimensional.%

Let \stack[X]{E^\infty} be the full subcategory of \stack[X]{H^\infty} consisting of all objects \Hfield E whose sheaf of sections is locally finitely generated over $X$, in the sense of Axiom \textsc{vi}. \stack[X]{E^\infty} is a complex tensor subcategory i.e. it is closed under $\otimes$, $\boldsymbol*$ and it contains the tensor unit: indeed, ${\Hfield[s]E \otimes_{\smooth} \Hfield[s]{E'}}$, which is a locally finitely generated sheaf of modules over $X$ because such are $\Hfield[s]E$ and $\Hfield[s]{E'}$, surjects (as a sheaf) onto $\Hfield[s]{(E\otimes E')}$, by Axiom \textsc i, so the latter will be locally finite too, as contended. Moreover, the pull-back functor $f^*: \stack{H^\infty}(Y) \to \stack{H^\infty}(X)$ carries \stack[Y]{E^\infty} into \stack[X]{E^\infty}. We obtain a smooth substack $\stack{E^\infty} \subset \stack{H^\infty}$ of additive complex tensor categories; it is clear that \stack{E^\infty} satisfies Axioms \textsc i\nobreakdash--\textsc{vi}.%

The objects of the subcategory $\stack[X]{E^\infty} \subset \stack[X]{H^\infty}$ will be referred to as \index{smooth Euclidean field}\textit{smooth Euclidean fields} over $X$.

\sezione{Construction of Equivariant Maps}\label{N.17}
Let \stack F denote an arbitrary stack of smooth fields, to be regarded as fixed throughout the present section.

The next lemma is to be used in combination with Lemma \ref{N.15}.\ref{N.iii4}.%
\begin{lemma}\label{N.iv2}%
Let \G\ be a (locally) transitive Lie groupoid, and let $X$ be its base manifold. Consider any representation $(E,\rho) \in \R[F]G$. Then $E \in \V[F]X$ i.e.\ $E$ is a locally trivial object of \stack[X]F.%
\end{lemma}
\begin{proof}
\index{groupoid!locally transitive}\index{groupoid!transitive}\index{locally transitive groupoid}\index{transitive groupoid}Local transitivity means that the mapping $(\s,\t): \G \to {X\times X}$ is a submersion. Fix a point $x \in X$. Since $(x,x)$ lies in the image of the map $(\s,\t)$, the latter admits a local smooth section ${U\times U} \to \G$ over some open neighbourhood of $(x,x)$. Let us consider the `restriction' $g: U \to \G$ of this section to $U \equiv {U\times \{x\}}$: $g$ will be a smooth map for which the identities $\s(g(u)) = u$ and $t(g(u)) = x$ hold for all $u \in U$.

Let $\pt \xto x X$ denote the map $\pt \mapsto x$. We have already noticed that, by the `dimension' Axiom \refequ[N.15]{N.iiiA5}, there is an isomorphism ${x^*E} \iso {\TU \oplus \cdots \oplus \TU}$ (a trivialization) in \stack[\pt]F. Now, it will be enough to pull $\rho$ back to $U$ along the smooth map $g$ and observe that there is a factorization of the map ${\t\circ g}$ as the collapse $c: U \to \pt$ followed by $x: \pt \to X$ in order to conclude that there is also a trivialization $E|_U \iso {\TU_U \oplus \cdots \oplus \TU_U}$ in \stack[U]F. Indeed, since $\rho$ is an isomorphism, one can form the following long invertible chain%
\begin{multline}\notag%
E|_U = {i_U^*E} = {({\s\circ g})^*E} \can {g^*\s^*E} \xto{\:g^*\rho\:} {g^*\t^*E} \can {({\t\circ g})^*E} =\\%
= {({x\circ c})^*E} \can {c^*(x^*E)} \iso {c^*({\TU \oplus \cdots \oplus \TU})} = {\TU_U \oplus \cdots \oplus \TU_U}%
\end{multline}%
(recall that the pull-back $c^*$ preserves direct sums).%
\end{proof}%

Let $i: S \into X$ be an \index{invariant submanifold, subset}\textit{invariant} immersed submanifold, viz.\ one whose image $i(S)$ is an invariant subset under the `tautological' action of \G\ on its own base. The pull-back of \G\ along $i$ makes sense and proves to be a Lie subgroupoid\footnote{In general, a \gm{{\em Lie subgroupoid}} is a Lie groupoid homomorphism $(\varphi,f)$ such that both $\varphi$ and $f$ are injective immersions.} $\iota: \G|_S \into \G$ of \G. (Observe that $\G|_S = \G^S = \s_{\G}^{-1}(S)$.) In the special case of an orbit immersion, $\G|_S$ will be a transitive Lie groupoid over $S$. Then the lemma says that for any $(E,\rho) \in {\Ob\,\R\G}$ the pull-back ${i_S^*E}$ is a locally trivial object of \stack[S]F, because the transitive Lie groupoid \R{\G|_{\mathnormal S}} acts on ${i_S^*E}$ via ${\iota_S^*\rho}$. In particular, when the orbit $S \into X$ is a submanifold, we can also write $E|_S = {i_S^*E} \in \V[F]S$.

\begin{nota}\label{N.iv1}%
The notion of Lie groupoid representation we have been working with so far is completely intrinsic. We were able to prove all results by means of purely formal arguments, involving only manipulations of commutative diagrams. For the purposes of the present section, however, we have to change our point of view.%

Let \G\ be a Lie groupoid. Consider a representation $(E,\rho) \in {\Ob\,\R\G}$, $\s^*E \xto\rho \t^*E$. Each arrow $g$ determines a linear map $\rho(g): E_{\s(g)} \to E_{\t(g)}$ defined via the commutativity of the diagram%
\begin{equazione}\label{N.iv1a}%
\begin{split}%
\xymatrix{[g^*\s^*E]_*\ar[d]^{[g^*\rho]_*}\ar[r]^-{[\can]_*} & [\s(g)^*E]_*\ar@{=}[r]^-{\text{def.}} & E_{\s(g)}\ar@{-->}[d]^{\rho(g)} \\ [g^*\t^*E]_*\ar[r]^-{[\can]_*} & [\t(g)^*E]_*\ar@{=}[r]^-{\text{def.}} & E_{\t(g)}}%
\end{split}%
\end{equazione}%
where the notation \refequ[N.11]{N.i10} is used. It is routine to check that the cocycle conditions \refequ[N.13]{N.equR2} and \refequ[N.13]{N.equR3} in the definition of representation imply that the correspondence $g \mapsto \rho(g)$ is multiplicative i.e.\ that $\rho(g'g) = {\rho(g') \circ \rho(g)}$ and $\rho(x) = \id$ for each point of the base manifold $X$.%

Next, consider any arrow $g_0$. Also, let $\zeta \in \sections E(U)$ be a section defined over a neighbourhood of $\s(g_0)$ in $X$. Recall that according to \refequ[N.11]{N.i19} $\zeta$ will determine the section ${\zeta\circ\s} \in \sections[{\G}]{(\s^*E)}(\G^U)$, defined over the open subset $\G^U = \s^{-1}(U)$ of the manifold of arrows \mca[1]{\G}; the morphism of sheaves of modules \sections{\rho} can be evaluated at ${\zeta\circ\s}$: $[\sections\rho\,(\G^U)]({\zeta\circ\s}) \in \sections{(\t^*E)}(\G^U)$. Axiom \refequ[N.15]{N.iiiA1} implies that there exists an open neighbourhood $\Gamma \subset \G^U$ of $g_0$ over which $[\sections\rho\,(\G^U)]({\zeta\circ\s})$ can be expressed as a finite linear combination, with coefficients in $\C^\infty(\Gamma)$, of sections of the form ${\zeta'_i\circ\t}$ with $\zeta'_i$, $i = 1, \cdots, d$ defined over $\t(\Gamma)$. Explicitly,%
\begin{equazione}\label{N.iv1b}%
\bigl[\sections\rho\,(\Gamma)\bigr]({\zeta \circ \s|_\Gamma}) \:= \:\txtsum{i=1}d{r_i \, ({\zeta'_i \circ \t})|_\Gamma}%
\end{equazione}%
with $r_1, \ldots, r_d \in \C^\infty(\Gamma)$ and $\zeta'_1, \ldots, \zeta'_d \in (\sections E)(\t(\Gamma))$. This equality can be evaluated at $g \in \Gamma$ in the abstract sense of \refequ[N.11]{N.i12}, also taking \refequ{N.iv1a} into account, to get a more intuitive expression%
\begin{equazione}\label{N.iv1c}%
{\reprho(g) \cdot \zeta(\s[g])} \:= \:\txtsum{i=1}d{r_i(g) \,\zeta'_i(\t[g])}\text.%
\end{equazione}%

To summarize: any \G\nobreakdash-action $(E,\rho)$ determines an operation $g \mapsto \rho(g)$ which assigns a linear isomorphism $E_x \xto{\rho(g)\!\!} E_{x'}$ to each arrow $x \xto g x'$ in such a way that the composition of arrows is respected; moreover, the operation enjoys a `smoothness property' whose technical formulation is synthesized in Equation~\refequ{N.iv1c}. Conversely, it is yet another exercise to recognize that such data determine an action of \G\ on $E$, by Axiom~\refequ[N.15]{N.iiiA3}. Therefore we see that for the representations whose type is a stack of smooth fields the intrinsic definition of \S\ref{N.13} is equivalent to a more concrete definition involving an operation $g \mapsto \rho(g)$ and a `smoothness condition' expressed pointwise.%
\end{nota}%

Let \G\ be a Lie groupoid over a manifold $X$. Consider any representation $(E,\rho) \in {\Ob\,\R\G}$. Fix an arbitrary point $x_0 \in X$. Using the remarks of the preceding note, the fact that the fibre $E_0 \bydef E_{x_0}$ is a finite dimensional vector space, by Axiom~\refequ[N.15]{N.iiiA5}, and the fact that the evaluation map \refequ[N.15]{N.iii1}%
$$%
(\sections E)_0 \to E_0\text, \quad \zeta \mapsto \zeta(x_0)%
$$%
is surjective, one sees at once that the operation%
\begin{equazione}\label{N.iv3}%
\rho_0: \G_0 \to \GL(E_0)\text, \quad g \mapsto \rho(g)%
\end{equazione}%
is a smooth representation of the Lie group $G = \G_0$ (= the isotropy group at $x_0$) on the finite dimensional vector space $E_0$.%

Now, suppose we are given a $G$\nobreakdash-equivariant linear map $A: E_0 \to F_0$, for some other \G\nobreakdash-action $(F,\sigma)$. Let $S \into X$ be the orbit through $x_0$; just to fix ideas, assume it is a submanifold. The theory of Morita equivalences of \S\ref{N.14} says that there exists a unique morphism $A': (E|_S,\rho|_S) \to (F|_S,\sigma|_S)$ in \R{\G|_{\mathnormal S}} such that $(A')_0 = A$, up to the standard canonical identifications. Actually, for any point $z \in S$ and any arrow $g \in \G(x_0,z)$ one has%
\begin{equazione}\label{N.iv4}%
(A')_z = {\sigma(g) \cdot A \cdot {\rho(g)}^{-1}}: E_z \to F_z\text.%
\end{equazione}%
Set $E' = E|_S$. As remarked earlier, since the groupoid $\G|_S$ is transitive it follows that the object $E'$ is locally trivial, by Lemma~\ref{N.iv2}. If the submanifold $S \into X$ is in addition closed then, since base manifolds of Lie groupoids are always paracompact, Lemma~\ref{N.15}.\ref{N.iii4} will yield a morphism $a: E \to F$ extending $A'$ and hence, a fortiori, $A$.%

\sottosezione{The averaging operator}%

We are now ready to describe an \gm{\index{averaging technique}averaging technique} which is of central importance in our work\nobreakdash---as the reader will see. We explain in detail how, starting from any (right-invariant) Haar system $\mu = \{\mu^x\}$ on a proper Lie groupoid \G\ over a manifold $M$, one can construct, for each pair of representations $R = (E,\rho), S = (F,\sigma) \in \R\G$ (of type \stack F), a linear operator
\begin{equazione}\label{N.ivB1}%
\Av\mu: \Hom_{\stack[M]F}(E,F) \to \Hom_{\R\G}(R,S)%
\end{equazione}%
called the \gm{averaging operator (of type \stack F)} associated with $\mu$, with the property that $\Av\mu(a) = a$ whenever $a$ already belongs to the subspace $\Hom_{\R\G}(R,S) \subset \Hom_{\stack[M]F}(E,F)$. This construction will be compatible with the restriction to an invariant submanifold of the base: namely, if $N \subset M$ is any such submanifold then, letting $\nu$ denote the Haar system induced by $\mu$ on the subgroupoid $\G|_N = \G^N \stackrel{\:\iota_N}{\into} \G$ (what we are saying makes sense because $N$ is invariant), the following diagram will commute%
\begin{equazione}\label{N.ivB2}%
\begin{split}%
\xymatrix{\Hom_{\stack[M]F}(E,F)\ar[d]^{i_N^*}\ar[r]^-{\Av\mu} & \Hom_{\R\G}(R,S)\ar[d]^{\iota_N^*} \\ \Hom_{\stack[N]F}(E|_N,F|_N)\ar[r]^(.49){\Av\nu} & \Hom_{\R{\G|_\mathnormal{N}}}({\iota_N^*R},{\iota_N^*S})\text.\!\!}%
\end{split}%
\end{equazione}%
Thus, in particular, if $a$ restricts to an invariant morphism over $N$ then $\Av\mu(a)|_N = a|_N$. Since $\mu$ will be fixed throughout the present discussion, we abbreviate $\Av\mu(a)$ into $\tilde a$ from now on.%

\sloppy%
We start from a very simple remark, valid even without assuming \G\ to be proper. Suppose that $\zeta \in \sections E(U)$ and $\eta_1, \ldots, \eta_n \in \sections F(U)$ are sections over some open subset of $M$, and moreover that $\eta_1, \ldots, \eta_n$ are local generators for \sections F over $U$; then for each $g_0 \in \G^U = \s^{-1}(U)$ there exists an open neighbourhood $g_0 \in \Gamma \subset \G^U$, along with smooth functions $\phi_1, \ldots, \phi_n \in \C^\infty(\Gamma)$, such that the identity%
\begin{equazione}\label{N.ivB3}%
{\sigma(g)^{-1} \cdot a_{\t(g)} \cdot \rho(g) \cdot \zeta(\s[g])} = \txtsum{j=1}n{\phi_j(g)\eta_j(\s[g])}%
\end{equazione}%
holds in the fibre $F_{\s(g)}$ for all $g \in \Gamma$. To see this, recall that, according to Note~\ref{N.iv1}, there are an open neighbourhood $\Gamma$ of $g_0$ in $\G^U$ and local smooth sections $\zeta'_1, \ldots, \zeta'_m$ of $E$ over $U' = \t(\Gamma)$, such that ${\rho(g) \zeta(\s[g])} = \txtsum{i=1}m{r_i(g) \zeta'_i(\t[g])}$ for some smooth functions $r_1, \ldots, r_m \in \C^\infty(\Gamma)$. For $i = 1, \ldots, m$, put $\eta'_i = {\sections a(U')(\zeta'_i)} \in \sections F(U')$. Since $\Gamma^{-1}$ is a neighbourhood of $g_0^{-1}$ we can assume\inciso{again by Note~\ref{N.iv1}, using the hypothesis that the $\eta_j$'s are generators}$\Gamma$ to be so small that for each $i = 1, \ldots, m$ there exist smooth functions $s_{1,i}, \ldots, s_{n,i} \in \C^\infty(\Gamma^{-1})$ with ${\sigma(g^{-1}) \eta'_i(\t[g])} = \txtsum{j=1}n{s_{j,i}(g^{-1}) \eta_j(\s[g])}$ ${\forall g} \in \Gamma$. Hence for all $g \in \Gamma$ we get%
\begin{multline}\notag%
{\sigma(g)^{-1} \cdot a_{\t(g)} \cdot \rho(g) \cdot \zeta(\s[g])} = {\sigma(g^{-1}) \cdot a_{\t(g)} \cdot \txtsum{i=1}m{r_i(g) \zeta'_i(\t[g])}} =\\%
= \txtsum{i=1}m{r_i(g) \sigma(g^{-1}) \eta'_i(\t[g])} = \txtsum{j=1}n{\left[\txtsum{i=1}m{r_i(g) s_{j,i}(g^{-1})}\right] \eta_j(\s[g])}\text,%
\end{multline}%
which is \refequ{N.ivB3} with $\phi_j(g) = \txtsum{i=1}m{r_i(g)s_{j,i}(g^{-1})}$, $j = 1, \ldots, n$.%

\fussy%
Let $\alpha = \sections a \in \Hom_{\smooth}(\sections E,\sections F)$. We can use the last remark to obtain a morphism $\tilde\alpha: \sections E \to \sections F$ of sheaves of modules over $M$, in the following way. Let $\zeta$ be a local smooth section of $E$, defined over an open subset $U \subset M$ so small that there exists a system $\eta_1, \ldots, \eta_n$ of local generators for $F$ over $U$ (such a system can always be found locally, because \stack F satisfies Axiom~\refequ[N.15]{N.iiiA6}). For each $g_0 \in \G^U = \s^{-1}(U)$, select an open neighbourhood $\Gamma(g_0)$, along with smooth functions $\phi^{g_0}_1, \ldots, \phi^{g_0}_n \in \C^\infty\bigl(\Gamma(g_0)\bigr)$, as in \refequ{N.ivB3}. Since the manifold of arrows of \G, and\inciso{consequently}its open submanifold $\G^U$, is paracompact (we are assuming \G\ proper now; cf.\ \S\ref{O.1.5}), there will be a smooth partition of unity $\{\theta_i\}$, $i \in I$ on $\G^U$ subordinated to the open cover $\{\Gamma(g)\}$, $g \in \G^U$. Then we put%
\begin{equazione}\label{N.ivB4}%
{\tilde\alpha(U)\mspace{1mu}\zeta} = \txtsum{j=1}n{\Phi_j\mspace{1mu}\eta_j}\text, \quad \text{where} \quad \Phi_j(u) = {\int_{\G^u} \txtsum{i\in I}{}{\theta_i(g)\mspace{1mu}\phi^i_j(g)} \:\mathrm{d}\mu^u(g)}%
\end{equazione}%
(note that the integrand $\txtsum{i\in I}{}{\theta_i\mspace{1mu}\phi^i_j}$ is a smooth function on $\G^U$ and hence $\Phi_j \in \C^\infty(U)$, $j = 1, \ldots, n$). Of course, many arbitrary choices are involved here, so one has to make sure that this definition is not ambiguous (however, as soon as \refequ{N.ivB4} is known to be independent of all these choices, it will certainly define a morphism of sheaves of modules over $M$). One can do this, in two steps, by introducing independently a certain bundle of linear maps $\{\lambda_x: E_x \to F_x\}$ over $M$ first and then checking that ${[\tilde\alpha(U)\mspace{1mu}\zeta]}(u) = {\lambda_u\bigl(\zeta(u)\bigr)}$ for all $u \in U$. Since the right-hand term in the last equality will not depend on any choice, Axiom~\refequ[N.15]{N.iiiA2} will imply at once that ${\tilde\alpha(U)\mspace{1mu}\zeta}$ is a well-defined section of $F$ over $U$. The same equality will furthermore yield the conclusion that $\tilde\alpha \in \Hom_{\smooth[M]}(\sections E,\sections F)$ is equal to \sections{\mspace{1mu}\tilde a} for a unique $\tilde a \in \Hom_{\stack[M]F}(E,F)$, by Axiom~\refequ[N.15]{N.iiiA3}. It should be clear how to proceed now, but let us carry out the details anyway, for completeness. If we look at \refequ{N.ivB3} with $\s(g) = x$ fixed, we immediately recognize that the map%
\begin{equazione}\label{N.ivB5}%
\G^x \to F_x\text, \qquad g \:\mapsto \:{\sigma(g)^{-1} \cdot a_{\t(g)} \cdot \rho(g) \cdot \zeta(x)}\text,%
\end{equazione}%
of the manifold $\G^x = \s^{-1}(x)$ into the finite dimensional vector space $F_x$, is of class $\C^\infty$ and hence continuous. Since for each $v \in E_x$ there is some local section $\zeta$ of $E$ about $x$ such that $v = \zeta(x)$, by Axiom~\refequ[N.15]{N.iiiA1}, we can write down the integral%
\begin{equazione}\label{N.ivB6}%
{a^\mu(x) \cdot v} \bydef {\int_{\G^x} {\sigma(g)^{-1} \cdot a_{\t(g)} \cdot \rho(g) \cdot v} \:\:\mathrm{d}\mu^x(g)}%
\end{equazione}%
for each $v \in E_x$. Clearly $v \mapsto {a^\mu(x) \cdot v}$ defines a linear map of $E_x$ into $F_x$, so we get our bundle of linear maps $\bigl\{a^\mu(x): E_x \to F_x\bigr\}$. It remains to check, for an arbitrary $u \in U$, the equality ${[\tilde\alpha(U)\mspace{1mu}\zeta]}(u) = {a^\mu(u) \cdot \zeta(u)}$ with ${\tilde\alpha(U)\mspace{1mu}\zeta}$ given by \refequ{N.ivB4}. The computation is straightforward:%
\begin{align*}%
{[\tilde\alpha(U)\mspace{1mu}\zeta]}(u) &\:= \:\txtsum{j=1}n{\Phi_j(u)\mspace{1mu}\eta_j(u)} \:= \:\txtsum{j=1}n{{\int_{\G^u} \txtsum{i\in I}{}{\theta_i(g)\mspace{1mu}\phi^i_j(g)} \:\mathrm{d}\mu^u(g)} \:\eta_j(u)}%
\\%
&\:= {\int_{\G^u} \txtsum{i\in I}{}{\theta_i(g) \:\txtsum{j=1}n{\phi^i_j(g)\mspace{1mu}\eta_j(\s[g])}} \:\mathrm{d}\mu^u(g)}%
\\%
&\:= {\int_{\G^u} \txtsum{i\in I}{}{\theta_i(g)} {\bigl[\sigma(g)^{-1} \cdot a_{\t(g)} \cdot \rho(g) \cdot \zeta(\s[g])\bigr]} \,\mathrm{d}\mu^u(g)}%
\\%
&\:= \:{a^\mu(u) \cdot \zeta(u)}\text.%
\end{align*}%

In conclusion, we define $\Av\mu(a)$ as the unique morphism $\tilde a: E \to F \in \stack[M]F$ such that $\sections{\tilde a} = \widetilde{(\sections a)}$. The linearity of $a \mapsto \Av\mu(a)$ follows now from \refequ{N.ivB6}, the relation $[{\tilde\alpha(U) \zeta}](u) = {a^\mu(u)\cdot \zeta(u)}$ and the faithfulness of $a \mapsto \sections a$. It remains to show that $\Av\mu(a)$ belongs to $\Hom_{\R\G}(R,S)$ and that $\Av\mu(a)$ equals $a$ when $a$ already belongs to $\Hom_{\R\G}(R,S)$; although the calculation is completely standard, we review it because of its importance. In order to prove that $\tilde a \equiv \Av\mu(a)$ is a morphism of \G\nobreakdash-actions, it will be enough (by Axiom~\refcnt[N.15]{N.iiiA2}) to check the identity ${\tilde a_{\t(g)} \circ \varrho(g)} = {\sigma(g) \circ \tilde a_{\s(g)}}$ or equivalently, letting $x = \s(g)$ and $x' = \t(g)$, the identity ${a^\mu(x')\circ \varrho(g)} = {\sigma(g)\circ a^\mu(x)}$ for each arrow $g$; the corresponding computation reads as follows:%
\begin{align*}%
{a^\mu(x')\circ \varrho(g)} &= {\int_{\G(x',\text-)} {\sigma(g')^{-1} a_{\t(g')} \varrho(g') \varrho(g)}\, \mathrm d\mu^{x'}(g')} && \text{by~\refequ{N.ivB6}}\\%
&= {\int_{\G(x,\text-)} {\sigma(g) \sigma(h)^{-1} a_{\t(h)} \varrho(h)}\, \mathrm d\mu^x(h)} && \text{by~right-invariance}\\%
&= {\sigma(g)\circ a^\mu(x)} && \text{by~\refequ{N.ivB6}~again.}%
\end{align*}%
Next, whenever $a$ is an element of $\Hom_{\R\G}(R,S)$, the computation%
\begin{align*}%
a^\mu(x) &= {\int_{\G(x,\text-)} {\sigma(g)^{-1} a_{\t(g)} \varrho(g)}\, \mathrm d\mu^{x}(g)} && \text{by~\refequ{N.ivB6}}\\%
&= {\int_{\G(x,\text-)} a_x\, \mathrm d\mu^x(g)} && \text{because~$a \in \Hom_{\R\G}(R,S)$}\\%
&= a_x && \text{because~$\mu$~is~normalized}%
\end{align*}%
proves the identity $\tilde a = a$.%

\sottosezione{Applications}%

For the reader's convenience and for future reference, it will be useful to collect the conclusions of the previous subsection into a single statement. As ever, \stack F will denote an arbitrary stack of smooth fields, for example the stack of smooth vector bundles or the stack of smooth Euclidean fields.%

\begin{proposizione}[Averaging Lemma]\label{N.ivB7}%
Let \G\ be a proper Lie groupoid over a manifold $M$, and let $\mu$ be a right-invariant Haar system on \G.%

Then for any given \G\nobreakdash-actions $R = (E,\varrho)$ and $S = (F,\sigma)$ of type \stack F, each morphism $a: E \to F$ in the category \stack[M]F determines a (unique) morphism $\tilde a = \Av\mu(a): R \to S \in \R[F]\G$ through the requirement that for each $x \in M$ the map $\tilde a_x: E_x \to F_x$ should be given by the formula%
\begin{equazione+}\label{N.ivB8}%
\tilde a_x(v) = {\int_{\G^x} {\sigma(g)^{-1}\cdot a_{\t(g)}\cdot \varrho(g)\cdot v} \:\mathrm d\mu^x(g)}\text. & \bigl(\forall v &\in E_x\bigr)%
\end{equazione+}%
In particular, $\tilde a = a$ for all \G\nobreakdash-equivariant $a$.%
\end{proposizione}%

\separazione%

\noindent We will now derive a series of useful corollaries, which enter as key ingredients in many proofs throughout \refsez{N.20}.%

\begin{corollario}[Isotropy Extension Lemma]\label{N.iv5}%
Let \G\ be a proper Lie groupoid over a manifold $M$, and let $x_0 \in M$ be any point.%

Let $R = (E,\varrho)$ and $S = (F,\sigma)$ be \G\nobreakdash-actions of type \stack F and put $E_0 \equiv E_{x_0}$ and $F_0 \equiv F_{x_0}$. Moreover, let $A: E_0 \to F_0$ be a $G$\nobreakdash-equivariant linear map, where $G \equiv \G_0$ denotes the isotropy group of \G\ at $x_0$.%

Then there exists a morphism $a: R \to S$ in \R[F]\G\ such that $a_0 \equiv a_{x_0} = A$.%
\end{corollario}%
\begin{proof}%
Apply Lemma~\refcnt[N.15]{N.iii4} and then the Averaging Lemma to the morphism $A^S: (E|_S,\varrho|_S) \to (F|_S,\sigma|_S) \in \R[F]{\G|_{\mathnormal S}}$ \refequ{N.iv4}, where $S = {\G\cdot x_0}$. The corollary will follow from the formula~\refequ{N.ivB8} written at $x=x_0$.%
\end{proof}%

\begin{corollario}[Existence of Invariant Metrics]\label{N.iv6}%
Let \G\ be a proper Lie groupoid over a manifold $M$. Let $R = (E,\varrho) \in \R\G$ be a representation. Then there exists a metric $m: {R\otimes R^*} \to \TU$ in \R\G.%
\end{corollario}%
\begin{proof}%
Choose any metric $\phi: {E\otimes E^*} \to \TU$ in \stack[M]F (such metrics exist because \stack F satisfies Axiom~\refcnt[N.15]{N.iiiA4} and $M$ is paracompact); also fix any right-invariant Haar system $\mu$ on \G. By applying the averaging operator we obtain a morphism $\tilde\phi = \Av\mu(\phi): {R\otimes R^*} \to \TU$ in \R\G. We contend that $\tilde\phi$ is an invariant metric on $R$. It suffices to prove that for each $x \in M$ the induced form $\tilde\phi_x: {E_x\otimes {E_x}^*} \to \nC$ is a Hilbert metric (i.e.\ Hermitian and positive definite). Formula~\refequ{N.ivB8} reads%
\begin{equazione+}\label{N.iv7}%
\tilde\phi_x(v,w) = {\int_{\G^x} \bigsca{\varrho(g)v}{\varrho(g)w}\hspace{0pt}_\phi \:\mathrm d\mu^x(g)}\text, & \bigl(\forall v,w \in E_x\bigr)%
\end{equazione+}%
whence our claim is evident.%
\end{proof}%

Let $R = (E,\varrho)$ be any \G\nobreakdash-action. By a {\em \G\nobreakdash-invariant} section of $E$, defined over an invariant submanifold $N$ of the base $M$ of \G, we mean any section $\zeta \in \Gamma(N;E|_N)$ which is at the same time a morphism $\TU \to R|_N$ in $\R{\G|_{\mathnormal N}}$.%
\begin{corollario}[Invariant Sections]\label{iv7}%
Let $S$ be a closed invariant submanifold of the base $M$ of a proper Lie groupoid \G. Let $R = (E,\varrho) \in \R\G$ be a representation.%

Then each \G\nobreakdash-invariant section $\xi$ of $E$ over $S$ can be extended to a global \G\nobreakdash-invariant section; in other words, there exists some \G\nobreakdash-invariant $\Xi \in \Gamma(M;E)$ such that $\Xi|_S = \xi$.%
\end{corollario}%
\begin{proof}%
Apply Lemma~\refcnt[N.15]{N.iii4} and the Averaging Lemma.%
\end{proof}%

In general, we shall say that a partial function $\varphi: S \to \nC$, defined on an arbitrary subset $S \subset M$, is {\em smooth} when for each $x \in M$ one can find an open neighbourhood $B$ of $x$ in $M$ and a smooth function $B \to \nC$ that restricts to $\varphi$ over ${B\cap S}$.%
\begin{corollario}[Invariant Functions]\label{iv8}%
Let $S$ be any invariant subset of the base manifold $M$ of a proper Lie groupoid \G. Suppose $\varphi: S \to \nR$ is a smooth invariant function (i.e.\ $\varphi({g\cdot s}) = \varphi(s)$ for all $g$). Then there exists a smooth invariant function $\Phi: M \to \nR$ extending $\varphi$ outside $S$.%
\end{corollario}%
\begin{proof}%
Apply the Averaging Lemma to any smooth function extending $\varphi$ outside $S$ (such an extension can be obtained by means of a partition of unity over $M$, because of the smoothness of $\varphi$).%
\end{proof}%

\sezione{Fibre Functors}\label{N.18}
Let \stack F be a stack of complex smooth fields, to be regarded as fixed once and for all. Let $M$ be a paracompact smooth manifold.

\begin{definizione}\label{xviii.1}
By a \index{fibre functor}\index{fibre functor!over a manifold|see{base}}\textit{fibre functor (of type \stack F) over $M$,} or \index{base}{\em with base $M$,} we mean a faithful complex tensor functor
\begin{equazione}\label{N.v0}%
\fibfunc: \Kt \longto \stack[M]F\text,%
\end{equazione}%
of some additive complex tensor category \Kt, with values into \stack[M]F. We do {\em not} assume \Kt\ to be rigid.%
\end{definizione}%

When a fibre functor \fibfunc\ is assigned over $M$, one can construct a groupoid $\tanngpd{\fibfunc}$ having the points of $M$ as objects. Under reasonable assumptions, it is possible to make $\tanngpd{\fibfunc}$ a topological groupoid over the (topological) space $M$; the choice of a topology is dictated by the idea that the objects of \Kt\ should give rise to \gm{continuous representations} of $\tanngpd{\fibfunc}$ and that, vice versa, continuity of these representations should be enough to characterize the topological structure. An improvement of the same idea leads one to study a certain {\em functional structure} on $\tanngpd{\fibfunc}$, in the sense of \textit{Bredon (1972),} p.~297, and the important related problem of determining sufficient conditions for this functional structure to be compatible with the groupoid operations. Another fundamental issue here is to understand whether one gets in fact a {\em manifold structure}\footnote{A manifold can be defined as a topological space endowed with a functional structure locally looking like the structure of smooth real valued functions on some $\nR^d$.} making $\tanngpd{\fibfunc}$ a Lie groupoid over $M$; if this proves to be the case, we say that the fibre functor \fibfunc\ is \index{fibre functor!smooth}\index{smooth fibre functor}\textit{smooth.}

Some notation is needed first of all. Let $x$ be a point of $M$. If $x$ also denotes the (smooth) map $\pt \to M$, $\pt \mapsto x$, one can consider the complex tensor functor `fibre at $x$' which was introduced in \S\ref{N.11}%
\begin{equazione}\label{N.v1}%
\stack[M]F \to \VectorSpaces\text, \quad E \mapsto E_x \bydef {(x^*E)}_*\text.%
\end{equazione}%
Let $\fibfunc_x$ be the composite complex tensor functor%
\begin{equazione}\label{N.v2}%
\Kt \xto{\;\fibfunc\,} \stack[M]F \xto{\:{(\text{-})}_x} \VectorSpaces\text, \quad R \mapsto \fibfunc_x(R) \bydef \smash{\bigl(\fibfunc(R)\bigr)}_x\text.%
\end{equazione}%
Define the \textit{complex,} resp.\ \index{Tannakian groupoid!T(omega)@\tannakian\fifu|emph}\index{T(omega)@\tannakian{\fifu} (Tannakian groupoid associated with a fibre functor)|emph}\textit{real, Tannakian groupoid of \fibfunc} in the following way: for $x, x' \in M$, put
\begin{equazione}\label{N.v3}%
\left\{\begin{aligned}%
\tanngpd[\nC]\fibfunc(x,x') \;&\bydef \;\Iso^\otimes(\fibfunc_x,\fibfunc_{x'})%
\\%
\tanngpd[\nR]\fibfunc(x,x') \;&\bydef \;\Iso^{\otimes,\boldsymbol*}(\fibfunc_x,\fibfunc_{x'})\text.%
\end{aligned}\right.%
\end{equazione}%
(Recall that the right-hand side of the second equal sign denotes the set of all the \index{self-conjugate|emph}self-conjugate tensor preserving natural isomorphisms $\fibfunc_x \isoto \fibfunc_{x'}$, that is to say, the subset of $\Iso^\otimes(\fibfunc_x,\fibfunc_{x'})$ consisting of those $\lambda$ which make the following square commutative for each object $R \in \Ob(\Kt)$:%
\begin{equazione}\label{N.v4}%
\begin{split}%
\xymatrix@C=27pt{\fibfunc_x(R)^*\ar[d]_{\text{can.}}^\can\ar[r]^-{\lambda(R)^*} & \fibfunc_{x'}(R)^*\ar[d]_{\text{can.}}^\can \\ \fibfunc_x(R^*)\ar[r]^-{\lambda(R^*)} & \fibfunc_{x'}(R^*)\text{.)}\mspace{-13mu}}%
\end{split}%
\end{equazione}%
Setting $({\lambda'\lambda})(R) = {\lambda'(R)\circ\lambda(R)}$ and $x(R) = \id$, one obtains two groupoids over the set of points of $M$, with inverse given by $\lambda^{-1}(R) = \lambda(R)^{-1}$. We may also express \refequ{N.v3} in short by writing $\tanngpd[\nC]\fibfunc = \Aut^\otimes(\fibfunc)$ and $\tanngpd[\nR]\fibfunc = \Aut^{\otimes,\boldsymbol*}(\fibfunc)$.%

\sloppy%
Let us investigate the relationship between the complex tannakian groupoid \tanngpd[\nC]{\fibfunc} and its subgroupoid \tanngpd[\nR]{\fibfunc} first. As a convenient notational device, we omit writing \fibfunc\ when we simply refer to the set of arrows of the tannakian groupoid; thus for instance \tanngpd[\nC]{} is the set of arrows of the groupoid \tanngpd[\nC]\fibfunc. We define a map $\tanngpd[\nC]{} \to \tanngpd[\nC]{}$, $\lambda \mapsto \overline\lambda$, which we call complex conjugation, by setting $\overline\lambda(R) = \lambda(R^*)^*$; more precisely, $\overline\lambda(R)$ is defined by imposing the commutativity of%
\begin{equazione}\label{N.v5}%
\begin{split}%
\xymatrix@C=35pt@M=5pt{\fibfunc_x(R^*)^*\ar[d]^{\lambda(R^*)^*}\ar[r]^-\can & \fibfunc_x(R^{**})\ar[r]^-{\fibfunc_x(\can)} & \fibfunc_x(R)\ar@{-->}[d]^{\overline\lambda(R)} \\ \fibfunc_{x'}(R^*)^*\ar[r]^-\can & \fibfunc_{x'}(R^{**})\ar[r]^-{\fibfunc_{x'}(\can)} & \fibfunc_{x'}(R)\text.\!\!}%
\end{split}%
\end{equazione}%
It is straightforward to check that $\lambda \in \Hom^\otimes(\fibfunc_x,\fibfunc_{x'})$ implies $\overline\lambda \in \Hom^\otimes(\fibfunc_x,\fibfunc_{x'})$ and that $\lambda \mapsto \overline\lambda$ is a groupoid homomorphism of \tanngpd[\nC]\fibfunc\ into itself, identical on objects; this endomorphism is moreover involutive viz.\ $\overline{\overline\lambda} = \lambda$. Then we can characterize the arrows belonging to the subgroupoid \tanngpd[\nR]\fibfunc\ as the fixed points of the involution $\lambda \mapsto \overline\lambda$:%
\begin{equazione}\label{N.v6}%
\tanngpd[\nR]{} = \{\lambda \in \tanngpd[\nC]{}: \overline\lambda = \lambda\}\text.%
\end{equazione}%

\fussy
Next, we endow the set \tanngpd{} = \tanngpd[\nC]{} or \tanngpd[\nR]{} with a topology. In order to do this, we need to introduce the notion of \gm{metric} in \stack[M]F. Let $E$ be an object of \stack[M]F. A \index{metric}\textit{metric on $E$,} or {\em supported by $E$,} is a Hermitian form $\phi: {E\otimes E^*} \to \TU$ in \stack[M]F such that for all $x \in M$ the induced Hermitian form $\phi_x$ on the fibre $E_x$
\begin{equazione}\label{N.v7}%
{E_x\otimes {E_x}^*} \can ({E\otimes E^*})_x \xto{\phi_x} \TU_x \can \nC%
\end{equazione}%
is positive definite (and hence turns $E_x$ into a complex Hilbert space of finite dimension).%

We start by defining a collection \index{R@{\Rf} (collection of representative functions)}\repfunc\ of complex valued functions on \tanngpd{}, which we may call the \gm{\index{representative function}representative functions}. (Whenever we need to distinguish between \tanngpd[\nC]{} and \tanngpd[\nR]{}, we can write $\repfunc(\nC)$ or $\repfunc(\nR)$ as the case may be.)

Choose an object $R \in \Ob(\Kt)$, and let $\phi$ be a metric on the object $\fifu(R)$ of \stack[M]F. Also fix a pair of global smooth sections $\zeta, \zeta' \in \sections{({\fifu R})}(M)$. You get a complex function%
\begin{multiriga}\label{N.v8}%
r_{R,\phi,\zeta,\zeta'}: \tanngpd{} \to \nC\text,\quad \lambda\: \mapsto\: \smash{\bigsca{\lambda(R)\cdot\zeta(\s[\lambda])}{\zeta'(\t[\lambda])}}_\phi%
\\%
\equiv\: \phi_{\t(\lambda)}\bigl({\lambda(R)\cdot\zeta(\s[\lambda])},\zeta'(\t[\lambda])\bigr)\text.%
\end{multiriga}%
Then put%
\begin{multiriga}\label{N.v9}%
\repfunc = \bigl\{r_{R,\phi,\zeta,\zeta'}:\:\: R \in \Ob(\Kt)\text,\:\: \text{$\phi$ metric on $\fifu(R)$ in \stack[M]F,}%
\\%
\zeta, \zeta' \in \sections{({\fifu R})}(M)\bigr\}\text.%
\end{multiriga}%
We endow \tanngpd{} with the coarsest topology making all the functions in \repfunc\ continuous. From now on in our discussion \tanngpd[\nC]{} and \tanngpd[\nR]{} will always be regarded as topological spaces, with this topology. Observe that \tanngpd[\nR]{} turns out to be a subspace of \tanngpd[\nC]{}; more explicitly, the topology on \tanngpd[\nR]{} induced by $\repfunc(\nR)$ coincides with the topology induced from \tanngpd[\nC]{} along the inclusion $\tanngpd[\nR]{} \subset \tanngpd[\nC]{}$.%

We now want to establish a few fundamental algebraic properties of the collection \repfunc\ of complex valued functions on \tanngpd{}. We are going to show that \repfunc\ is a complex algebra of functions, and moreover that $\repfunc(\nR)$ is closed under taking the complex conjugate. Both assertions are immediate consequences of the following identities:%
\begin{elenco}%
\item For all smooth functions $a \in \C^\infty(M)$,%
\begin{equazione}\label{N.v10i}%
{({a\circ\s}) r_{R,\phi,\zeta,\zeta'}} = r_{R,\phi,{a\zeta},\zeta'} \quad \text{and} \quad {({a\circ\t}) r_{R,\phi,\zeta,\zeta'}} = r_{R,\phi,\zeta,{\overline a\zeta'}}\text;%
\end{equazione}%
in particular, if $c \in \nC$ is constant, $r_{R,\phi,{c\zeta},\zeta'} = {c\, r_{R,\phi,\zeta,\zeta'}} = r_{R,\phi,\zeta,{\overline c\zeta'}}$.%
\item If we let $\tau$ denote the metric on $\fibfunc(\TU)$ corresponding to the trivial metric ${\TU\otimes\TU^*} \can {\TU\otimes\TU} \can \TU$ on the object $\TU$ of \stack[M]F, and $1 \in \sections{(\fibfunc\TU)}(M)$ correspond to the ``unity section'' of $\TU \in \stack[M]F$ under the iso $\upsilon: \TU \isoto \fibfunc(\TU)$, then%
\begin{equazione}\label{N.v10ii}%
\text{`unity constant function'} = r_{\TU,\tau,1,1}\text.%
\end{equazione}%
\item For any choice of a direct sum $R \into {R\oplus S} \infrom S$ in \Kt,%
\begin{equazione}\label{N.v10iii}%
{r_{R,\phi,\zeta,\zeta'} + r_{S,\psi,\eta,\eta'}} = r_{{R\oplus S},{\phi\oplus\psi},{\zeta\oplus\eta},{\zeta'\oplus\eta'}}\text,%
\end{equazione}%
where ${\zeta\oplus\eta} \in \sections{(\fibfunc({R\oplus S}))}(M)$ etc.\ are obtained by setting $\fibfunc(R) \oplus \fibfunc(S) = \fibfunc({R\oplus S})$.%
\item Allowing the obvious (canonical) identifications,%
\begin{equazione}\label{N.v10iv}%
{r_{R,\phi,\zeta,\zeta'} r_{S,\psi,\eta,\eta'}} = r_{{R\otimes S},{\phi\otimes\psi},{\zeta\otimes\eta},{\zeta'\otimes\eta'}}\text.%
\end{equazione}%
(For instance, ${\zeta\otimes\eta}$ here denotes really the global section of $\fibfunc({R\otimes S})$ corresponding to the ``true'' ${\zeta\otimes\eta}$ in the iso $\tau_{R,S}: {\fifu(R)\otimes\fifu(S)} \isoto \fibfunc({R\otimes S})$.)%
\item Allowing again some loose notation,%
\begin{equazione}\label{N.v10v}%
\overline{r_{R,\phi,\zeta,\zeta'}} = {r_{R^*,\phi^*,\zeta^*,\zeta'^*} \circ \text{``$\lambda \mapsto \overline\lambda$''}}\text.%
\end{equazione}%
In particular, since the complex conjugation ``$\lambda \mapsto \overline\lambda$'' restricts to the identity on \tanngpd[\nR]{}, it follows that $\overline{r_{R,\phi,\zeta,\zeta'}} = r_{R^*,\phi^*,\zeta^*,\zeta'^*}$ in $\repfunc(\nR)$.%
\end{elenco}%
Notice that from the fact that $\repfunc(\nR)$ is closed under complex conjugation it follows immediately that the real and imaginary parts of any function in $\repfunc(\nR)$ will belong to $\repfunc(\nR)$ as well. Thus, if we let $\nR[\repfunc] \subset \repfunc(\nR)$ denote the subset of all the real valued functions, we can express $\repfunc(\nR) = {\nC\otimes{\nR[\repfunc]}}$ as the complexification of a real functional algebra.%

For the rest of the section\inciso{and for the purposes of the present thesis}we will only be interested in studying the real tannakian groupoid \tanngpd[\nR]\fibfunc. So from now on we forget about \tanngpd[\nC]\fibfunc\ and simply write \tanngpd\fibfunc\ for \tanngpd[\nR]\fibfunc. There is one further piece of structure we want to consider on \tanngpd\fibfunc, besides the topology.%

Let the sheaf of continuous (real valued) functions on an arbitrary topological space $T$ be denoted by \continuous[T]. Then recall that according to {\em Bredon (1972),} a ``\index{functionally structured space}functionally structured space'' is a topological space $T$, endowed with a \index{functional structure}sheaf of real algebras of continuous functions on $T$\nobreakdash---in other words, a subsheaf of algebras of \continuous[T]. A morphism of such ``functionally structured spaces'' is then defined as a continuous mapping such that the pullback of continuous functions along the mapping is compatible with the functional structures. For more details, we refer the reader to {\em loc.~cit.,} p.~297. We adopt this point of view in order to obtain a natural surrogate on \tanngpd[\nR]{} of the notion of \gm{smooth function}, drawing on the intuition that the representative functions should be regarded as the prototype \gm{smooth functions}.

It is obvious that if we start from the idea that the (real) representative functions are ``{smooth}'' then so we have to regard any function obtained by composing them with a smooth function $f: \nR^d \to \nR$. Define \index{R infinity@\ensuremath{\Rf^\infty} (canonical functional structure on the Tannakian groupoid)}\index{canonical functional structure|see{standard \ensuremath{\C^\infty}-structure}}\index{standard C infinity-structure@standard \ensuremath{\C^\infty}-structure \ensuremath{\Rf^\infty}|emph}$\repfunc^\infty$ to be the sheaf, of continuous real valued functions on the space \tanngpd{} = \tanngpd[\nR]{}, generated by the presheaf
\begin{multiriga}\label{N.v11}%
\Omega \:\mapsto \:\bigl\{f(r_1|_\Omega,\ldots,r_d|_\Omega):\:\: \text{$f: \nR^d \to \nR$ of class $\C^\infty$,}%
\\%
r_1, \ldots, r_d \in \nR[\repfunc]\bigr\}\text.%
\end{multiriga}%
In other words, $\repfunc^\infty$ is the smallest subsheaf of \continuous[\tanngpd{}] containing \refequ{N.v11} as a sub-presheaf. The expression $f(r_1|_\Omega,\ldots,r_d|_\Omega)$ denotes of course the function $\lambda \mapsto f\bigl(r_1(\lambda),\ldots,r_d(\lambda)\bigr)$, $\lambda \in \Omega$. Since \refequ{N.v11} is evidently a presheaf of \nR\nobreakdash-algebras of continuous functions on \tanngpd{}, $\repfunc^\infty$ will be a sheaf of such algebras and hence the pair $(\tanngpd{},\repfunc^\infty)$ will constitute a functionally structured space.%

Of course, we would like to say that the functional structure $\Rf^\infty$ on \tannakian{} is \index{C infinity-structured groupoid@\ensuremath{\C^\infty}-structured groupoid|emph}\index{functionally structured groupoid}compatible with the groupoid structure of \tannakian\fifu. This means that the structure maps of \tannakian{\fifu} should be all morphisms of functionally structured spaces, the base $M$ being regarded as such a space by means of its own sheaf of smooth real valued functions; in particular, the structure maps should be all continuous. What we are saying is not very precise, of course, unless we turn the space of composable arrows itself into a functionally structured space. Let us begin by observing that if $(X,\mathscr F)$ and $(Y,\mathscr G)$ are any functionally structured spaces then so is their Cartesian product endowed with the sheaf ${\mathscr F\otimes\mathscr G}$ locally generated by the functions $({\varphi\otimes\psi})(x,y) = {\varphi(x) \psi(y)}$. Then one can repeat the foregoing procedure to obtain, on ${X\times Y}$, a sheaf $({\mathscr F\otimes\mathscr G})^\infty$ \textit{of class $C^\infty$,} i.e.\ closed under composition with arbitrary smooth functions as in \refequ{N.v11}. Any subspace $S \subset {X\times Y}$ may be finally regarded as a functionally structured space by endowing it with the induced sheaf $({\mathscr F\otimes\mathscr G})^\infty|_S \bydef {i_S}^*[({\mathscr F\otimes\mathscr G})^\infty]$, where $i_S$ denotes the inclusion mapping of $S$ into ${X\times Y}$. (Recall that if $f: S \to T$ is any continuous mapping into a functionally structured space $(T,\mathscr T)$ then ${f^*\mathscr T}$ is the functional sheaf on $S$ associated with the presheaf
\begin{displaymath}%
U \mapsto \inductivelim{V\supset f(U)}{\mathscr T(V)}\text{.~}\bigr)%
\end{displaymath}%
Notice that in case $X$ and $Y$ are smooth manifolds and $S \subset {X\times Y}$ is a submanifold, one recovers the correct functional structures: $({\smooth[X] \otimes \smooth[Y]})^\infty = \smooth[{X\times Y}]$ and $\smooth[{X\times Y}]|_S^{} = \smooth[S]$. It is therefore perfectly reasonable to endow the space of composable arrows $\mca[2]{\tannakian{}} = \mca{\tannakian{}}$ with the functional structure $\smash{\mca[2]\Rf}^{,\infty} \bydef ({\Rf^\infty \otimes \Rf^\infty})^\infty|_{\mca[2]{\tannakian{}}}$ and to call the composition map $\c: \mca[2]{\tannakian{}} \to \tannakian{}$ ``smooth'' whenever it is a morphism of such functionally structured space into $(\tannakian{},\Rf^\infty)$.%

Later on we will show that \tannakian\fifu\ is actually a functionally structured groupoid in the two cases of major interest for us, namely when \fifu\ is the standard fibre functor \forget\G\ associated with a proper Lie groupoid (\refsez{N.20}) or when \fifu\ is a \gm{classical} fibre functor (\refsez{N.21}). However, we can already very easily check the ``smoothness'' (in particular, the continuity) of some of the structure maps:%
\newline%
\textsl{(a)~The source map $\s: \tannakian{} \to M$.}~First of all observe that for an arbitrary $a \in \C^\infty(M)$ we have ${a\circ\s} \in \Rf$, by \refequ{N.v10i} and \refequ{N.v10ii}. Let $U \subset M$ be open. For each $u \in U$ there exists $f_u \in \C^\infty(M)$ with $\support{f_u} \subset U$ and $f_u(u) = 1$. Since ${f_u\circ\s} \in \Rf$, the subset $({f_u\circ\s})^{-1}(\nC_{\scriptscriptstyle\neq0}) \subset \tannakian{}$ must be open. Now $({f_u\circ\s})^{-1}(\nC_{\scriptscriptstyle\neq0}) = \s^{-1}\bigl(\smash{f_u}^{-1}(\nC_{\scriptscriptstyle\neq0})\bigr) \subset \s^{-1}(U)$, so $\s^{-1}(U)$ can be expressed as a union of open subsets of \tannakian{} and therefore it is open. This shows that \s\ is continuous; since ${a\circ\s} \in \nR[\Rf]$ whenever $a$ is real valued, it also follows that \s\ is a morphism of functionally structured spaces.%
\newline%
\textsl{(b)~The target map $\t: \tannakian{} \to M$.}~The discussion here is entirely analogous, starting from the other identity ${a\circ\t} = r_{\TU,\tau,1,\overline a1} \in \Rf$.%
\newline%
\textsl{(c)~The unit section $\u: M \to \tannakian{}$.}~This time let $r = r_{R,\phi,\zeta,\zeta'} \in \Rf$ be given; we must show that ${r\circ\u} \in \C^\infty(M)$. This is trivial because%
$$%
({r\circ\u})(x) = \smash{\bigsca{x(R)\cdot\zeta(x)}{\zeta'(x)}}_\phi = \smash{\bigsca{\zeta(x)}{\zeta'(x)}}_\phi = \smash{\scalare\zeta{\zeta'}}_\phi(x)\text.%
$$%

Finally, let us remark that, as a consequence of the existence of metrics on any object of \stack[M]F (because \stack F is a stack of smooth fields and $M$ admits partitions of unity), \textit{the space \tannakian{} of arrows of \tannakian\fifu\ is always Hausdorff.} Indeed, let $\mu \neq \lambda \in \tannakian{}$. We can assume $\s(\mu) = x = \s(\lambda)$ and $\t(\mu) = x' = \t(\lambda)$ otherwise we are immediately done by using the Hausdorffness of $M$ and the continuity of either the source or the target map. Then there exists $R \in \Ob(\Kt)$ with $\mu(R) \neq \lambda(R)$. Choose any metric $\phi$ on $\fifu(R)$ (there is at least one): since $\phi_{x'}$ is in particular non-degenerate on $E_{x'}$, there will be {\em global}\inciso{again, because of the existence of partitions of unity}sections $\zeta, \zeta' \in \sections{(\fifu R)}(M)$ with%
$$%
z_\mu \,= \,\smash{\bigl\langle{\mu(R)\cdot\zeta(x)},\zeta'(x')\bigr\rangle}_\phi \,\neq \,\smash{\bigl\langle{\lambda(R)\cdot\zeta(x)},\zeta'(x')\bigr\rangle}_\phi \,= \,z_\lambda\text.%
$$%
Let $D_\mu, D_\lambda \subset \nC$ be disjoint open disks about $z_\mu, z_\lambda$ respectively. Then, setting $r = r_{R,\phi,\zeta,\zeta'}^{}$, the inverse images $r^{-1}(D_\mu)$ and $r^{-1}(D_\lambda)$ will be disjoint open neighbourhoods of $\mu$ and $\lambda$ in \tannakian{}.%

\sezione{Properness}\label{N.19}
We shall say that a metric $\phi$ on the object $\fifu(R)$, $R \in \Ob(\Kt)$ of \stack[M]F is \index{metric!omega-invariant@\fifu-invariant|emph}\index{omega-invariant metric@\fifu-invariant metric|emph}\textit{\fifu\nobreakdash-invariant,} when there exists a Hermitian form $m: {R\otimes R^*} \to \TU$ in \Kt\ such that $\phi$ coincides with the induced Hermitian form
\begin{equazione}\label{N.v12}%
{\fifu(R)\otimes \fifu(R)^*} \can \fifu({R\otimes R^*}) \xto{\fifu(m)} \fifu(\TU) \can \TU\text.%
\end{equazione}%
We express this in short by writing $\phi = {\fifu_*^{}m}$. Note that by the faithfulness of \fifu\ there is at most one such $m$.%

\begin{definizione}\label{xix.1}%
A fibre functor $\fifu: \Kt \longto \stack[M]F$ will be called \index{fibre functor!proper|emph}\index{proper fibre functor|emph}\textit{proper} if
\begin{elenco}%
\item the continuous mapping $\,(\s,\t): \tannakian{} \to {M\times M}\,$ is proper, and%
\item for every object $R \in \Ob(\Kt)$, the object $\fifu(R)$ of \stack[M]F supports an \fifu\nobreakdash-invariant metric.%
\end{elenco}%
We can express the second condition more succinctly by saying that \gm{there are enough \fifu\nobreakdash-invariant metrics}.%
\end{definizione}%

\begin{esempio}\label{N.v13}
Let \fifu\ be the standard functor \index{omega(G)@\forget[T]\G, \forget[\infty]{\G} (forgetful functor)}\index{standard fibre functor!type T@(of type \stack T) \forget[T]\G}$\forget\G: \R\G \longto \stack[M]F$, of type \stack F, associated with a proper Lie groupoid \G\ over $M$. Then \fifu\ is a proper fibre functor.

In order to see this, observe (cfr.\ also \S\ref{N.20}) that there is an obvious \index{envelope homomorphism@envelope homomorphism \envelope[T]\G, \envelope[\infty]\G}\index{pi(G)@\envelope[T]\G, \envelope[\infty]{\G} (envelope homomorphism)}homomorphism of groupoids
\begin{equazione}\label{N.v14}
\envehom: \G \longto \tannakian\G \bydef \tannakian{\forget\G}\text,%
\end{equazione}%
identical on the base, called the \gm{\stack F-envelope homomorphism} of \G\ and defined by setting $\envelope{g}(R) = \varrho(g)$ for each object $R = (E,\varrho)$ of \R\G; the notation $\varrho(g)$ was introduced in \S\ref{N.17}. The mapping $\mca[1]\envehom: \mca[1]\G \to \mca[1]{\tannakian{}}$ is continuous. Indeed, if we fix any representative function $r = r_{R,\phi,\zeta,\zeta'} \in \repfunc$, let us say with $R = (E,\varrho)$, and a small open subset $\Gamma \subset \G$ on which we have, for $\varrho$ acting on $\zeta$, the sort of expression
\begin{displaymath}%
{\varrho(g) \cdot \zeta(\s[g])} = \txtsum{i=1}\ell{r_i(g) \zeta'_i(\t[g])}\text, \quad r_i \in \C^\infty(\Gamma)%
\end{displaymath}%
we derived in \S\ref{N.17}, then for all $g \in \Gamma$ we obtain%
\begin{displaymath}%
({r\circ\envehom})(g) = \smash{\bigsca{\envehom(g)(R) \cdot \zeta(\s[g])}{\zeta'(\t[g])}}_\phi = \txtsum{i=1}\ell{r_i(g){\bigl(\smash{\scalare{\zeta'_i}{\zeta}}_\phi \circ \t\bigr)}(g)}\text.%
\end{displaymath}%
Therefore, we conclude that ${r\circ\envehom} \in \C^\infty(\G)$ and hence, in particular, that ${r\circ\envehom}$ is continuous. Note that in fact this argument shows that the map \envehom\ is a morphism of functionally structured spaces, of $(\G,\smooth[\G])$ into $(\tannakian{},\repfunc^\infty)$. We will prove in \S\ref{N.20} that the envelope mapping \envehom\ is also surjective; the properness of $\,(\s,\t): \tannakian{} \to {M\times M}\,$ is now a trivial consequence of this fact and the properness of $\,(\s,\t): \G \to {M\times M}$. The existence of enough invariant metrics was established in \S\ref{N.17} as a corollary to the Averaging Lemma.%
\end{esempio}%

Back to general notions, it turns out that in order to characterize the topology of \tanngpd{} the \fifu\nobreakdash-invariant metrics are (for \fifu\ proper) as good as the generic, `not necessarily invariant' ones. More exactly, let $\repfunc' \subset \repfunc$ be the set of all the representative functions $r_{R,\phi,\zeta,\zeta'}$ with $\phi = {\fifu_*^{}m}$ an \fifu\nobreakdash-invariant metric on $\fifu(R)$. Note that $\repfunc'$ is a subalgebra of $\repfunc$, closed under complex conjugation; this follows from the identities proved above, by observing that ${\fifu_*^{}m\otimes\fifu_*^{}n} = \fifu_*^{}({m\otimes n})$ and so on. Then we claim that%
\begin{lemma}\label{N.v15}%
The topology on \tanngpd{} is also the coarsest making all the functions in $\repfunc'$ continuous.%
\end{lemma}%
\begin{proof}%
Recall that the topology on \tanngpd{} was defined as the coarsest making all the functions belonging to $\repfunc$ continuous. We have already observed that $\repfunc'$ is an algebra of continuous complex functions on \tanngpd{}, closed under conjugation. Moreover, it separates points, because of the existence of enough \fifu\nobreakdash-invariant metrics, cf.\ the argument used to prove Hausdorffness of \tanngpd{}. Henceforth, for every open subset $\Omega \subset \tanngpd{}$ with compact closure $\overline\Omega$, the involutive subalgebra $\repfunc'_{\overline\Omega} \subset \C^0(\overline\Omega;\nC)$, formed by the restrictions to $\overline\Omega$ of elements of $\repfunc'$, is sup-norm dense in $\C^0(\overline\Omega)$ and a fortiori in $\repfunc_{\overline\Omega} = \left\{r|_{\overline\Omega}: r \in \repfunc\right\}$, as a consequence of the Stone\nobreakdash--Weierstrass theorem.%

This remark applies in particular to $\Omega = \tanngpd{}|_{U\times U'}$ where $U$ and $U'$ are open subsets of $M$ with compact closure. (Here is where we use the properness of $\tanngpd{} \xto{\scriptscriptstyle(\s,\t)} {M\times M}$.) Note that the subset $\tanngpd{}|_{U\times U'}$ is also open in the space $\tanngpd{}'$ = \gm{\tanngpd[\nR]{} with the topology generated by $\repfunc'$} because $\tanngpd{}' \xto{\scriptscriptstyle(\s,\t)} {M\times M}$ is clearly still continuous. Since the subsets $\tanngpd{}|_{U\times U'}$ cover \tanngpd{}, we are now reduced to showing that the identity mappings%
$$%
\tanngpd{}|_{U\times U'} \,\xto= \,\tanngpd{}'|_{U\times U'}%
$$%
are homeomorphisms.%

To simplify the notation, we reformulate our claim as follows: given a subset $\Omega \subset \tanngpd[\nR]{}$, open in both \tanngpd{} and $\tanngpd{}'$ and with compact closure in \tanngpd{}, show that the identity mapping $\Omega' \xto= \Omega$ is continuous (here $\Omega'$ denotes of course the open subset, viewed as a subspace of $\tanngpd{}'$). Notice that the topology on $\Omega$ generated by the collection of functions $\repfunc_\Omega = \left\{r|_\Omega: r \in \repfunc\right\}$ coincides with the subspace topology induced from \tanngpd{}. Then, let $r \in \repfunc$ be fixed; since $\overline\Omega$ is compact in \tanngpd{}, the restriction $r|_\Omega$ will be, as remarked at the beginning, a uniform limit of continuous functions on $\Omega'$ and hence itself a continuous function on $\Omega'$.%
\end{proof}%
We shall make implicit use of the lemma throughout the rest of the present subsection.%

Another easy, although important observation is that all $\lambda \in \tanngpd[\nR]{}$ will act unitarily under any \fifu\nobreakdash-invariant metric. More precisely, for any object $R \in \Ob(\Kt)$ and any \fifu\nobreakdash-invariant metric $\phi$ on $\fifu(R)$, the linear isomorphism $\lambda(R)$ will preserve the inner product $\smash{\scalare{}{}}_\phi$:%
\begin{equazione}\label{N.v16}%
\smash{\bigl\langle{\lambda(R)\cdot v},{\lambda(R)\cdot v'}\bigr\rangle}_\phi \,= \,\smash{\scalare{v}{v'}}_\phi\text.%
\end{equazione}%
We use this observation to prove the following%
\begin{proposizione}\label{N.v17}%
Let \fifu\ be a proper fibre functor. Then \tanngpd{\fifu} is a (Hausdorff, proper) topological groupoid.%
\end{proposizione}%
\begin{proof}%
We must show that the inverse and composition maps of \tanngpd{\fifu} are continuous.%

\textsl{a)}~Continuity of the inverse map $\i: \tanngpd{} \to \tanngpd{}$. It must be proved that the composite ${r\circ\i}$ is continuous on \tanngpd{}, for any $r = r_{R,\phi,\zeta,\zeta'} \in \repfunc$ with $\phi$ an \fifu\nobreakdash-invariant metric on $\fifu(R)$. This is immediate, because%
\begin{align*}%
({r_{R,\phi,\zeta,\zeta'}^{} \circ \i})(\lambda) \;&= \;\smash{\bigsca{{\lambda(R)}^{-1} \cdot \zeta(\t\lambda)}{\zeta'(\s\lambda)}}_\phi \;= \;\smash{\bigsca{\zeta(\t\lambda)}{\lambda(R) \cdot \zeta'(\s\lambda)}}_\phi%
\\%
&= \;\overline{\bigsca{\lambda(R) \cdot \zeta'(\s\lambda)}{\zeta(\t\lambda)}\mspace{0mu}_\phi} \;= \;\overline{r_{R,\phi,\zeta',\zeta}}\,(\lambda)\text,%
\end{align*}%
in view of \refequ{N.v16}.%

\textsl{b)}~Continuity of composition $\c: \mca{\tanngpd{}} \to \tanngpd{}$ (the domain of the map being topologized as a subspace of the cartesian product ${\tanngpd{} \times \tanngpd{}}$). We make a technical observation first.%

\sloppy%
Fix $\lambda \in \tanngpd{}$, let us say $\lambda: x \to x'$. Let $R \in {\Ob\:\Kt}$ and let $\phi$ be any \fifu\nobreakdash-invariant metric on $E = \fifu(R)$. Fix a local $\phi$\nobreakdash-orthonormal system $\zeta'_1, \ldots, \zeta'_d \in \sections{(\fifu R)}(U')$ for $E$ about $x'$ as in \refequ[N.15]{N.iii7}; hence, in particular,%
\begin{equazione}\label{N.v18}%
E_{x'} = \mathrm{Span}\,\{\zeta'_1(x'), \ldots, \zeta'_d(x')\}\text.%
\end{equazione}%
Since $M$ is paracompact, it is no loss of generality to assume that for every $i = 1, \ldots, d$ $\zeta'_i$ is the restriction to $U'$ of a global section $\zeta_i$ of $\fifu(R)$. Let $\zeta \in \sections{(\fifu R)}(M)$ be another global section. Consider an open neighbourhood $\Omega$ of $\lambda$ in \tanngpd{} such that $\t(\Omega) \subset U'$. Also let $\Phi_i \in \C^0(\Omega;\nC)$ $(i = 1, \ldots, d)$ be a list of continuous complex functions on $\Omega$. Then the norm function%
\begin{equazione}\label{N.v19}%
\mu \;\mapsto \;\modulo{{\mu(R) \cdot \zeta({\s\mu})} - \txtsum{i=1}d{\Phi_i^{}(\mu) \zeta'_i({\t\mu})}}%
\end{equazione}%
is certainly continuous on $\Omega$: indeed, its square is%
\begin{align*}%
&\phantom= \bigmod{\mu(R) \zeta(\s\mu)}^2 - 2\displaysum{i}{}{\parteRe{\left[\overline{\Phi_i^{}(\mu)} \bigsca{\mu(R) \zeta(\s\mu)}{\zeta'_i(\t\mu)}\right]}} + \modulo{\txtsum{i=1}{d}{\Phi_i^{}(\mu) \zeta'_i(\t\mu)}}^2%
\\%
&= \bigmod{\zeta(\s\mu)}^2 - 2\displaysum{i}{}{\parteRe{\left[\overline{\Phi_i(\mu)} \bigsca{\mu(R) \zeta(\s\mu)}{\zeta_i(\t\mu)}\right]}} + \txtsum{i=1}d{\bigmod{\Phi_i(\mu)}^2}%
\end{align*}%
(because $\mu(R)$ is unitary \refequ{N.v16} and the vectors $\zeta'_i(\t\mu)$, $i = 1, \ldots, d$ form an orthonormal system in $E_{\t(\mu)}$). Now, make $\Phi_i(\mu) = \bigsca{\mu(R) \zeta(\s\mu)}{\zeta_i(\t\mu)}$ in \refequ{N.v19} and evaluate the function you get at $\mu = \lambda$: the result will be zero, because the vectors $\zeta_i(x')$, $i = 1, \ldots, d$ constitute an orthonormal {\em basis} of $E_{x'}$. Hence, by the just observed continuity, for each $\varepsilon > 0$ there will be an open neighbourhood of $\lambda$ in \tanngpd{}, let us call it $\Omega^\varepsilon(\lambda)$, over which the following estimate holds%
\begin{equazione}\label{N.v20}%
\modulo{{\mu(R) \cdot \zeta(\s\mu)} - \txtsum{i=1}d{r_{R,\phi,\zeta,\zeta_i}^{}(\mu) \zeta_i^{}(\t\mu)}} < \varepsilon\text.%
\end{equazione}%

\fussy%
With this preliminary observation at hand it is easy to show continuity of the composition of arrows. Indeed, consider an arbitrary object $R \in {\Ob\:\Kt}$, an arbitrary \fifu\nobreakdash-invariant metric $\phi$ on $\fifu(R)$, and arbitrary global sections $\zeta, \eta \in \sections{(\fifu R)}(M)$. We have to check the continuity of the function%
\begin{equazione}\label{N.v21}%
(\mu',\mu) \;\mapsto \;({r_{R,\phi,\zeta,\eta}^{} \circ \c})(\mu',\mu) \,= \,\smash{\bigsca{\mu'(R) \cdot \mu(R) \cdot \zeta(\s\mu)}{\eta(\t\mu')}}_\phi%
\end{equazione}%
on the space of composable arrows \mca[2]{\tanngpd{}}. Let $x \xto\lambda x' \xto{\lambda'} x''$ be an arbitrary pair of composable arrows, which we regard as fixed. Choose a local $\phi$\nobreakdash-orthonormal system about $x'$ as before. Then, by the estimate \refequ{N.v20} and our remark \refequ{N.v16} that $\mu'(R)$ is unitary, for all $(\mu',\mu)$ close enough to $(\lambda',\lambda)$, let us say for $\mu \in \Omega^\varepsilon(\lambda)$, the function \refequ{N.v21} will differ from the function%
$$%
\txtsum{i=1}d{r_{R,\phi,\zeta,\zeta_i}^{}(\mu) \,\smash{\bigsca{\mu'(R)\cdot\zeta_i^{}(\s\mu')}{\eta(\t\mu')}}_\phi} \:= \:\txtsum{i=1}d{r_{R,\phi,\zeta,\zeta_i}^{}(\mu) \,r_{R,\phi,\zeta_i,\eta}^{}(\mu')}%
$$%
up to ${\varepsilon\norma\eta}$, where \norma\eta\ is a positive bound for the norm of $\eta$ in a neighbourhood of $x''$. This proves the desired continuity, because the last function is certainly continuous on ${\tanngpd{} \times \tanngpd{}}$ and hence on \mca[2]{\tanngpd{}}.%
\end{proof}%

\sezione{Reconstruction Theorems}\label{N.20}
When applying the formal apparatus of \refsez{N.18} to the standard fibre functor \forget[F]{\G} associated with a Lie groupoid \G, we prefer to use the alternative notation \index{T(G)@\tannakian[T]\G, \tannakian[\infty]{\G} (Tannakian groupoid associated with a Lie groupoid)|emph}\index{Tannakian groupoid!T(G) type T@(of type \stack T) \tannakian[T]\G|emph}\tannakian[F]{\G} for the real Tannakian groupoid \xtannakian{\forget[F]\G;\nR} and refer to the latter as the \index{envelope|see{Tannakian groupoid}}\textit{(real) \stack F-envelope} of \G. If explicit mention of type is not necessary, we normally just write \tannakian\G.

The \index{envelope homomorphism@envelope homomorphism \envelope[T]\G, \envelope[\infty]\G|emph}\index{pi(G)@\envelope[T]\G, \envelope[\infty]{\G} (envelope homomorphism)|emph}\textit{\stack F-envelope homomorphism} associated with a Lie groupoid \G\ is the groupoid homomorphism $\envelope{}: \G \to \tannakian\G$, or, more pedantically,
\begin{equazione}\label{O.esp86}%
\envelope[F]\G: \G \longto \tannakian[F]\G%
\end{equazione}%
defined by the formula ${\envelope{g}(E,\varrho)} \bydef \varrho(g)$. (Having a look at Note \refcnt[N.17]{N.iv1} one more time might be useful at this point.) The study of properties of the envelope homomorphism \envelope\G\ for proper \G\ will constitute our main concern in this section.

Let ${M/\G}$ be the topological space obtained by endowing the set of \index{orbit}orbits $\{{\G\cdot x}| x \in M\}$ with the quotient topology induced by the \index{orbit map, space}\textit{orbit map}
\begin{equazione}\label{O.esp91}
\mathit o: M \to {M/\G}
\end{equazione}
(the map sending a point $x$ to the respective \G\nobreakdash-orbit $\mathit o(x) = {\G\cdot x}$). Note that the map $\mathit o$ is open: indeed, if $U \subset M$ is an open subset then so is $\mathit o^{-1}(\mathit o(U)) = \t(\s^{-1}(U))$ because \t\ is an open map\nobreakdash---actually, a submersion. Furthermore, ${M/\G}$ is a locally compact Hausdorff space. Indeed, suppose $\G(x,x')$ empty. Properness of \G, applied to some sequence of balls ${B_i\times {B_i}'}$ shrinking to the point $(x,x')$, will yield open balls $B, B' \subset M$ at $x, x'$ such that $(\s,\t)^{-1}({B\times B'})$ is empty, in other words, such that ${\mathit o(B)\cap\mathit o(B')} = \varnothing$, as contended. In particular, every orbit ${\G\cdot x} = \mathit o^{-1}\{\mathit o(x)\}$ is a closed subset of $M$.%
\begin{theorem}\label{O.thm3}%
Let \stack F be any stack of smooth fields. Let \G\ be a proper Lie groupoid. Then the \stack F\nobreakdash-envelope homomorphism $\envelope[F]\G: \G \to \tannakian[F]\G$ is full (i.e.\ surjective, as a mapping of the spaces of arrows).%
\end{theorem}%
\begin{proof}%
To begin with, let us prove that $\G(x,x')$ empty implies $\tannakian\G(x,x')$ empty. Put $S = {{\G x}\cup {\G x'}}$ and let $\varphi: S \to \nC$ be the function which takes the value $1$ over the orbit ${\G x}$ and the value $0$ over the orbit ${\G x'}$; $\varphi$ is well-defined because ${{\G x}\cap {\G x'}} = \varnothing$. $S$ is an invariant submanifold of $M$. Since $S$ is the union of two disjoint closed subsets of $M$, it is also a closed submanifold. Moreover, $\varphi$ is equivariant with respect to the trivial representation of \G, i.e.\ $\varphi({g\cdot s}) = \varphi(s)$. Corollary \refcnt[N.17]{iv8} says that there is some smooth invariant function $\Phi: M \to \nC$, extending $\varphi$, equivalently, some smooth function $\Phi$ on $M$, constant along the \G\nobreakdash-orbits and with $\Phi(x) = 1$, $\Phi(x') = 0$. By setting $b_z \bydef {\Phi(z)\mspace{1.5mu}\id}$, one gets an endomorphism $b$ of the trivial representation with $b_x = \id$ and $b_{x'} = 0$. Now, suppose there exists some $\lambda \in \tannakian\G(x,x')$: then, by the naturality of $\lambda$, one gets a commutative square%
\begin{displaymath}%
\xymatrix@R=17pt{\nC\ar[d]^\id\ar[r]^\lambda & \nC\ar[d]^0 \\ \nC\ar[r]^\lambda & \nC}%
\end{displaymath}%
which contradicts the invertibility of $\lambda(\TU)$.%

In order to finish the proof of the theorem, it will be sufficient to prove surjectivity of all isotropy homomorphisms induced by \envelope{}, because%
\begin{displaymath}%
\xymatrix@C=31pt{\G|_x\ar[d]^{g\,\text-}_\iso\ar[r]^-{\envelope{}_x} & {\tannakian\G}|_x\ar[d]^{\envelope g\,\text-}_\iso \\ \G(x,x')\ar[r]^-{\envelope{}_{x,x'}} & {\tannakian\G}(x,x')}%
\end{displaymath}%
commutes for all $g \in \G(x,x')$. More explicitly, it will be sufficient to prove that $\envelope{}_x: \G|_x \to {\tannakian\G}|_x$ is an epimorphism of groups, for every $x \in M$. This follows immediately from Proposition \refcnt[O.1.4+.3.3]{O.prp5} and Corollary \refcnt[N.17]{N.iv5}.%
\end{proof}%
\separazione%
We continue to work with an arbitrary stack of smooth fields.%

\begin{definizione}\label{O.dfn3}
A Lie groupoid \G\ will be said to be \index{reflexive|emph}\index{groupoid!reflexive|emph}\textit{\stack F-reflexive,} or \index{self-dual|see{reflexive}}\index{groupoid!self-dual|see{reflexive}}\textit{self-dual relative to \stack F,} if its \stack F-envelope homomorphism $\envelope[F]\G: \G \to \tannakian[F]\G$ is an isomorphism of topological groupoids.
\end{definizione}

It turns out, for proper Lie groupoids, that the requirement that the continuous mapping $\mca[1]{\envelope{}}: \mca[1]\G \to \mca[1]{\tannakian\G}$ should be open is superfluous; more precisely, one has the following statement:%
\begin{theorem}\label{O.thm1}%
Let \G\ be a proper Lie groupoid. Let \stack F be any stack of smooth fields. Then \G\ is \stack F\nobreakdash-reflexive if and only if the homomorphism \envelope[F]\G\ is faithful (i.e.\ injective, as a mapping of the spaces of arrows).%
\end{theorem}%
\begin{proof}%
The assertion that injectivity implies bijectivity, or, to say the same thing differently, that faithfulness implies full faithfulness, is an immediate consequence of Theorem \refcnt{O.thm3} above.%

As to the statement that the mapping \envelope{} is open, we have to show that whenever $\Gamma$ is an open subset of \mca[1]{\G} and $g_0$ a point of $\Gamma$, the image \envelope{\Gamma} is a neighbourhood of \envelope{g_0} in \tannakian{\G}.%

To fix ideas, suppose $g_0 \in \G(x_0,x_0')$. Let us start by observing that, as in the proof of Proposition~\refcnt[O.1.4+.3.3]{O.prp5}, it is possible to find a representation $R = (E,\varrho) \in {\Ob\,\R\G}$ whose associated $x_0$\nobreakdash-th isotropy homomorphism $\varrho_0: G_0 \to \GL(E_0)$ is injective (same notation as in Eq.~\refequ[N.17]{N.iv3}); for such an $R$, the map $\G(x_0,x_0') \to \Lis(E_{x_0},E_{x_0'})$, $g \mapsto \varrho(g)$ is also injective. We regard $R$ as fixed once and for all. Moreover, we choose an arbitrary Hilbert metric $\phi$ on $E$. As we know from \refsez{N.15}, there are local $\phi$\nobreakdash-orthonormal frames for $E$%
\begin{equazione}\label{xx3}%
\left\{\begin{aligned}%
\zeta_1, \ldots, \zeta_d &\in \sections E(U) &\quad &\text{about $x_0$ and}%
\\[.4\medskipamount]%
\zeta_1', \ldots, \zeta_d' &\in \sections E(U') &\quad &\text{about $x_0'$;}%
\end{aligned}\right.%
\end{equazione}%
their cardinality turns out to be the same because $E_{x_0} \iso E_{x_0'}$. Since $M$ is paracompact, it is no loss of generality to assume that the $\zeta_i$ and the $\zeta_{i'}'$ are (restrictions of) global sections. Finally, we select any compactly supported smooth functions $a, a': M \to \nC$ with $\support a \subset U$ and $\support{a'} \subset U'$, such that $a(x) = 1 \aeq x = x_0$ and $a'(x') = 1 \aeq x' = x_0'$.%

Let us put, for all $1 \leqq i, i' \leqq d$,%
\begin{equazione+}\label{xx4}%
\varrho_{i,i'} = {r_{i,i'}\circ \envelope{}} \bydef {r_{R,\phi,\zeta_i,\zeta_{i'}'} \circ \envelope{}}: \G \to \nC\text, & \text{[using notation \refequ[N.18]{N.v8}]}%
\end{equazione+}%
and for $i=0$ and $0\leqq i'\leqq d$, resp.\ $0\leqq i\leqq d$ and $i'=0$,\footnote{For $i=i'=0$ either choice will do; for $d=0$ there are obvious modifications which we leave to the reader. The only thing that really matters is that both ${a\circ\s}$ and ${a'\circ\t}$ should occur in the intersection (\ref{xx6}) at least once.}%
\begin{equazione}\label{xx5}%
\left\{\begin{aligned}%
\varrho_{0,i'} = {r_{0,i'}\circ \envelope{}} &\bydef {a\circ \s_{\G}} = {\bigl({a\circ \s_{\mspace{1.8mu}\tannakian\G}}\bigr)\circ \envelope{}}: \G \to \nC\text, \quad\text{resp.}%
\\[.3\medskipamount]%
\varrho_{i,0} = {r_{i,0}\circ \envelope{}} &\bydef {a'\circ \t_{\G}} = {\bigl({a'\circ \t_{\mspace{2.2mu}\tannakian\G}}\bigr) \circ \envelope{}}: \G \to \nC\text.%
\end{aligned}\right.%
\end{equazione}%
Also, put $z_{i,i'} = \varrho_{i,i'}(g_0) \in \nC$. We claim that, as a consequence of properness, there exist open disks $D_{i,i'} \subset \nC$ centred at $z_{i,i'}$ such that%
\begin{equazione}\label{xx6}%
\displaycap{0\leqq i,i'\leqq d}{}{{\varrho_{i,i'}}^{-1}(D_{i,i'})} \subset \Gamma\text.%
\end{equazione}%
Before we go into the proof of this claim, let us show how the statement that \envelope{\Gamma} is a neighbourhood of \envelope{g_0} follows from \refequ{xx6}. Since, by Theorem~\refcnt{O.thm3}, \envelope{} is surjective as a mapping of \mca[1]\G\ into \mca[1]{\smash{\tannakian\G}}, we have%
\begin{multline}\notag%
\displaycap{}{}{{r_{i,i'}}^{-1}\bigl(D_{i,i'}\bigr)} = {\envelope{}\envelope{}^{-1}\left(\smash[b]{\displaycap{}{}{{r_{i,i'}}^{-1}\bigl(D_{i,i'}\bigr)}}\right)} = {\envelope{}\left(\smash[b]{\displaycap{}{}{\envelope{}^{-1}{r_{i,i'}}^{-1}\bigl(D_{i,i'}\bigr)}}\right)}%
\\%
= \xenvelope{\smash[b]{\displaycap{}{}{{\varrho_{i,i'}}^{-1}\bigl(D_{i,i'}\bigr)}}} \subset \envelope{\Gamma}\text. \qquad\text{(by the inclusion \refequ{xx6})}%
\end{multline}%
Now we are done, because $g_0 \in {r_{i,i'}}^{-1}\bigl(D_{i,i'}\bigr)$ and $r_{i,i'} \in \C^0(\mca[1]{\tannakian\G};\nC)$ for all $0\leqq i,i'\leqq d$.%

In order to prove our claim \refequ{xx6}, let us consider, for each $0\leqq i,i'\leqq d$, a decreasing sequence of open disks%
\begin{equazione}\label{xx7}%
\cdots \subset D_{i,i'}^{\ell+1} \subset D_{i,i'}^\ell \subset \cdots \subset D_{i,i'}^0 \subset \nC%
\end{equazione}%
centred at $z_{i,i'}$ and whose radius $\delta_{i,i'}^\ell$ tends to zero. If we make the innocuous assumption $\delta_{i,i'}^0 = 1$ then it will follow from our hypotheses on the functions $a, a'$ that the sets%
\begin{equazione+}\label{xx8}%
\Sigma^\ell \bydef {\displaycap{0\leqq i,i'\leqq d}{}{{r_{i,i'}}^{-1}\Bigl(\overline{D_{i,i'}^\ell}\Bigr)} - \Gamma} & \bigl(\ell = 1, 2, \ldots\bigr)%
\end{equazione+}%
are closed subsets of the compact space $\G(K,K')$, where $K = \support a$ and $K' = \support{a'}$. The sets $\Sigma^\ell$ form a decreasing sequence. Their intersection \txtcap{\ell=1}{\infty}{\Sigma^\ell} has to be empty because of the faithfulness of $g \mapsto \varrho(g)$ on $\G(x_0,x_0')$ and our hypotheses on $a$, $a'$. Hence, by compactness, there will be some $\ell$ such that $\Sigma^\ell = \varnothing$. This proves the claim, and therefore, the theorem.%
\end{proof}%

\begin{nota}\label{xx.16}
{\em(The present remark will be used nowhere else and therefore it may be skipped without consequences. You should read \refsez[O.2.3]{O.2.4} first, anyway.)}

Observe that whenever \G\ and \H\ are \index{Morita equivalence}Morita equivalent Lie groupoids, one of them is \stack F\nobreakdash-reflexive if and only if the other is. Indeed, by naturality of the envelope transformation $\envelope[F]{\text-}: \Id \to \tannakian[F]{\text-}$, one gets a commutative square of topological groupoid homomorphisms%
\begin{equazione}\label{xx.17}%
\begin{split}%
\xymatrix@C=50pt{\G\ar[d]^\varphi_{\text{Morita~eq.}}\ar[r]^-{\envelope\G} & {\tannakian\G}\ar[d]^{\tannakian\varphi} \\ \H\ar[r]^-{\envelope\H} & {\tannakian\H}}%
\end{split}%
\end{equazione}%
in which both $\varphi$ and $\tannakian\varphi$ are fully faithful. It follows immediately that \envelope\G\ is fully faithful if and only if the same is true of \envelope\H. With a bit more work, it can be shown that \envelope\G\ is an open map if and only if \envelope\H\ is so (use the simplifying assumption that $\mca[0]\varphi: \mca[0]\G \to \mca[0]\H$ is a surjective submersion).%
\end{nota}%
\separazione%
By definition, a Lie groupoid \G\ is \stack F\nobreakdash-reflexive if and only if one can solve topologically the problem of reconstructing \G\ from its representations of type \stack F (that is to say one can recover \G\ up to isomorphism of topological groupoids from such representations). In the case of Lie groups, a topological solution provides a completely satisfactory answer because the smooth structure of any Lie group is uniquely determined by the topology of the group itself. However, in the present more general context it is not evident a priori that the notion of reflexivity we introduced above is as strong as to settle the smoothness problem mentioned at the beginning of \refsez{N.18}, think e.g.\ of $\G = M$ a smooth manifold. More precisely, we consider the following question: does reflexivity of \G, in the foregoing purely topological sense, actually imply that the functionally structured space $(\mca[1]{\smash{\tannakian\G}},\Rf^\infty)$ defined in \refsez{N.18} is a smooth manifold and the envelope map $\mca[1]{\envelope{}}: \mca[1]\G \to \mca[1]{\smash{\tannakian\G}}$ a diffeomorphism? The answer proves to be affirmative, as we shall now see.%

Let \G\ be an arbitrary Lie groupoid. Choose an arrow $g_0 \in \G(x_0,x_0')$ and a representation $R = (E,\varrho)$ of \G\ first of all. Then choose an arbitrary metric $\phi$ on $E$ and global sections $\zeta_1, \ldots, \zeta_d$, resp.\ $\zeta_1', \ldots, \zeta_d'$, forming a local $\phi$\nobreakdash-orthonormal frame for $E$ about $x_0$, resp.\ $x_0'$, as in the proof of Theorem~\refcnt{O.thm1}. These data determine a smooth mapping%
\begin{multiriga}\label{xx9}%
\varrho^{\zeta_1\ldots,\zeta_d}_{\zeta_1',\ldots,\zeta_d'}: \mca[1]\G \longto {M\times M\times \mathit M(d;\nC)}\text, \\ \text{as follows:}\quad g \mapsto \bigl(\s(g);\t(g);\varrho_{1,1}(g),\ldots,\varrho_{i,i'}(g),\ldots,\varrho_{d,d}(g)\bigr)%
\end{multiriga}%
(the functions $\varrho_{i,i'}$ are those defined in \refequ{xx4}; $\mathit M(d;\nC) = \End(\nC^d)$ is the space of ${d\times d}$ complex matrices).%

If the envelope homomorphism $\envelope\G: \G \to \tannakian\G$ of the Lie groupoid \G\ is faithful, it follows from Lemma \refcnt[O.1.4+.3.3]{O.lem8} that for every point $x$ of the base manifold $M$ of \G\ there exists a representation $(E,\varrho) \in {\Ob\,\R\G}$ such that $\kernel{\varrho_x}$ is a discrete subgroup of the isotropy group $G_x = \G|_x$. Consequently, for an arbitrary arrow $g_0 \in \G(x_0,x_0')$ there will exist $(E,\varrho) \in {\Ob\,\R\G}$ such that the map $\G(x_0,x_0') \to \Lis(E_{x_0},E_{x_0'})$, $g \mapsto \varrho(g)$ is injective on some open neighbourhood of $g_0$ in $\G(x_0,x_0')$. Then the following lemma applies:%
\begin{lemma}\label{O.lem1}%
Let \G\ be a Lie groupoid. Fix an arrow $g_0 \in \G(x_0,x_0')$ and let $(E,\varrho) \in {\Ob\,\R\G}$ be a representation. Suppose the map ${g\mapsto\varrho(g)}: \G(x_0,x_0') \to \Lis(E_{x_0},E_{x_0'})$ is injective on some open neighbourhood of $g_0$ in $\G(x_0,x_0')$.%

Then the smooth mapping $\varrho^\zeta_{\zeta'}: \mca[1]\G \to {M\times M\times \mathit M(d;\nC)}$ \refequ{xx9} is an immersion at $g_0$, for any choice of a metric and of related orthonormal frames $\zeta = \{\zeta_1,\ldots,\zeta_d\}$, $\zeta' = \{\zeta_1',\ldots,\zeta_d'\}$.%
\end{lemma}%
\begin{proof}%
Let $M$ be the base manifold of \G. Fix open balls $U, U' \subset M$, centred at $x_0, x_0'$ respectively and so small that the sections $\zeta_1, \ldots, \zeta_d$, resp.\ $\zeta_1', \ldots, \zeta_d'$ form a local orthonormal frame for $E$ over $U$, resp.\ $U'$. Since the source map \s\ of \G\ is a submersion, one can always choose $U$ also so small that there exists a local trivialization $\Gamma \iso {U\times B} \xto\pr U$ for \s\ in a neighbourhood $\Gamma$ of $g_0$ in \mca[1]\G, where $B$ is an open euclidean ball. It is no loss of generality to assume $\t(\Gamma) \subset U'$. Then we obtain, for the restriction of the mapping $\varrho^\zeta_{\zeta'} = \varrho^{\zeta_1\ldots,\zeta_d}_{\zeta_1',\ldots,\zeta_d'}$ to $\Gamma$, a ``coordinate expression'' of the following form%
\begin{equazione}\label{xx.10}%
{U\times B} \to {U\times U'\times \mathit M(d;\nC)}\text, \quad (u,b) \mapsto \bigl(u,u'(u,b),\boldsymbol\varrho(u,b)\bigr)%
\end{equazione}%
where $\boldsymbol\varrho(g) \in \mathit M(d;\nC)$ denotes the matrix $\{\varrho_{i,i'}(g)\}_{1\leqq i,i'\leqq d}$. The differential of the mapping \refequ{xx.10} at, let us say, $g_0 = (x_0,0)$ reads%
\begin{equazione}\label{xx.11}%
\begin{split}\left(%
\begin{matrix}%
\mathit{Id} & 0 \\ * & {\D_2u'}(x_0,0) \\ * & {\D_2\boldsymbol\varrho}(x_0,0)%
\end{matrix}%
\right)\end{split}%
\end{equazione}%
and it is therefore injective if and only if such is the differential of the partial map ${b\mapsto \bigl(u'(x_0,b),\boldsymbol\varrho(x_0,b)\bigr)}: B \to {U'\times \mathit M(d;\nC)}$ at the origin of $B$.%

We are now reduced to showing that the restriction%
\begin{displaymath}%
\varrho^\zeta_{\zeta'}: \G(x_0,\text-) \longto {M\times \GL(d)} = {\{x_0\}\times M\times \GL(d;\nC)}%
\end{displaymath}%
is an immersion at $g_0$. Let $G_0 = \G|_{x_0}$ be the isotropy group at $x_0$ and choose, in the vicinity of $g_0$, a local (equivariant) trivialization $\G(x_0,S) \iso {S\times G_0}$ for the principal $G_0$\nobreakdash-bundle $\t^{x_0}: \G(x_0,\text-) \to {\G x}$; we can assume that $S$ is a submanifold of $U'$ and that in this local chart $g_0 = (x_0',e_0)$, where $e_0$ stands for the neutral element of $G_0$. We then obtain a new coordinate expression for the restriction of $\varrho^\zeta_{\zeta'}$ to $\G(x_0,\text-)$, namely%
\begin{equazione}\label{xx.12}%
{S\times G_0} \to {U'\times \GL(d;\nC)}\text, \quad (s,g) \mapsto \bigl(s,\boldsymbol\varrho(s,g)\bigr)\text.%
\end{equazione}%
Since its first component is the inclusion of a submanifold, this map will be an immersion at $g_0 = (x_0',e_0)$ provided the partial map $g \mapsto \boldsymbol\varrho(x_0',g)$ is an immersion at $e_0$. The latter corresponds to the diagonal of the square%
\begin{displaymath}%
\xymatrix@R=19pt{G_0\ar[d]^{g_0\,\text-}_\iso\ar[r]^-\varrho & \Aut(E_{x_0})\ar[d]^{\rho(g_0)\,\text-}_\iso \\ \G(x_0,x_0')\ar[r]^-\varrho & \Lis(E_{x_0},E_{x_0'})\text,\!\!}%
\end{displaymath}%
so our problem reduces to proving that the homomorphism $\varrho: G_0 \to \GL(E_{x_0})$ is immersive. By hypothesis, this is injective in an open neighbourhood of $e_0$ and hence our claim follows at once.%
\end{proof}%

We are now ready to establish our previous claims about the functional structure $\Rf^\infty$ on the Tannakian groupoid \tannakian\G. Let \G\ be any \stack F\nobreakdash-reflexive Lie groupoid (\stack F an arbitrary stack of smooth fields, as ever).%

Fix an arrow $\lambda_0 \in \mca[1]{\smash{\tannakian\G}}$. Our first task will be to find some open neighbourhood $\Omega$ of $\lambda_0$ such that $(\Omega,\Rf^\infty_\Omega)$ turns out to be isomorphic, as a functionally structured space, to a smooth manifold $(X,\smooth[X])$. Since we are working under the hypothesis that \G\ is reflexive, there is a unique $g_0 \in \mca[1]\G$ such that $\lambda_0 = \envelope{g_0}$. By Lemma~\refcnt{O.lem1} and the comments preceding it, we can find, for a conveniently chosen $(E,\varrho) \in {\Ob\,\R\G}$, an open neighbourhood $\Gamma$ of $g_0$ in \mca[1]\G\ such that the smooth map $\varrho^\zeta_{\zeta'}: \mca[1]\G \to {M\times M\times \mathit M(d;\nC)}$ \refequ{xx9} induces a diffeomorphism of $\Gamma$ onto a submanifold $X \subset {M\times M\times \mathit M(d;\nC)}$. Notice that the same data which determine the map \refequ{xx9} also determine a map of functionally structured spaces%
\begin{multiriga}\label{xx.13}%
r^\zeta_{\zeta'} = r^{\zeta_1\ldots,\zeta_d}_{\zeta_1',\ldots,\zeta_d'}: \mca[1]{\smash{\tannakian\G}} \longto {M\times M\times \mathit M(d;\nC)}\text, \\ \lambda \mapsto \left(\s(\lambda);\t(\lambda);\{r_{i,i'}(\lambda)\}_{1\leqq i,i'\leqq d}\right)\text,%
\end{multiriga}%
where we put $r_{i,i'} = r_{R,\phi,\zeta_i,\zeta'_{i'}} \in \Rf$ \refequ[N.18]{N.v9}. From the reflexivity of \G\ again, it follows that the envelope map \envelope{} induces a homeomorphism between $\Gamma$ and the open subset $\Omega \bydef \envelope\Gamma$ of \mca[1]{\smash{\tannakian\G}}. The following diagram%
\begin{equazione}\label{xx.14}%
\begin{split}%
\xymatrix@C=60pt{\Gamma\ar[dr]_{\envelope{}|_\Gamma}^(.6){\:\iso\text{~homeo}}\ar[rr]^-{\varrho^\zeta_{\zeta'}|_\Gamma}_-{\iso\text{~diffeo}} && X\subset{M\times M\times \mathit M(d;\nC)} \\ & \Omega\ar[ur]!DL_(.7){r^\zeta_{\zeta'}|_\Omega} &}%
\end{split}%
\end{equazione}%
is clearly commutative. We contend that the map $r^\zeta_{\zeta'}|_\Omega$ provides the desired isomorphism of functionally structured spaces. Explicitly, this means that an arbitrary function $f: X' \to \nC$ belongs to $\C^\infty(X')$ if and only if its pullback $h = {f\circ r^\zeta_{\zeta'}}$ belongs to $\Rf^\infty(\Omega')$, for each fixed pair of corresponding open subsets $\Omega' \subset \Omega$, $X' \subset X$. Note that since the problem is local, we can make the simplifying assumption $\Omega' = \Omega$, $X' = X$. Thus, suppose $f \in \C^\infty(X)$ first; because of the local character of the problem again, it is not restrictive to assume that $f$ admits a smooth extension $\tilde f \in \C^\infty\bigl({M\times M\times \mathit M(d)}\bigr)$. Then $h$ coincides with the restriction to $\Omega$ of a global function $\tilde h = {\tilde f\circ r^\zeta_{\zeta'}}: \mca[1]{\tannakian{}} = \mca[1]{\smash{\tannakian\G}} \to \nC$ belonging to $\Rf^\infty(\mca[1]{\tannakian{}})$ because \refequ{xx.13} is a map of functionally structured spaces. Conversely, suppose $h = {f\circ r^\zeta_{\zeta'}} \in \Rf^\infty(\Omega)$. We know, from Example~\refcnt[N.19]{N.v13}, that the envelope map \envelope{} is a morphism of functionally structured spaces. Hence the composite ${h\circ\envelope{}}$ will belong to $\C^\infty(\Gamma)$. Since ${h\circ\envelope{}} = {f\circ r^\zeta_{\zeta'} \circ \envelope{}} = {f\circ \varrho^\zeta_{\zeta'}}$ and $\varrho^\zeta_{\zeta'}|_\Gamma$ is a diffeomorphism of $\Gamma$ onto $X$, it follows that $f \in \C^\infty(X)$, as contended.%

We have therefore proved that if a Lie groupoid \G\ is \stack F\nobreakdash-reflexive then the space $(\mca[1]{\smash{\tannakian[F]\G}},\Rf^\infty)$ is actually a (Hausdorff) smooth manifold. There is little work left to be done by now:%
\begin{proposizione}\label{O.prp10}%
Let \stack F be an arbitrary stack of smooth fields and let \G\ be a Lie groupoid. Suppose \G\ is \stack F\nobreakdash-reflexive.%

Then the Tannakian groupoid \tannakian[F]\G, endowed with its canonical functional structure $\Rf^\infty$, turns out to be a Lie groupoid; moreover, the \stack F\nobreakdash-envelope homomorphism%
\begin{equazione}\label{xx.15}%
\envelope[F]\G: \G \xto{\mspace{10mu}\iso\mspace{10mu}} \tannakian[F]\G%
\end{equazione}%
turns out to be an isomorphism of Lie groupoids.%
\end{proposizione}%
\begin{proof}%
We know from the foregoing discussion that $(\mca[1]{\tannakian{}},\Rf^\infty)$ is a smooth manifold. Then all we have to show now, clearly, is that the envelope map $\envelope{}: \mca[1]\G \to \mca[1]{\tannakian{}}$ is a diffeomorphism. Equivalently, we have to show that \envelope{} is an isomorphism of functionally structured spaces between $(\mca[1]\G,\smooth[{\mca[1]\G}])$ and $(\mca[1]{\tannakian{}},\Rf^\infty)$. This follows immediately, locally, from the commutativity of the triangles \refequ{xx.14} and the previously established fact that both $\varrho^\zeta_{\zeta'}|_\Gamma$ and $r^\zeta_{\zeta'}|_\Omega$ are functionally structured space isomorphisms onto $(X,\smooth[X])$.%
\end{proof}%
\separazione%
Let us pause for a moment to summarize our current knowledge of the \stack F\nobreakdash-envelope $\envelope[F]\G: \G \to \tannakian[F]\G$ of an arbitrary {\em proper} Lie groupoid \G. First of all, we know that \envelope\G\ is faithful \textit{(Thm.~\refcnt{O.thm3}).} We have also ascertained that \tannakian\G\ is a topological groupoid \textit{(Ex.~\refcnt[N.19]{N.v13} and Prop.~\refcnt[N.19]{N.v17}).} Moreover, it has been established that \envelope\G\ is necessarily an isomorphism of topological groupoids in case \envelope\G\ is faithful \textit{(Thm.~\refcnt{O.thm1});} whenever this happens to be true, one can completely solve the reconstruction problem for \G\ \textit{(Prop.~\refcnt{O.prp10}).} Now observe that faithfulness of \envelope\G\ is equivalent to the following property: \textsl{if $g \neq \u(x)$ in the isotropy group $\G|_x$ then there exists a representation $(E,\varrho) \in {\Ob\,\R\G}$ such that $\varrho(g) \neq \id \in \Aut(E_x)$.} We can therefore conclude by saying that an arbitrary proper Lie groupoid can be recovered from its representations of type \stack F if and only if such representations are \gm{enough} in the sense of the above-mentioned property.%

The final part of the present section will be devoted to showing that \textsl{any proper Lie groupoid admits enough representations of type \stack{E^\infty}} (= smooth Euclidean fields, cfr.\ \refsez{N.16}). By the foregoing remarks, this will immediately imply the general reconstruction theorem we were striving for. Recall that our approach via smooth Euclidean fields is motivated by the impossibility to obtain that result by using representations of type \stack{V^\infty} (smooth vector bundles), as illustrated by the examples discussed in \refsez{O.3.1+.4.1}.%
\separazione%
We begin with some preliminary remarks of a purely topological nature. Let \G\ be a proper Lie groupoid and let $M$ denote the base manifold of \G. Recall that a subset $S \subset M$ is said to be \index{invariant submanifold, subset}\textit{invariant} when $s \in S$ implies ${g\cdot s} \in S$ for all arrows $g \in \mca[1]\G$. If $S$ is an arbitrary\inciso{viz., not necessarily invariant}subset of $M$, we let ${\G\cdot S}$ denote the \index{saturation}\textit{saturation} of $S$, that is to say the smallest invariant subset of $M$ containing $S$, so that $S$ is invariant if and only if ${\G\cdot S} = S$; note that the saturation of an open subset is also open. \textsl{Now let $V$ be any open subset with compact closure: we contend that ${\G\cdot\overline V} = \overline{\G\cdot V}$.} The direction `$\subset$' of this equality is valid even for a non-proper Lie groupoid; it follows for instance from the existence of local bisections. To check the opposite inclusion, one can resort to the well-known fact that the orbit space\footnote{The quotient of $M$ associated with the equivalence $x \thicksim {g\cdot x}$. We will indicate by $\mathit o$ the map (of $M$ into this quotient) which sends $x$ to its equivalence class.} of a proper Lie groupoid is Hausdorff and then use the compactness of $\overline V$; in detail: since the image of the compact set ${\overline V}$ under the continuous mapping $\mathit o: M \to M/\G$ is a compact and hence closed subset of the Hausdorff space $M/\G$, the inverse image ${\G\cdot\overline V} = \mathit o^{-1}\left(\mathit o\left(\overline V\right)\right)$ must be closed as well. Next, let $U$ be an invariant open subset of $M$. From the equality we have just proved, it follows immediately that \textsl{$U$ coincides with the union of all its open invariant subsets $V$, $\overline V \subset U$.} Indeed, since any given point $u_0 \in U$ admits an open neighbourhood $W$ with compact closure contained in $U$, one has%
$$%
u_0 \in {\G\cdot W} = V \subset \overline V = \overline{\G\cdot W} = {\G\cdot\overline W} \subset {\G\cdot U} = U\text.%
$$%
The latter remark applies to the construction of \G\nobreakdash-invariant partitions of unity on $M$; for our purposes it will be enough to illustrate a special case of this construction. Consider an arbitrary point $x_0 \in M$ and let $U$ be an open {\em invariant} neighbourhood of $x_0$. Choose another open neighbourhood $V$ of $x_0$, invariant and with closure contained in $U$. The orbit ${\G\cdot x_0}$ and the set-theoretic complement ${\complement V}$ are invariant disjoint closed subsets of $M$, so \textit{Corollary~\refcnt[N.17]{iv8}} provides us with an invariant \index{cut-off function|emph}function $\Phi \in \C^\infty(M;\nR)$ such that $\Phi(x_0) = 1$ and $\Phi = 0$ outside $V$.%

We are now ready to establish a basic extension property enjoyed by the representations of type \stack{E^\infty} of proper Lie groupoids; our \gm{main theorem} below will be essentially a consequence of this property and of Zung's results on local linearizability. Our goal will be achieved by means of an obvious cut-off technique which is of course not available when one limits oneself to representations on vector bundles.%

Since throughout the subsequent discussion the type \stack F = \stack{E^\infty} is fixed, we agree to systematically suppress any reference to type. Let \G\ be an arbitrary {\em proper} Lie groupoid and let $M$ denote its base as usual. Let $U \subset M$ be a {\em \G\nobreakdash-invariant} open neighbourhood of a point $x_0 \in M$, and suppose we are given a {\em partial} representation $(\E_U,\varrho_U) \in \R{\G|_{\mathnormal U}}$. We know from \refsez{N.17} that there is an induced Lie group representation%
\begin{equazione}\label{N.vi1}%
\varrho_{U,0}: G_0 \longto \GL(\E_{U,0})%
\end{equazione}%
of the isotropy Lie group $G_0 = \G|_{x_0}$ on the vector space $\E_{U,0} = \smash{\left(\E_U\right)}_{x_0}$. \textsl{We contend that one can construct a global representation $(\E,\varrho) \in \R\G$ for which it is possible to exhibit an isomorphism of $G_0$\nobreakdash-spaces $\E_0 \bydef \E_{x_0} \iso \E_{U,0}$.} (The $G_0$\nobreakdash-space structure on $\E_0$ comes from the induced representation%
\begin{equazione}\label{N.vi2}%
\varrho_0: G_0 \longto \GL(\E_0)\text,%
\end{equazione}%
that on $\E_{U,0}$ from \refequ{N.vi1}.)%

To begin with, let us fix any invariant smooth function $a \in \C^\infty(M)$ with $a(x_0) = 1$ and $\support a \subset U$; such functions always exist\inciso{as we saw before}in view of the properness of \G. Let $V \subset M$ denote the open subset consisting of all $x$ such that $a(x)\neq 0$. One can define the following bundle $\{\E_x\}$ of Euclidean spaces over $M$:%
\begin{equazione}\label{N.vi3}%
\E_x =%
\begin{cases}%
\E_{U,x} & \text{if $x \in V$}\\%
\{0\} & \text{otherwise.}%
\end{cases}%
\end{equazione}%
Let \f E be the smallest sheaf of sections of the bundle $\{\E_x\}$ which contains the following presheaf%
\begin{equazione}\label{N.vi4}%
W \;\mapsto \;\bigl\{{a\zeta}: \zeta \in \f{(E_{\mathnormal U})}({U\cap W})\bigr\}\text.%
\end{equazione}%
(Here of course ${a\zeta}$ is to be understood as the appropriate ``prolongation by zero'' of the indicated section; note that since $M$ admits partitions of unity \refequ{N.vi3} actually equals \f E.) One verifies immediately that these data define a smooth Euclidean field \E\ over $M$. Next, introduce $\varrho$ by putting%
\begin{equazione}\label{N.vi5}%
\varrho(g) =%
\begin{cases}%
\varrho_U(g) & \text{for $g \in \G|_V$}\\%
0 & \text{otherwise.}%
\end{cases}%
\end{equazione}%
This law must be understood as describing a bundle $\bigl\{\varrho(g): ({\s^*\E})_g \isoto ({\t^*\E})_g\bigr\}$ of linear isomorphisms indexed over the manifold \G. The compatibility of this family of maps with the composition of arrows, amounting to the equalities $\varrho({g'g}) = {\varrho(g') \varrho(g)}$ and $\varrho(x) = \id$, is clear. Now, $\varrho$ will be an action of \G\ on \E\ provided it is a morphism ${\s^*\E} \to {\t^*\E}$ of Euclidean fields over \G: this is obvious, because for suitable functions $r_i \in \C^\infty$ one has%
$$%
{\varrho(g) {a\zeta}(\s[g])} = {a(\s[g]) {\varrho(g) \zeta(\s[g])}} = {a(\t[g]) \txtsum{i=1}\ell{r_i(g) \zeta'_i(\t[g])}} = \txtsum{i=1}\ell{r_i(g) {a\zeta'_i}(\t[g])}\text,%
$$%
in view of the \G\nobreakdash-invariance of $a$. Hence $(\E,\varrho) \in \R\G$. Finally, the identity $\E_0 = \E_{x_0} \bydef \E_{U,x_0} = \E_{U,0}$ provides the required $G_0$\nobreakdash-equivariant isomorphism.%

\begin{theorem}[General Reconstruction Theorem, Main Theorem]\label{xx.18}\index{main theorem}
Each proper Lie groupoid is \stack{E^\infty}\nobreakdash-reflexive.
\end{theorem}
\begin{proof}
Let \G\ be any such groupoid and fix a point $x_0$ of its base manifold $M$. We need to show the existence of a Euclidean representation $(\E,\varrho) \in {\Ob\,\R\G}$ inducing a faithful isotropy representation $\varrho_0: G_0 \into \GL(\E_0)$ \refequ{N.vi2} (we freely use the notation above). By the previously established extension property of Euclidean representations, it will be enough to find a partial representation $(E_U,\varrho_U) \in {\Ob\,\R{\G|_{\mathnormal U}}}$ defined over some invariant open neighbourhood $U$ of $x_0$ and with $\varrho_{U,0}: G_0 \into \GL(E_{U,0})$ \refequ{N.vi1} injective.%

It was observed in \refsez{N.4} that Zung's Local Linearizability Theorem yields the existence of \textsl{(a)} a smooth representation $G_0 \to \GL(V)$ on some (real) finite dimensional vector space \textsl{(b)} an embedding of manifolds $V \stackrel i\into M$ such that $0 \mapsto x_0$ and such that $U \bydef {\G\cdot i(V)}$ is an open subset of $M$ \textsl{(c)} a Morita equivalence ${G_0\ltimes V} \xto\iota \G|_U$ inducing $V \stackrel i\into U$ at the level of base manifolds. Note that the isotropy of ${G_0\ltimes V}$ at $0$ equals $G_0$ and that the equivalence $\iota$ induces an automorphism $\iota_0 \in \Aut(G_0)$ (which can be assumed to be the identity, just to fix ideas).%

Now let $\Phi: G_0 \into \GL(\boldsymbol E)$ be any faithful representation on a finite dimensional complex vector space. One has an induced faithful representation $\widetilde\Phi$ of ${G_0\ltimes V}$ on ${V\times\boldsymbol E}$ (cfr.\ the end of \refsez{O.4.3}). By the theory of \refsez{N.14}, there exists some representation $(E_U,\varrho_U) \in {\Ob\,\R{\G|_{\mathnormal U}}}$ such that $\iota^*(E_U,\varrho_U) \iso ({V\times\boldsymbol E},\widetilde\Phi)$; this is precisely the one we are looking for, because $\varrho_{U,0}: G_0 \into \GL(E_{U,0}) \iso \GL(\boldsymbol E)$ must coincide with $\Phi$.%
\end{proof}%

\capitolo{Classical Fibre Functors}\label{5}

In the present chapter we will again occupy ourselves with the study of the abstract notion of fibre functor. However, we shall be exclusively interested in fibre functors which take values in the category of smooth vector bundles over a manifold, in other words fibre functors of the form $\fifu: \Kt[C] \to \V[\infty]M$ or, equivalently, of type \stack{V^\infty}. Moreover, since in all examples of such functors we have in mind the tensor category \Kt[C] invariably turns out to be rigid, we shall make the assumption that \Kt[C] is rigid even though this is not indispensable; note that in this case $\End^\otimes(\fifu) = \Aut^\otimes(\fifu)$ ie $\lambda$ tensor preserving implies $\lambda$ invertible, see, for instance, \cite{DeMi82} \mbox{Prop.~1.13.} We shall use the adjective `{classical}' to refer to fibre functors of this sort.

Section \ref{N.21} is devoted to the study of some general properties of classical fibre functors. To start with, the Tannakian groupoid \tannakian{\fifu} associated with a classical fibre functor \fifu\ proves to be a $\C^\infty$-structured groupoid, that is to say all the structure maps of \tannakian{\fifu} turn out to be morphisms of functionally structured spaces; compare \refsez{N.18}. This allows us to introduce the category \R[\infty]{\tannakian\fifu} of $\C^\infty$-representations of the $\C^\infty$-structured groupoid \tannakian\fifu, along with an ``{evaluation}'' functor
$$
\ev: \Kt[C] \longto \R[\infty]{\tannakian\fifu}\text.
$$
The latter is in fact a tensor functor, by which the category \Kt[C] is put in relation to \R[\infty]{\tannakian\fifu}; we shall say more about this functor in \refsez{O.3.5}. Finally, we observe that a classical fibre functor \fifu\ which admits enough \fifu-invariant metrics (in the sense of Definition \refcnt[N.19]{xix.1}) is proper\nobreakdash---in other words, so is the corresponding map $(\s,\t): \tannakian{\fifu} \to {M\times M}$.

Section \ref{O.2.1} deals with the technical notion of tame submanifold, and is preliminary to \refsez[O.2.2+N.23]{O.2.4}. However, in order to read the latter sections a thorough understanding of \refsez{O.2.1} is not really necessary: it is actually enough to know what tame submanifolds are and the statements of Propositions \refcnt[O.2.1]{O.prp101}, \refcnt[O.2.1]{O.prp1}; one may skip what remains of \refsez{O.2.1} at first reading.

Section \ref{O.2.2+N.23} provides, for the Tannakian groupoid \tannakian{\fifu} associated with a classical fibre functor $\fifu: \Kt[C] \to \V[\infty]M$, an alternative characterization of the property of smoothness in terms of what we call representative charts. Such charts arise from the objects of the category \Kt, and their definition involves tame submanifolds of linear groupoids $\GL(E)$ over the manifold $M$.

Sections \ref{O.2.3}--\ref{O.2.4} are devoted to morphisms of fibre functors. For each morphism between two classical fibre functors there exists a corresponding homomorphism between the associated Tannakian groupoids, which turns out to be ``{smooth}'' ie a homomorphism of $\C^\infty$-structured groupoids. In \refsez{O.2.4} we introduce, as a special case, the notion of weak equivalence; the alternative characterization of smoothness provided in \refsez{O.2.2+N.23} is here put to work to show that the property of smoothness is, for classical fibre functors, invariant under weak equivalence. Finally, the homomorphism associated with a weak equivalence of smooth classical fibre functors is proved to be a Morita equivalence.

\sezione{Basic Definitions and Properties}\label{N.21}
In this section we study general properties of classical fibre functors. Let us begin by giving a precise definition:

\begin{definizione}\label{xxi1}
We shall call a fibre functor $\fifu: \Kt[C] \to \stack[M]F$ \index{classical fibre functor}\index{fibre functor!classical}\textit{classical} if it meets the following requirements:
\begin{elenco}
\item the domain tensor category \Kt\ is \index{rigid tensor category}\index{tensor category!rigid}rigid;
\item for every $R \in \Ob(\Kt)$, $\fifu(R)$ is a \index{locally trivial object}locally trivial object of \stack[M]F.
\end{elenco}
\end{definizione}
Observe that since the type \stack F is a stack of smooth fields, $\fifu(R)$ in \textsl{ii)} will actually belong to ${\Ob\,\V[F]M}$ ie it will be a locally trivial object of \stack[M]F of locally finite rank (cfr \refsez{N.11}). Since \V[F]M is equivalent to the category \stack[M]{V^\infty} of smooth vector bundles of locally finite rank over $M$ (recall that the base $M$ is always paracompact), it follows that the theory of classical fibre functors essentially reduces to just one type \stack F = \stack{V^\infty}. Because of this, for the rest of the present chapter\inciso{actually, for the rest of the present work}we shall omit any reference to type. So, for instance, we will write \index{V(X;k)@\V[\infty]{X;k}, \stack[X]{V^\infty} (category of vector bundles)}\stack[M]{V^\infty} or \V[\infty]{M} at all places where we would otherwise write \stack[M]F.

The pivotal fact of classical fibre functor theory is that for such fibre functors one has local formulas analogous to \refequ[N.17]{N.iv1c}. Namely, let $\fifu: \Kt \to \V[\infty]{M}$ be a classical fibre functor. Let an object $R \in \Ob(\Kt)$ and an arrow $\lambda_0 \in \tannakian{} \equiv \mca[1]{\smash{\tannakian\fifu}}$ be given. Choose, on $E \equiv \fifu(R)$, an arbitrary Hilbert metric $\phi$, whose existence is guaranteed by the paracompactness of $M$. By the local triviality assumption on $E$, it will be possible to find a local \index{orthonormal frame}$\phi$-orthonormal frame ${\zeta_1}', \ldots, {\zeta_d}' \in \sections E(U')$ about ${x_0}' \equiv \t(\lambda_0)$ such that $E_{u'} = \mathrm{Span}\,\bigl\{{\zeta_1}'(u'), \ldots, {\zeta_d}'(u')\bigr\}$ for all $u' \in U'$. (Note that here one really needs local triviality of $E$ within \stack F, in the sense of \refsez{N.11}, and not just the hypothesis that \sections E is locally free as a sheaf of modules over $M$.) Then for any given local section $\zeta \in \sections E(U)$, defined in a neighbourhood $U$ of $x_0 \equiv \s(\lambda_0)$, one gets, by letting $\Omega \equiv {\s^{-1}(U) \cap t^{-1}(U')} \subset \tannakian{}$,
\begin{equazione+}\label{xxi2}
{\lambda(R) \cdot \zeta(\s[\lambda])} = \txtsum{i'=1}d{r_{R,\phi,\zeta,{\zeta_{i'}}'}(\lambda) {\zeta_{i'}}'(\t[\lambda])}\text, & (\forall\lambda \in \Omega)
\end{equazione+}
where $r_{R,\phi,\zeta,{\zeta_{i'}}'} \in \Rf^\infty(\Omega)$ denotes\inciso{as in \refsez{N.18}}the representative function $\lambda \mapsto \bigsca{\lambda(R) \cdot \zeta(\s[\lambda])}{{\zeta_{i'}}'(\t[\lambda])}\hspace{0pt}_\phi$.

We shall immediately put this basic remark to work in the proof of the following

\begin{proposizione}\label{xxi3}
For every classical fibre functor $\fifu: \Kt[C] \to \V[\infty]{M}$, the \index{Tannakian groupoid!T(omega)@\tannakian\fifu}\index{T(omega)@\tannakian{\fifu} (Tannakian groupoid associated with a fibre functor)}Tannakian groupoid \tannakian{\fifu} is a \index{C infinity-structured groupoid@\ensuremath{\C^\infty}-structured groupoid}$\C^\infty$-structured groupoid (with respect to the standard \index{standard C infinity-structure@standard \ensuremath{\C^\infty}-structure \ensuremath{\Rf^\infty}}$\C^\infty$-structure $\Rf^\infty$ defined in \refsez{N.18}). \tannakian{\fifu} is, in particular, a topological groupoid for every classical \fifu.
\end{proposizione}
\begin{proof}
Let us take an arbitrary representative function $r = r_{R,\phi,\zeta,\zeta'}: \tannakian{} \to \nC$ on the space $\tannakian{} \equiv \mca[1]{\smash{\tannakian\fifu}}$, as in \refequ[N.18]{N.v8}. We shall regard $r$ as fixed throughout the entire proof.

\sloppy
To begin with, we consider the composition map $\mca[2]{\tannakian{}} = \mca{\tannakian{}} \xto\c \tannakian{}$. Our goal is to show that the function ${r\circ\c}$ is a global section of the sheaf $\smash{\mca[2]\Rf}^{,\infty} \equiv ({\Rf^\infty \otimes \Rf^\infty})^\infty|_{\mca[2]{\tannakian{}}}$. (Review, if necessary, the discussion about functionally structured groupoids in \refsez{N.18}.) Fix any pair of composable arrows $({\lambda_0}',\lambda_0) \in \mca[2]{\tannakian{}}$. There will be some $\phi$-orthonormal frame ${\zeta_1}', \ldots, {\zeta_d}' \in \sections{(\fifu R)}(U')$ about ${x_0}' \equiv \t(\lambda_0)$, such that Eq.~\refequ{xxi2} above holds for all $\lambda \in \Omega' \equiv \t^{-1}(U')$. Then, for every pair $(\lambda',\lambda)$ belonging to the open subset $\Omega'' \equiv {\s^{-1}(U') {_{\s}\times_{\t}} t^{-1}(U')} \subset \mca[2]{\tannakian{}}$, one gets the identity
\begin{multline}\notag
({r\circ\c})(\lambda',\lambda) = r({\lambda'\circ\lambda}) = \smash{\bigsca{\lambda'(R) \cdot \lambda(R) \cdot \zeta(\s[\lambda])}{\zeta'(\t[\lambda'])}}_\phi\\
= \txtsum{i'=1}d{r_{R,\phi,\zeta,{\zeta_{i'}}'}(\lambda) r_{R,\phi,{\zeta_{i'}}',\zeta'}(\lambda')}\text, \qquad \text{by \refequ{xxi2}}
\end{multline}
which expresses $({r\circ\c})|_{\Omega''}$ in the desired form, namely as an element of $\smash{\mca[2]\Rf}^{,\infty}(\Omega'')$.

\fussy
Next, consider the inverse map $\tannakian{} \xto\i \tannakian{}$. Fix any $\lambda_0 \in \tannakian{}$. In a neighbourhood $U$ of $x_0 = \s(\lambda_0)$ it will be possible to find a trivializing $\phi$\nobreakdash-orthonormal frame $\zeta_1, \ldots, \zeta_d \in \sections{(\fifu R)}(U)$. One can write down \refequ{xxi2} for each $\zeta_i$ ($i = 1, \ldots, d$):%
\begin{equazione+}\label{xxi4}%
{\lambda(R)\cdot \zeta_i(\s[\lambda])} = \txtsum{i'=1}d{r_{R,\phi,\zeta_i,\zeta'_{i'}}(\lambda) \zeta'_{i'}(\t[\lambda])}\text. & (\lambda \in \Omega = {\s^{-1}(U)\cap t^{-1}(U')}\text.)%
\end{equazione+}%
Letting $\{r'_{i',i}(\lambda): 1\leqq i',i\leqq d\}$ denote the inverse of the matrix $\{r_{R,\phi,\zeta_i,\zeta'_{i'}}(\lambda): 1\leqq i',i\leqq d\}$ for each $\lambda$ (this makes sense because $\lambda(R)$ is a linear iso), we see from the standard formula involving the inverse of the determinant that $r'_{i',i} \in \Rf^\infty(\Omega)$ for all $1\leqq i',i\leqq d$. If we now put $a'_{i'} = \scalare{\zeta}{\zeta'_{i'}}_\phi \in \C^\infty(U')$ for all $i' = 1, \ldots, d$ and $a_i = \scalare{\zeta_i}{\zeta'}_\phi \in \C^\infty(U)$ for all $i = 1, \ldots, d$, we obtain the following expression for $({r\circ\i})|_\Omega$%
\begin{align*}%
({r\circ\i})(\lambda) &= r(\lambda^{-1}) = \smash{\bigsca{\lambda(R)^{-1}\cdot \zeta(\t[\lambda])}{\zeta'(\s[\lambda])}}_\phi =%
\\%
&= \txtsum{i'=1}d{a'_{i'}(\t[\lambda]) \smash{\bigsca{\lambda(R)^{-1}\cdot \zeta'_{i'}(\t[\lambda])}{\zeta'(\s[\lambda])}}_\phi} =%
\\%
&= \txtsum{i'=1}d{\txtsum{i=1}d{a'_{i'}(\t[\lambda]) r'_{i',i}(\lambda) a_i(\s[\lambda])}}\text,%
\end{align*}%
which clearly shows membership of $({r\circ\i})|_\Omega$ in $\Rf^\infty(\Omega)$.%

The ``smoothness'' of the remaining structure maps was already proved in \S\ref{N.18} for an arbitrary fibre functor.%
\end{proof}%
\separazione%
By exploiting the categorical equivalence $\V{M} \xto\iso \stack[M]{V^\infty}$, $E \mapsto \tilde E$ \refequ[N.12]{xii.33}, one can make sense of the expression \index{linear groupoid@linear groupoid \ensuremath{\GL(E)}}\index{GL(E)@\ensuremath{\GL(E)} (linear groupoid)}$\GL(E)$ for every $E \in {\Ob\,\V{M}}$ simply by regarding $\GL(E)$ as short for $\GL(\tilde E)$. If $\fifu: \Kt \to \V{M}$ is a classical fibre functor, each object $R \in \Ob(\Kt)$ will determine a \index{evaluation representation@evaluation representation \ensuremath{\ev_R}|emph}\index{ev R@\ensuremath{\ev_R} (evaluation representation)|emph}homomorphism of functionally structured groupoids
\begin{equazione}\label{xxi5}%
\ev_R: \tannakian\fifu \longto \GL({\fifu R})\text, \quad\lambda \mapsto \lambda(R)%
\end{equazione}%
(note that if $\phi$ is any Hilbert metric on $E = \omega(R)$, the functions $q_{\phi,\zeta,\zeta'}: \mca[1]{\smash{\GL(E)}} \to \nC$, $\mu \mapsto \smash{\bigsca{\mu\cdot \zeta(\s[\mu])}{\zeta'(\t[\mu])}}_\phi$ will provide suitable local coordinate systems for the manifold $\mca[1]{\smash{\GL(E)}}$), which may be thought of as a ``smooth'' representation of \tannakian\fifu.%

It is worthwhile mentioning the following universal property, which characterizes the functional structure (and topology) we endowed the Tannakian groupoid with. Let \fifu\ be a classical fibre functor. Then for any functionally structured space $(Z,\mathscr F)$, a mapping $f: Z \to \tannakian{} = \mca[1]{\smash{\tannakian\fifu}}$ is a morphism of $(Z,\mathscr F)$ into $(\tannakian{},\Rf^\infty)$ (or simply, a continuous mapping of $Z$ into \tannakian{}) if and only if such is ${\ev_R\circ f}$ for every $R \in {\Ob\,\Kt}$. The `only if' direction is clear because of the foregoing remarks about the ``smoothness'' of $\ev_R$. Conversely, consider any representative function $r = r_{R,\phi,\zeta,\zeta'}: \tannakian{} \to \nC$; if $q_{\phi,\zeta,\zeta'}: \mca[1]{\smash{\GL({\fifu R})}} \to \nC$ is the smooth function defined above then one has ${r\circ f} = {q_{\phi,\zeta,\zeta'}\circ \ev_R\circ f} \in \mathscr F(Z)$, because by assumption ${\ev_R\circ f}$ is a morphism of $(Z,\mathscr F)$ into the smooth manifold $\mca[1]{\smash{\GL({\fifu R})}}$. The equivalence is now proven.%

In a manner entirely analogous to \refsez{O.3.1+.4.1}, one can define the complex tensor category \index{R(T;k)@\R[\infty]{\mathcal T;k} (category of smooth representations on vector bundles)|emph}\R[\infty]{\tannakian\fifu;\nC} of all ``{\index{C infinity representation@\ensuremath{C^\infty}-representation|emph}\index{representation!C infinity or smooth@\ensuremath{C^\infty}- or smooth|emph}\index{smooth representation|see{\ensuremath{C^\infty}-representation}}smooth}'' representations of the functionally structured groupoid \tannakian{\fifu} on smooth complex vector bundles over the base manifold $M$ of \fifu. Precisely, any such representation will consist of a complex vector bundle $E \in {\Ob\,\stack[M]{V^\infty}}$ and a homomorphism $\varrho: \tannakian\fifu \to \GL(E)$ of functionally structured groupoids over $M$ ($\varrho$ identical on $M$). Then one has the complex tensor functor
\begin{equazione}\label{xxi6}%
\ev: \Kt \longto \R[\infty]{\tannakian\fifu;\nC}\text, \quad R \mapsto (\widetilde{\fifu(R)},\ev_R)
\end{equazione}
(the so-called ``\index{evaluation functor@evaluation functor \ensuremath{\ev}|emph}\index{ev@\ensuremath{\ev} (evaluation functor)|emph}evaluation functor''). The parallel with the situation depicted in \S\ref{O.1.3} leads us to formulate the problem of determining whether or not the functor \refequ{xxi6} is in general\inciso{for an arbitrary classical fibre functor}a categorical equivalence. The answer is known to be yes, actually in the strong form of an isomorphism of categories, for a large class of examples: see \S\ref{O.3.5}, Proposition \refequ[O.3.5]{O.prp12} and related comments.

\separazione

We conclude this introductory section with an observation about proper classical fibre functors (cfr.\ \S\ref{N.19}). We intend to show that, in the classical case, existence of enough invariant metrics is sufficient to ensure properness and hence that the first condition of Definition~\ref{N.19}.\ref{xix.1} is actually redundant for any classical fibre functor.%

Notice first of all that each Hilbert metric $\phi$ on a complex vector bundle $E \in {\Ob\,\stack[M]{V^\infty}}$ determines a subgroupoid $\mathit U(E,\phi) \subset \GL(E)$, consisting of all {\em$\phi$\nobreakdash-unitary} linear isomorphisms between the fibres of $E$; more explicitly, the arrows $x \to x'$ in $\mathit U(E,\phi)$ are the unitary isomorphisms of $(E_x,\phi_x)$ onto $(E_{x'},\phi_{x'})$. Clearly, $\mathit U(E,\phi)$ is a proper Lie groupoid over the manifold $M$, embedded into $\GL(E)$. When there is no danger of ambiguity about the metric, we will just suppress $\phi$ from the notation.%

From our elementary remark \refequ[N.19]{N.v16} it follows that for any \fifu\nobreakdash-invariant Hilbert metric $\phi$ on $\fifu(R)$ ($R \in {\Ob\,\Kt}$) the evaluation homomorphism $\ev_R$ \refequ{xxi5} must factor through the subgroupoid $\mathit U({\fifu R}) \into \GL({\fifu R})$. Hence one may view $\ev_R$ as a ``smooth'' homomorphism%
\begin{equazione}\label{xxi7}%
\ev_R: \tannakian\fifu \longto \mathit U({\fifu R})\text, \quad\lambda \mapsto \lambda(R)\text.%
\end{equazione}%

\begin{proposizione}\label{xxi8}
Let $\fifu: \Kt \to \V[\infty]{M}$ be a classical fibre functor. Suppose there are enough \index{metric!omega-invariant@\fifu-invariant}\index{omega-invariant metric@\fifu-invariant metric}\fifu-invariant metrics (cfr \S\ref{N.19}, Definition~\ref{xix.1}). Then \fifu\ is \index{fibre functor!proper}\index{proper fibre functor}proper; in particular, \tannakian{\fifu} is a proper groupoid.
\end{proposizione}
\begin{proof}
Let us assign, to each object $R \in {\Ob\,\Kt}$, an \fifu\nobreakdash-invariant metric $\phi_R$ on $\fifu(R)$ once and for all. We shall simply write $\mathit U({\fifu R})$ in place of $\mathit U(\fifu(R),\phi_R)$.%

Let $K$ be an arbitrary compact subset of the base manifold $M$. We have to show that $\tannakian{}|_K = (\s,\t)^{-1}({K\times K})$ is a compact subset of the topological space $\tannakian{} = \mca[1]{\smash{\tannakian\fifu}}$. Consider the auxiliary space%
\begin{equazione}\label{xxi9}%
Z_K \bydef \displayprod{R\in{\Ob\,\Kt}}{}{\mathit U({\fifu R})|_K}%
\end{equazione}%
(product of topological spaces) and observe that $Z_K$ is compact because the same is true of each factor $\mathit U({\fifu R})|_K$. There is an obvious continuous injective map $e: \tannakian{}|_K \into Z_K$ given by $\lambda \mapsto \{\lambda(R)\}_{R\in{\Ob\,\Kt}}$. We claim that this map is actually a topological embedding of $\tannakian{}|_K$ onto a closed subset of $Z_K$: this will entail the required compactness of $\tannakian{}|_K$.%

\textit{The map $e$ is an embedding.} This will be implied at once by the following extension property of representative functions: for every $r = r_{R,\phi,\zeta,\zeta'} \in \Rf$ \refequ[N.18]{N.v9}, there exists a continuous function $h: Z_K \to \nC$ such that $r = {h\circ e}$ on $\tannakian{}|_K$. In order to obtain such an extension of $r$, note simply that on $\tannakian{}|_K$ one has $r_{R,\phi,\zeta,\zeta'} = {({q_{\phi,\zeta,\zeta'}\circ \pi_R})\circ e}$, where $\pi_R: Z_K \to \mathit U({\fifu R})|_K$ is the $R$\nobreakdash-th projection and $q_{\phi,\zeta,\zeta'}$ is the (restriction to $\mathit U({\fifu R})|_K$ of) the smooth function $\GL({\fifu R}) \to \nC$, $\mu \mapsto \smash{\bigsca{\mu\cdot \zeta(\s[\mu])}{\zeta'(\t[\mu])}}_\phi$.%

\textit{The image of $e$ is a closed subset of $Z_K$.} It is sufficient to observe that the conditions expressing membership of $\mu = \{\mu_R\}_{R\in{\Ob\,\Kt}} \in \txtprod{}{}{\mathit U({\fifu R})|_K}$ in the image of $e$\inciso{namely that $\s(\mu_R) = \s(\mu_S)$ and $\t(\mu_R) = \t(\mu_S)$ $\forall R, S \in {\Ob\,\Kt}$, naturality of $\mu$ and its being tensor preserving and self-conjugate}are each stated in terms of a huge number of identities which involve only the coordinates $\mu_R = \pi_R(\mu)$ in a continuous way.%
\end{proof}%

\begin{nota}\label{xxi10}%
A very marginal comment about proper classical fibre functors, improving, in the classical case, Lemma \ref{N.19}.\ref{N.v15}: \textit{for any proper classical \fifu, the equality ${\Rf = \Rf'}$ holds.} In order to see this, notice first of all that if $U$ is any open subset of $M$ on which $E|_U$ ($E = \fifu(R)$) trivializes then we can find $a\in \Aut(E|_U)$ such that $\phi_u(v,v') = \phi_{R,u}(v,{a_u\cdot v'})$ for all $u \in U$ ($\phi$ an arbitrary metric on $E$, $\phi_R$ as in the proof of the preceding proposition, $v, v' \in E_u$). Now, if we put $\xi'_U = {a(U)\zeta'_U}$ where $\zeta'_U$ is the restriction to $U$ of $\zeta'$, we get $r_{R,\phi,\zeta,\zeta'} = r_{R,\phi_R,\zeta,\xi'_U}$ on $\t^{-1}(U) \subset \tannakian{}$. We can use a partition of unity over all such $U$'s to obtain a global section $\xi'$ with the property that $r_{R,\phi,\zeta,\zeta'} = r_{R,\phi_R,\zeta,\xi'} \in \Rf'$.%
\end{nota}%

\sezione{Tame Submanifolds of a Lie Groupoid}\label{O.2.1}
Let \G\ be a Lie groupoid over a manifold $M$.

\begin{definizione}\label{xxii1}
A submanifold $\Sigma$ of the manifold of arrows \mca[1]{\G} will be said to be \index{principal submanifold}\textit{principal} if it can be covered with local parametrizations (viz inverses of local charts or, equivalently, open embeddings) of the form
\begin{equazione}\label{O.equ18}
\left\{\begin{array}{l}
{Z\times A} \into \Sigma
\\[\smallskipamount]
(z,a) \mapsto {\tau(z) \cdot \eta(a)}\text,
\end{array}\right.
\end{equazione}
where $Z$ is a submanifold of $M$, $\tau: Z \to \G(x,\text-)$ is, for some point $x \in M$, a smooth section to the target map of the groupoid, $\eta: H \into G_x$ is a Lie subgroup of the $x$-th isotropy group $G_x$ of \G\ and $A$ is an open subset of $H$ such that $\eta$ restricts to an embedding of $A$ into $G_x$.
\end{definizione}

Note that the image $\Sigma = {\tau(Z) \cdot \eta(A)}$ of a map of the form \refequ{O.equ18} is always a submanifold of \mca[1]{\G} and that the same map induces a smooth isomorphism of ${Z\times A}$ onto $\Sigma$. So, in particular, it makes sense to use such maps as local parametrizations. (Details can be found in Note \ref{O.sub1} below.)

Note also that any principal submanifold of \mca[1]{\G} admits an open cover by local parametrizations of type \refequ{O.equ18} with the additional property that the Lie group $H$ is connected and $A$ is an open neighbourhood in $H$ of the neutral element $e$. (Indeed, let $\sigma \in \Sigma$ be a given point and choose a local parametrization ${\tau\cdot\eta}$ of the form \refequ{O.equ18}. Suppose $\sigma = (z,a) \in {Z\times A}$ in this local chart. Replacing $A$ with ${a^{-1}A}$ and $\tau$ with ${\tau\cdot\eta(a)}$ accomplishes the reduction to the situation where $A$ is a neighbourhood of $e$ and $\sigma = (z,e)$; intersecting with the connected component of $e$ in $H$ finishes the job.)

\begin{lemma}\label{O.lem5}
Let $\varphi: \G \to \G'$ be a Lie groupoid homomorphism, inducing an immersion $f: M \to M'$ at the level of manifolds of objects. Assume that $\Sigma$ and $\Sigma'$ are principal submanifolds of \G\ and $\G'$ respectively, with the property that $\varphi$ maps $\Sigma$ injectively into $\Sigma'$.

Then $\varphi$ restricts to an immersion of $\Sigma$ into $\Sigma'$.
\end{lemma}
\begin{proof}
Fix any point $\sigma_0 \in \Sigma$ and let $x_0 \equiv \s(\sigma_0)$, $z_0 \equiv \t(\sigma_0)$. Choose local parametrizations ${\tau\cdot\eta}: {Z\times A} \into \Sigma$ and ${\tau'\cdot\eta'}: {Z'\times A'} \into \Sigma'$ of type \refequ{O.equ18} with, let us say, $\sigma_0 = (z_0,e) \in {Z\times A}$ and $\varphi(\sigma_0) = (f(z_0),e') \in {Z'\times A'}$, where $e$, resp.\ $e'$ is the neutral element of the Lie subgroup $\eta: H \into G_{x_0}$, resp.\ $\eta': H' \into G'_{f(x_0)}$. As remarked above, the Lie groups $H$ and $H'$ can be assumed to be connected. Let the domain of the first parametrization shrink around the point $(z_0,e)$ until the smooth injection $\varphi: \Sigma \into \Sigma'$ admits a local representation relative to the chosen parametrizations, namely
$$
\xymatrix@C=30pt@R=19pt{\Sigma\ar[r]^-\varphi & \Sigma' \\ {Z\times A}\ar@{_{(}->}[u]_{\tau\cdot\eta}\ar@{-->}[r]^-{\tilde\varphi} & {Z'\times A'}\text.\!\!\ar@{_{(}->}[u]_{\tau'\cdot\eta'}}
$$
$\tilde\varphi$ will be a smooth injective map, of the form $(z,a) \mapsto \bigl(z'(z,a),a'(z,a)\bigr)$. Note that $z'(z,a) = f(z)$ so that, in particular, $f$ maps $Z$ into $Z'$; this follows by comparing the target of the two sides of the equality
$$
{\tau'(z') \cdot \eta'(a')} = {\varphi(\tau(z)) \cdot \varphi(\eta(a))}\text.
$$
Since the restriction of $f$ to $Z$ is an immersion of $Z$ into $Z'$, the mapping $\tilde\varphi$ is immersive at $(z_0,e)$ if and only if the corresponding partial map $a \mapsto a'(z_0,a)$ is immersive at $e \in A$. Now, consider the following huge commutative diagram, where we put $x_0' \equiv f(x_0)$ and $z_0' \equiv f(z_0)$:
$$
\xymatrix@C=38pt{G_{x_0}\ar@{=}[r] & G_{x_0}\ar[r]^-\varphi & G'_{f(x_0)}\ar@{=}[r] & G'_{f(x_0)} \\ & \G(x_0,z_0)\ar[u]_{\tau(z_0)^{-1}\cdot}\ar[r]^-\varphi & \G'(x_0',z_0')\ar[u]_{\varphi(\tau(z_0))^{-1}\cdot} & \\ A\ar@{^{(}->}[uu]^\eta\ar@{=}[r] & {\{z_0\}\times A}\ar@{_{(}->}[u]_{\tau\cdot\eta}\ar[r] & {\{z_0'\}\times A'}\ar@{_{(}->}[u]_{\tau'\cdot\eta'}\ar@{=}[r] & A'\ar@{^{(}->}[uu]^{\eta'}}
$$
[the rectangle on the right commutes because $\varphi(\tau(z_0)) = \varphi(\sigma_0) = \tau'(f(z_0)) = \tau'(z_0')$]. The commutativity of the outer rectangle entails that the bottom map in this diagram, namely $a \mapsto a'(z_0,a)$, coincides with the restriction to $A$ of a (necessarily unique) Lie group homomorphism $\zeta: H \to H'$; the same map is therefore an immersion, because a Lie group homomorphism which is injective in a neighbourhood of $e$ must be immersive, see eg \mbox{Br\"ocker} and \mbox{tom Dieck} \cite{Broecker&tomDieck'85}, p.~27. The proof of the existence of the homomorphism of Lie groups $\zeta$ is deferred to Note \refcnt{O.sub4} below.
\end{proof}

\sloppy
\begin{definizione}\label{O.dfn1}
A submanifold $\Sigma$ of the arrow manifold of a Lie groupoid \G\ will be said to be \index{tame submanifold|emph}\textit{tame} if the following conditions are satisfied:
\begin{elenco}
\item the source map of \G\ restricts to a submersion of $\Sigma$ onto an open subset of the base manifold $M$ of \G;
\item for each point $x \in M$, the corresponding source fibre $\Sigma(x,\text-) \equiv {\Sigma\cap \G(x,\text-)}$ is a principal submanifold.
\end{elenco}
Note that from the first condition it already follows that the source fibre $\Sigma(x,\text-)$ is a submanifold (of $\Sigma$ and hence) of \mca[1]\G.
\end{definizione}

\fussy
\begin{proposizione}\label{O.prp101}
Let $\varphi: \G \to \G'$ be a Lie groupoid homomorphism, inducing an immersion $f: M \to M'$ at the level of base manifolds. Suppose that $\Sigma$, resp.\ $\Sigma'$ is a tame submanifold of \G, resp.\ $\G'$ and that $\varphi$ maps $\Sigma$ injectively into $\Sigma'$.

Then $\varphi$ restricts to an immersion of $\Sigma$ into $\Sigma'$.
\end{proposizione}
\begin{proof}
Fix $\sigma_0 \in \Sigma$, and put $x_0 = \s(\sigma_0)$. Choose local parametrizations ${U\times B} \into \Sigma$ at $\sigma_0 \iso (x_0,0) \in {U\times B}$, and ${U'\times B'} \into \Sigma'$ at $\varphi(\sigma_0) \iso (f(x_0),0) \in {U'\times B'}$, locally trivializing the respective source map\inciso{which is a submersion because of Condition \textsl{i)} of Definition \refcnt{O.dfn1}}over the open subsets $U \subset M$, $U' \subset M'$. (Here $B$ and $B'$ are open balls.) This means, for instance, that the first parametrization makes the diagram%
$$%
\xymatrix{{U\times B}\,\ar[dr]_{\pr}\ar@{^{(}->}[rr] & & \Sigma\ar[dl]^{\s} \\ & U &}%
$$%
commute. If the domain of the first parametrization is made to be conveniently small around the center $(x_0,0)$, the mapping $\varphi: \Sigma \into \Sigma'$ will induce a smooth and injective local expression%
$$%
\xymatrix{\Sigma\ar[r]^-{\varphi} & \Sigma' \\ {U \times B}\ar@{_{(}->}[u]\ar@{-->}[r] & {U' \times B'}\ar@{_{(}->}[u]}%
$$%
of the form $(x,b) \mapsto (x'(x,b),b'(x,b)) = (f(x),b'(x,b))$, so that, in particular, $f$ will map $U$ into $U'$. Since $f: U \to U'$ is then an immersion by assumption, the above local expression is an immersive map at $(x_0,0)$ if and only if the partial map $b \mapsto b'(x_0,b)$ is immersive at $0\in B$. At this point we can use Lemma \ref{O.lem5} to conclude the proof.
\end{proof}%

In particular, it follows that when a homomorphism $\varphi$ of Lie groupoids (let us say over the same manifold $M$ and with $f = \id$) induces a homeomorphism between two tame submanifolds $\Sigma$ and $\Sigma'$, then it restricts in fact to a diffeomorphism of $\Sigma$ onto $\Sigma'$. This will be for us the most useful property of tame submanifolds, and we shall make repeated application of it in the subsequent sections. Actually, the motivation for introducing the concept of tame submanifold was precisely to ensure this kind of automatic ``differentiability out of continuity''.%

\begin{nota}\label{O.sub1}%
Let $S = {\G m}$ be the $m$-th orbit. As a notational convention, we shall use the letter $S$ when we think of this orbit as a manifold, endowed with the unique differentiable structure that turns the target map%
\begin{equazione}\label{O.equ1}%
\t: \G(m,\text-) \to S%
\end{equazione}%
into a principal bundle with fibre the Lie group $G_m$ (acting on the manifold $\G(m,\text-)$ from the right, in the obvious way); \refequ{O.equ1} is in particular a fibre bundle, which is in fact equivariantly locally trivial. The inclusion $S \into M$ is an injective immersion, although not in general an embedding of manifolds. See also \textit{Moerdijk and Mr\v cun (2003),} \cite{Moerdijk&Mrcun'03}~pp.~115--117.%

To begin with, we show that the inclusion map is an embedding of the manifold $Z$ into $S$. Of course, $Z$ is a submanifold of $M$ and we have the inclusion $Z \subset {\G m}$, but from this fact we cannot a priori conclude that $Z$ embeds into $S$, not even that the inclusion map $Z \into S$ is continuous; the reason why we can do away with this difficulty is that over $Z$ there exists, by assumption, a smooth section $\tau$ to the target map $\G(m,\text-) \to M$. (Incidentally, observe that any such $\tau: Z \to \G(m,\text-)$ is an embedding of manifolds. Clearly, it will be enough to see that $\tau$ is an embedding of $Z$ into \G. Since $\tau$ is a smooth section over $Z$ to $\t: \G \to M$, it is an injective immersion; moreover, for any open subset $U$ of $M$ we have%
$$%
\tau({Z\cap U}) = {\tau(Z) \cap \t^{-1}(U)}\text{.)}%
$$%
Now, from the existence of $\tau$ it follows immediately that the inclusion ${\t \circ \tau}$ of $Z$ into $S$ is a smooth mapping; moreover, we have that this is actually an injective immersion, because on composing it with $S \into M$ one obtains the embedding $Z \into M$. It only remains to notice that if $U$ is open in $M$ then ${Z\cap U}$ coincides with ${Z\cap W}$ where $W = {\t \, \G(m,U)}$ is open in $S$.%

Next, we show that%
\begin{lemma}\label{O.lem2}%
For every $z_0 \in Z$, there is a local trivialization of the principal bundle \refequ{O.equ1}, of the form%
$$%
\G(m,W) \iso {W\times G_m}%
$$%
over an open neighbourhood $W$ of $z_0$ in $S$, such that its unit section agrees with $\tau$ on ${Z\cap W}$. (Recall that the unit section of such a local trivialization is the mapping that corresponds to $W \into {W\times G_m}, w \mapsto (w,1_m)$.)%
\end{lemma}%
\begin{proof}%
Since $Z$ embeds as a submanifold of $S$, it is possible to find an open neighbourhood $W$ of $z_0$ in $S$ diffeomorphic to a product of manifolds%
$$%
W \iso {(W \cap Z) \times B}\text{,} \quad z_0 \iso (z_0,0)\text{,}%
$$%
where $B$ is an open euclidean ball. Moreover, it is clearly not restrictive to assume that the principal bundle \refequ{O.equ1} can be trivialized over $W$. Then, after having fixed one such trivialization, we can take the composite mapping%
$$%
W \iso {(W \cap Z) \times B} \xto{\:\pr\:} {W \cap Z} \xto{\:\tau\:} \G(m,W) \iso {W\times G_m} \xto{\:\pr\:} G_m\text{,}%
$$%
which we denote by $\theta: W \to G_m$, and use it to produce an equivariant change of charts and hence a new local trivialization for \refequ{O.equ1}, namely%
$$%
{W\times G_m} \isoto {W\times G_m} \iso \G(m,W)\text, \quad (w,g) \mapsto (w,{\theta(w)g})\text{,}%
$$%
whose unit section is immediately seen to agree with $\tau$ on ${Z \cap W}$.
\end{proof}%

Our aim was to prove that $\Sigma = {\tau(Z)\cdot\eta(A)}$ is a submanifold of \G\ and that ${\tau\cdot\eta}$ is a smooth isomorphism between ${Z\times A}$ and $\Sigma$. Thus, fix $\sigma_0 \in \Sigma$, an let $z_0 = \t(\sigma_0)$; the latter is a point of $Z$. Fix also a trivializing chart for the principal bundle \refequ{O.equ1} as in the statement of Lemma \ref{O.lem2}; then%
$$%
\xymatrix{{W\times G_m}\ar@{->}[r]_-{\iso}^-{\text{diffeo.}} & \G(m,W) \\ {({Z\cap W})\times A}\ar@{_{(}->}[u]_{\text{embed.}}\ar@{..>}[r]^-{\text{biject.}} & {\Sigma\cap \G(m,W)}\ar@{_{(}->}[u]_{\text{set-th. incl.}}}%
$$%
commutes, where on the left we have the obvious embedding of manifolds, and the bottom map is $(z,a) \mapsto {\tau(z)\cdot\eta(a)}$, the restriction of ${\tau\cdot\eta}$. (The diagram commutes precisely because the unit section of the chart agrees with $\tau$ over ${Z\cap W}$.) It is then clear that ${\Sigma\cap \G(m,W)}$ is a submanifold of the open neighbourhood  $\G(m,W)$ of $\sigma_0$ in $\G(m,\text-)$, and that ${\tau\cdot\eta}$ restricts to a diffeomorphism of ${(Z\cap W) \times A}$ onto this submanifold.%

Henceforth, $\Sigma$ is a submanifold of $\G(m,\text-)$ and ${\tau\cdot\eta}$ is a bijective local diffeomorphism between ${Z\times A}$ and $\Sigma$. (Note that the statement that $Z \into S$ is an embedding is really used here.)%
\end{nota}%

\begin{nota}\label{O.sub4}%
\textit{Assume that a commutative rectangle%
$$%
\xymatrix@M=4pt{A\ar[d]\ar@{^{(}->}[r] & H\ar@{..>}[d]^{\exists!\zeta}\ar@{^{(}->}[r]^{\eta} & G\ar[d]^{\varphi} \\ A'\ar@{^{(}->}[r] & H'\ar@{^{(}->}[r]^{\eta'} & G'}%
$$%
is given, where $G$, $G'$ are Lie groups, $\varphi$ is a Lie group homomorphism, $\eta: H \into G$ and $\eta': H' \into G'$ are Lie subgroups with $H$ connected, $A \subset H$, $A' \subset H'$ are open neighbourhoods of the unit elements $e$, $e'$ of $H$, $H'$ respectively, and $A \to A'$ is a smooth mapping. Then there exists a unique Lie group homomorphism $\zeta: H \to H'$ which fits in the diagram as indicated.}%

Indeed, since $A$ is an open neighbourhood of $e$ in $H$ and $H$ is connected, $A$ generates $H$ as a group, see \textit{Br\"ocker and tom Dieck (1995),} \cite{Broecker&tomDieck'85}~p.~10. So ${\varphi\eta}(A)$ generates ${\varphi\eta}(H)$, and therefore ${\varphi\eta}(H) \subset \eta'(H')$ because ${\varphi\eta}(A) \subset \eta'(A') \subset \eta'(H')$. Since $\eta': H' \to \eta'(H')$ is a bijective homomorphism of groups, there exists a unique group-theoretic solution $\zeta: H \to H'$ to the problem ${\eta'\circ\zeta} = {\varphi\circ\eta}$. The restriction of $\zeta$ to $A$ coincides with the given smooth map $A \to A'$, thus $\zeta$ is smooth in a neighbourhood of $e$; since left translations are Lie group automorphisms, the commutativity of%
$$%
\xymatrix{H\ar[d]_{\iso}^{h\cdot}\ar[r]^-{\zeta} & H'\ar[d]_{\iso}^{\zeta(h) \cdot} \\ H\ar[r]^-{\zeta} & H'}%
$$%
shows that $\zeta$ is smooth in the neighbourhood of any $h \in H$, and hence globally smooth, in other words a Lie group homomorphism.%
\end{nota}%

\sottosezione{Tameness and Morita equivalence}

There is still one fundamental point we need to discuss, for the treatment of weak equivalences of classical fibre functors in Section \ref{O.2.4} below. Namely, suppose one is given a \index{Morita equivalence}Morita equivalence of Lie groupoids $\varphi: \G \to \G'$ such that at the level of manifolds of objects it is given by a submersion $\varphi: M \to M'$. Let $\Sigma$ be a subset of the manifold of arrows of \G, and assume that every point of $\Sigma$ has an open neighbourhood $\Gamma$ in \G\ with
\begin{equazione}\label{O.equ3}%
{\varphi^{-1}(\Sigma')\cap \Gamma} \subset \Sigma\text{,}%
\end{equazione}%
where we put $\Sigma' = \varphi(\Sigma)$; note that this is equivalent to saying that%
$$%
{\forall \gamma} \in \Gamma, \quad \gamma \in \Sigma \; \aeq \; \varphi(\gamma) \in \Sigma'\text{.}%
$$%
Then one has what follows%
\begin{elenco}%
\item[1.]{\em$\Sigma$ is a submanifold of \G\ if and only if $\Sigma'$ is a submanifold of $\G'$;}%
\item[2.]{\em$\Sigma$ is a submanifold of \G\ verifying Condition \textsl{i)} of Definition \ref{O.dfn1} if and only if the same is true of $\Sigma'$ in $\G'$;}%
\item[3.]{\em for every $m \in M$, the restriction $\varphi: \Sigma(m,\text-) \to \Sigma'(\varphi(m),\text-)$ is an open mapping between topological subspaces of the manifolds $\G$ and $\G'$;}%
\item[4.]{\em for every $m \in M$, the fibre $\Sigma(m,\text-)$ is a principal submanifold of \G\ if and only if its image $\varphi(\Sigma(m,\text-))$ is a principal submanifold of $\G'$.}%
\end{elenco}%
Before we start with the proofs, let us show how these statements 1-4 may be used to derive the following main result%
\begin{proposizione}\label{O.prp1}%
Let $\varphi: \G \longto \G'$ be a Morita equivalence of Lie groupoids inducing a submersion at the level of base manifolds. Let $\Sigma$ be a subset of the manifold of arrows of \G\ which satisfies condition (\ref{O.equ3}) above, and put $\Sigma' = \varphi(\Sigma)$. Then $\Sigma$ is a tame submanifold of \G\ if and only if $\Sigma'$ is a tame submanifold of $\G'$.%
\end{proposizione}%
\begin{proof}%
($\Leftarrow$) Suppose $m \in M$ is given: we must show that $\Sigma(m,\text-)$ is a principal submanifold of \G. Because of Statement 3, $\varphi(\Sigma(m,\text-))$ is an open subset of the subspace $\Sigma'(\varphi(m),\text-) \subset \G'$. Since the latter is by assumption a principal submanifold of $\G'$, it follows that the open subset $\varphi(\Sigma(m,\text-))$ is a principal submanifold of $\G'$ as well, and hence, by Statement 4, that $\Sigma(m,\text-)$ is a principal submanifold of \G.%

($\seq$) Fix $m' \in M'$. According to Statement 3, we have the open covering%
$$%
\Sigma'(m',\text-) = {\bigcup_{m \in \varphi^{-1}(m')} \varphi(\Sigma(m,\text-))}\text{,}%
$$%
and every open set belonging to this covering is a principal submanifold of $\G'$, by Statement 4 and the assumption. Hence the whole submanifold $\Sigma'(m',\text-) \subset \G'$ is a principal submanifold of $\G'$.%
\end{proof}%

Now we come to the proofs of Statements 1 to 4:%

\vskip 3pt%
{\noindent \textsc{Proof of Statement 1.}} Recall from Note \ref{O.sub8}, (\ref{O.equ32}) below that, up to diffeomorphism, one has for the morphism $\varphi$ a canonical decomposition%
$$%
\xymatrix{\Gamma\ar[d]^{(\s,\t)}\ar[r]^-{\iso} & {\Gamma' \times B \times C}\ar[d]\ar[r]^-{\pr} & \Gamma'\ar[d]^{(\s',\t')} \\ {U \times V}\ar[r]^-{\iso \times \iso} & {U' \times B \times V' \times C}\ar[r]^-{\pr \times \pr} & {U' \times V'}}$$%
in a neighbourhood $\Gamma$ of every point of $\Sigma$, with $\Gamma$ verifying condition (\ref{O.equ3}). We have that ${\Sigma'\cap\Gamma'}$ is a submanifold of $\Gamma'$ if and only if ${({\Sigma'\cap\Gamma'})\times A}$ is a submanifold of ${\Gamma'\times A}$, where $A = {B\times C}$. Thus, since ${({\Sigma'\cap \Gamma'})\times A} = \pr^{-1}({\Sigma'\cap\Gamma'})$ corresponds to%
$$%
{\varphi^{-1}({\Sigma'\cap\Gamma'}) \cap \Gamma} = {\varphi^{-1}(\Sigma')\cap \Gamma} = {\Sigma\cap\Gamma}%
$$%
in the diffeomorphism $\Gamma \iso {\Gamma'\times B\times C}$, this is in turn equivalent to saying that ${\Sigma\cap\Gamma}$ is a submanifold of $\Gamma$. Thus we see that $\Sigma$ is a submanifold of $\G$ if and only if $\Sigma'$ is a submanifold of $\G'$.%

\vskip 3pt%
{\noindent \textsc{Proof of Statement 2.}} From the previous diagram, we get that, up to diffeomorphism, $\s: \Gamma \to U$ corresponds to ${\s'\times\pr}: {\Gamma'\times {B\times C}} \to {U'\times B}$, so it restricts to a submersion ${\Sigma\cap\Gamma} \to U$ if and only if ${\s'\times\pr}$ restricts to a submersion ${({\Sigma'\cap\Gamma'})\times {B\times C}} \to {U'\times B}$; and this is in turn true if and only if $\s': {\Sigma'\cap\Gamma'} \to U'$ is a submersion.%

\vskip 3pt%
{\noindent \textsc{Proof of Statement 3.}} Fix a point $\sigma_0 \in \Sigma(m,\text-)$ and an open neighbourhood of that point in \G. Then from Note~\ref{O.sub8} below, we have for the restriction of $\varphi$ to $\Sigma$ a canonical local decomposition%
$$%
\xymatrix{{\Sigma\cap\Gamma}\ar[d]^{\s}\ar[r]^-{\iso} & {({\Sigma'\cap\Gamma'})\times B\times C}\ar[d]^{\s'\times\id}\ar[r]^-{\pr} & {\Sigma'\cap\Gamma'}\ar[d]^{\s'} \\ U\ar[r]^-{\iso} & {U '\times B}\ar[r]^-{\pr} & U'}%
$$%
at $\sigma_0 = (\sigma'_0,0,0)$, where $\Gamma$ can be choosen as small as one likes around $\sigma_0$, simply by taking a smaller $\Gamma' = \varphi(\Gamma)$ at $\sigma_0' = \varphi(\sigma_0)$ and reducing the radius of the open balls $B$, $C$; in particular, $\Gamma$ can be chosen so small that it fits in the previously assigned open neighbourhood of $\sigma_0$ in \G.%

It is immediate to recognize that $\varphi({\Sigma(m,\text-)\cap \Gamma}) = {\Sigma'(\varphi(m),\text-) \cap \Gamma'}$, where the latter is clearly an open subset of the subspace $\Sigma'(\varphi(m),\text-)$ of $\G'$. Indeed, in the left-hand square of the preceding diagram, the \s-fibre above $m \in U$, namely%
$$%
({\Sigma\cap\Gamma})(m,\text-) = {\Sigma(m,\text-) \cap \Gamma}\text{,}%
$$%
corresponds to the ${\s'\times\pr}$-fibre above $(\varphi(m),0)$, namely%
$$%
{({\Sigma'\cap\Gamma'})(\varphi(m),\text-) \times {0\times C}}\text.%
$$%
The latter is mapped by the projection $\pr$ onto%
$$%
{({\Sigma'\cap\Gamma'})(\varphi(m),\text-)} = {\Sigma'(\varphi(m),\text-) \cap \Gamma'}\text{,}%
$$%
hence $\varphi$ maps ${\Sigma(m,\text-) \cap \Gamma}$ onto ${\Sigma'(\varphi(m),\text-) \cap \Gamma'}$, as contended.%

\vskip 3pt%
{\noindent \textsc{Proof of Statement 4.}} This will be based on the following lemma:%
\begin{lemma}\label{O.lem4}
Let $\varphi: \G \to \G'$ be a fully faithful homomorphism of Lie groupoids and let $\varphi: M \to M'$ be the map induced on base manifolds. Suppose that $\Sigma \subset \G$ and $\Sigma' = \varphi(\Sigma) \subset \G'$ are submanifolds. Suppose also that a commutative diagram%
\begin{equazione}\label{O.equ24}%
\xymatrix{\Sigma\ar[d]^{\t}\ar[r]^-{\iso} & {\Sigma'\times C}\ar[d]^{\t'\times\id}\ar[r]^-{\pr} & \Sigma'\ar[d]^{\t'} \\ V\ar[r]^-{\iso} & {V'\times C}\ar[r]^-{\pr} & V'}%
\end{equazione}%
is given, where $V \subset M$ and $V' \subset M'$ are open subsets, $C$ is an open ball and the $\iso$'s are diffeomorphisms such that the top row coincides with $\varphi$ (arrows) and the bottom one with $\varphi$ (objects). Let $\sigma_0 \in \Sigma$ be a point with $\sigma_0 \iso (\sigma'_0,0) \in {\Sigma'\times C}$.%

Then $\Sigma$ admits a local parametrization of type (\ref{O.equ18}) at $\sigma_0$ if and only if $\Sigma'$ admits such a parametrization at $\sigma'_0$.%
\end{lemma}%
\begin{proof}%
Notation: let $z_0 = \t(\sigma_0) \in V$ and $z_0' = \t'(\sigma_0') = \varphi(z_0) \in V'$. Observe that from (\ref{O.equ24}) it follows that $z_0$ corresponds to $(z_0',0)$ in the diffeomorphism $V \iso {V'\times C}$, because $\sigma_0$ corresponds to $(\sigma_0',0)$ in $\Sigma \iso {\Sigma'\times C}$.%

($\Leftarrow$) Suppose that $\Sigma'$ admits a type (\ref{O.equ18}) local parametrization ${\sigma'\cdot\eta'}: {Z'\times A'} \into \Sigma'$ at $\sigma'_0 \iso (z'_0,e') \in {Z'\times A'}$. It is clearly no loss of generality to assume that the whole $\Sigma'$ is the image of this local parametrization. $Z' = \t'(\sigma'(Z')) \subset \t'(\Sigma') \subset V'$ is a submanifold, because so is $Z' \subset M'$. Write the diffeomorphism $V \iso {V'\times C}$ as $v \mapsto (\varphi(v),c(v))$ and let $Z \subset V$ be the submanifold corresponding to ${Z'\times C}$. Define $\sigma: Z \to \Sigma$ as $\sigma(z) = (\sigma'(\varphi(z)),c(z)) \in {\Sigma'\times C} \iso \Sigma$, and $\eta$ by%
\begin{equazione}\label{O.equ27}%
\xymatrix{\G(m,m)\ar[rr]^-\varphi_-\iso & & \G'(m',m') \\ & H'\ar@{^{(}-->}[ul]^{\eta}\ar@{_{(}->}[ur]_{\eta'} &}%
\end{equazione}%
so that $\sigma$ is clearly a smooth \t-section%
\begin{align*}%
\t(\sigma(z)) &\iso {\bigl({\t'\times\id}\bigr)\bigl(\sigma'(\varphi(z)),c(z)\bigr)} = \bigl(\t'(\sigma'(\varphi(z))),c(z)\bigr) = (\varphi(z),c(z))%
\\%
&\iso z%
\end{align*}%
with $\sigma(z_0) \iso \bigl(\sigma'(\varphi(z_0)),c(z_0)\bigr) = (\sigma_0',0) \iso \sigma_0$, and $\eta: H \into G_m$ is a Lie subgroup, where we put $H = H'$. Let $A = A'$. It is immediate to calculate that the image of ${\sigma\cdot\eta}: {Z\times A} \into \G$ is the whole $\Sigma$: thus we have constructed a global parametrization of $\Sigma$ at $\sigma_0$.%

($\seq$) In the other direction, suppose we are given a local parametrization ${\sigma\cdot\eta}: {Z\times A} \into \Sigma$ of type (\ref{O.equ18}) such that $\sigma_0 \in \Sigma$ corresponds to $(z_0,e) = (\t(\sigma_0),e) \in {Z\times A}$. Clearly, $Z = \t(\sigma(Z)) \subset \t(\Sigma) \subset V$ is a submanifold since so is $Z \subset M$.%

To begin with, observe that it is not restrictive to assume that the submanifold $Z \subset V$ corresponds to ${Z'\times C}$ under the diffeomorphism $V \iso {V'\times C}$, where of course $Z' = \varphi(Z)$. Precisely, the diffeomorphism $\Sigma \iso {\Sigma'\times C}$, that identifies $\sigma_0$ with $(\sigma_0',0)$, allows one to choose a smaller open neighbourhood $(\sigma_0',0) \in {\Sigma_0'\times C_0} \subset {\Sigma'\times C}$ such that $\Sigma_0 \iso {\Sigma_0'\times C_0}$ is contained in the domain of the local chart $({\sigma\cdot\eta})^{-1}$. From the commutativity of the diagram%
$$%
\xymatrix@C=50pt{{Z\times A}\ar[d]^\pr & \Sigma_0\ar@{_{(}->}[l]_-{({\sigma\cdot\eta})^{-1}}^-{\text{open~emb.}}\ar[d]^{\t}\ar[r]^-{\iso} & {\Sigma_0'\times C_0}\ar[d]^{\t'\times\id} \\ Z & \t(\Sigma_0)\ar@{_{(}->}[l]_-{\text{inclusion}}\ar[r]^-{\iso} & {\t'(\Sigma_0') \times C_0}}%
$$%
it follows at once that $Z_0 = \t(\Sigma_0) \subset Z$ is an open subset such that $V \iso {V'\times C}$ induces a bijection $Z_0 \iso {Z_0'\times C_0}$, where $Z_0' = \t'(\Sigma_0')$. Since it is compatible with the aims of the present proof to replace $C$ with a smaller $C_0$ centered at $0$, we can work with the smaller local parametrization obtained by restricting $\sigma$ to the open subset $Z_0$ of $Z$.%

Secondly, the \t-section $\sigma: Z \to \Sigma$ induces, by means of the diffeomorphisms $Z \iso {Z'\times C}$ and $\Sigma \iso {\Sigma'\times C}$, a smooth mapping ${Z'\times C} \to {\Sigma'\times C}$ of the form $(z',c) \mapsto (\sigma'(z',c),c)$; indeed%
\begin{align*}%
(z',c) &\iso z = \t(\sigma(z)) \iso ({\t'\times\id})\bigl(\sigma'(z',c),c(z',c)\bigr)%
\\%
&= \bigl(\t'(\sigma'(z',c)),c(z',c)\bigr)\text,%
\end{align*}%
hence it follows $\t'(\sigma'(z',c)) = z'$ and $c(z',c) = c$. We claim that it is no loss of generality to assume that it actually is of the form $(z',c) \mapsto (\sigma'(z'),c)$, ie that $\sigma'$ does not really depend on the variable $c$. Indeed, define $\tau: Z \to \Sigma$ as $\tau(z) = \bigl(\sigma'(\varphi(z),0),c(z)\bigr) \in {\Sigma'\times C} = \Sigma$; such a map is also a smooth \t-section%
\begin{align*}%
\t(\tau(z)) &\iso ({\t'\times\id})\bigl(\sigma'(\varphi(z),0),c\bigr) = \bigl(\t'(\sigma'(z',c)),c\bigr)%
\\%
&= (\varphi(z),c) \iso z%
\end{align*}%
with $\tau(z_0) = \bigl(\sigma'(z_0',0),0\bigr) = \sigma(z_0) = \sigma_0$. Then we can apply Lemma~\refcnt{O.lem3} below, the `Reparametrization Lemma', to obtain a new type (\ref{O.equ18}) local parametrization of $\Sigma$ at $\sigma_0$, for which such an assumption holds as well. Then we can introduce a smooth $\t'$-section $\sigma': Z' \to \Sigma'$ such that $\sigma'(z'_0) = \sigma'_0$, by setting $\sigma'(z') = \sigma'(z',0)$; also, we define $\eta'$ by means of (\ref{O.equ27}) and put $H' = H$ and $A' = A$. Thus, from the simplifying assumption above, it follows that $\sigma'(\varphi(z)) = \varphi(\sigma(z))$ for every $z \in Z$, and therefore that the image of ${\sigma'\cdot\eta'}: {Z'\times A'} \into \G'$ coincides with $\varphi(\image\,{\sigma\cdot\eta})$. But $\image{\sigma\cdot\eta} \subset \Sigma$ is an open subset, and $\varphi: \Sigma \to \Sigma' = \varphi(\Sigma)$ is an open mapping, whence $\image{\sigma'\cdot\eta'}$ is an open subset of $\Sigma'$. This concludes the proof.%
\end{proof}%

\begin{nota}\label{O.sub8}%
Fix a point $\sigma_0 \in \Sigma$. Since $f$ is a submersion, one can choose open neighbourhoods $U$ and $V$ of $\s(\sigma_0)$ and $\t(\sigma_0)$ in $M$ respectively, so small that, up to diffeomorphism, ${f|_U}$ becomes an open projection $U \iso {U'\times B} \xto{\pr} U'$ ($U'$ is an open subset of $M'$ and $B$ is an open ball; moreover, we shall assume that $\s(\sigma_0)$ corresponds to $(f(\s(\sigma_0)),0)$ in the diffeomorphism $U \iso {U'\times B}$), and ${f|_V}$ becomes an open projection $V \iso {V'\times C} \xto{\pr} V'$ ($V'$ is an open subset of $M'$, and $C$ is an open ball; also, $\t(\sigma_0)$ corresponds to $(f(\t(\sigma_0)),0)$ in the diffeomorphism $V \iso {V'\times B}$). Since $\varphi$ is a Morita equivalence, we have the following pullback in the category of differentiable manifolds of class $C^\infty$%
$$%
\xymatrix{\G(U,V)\ar[d]^{(\s,\t)}\ar[r]^-{\varphi} & \G'(U',V')\ar[d]^{(\s',\t')} \\ {U\times V}\ar[r]^-{f\times f} & {U'\times V'}}%
$$%
which has therefore, up to diffeomorphism, the following aspect%
$$%
\xymatrix@C=45pt{\G(U,V)\ar[r]^-{\iso}_-{\text{diffeo.}}\ar[d]^{(\s,\t)} & {\G'(U',V')\times {B\times C}}\ar[d]^{(\s',\t')\times\id\times\id}\ar[r]^-{\pr} & \G'(U',V')\ar[d]^{(\s',\t')} \\ {U\times V}\ar[r]^-{\iso\,\times\,\iso} & {{U'\times B}\times {V'\times C}}\ar[r]^-{\pr\times\pr} & {U'\times V'}\text,\!\!}%
$$%
where the top composite arrow coincides with $\varphi$ and the bottom one with ${f\times f}$. Next, take an open neighbourhood $\Gamma$ of $\sigma_0$ in \G\ such that the relation \refequ{O.equ3} holds. Then the same relation is clearly also satisfied by any smaller open neighbourhood of $\sigma_0$ in \G, hence it is no loss of generality to assume that $\Gamma$ is contained in $\G(U,V)$ and that it corresponds to a product ${\Gamma'\times {B_0\times C_0}}$ (with $\Gamma' = \varphi(\Gamma)$ necessarily open in $\G'(U',V')$, because $\varphi: \G(U,V) \to \G'(U',V')$ is open as it is clear from the latter diagram, and with $B_0 \subset B$, $C_0 \subset C$ open balls centered at $0$ of smaller radius) in the diffeomorphism%
$$%
\G(U,V) \iso {\G'(U',V') \times {B\times C}}\text.%
$$%
Then, by our choice of $\Gamma$ we obtain a commutative diagram%
\begin{equazione}\label{O.equ32}%
\begin{split}%
\xymatrix@C=45pt{\Gamma\ar[r]^-{\iso}_-{\text{diffeo.}}\ar[d]^{(\s,\t)} & {\Gamma'\times {B_0\times C_0}}\ar[d]^{(\s',\t')\times\id\times\id}\ar[r]^-{\pr} & \Gamma'\ar[d]^{(\s',\t')} \\ {U_0\times V_0}\ar[r]^-{\iso\,\times\,\iso} & {{U'\times B_0}\times {V'\times C_0}}\ar[r]^-{\pr\times\pr} & {U'\times V'}}%
\end{split}%
\end{equazione}%
where the top composite arrow coincides with $\varphi$ and the bottom one with ${f\times f}$. Finally, by pasting the following commutative diagram%
$$%
\xymatrix@C=45pt{{U_0\times V_0}\ar[d]^{\pr}\ar[r]^-{\iso\,\times\,\iso} & {{U'\times B_0}\times {V'\times C_0}}\ar[d]^{\pr}\ar[r]^-{\pr\times\pr} & {U'\times V'}\ar[d]^{\pr} \\ V_0\ar[r]^{\iso} & {V'\times C_0}\ar[r]^{\pr} & V'}%
$$%
to the former one along the common edge, we obtain%
\begin{equazione}\label{O.equ30}%
\begin{split}%
\xymatrix@C=45pt{\Gamma\ar[r]^-{\iso}\ar[d]^{\t} & {\Gamma'\times {B_0\times C_0}}\ar[d]^{\t'\times\pr}\ar[r]^-{\pr} & \Gamma'\ar[d]^{\t'} \\ V_0\ar[r]^-{\iso} & {V'\times C_0}\ar[r]^-{\pr} & V'}%
\end{split}%
\end{equazione}%
and then, since property \refequ{O.equ3} holds for $\Gamma$,%
\begin{equazione}\label{O.equ31}%
\begin{split}%
\xymatrix{{\Sigma \cap \Gamma}\ar[r]^-\iso\ar[d]^\t & {({\Sigma' \cap \Gamma'})\times {B_0\times C_0}}\ar[d]^{\t'\times\pr}\ar[r]^-\pr & {\Sigma' \cap \Gamma'}\ar[d]^{\t'} \\ V_0\ar[r]^-\iso & {V'\times C_0}\ar[r]^-\pr & V'\text.\!\!}%
\end{split}%
\end{equazione}%
Both in \refequ{O.equ30} and in \refequ{O.equ31}, the top composite arrow coincides with the restriction of $\varphi$ and the bottom one with the restriction of $f$. Of course, one has analogous diagrams with source maps replacing target maps.%
\end{nota}%

\begin{nota}\label{O.sub2}%
Here we shall state and prove the Local Reparametrization Lemma, which was needed in the proof of Lemma~\refcnt{O.lem4}.%
\begin{lemma}[Local Reparametrization]\label{O.lem3}%
Let $\G \rightrightarrows M$ be a Lie groupoid. Suppose we are given: a point $m \in M$, a smooth \t-section $\tau: Z \to \G(m,\text-)$ defined over a submanifold $Z \subset M$, a Lie subgroup $\eta: H \into G_m$ and an open neighbourhood $A$ of the unit $e$ in $H$ such that the restriction of $\eta$ is an embedding. Let $\Sigma = {\tau(Z)\cdot\eta(A)}$ be the image of the mapping of type \refequ{O.equ18} obtained from these data.%

Let $\sigma_0 \iso (z_0,e) \in {Z\times A}$ be a given point in $\Sigma$, and suppose that $\sigma: Z \to \Sigma$ is any other smooth \t-section such that $\sigma(z_0) = \sigma_0 = \tau(z_0)$.%

Then there exists a smaller open neighbourhood ${Z_0\times A_0}$ of the point $(z_0,e)$ in ${Z\times A}$ such that%
$$%
{\sigma\cdot\eta}: {Z_0\times A_0} \into \Sigma%
$$%
is still a local parametrization for $\Sigma$ at $\sigma_0$.%
\end{lemma}%
\begin{proof}%
If we consider the composite ${({\tau\cdot\eta})^{-1} \circ \sigma}: Z \to \Sigma \to {Z\times A}$, we get smooth coordinate maps $z \mapsto (\zeta(z),\alpha(z))$, characterized by the equation $\sigma(z) = {\tau(\zeta(z)) \cdot \eta(\alpha(z))}$. Comparing the target of the sides of this equation we get $\zeta(z) = z$. Thus $\sigma$ is completely determined by the smooth mapping $\alpha: Z \to A$ via the relation $\sigma(z) = {\tau(z)\cdot \eta(\alpha(z))}$.%

Now, we choose a smaller open neighbourhood $A_0 \subset A$ of the unit $e$ such that ${A_0\cdot A_0} \subset A$, which exists by continuity of the multiplication of $H$, and next an open neighbourhood $Z_0$ of $z_0$ in $Z$ such that $\alpha(Z_0) \subset A_0$; this is possible because $\alpha(z_0) = e$, which follows from $\sigma(z_0) = \tau(z_0) = {\tau(z_0)\cdot \eta(e)}$. It is then clear that ${\sigma\cdot\eta}$ maps ${Z_0\times A_0}$ into $\Sigma$: indeed, ${\forall (z,a)} \in {Z_0\times A_0}$, ${\sigma(z)\cdot \eta(a)} = {({\tau(z)\cdot \eta(\alpha(z))}) \cdot \eta(a)} = {\tau(z)\cdot \eta({\alpha(z)\cdot a})}$, and this is clearly an element of ${\tau(Z_0) \cdot \eta({A_0\cdot A_0})} \subset {\tau(Z)\cdot\eta(A)} = \Sigma$.%

If again we compose ${({\tau\cdot\eta})^{-1} \circ ({\sigma\cdot\eta})}: {Z_0\times A_0} \to \Sigma \to {Z\times A}$, we get smooth coordinate maps $(z,a) \mapsto \bigl(\zeta(z,a),\alpha(z,a)\bigr)$, characterized by the relation ${\sigma(z)\cdot\eta(a)} = {\tau(\zeta(z,a)) \cdot \eta(\alpha(z,a))}$. Taking the target yields $\zeta(z,a) = z$, thus we have a smooth mapping ${Z_0 \times A_0} \to {Z \times A}$ of the form $(z,a) \mapsto (z,\alpha(z,a))$ characterized by the equation ${\sigma(z) \cdot \eta(a)} = {\tau(z) \cdot \eta(\alpha(z,a))}$. (So, in particular, $\alpha(z,e) = \alpha(z)$ and $\alpha(z_0,e) = e$.)%

To conclude, it will be enough to observe that this mapping has invertible differential at $(z_0,e) \in {Z_0 \times A_0}$, because if that is the case then the mapping induces a local diffeomorphism of an open neighbourhood of $(z_0,e)$ in ${Z_0 \times A_0}$ (which can be assumed to be ${Z_0 \times A_0}$ itself, up to shrinking) onto an open neighbourhood of $(z_0,e) \in {Z \times A}$, so that if we then compose back with ${\tau \cdot \eta}$ we see that ${\sigma \cdot \eta}$ is a diffeomorphism of ${Z_0 \times A_0}$ onto an open subset of $\Sigma$. To see the invertibility of the differential, it will be sufficient to prove that the partial map $a \mapsto \alpha(z_0,a)$ has invertible differential at $e \in A_0$. But from the characterizing equation (setting $z = z_0$)$$\alpha(z_0,a) = {\eta^{-1}(\tau^{-1}(z_0) \sigma(z_0)) \cdot a} = {\eta^{-1}(1_m) \cdot a} = a$$we see at once that this differential is in fact the identity.%
\end{proof}%
\end{nota}%

\begin{nota}\label{O.sub3}%
We include here a discussion of tame submanifolds in connection with embeddings of Lie groupoids, parallel to the one concerning Morita equivalences. Suppose one is given such an embedding, ie a Lie groupoid homomorphism $\iota: \G \into \G'$ such that the mapping $\iota$ itself and the mapping $i: M \into M'$ induced on bases are embeddings of manifolds. Let $\Sigma$ be a subset of \G, and put $\Sigma' = \iota(\Sigma) \subset \G'$. The following statements hold%
\begin{elenco}%
\item\textsl{$\Sigma$ is a submanifold of \G\ if and only if $\Sigma'$ is a submanifold of $\G'$, in which case the restriction $\iota: \Sigma \to \Sigma'$ is a diffeomorphism;}%
\item\textsl{$\Sigma$ is a principal submanifold of $\G$ if and only if $\Sigma'$ is a principal submanifold of $\G'$;}%
\item\textsl{in case $i: M \into M'$ is an open embedding, $\Sigma$ is a tame submanifold of $\G$ if and only if $\Sigma'$ is a tame submanifold of $\G'$.}%
\end{elenco}%
Note that, as a special case, we get invariance of tame submanifolds under isomorphisms of Lie groupoids.%
\end{nota}%

\sezione[Smoothness, Representative Charts]{Smoothness and Representative Charts}\label{O.2.2+N.23}
In \refsez{N.21} we discussed some general properties of classical fibre functors, which hold quite apart from the eventuality that the canonical $\C^\infty$-structure on the space of arrows of the Tannakian groupoid might prove not to be a smooth manifold structure. On the contrary, in the present section we turn our attention specifically to the problem of finding effective criteria to decide whether a given classical fibre functor is ``{smooth}'' in the sense illustrated at the beginning of \refsez{N.18}. Such criteria will be employed in \refsez{O.3.5}; they involve the technical notion of tame submanifold introduced in the preceding section.

To motivate our definitions (which may appear rather artificial at first glance) let us consider a smooth classical fibre functor \fifu\ over a manifold $M$. Recall that \fifu\ being \index{smooth fibre functor|emph}\index{fibre functor!smooth|emph}smooth means by definition that the standard $\C^\infty$-structure $\Rf^\infty$ on the space \mca[1]{\tannakian\fifu} turns \tannakian{\fifu} into a Lie groupoid over $M$; compare \refsez{N.18}. Consider any classical representation $\varrho: \tannakian{\fifu} \to \GL(E)$ on a smooth vector bundle $E$; we know from Lemma \refcnt[N.20]{O.lem1} that if the map $\lambda \mapsto \varrho(\lambda)$ is injective in the vicinity of $\lambda_0$ within the subspace $\tannakian\fifu(x_0,{x_0}')$ [$x_0 \equiv \s(\lambda_0)$, ${x_0}' \equiv \t(\lambda_0)$] of \mca[1]{\tannakian\fifu}, the same map must be an immersion, into the manifold of arrows of $\GL(E)$, of some open neighbourhood $\Omega \subset \tannakian{}$ of $\lambda_0$ and therefore it must induce, provided $\Omega$ is chosen small enough, a diffeomorphism of $\Omega$ onto a submanifold $\varrho(\Omega)$ of $\GL(E)$. When, in particular, $\varrho = \ev_R$ for some $R \in \Ob(\Kt)$, we agree to write $R(\Omega)$ for the submanifold [of the manifold of arrows of $\GL({\fifu R})$] that corresponds to $\Omega$, namely we put
\begin{equazione}\label{xxiii.1}
R(\Omega) \bydef \ev_R(\Omega)\text.
\end{equazione}
It is not exceedingly difficult to see that the submanifolds of $\GL(E)$ of the form $\varrho(\Omega)$, for all $\varrho$ and $\Omega$ such that $\varrho$ induces a diffeomorphism of $\Omega$ onto $\varrho(\Omega)$, are necessarily tame submanifolds of $\GL(E)$, cfr Lemma \refcnt[O.3.5]{O.lem7} below. It will be convenient to have a name for the local diffeomorphisms of the above-mentioned type:

\begin{definizione}\label{xxiii.2}
We shall call \index{representative chart|emph}\textit{representative chart} any pair $(\Omega,R)$ consisting of an open subset $\Omega$ of the space of arrows of \tannakian{\fifu} and an object $R \in \Ob(\Kt)$, such that $\ev_R: \tannakian{\fifu} \to \GL({\fifu R})$ restricts to a homeomorphism of $\Omega$ onto a tame submanifold $R(\Omega)$ of the linear groupoid $\GL({\fifu R})$.
\end{definizione}
Note that this definition has been formulated so that it makes sense for an arbitrary classical fibre functor \fifu; when \fifu\ is smooth and $(\Omega,R)$ is a representative chart, the map $\lambda \mapsto \lambda(R)$ induces a diffeomorphism of $\Omega$ onto the submanifold $R(\Omega)$ of $\GL({\fifu R})$: this justifies our definition.

Observe that if $R$ and $S$ are two isomorphic objects of \Kt[C] then $(\Omega,R)$ is a representative chart of \tannakian{\fifu} if and only if the same is true of $(\Omega,S)$ (see Note \refcnt{O.sub5} below). Moreover, if $(\Omega,R)$ is a representative chart of \tannakian\fifu, the same is obviously true of $(\Omega',R)$ for each open subset $\Omega' \subset \Omega$.

We know from Lemma \refcnt[O.1.4+.3.3]{O.lem8} that if a classical fibre functor \fifu\ is smooth then for each $\lambda_0$ there exists some $R \in \Ob(\Kt)$ such that the map $\lambda \mapsto \lambda(R)$ is injective in a neighbourhood of $\lambda_0$ within the subspace $\tannakian\fifu(\s[\lambda_0],\t[\lambda_0])$ of \mca[1]{\tannakian\fifu}. Now, as remarked before, this implies that $\lambda_0$ lies in the domain $\Omega$ of a representative chart $(\Omega,R)$: thus we see that \textsl{for any smooth classical fibre functor, the domains of representative charts form an open covering of the space of arrows of the corresponding Tannakian groupoid.}

\pagebreak

Next, let us consider an arbitrary representative chart $(\Omega,R)$ of \tannakian\fifu, for a smooth \fifu. Let $S$ be an arbitrary object of \Kt. By choosing direct sum representatives conveniently, we may suppose that $\fifu({R\oplus S}) = {{\fifu R} \oplus {\fifu S}}$. The evaluation map $\ev_{R\oplus S}$ will yield a one-to-one correspondence between $\Omega$ and the subspace $({R\oplus S})(\Omega)$ of $\GL({{\fifu R} \oplus {\fifu S}})$: indeed, since $\lambda({R\oplus S}) = {\lambda(R) \oplus \lambda(S)}$ for all $\lambda \in \tannakian{\fifu}$, it is clear that the map $\lambda \mapsto \lambda({R\oplus S})$ factors through the submanifold ${\GL({\fifu R}) \times_M \GL({\fifu S})} \into \GL({{\fifu R}\oplus {\fifu S}})$ (cfr Note \refcnt{O.sub6} below) as the map $\lambda \mapsto \bigl(\lambda(R),\lambda(S)\bigr)$ (the latter is evidently injective, because so is $\lambda \mapsto \lambda(R)$, by hypothesis). We contend that $\ev_{R\oplus S}$ actually induces a homeomorphism of $\Omega$ onto the respective image; since $\ev_{R\oplus S}$ is immersive (by Lemma \refcnt[N.20]{O.lem1}), our contention will imply at once that $({R\oplus S})(\Omega)$ is a submanifold of $\GL({{\fifu R} \oplus {\fifu S}})$ and that $\ev_{R\oplus S}$ yields a diffeomorphism between $\Omega$ and this submanifold. Now, let $\Omega' \subset \Omega$ be a given open subset; fix any open subset $\Lambda' \subset \GL({\fifu R})$ such that ${R(\Omega) \cap \Lambda'} = R(\Omega')$ (such $\Lambda'$ exist because $\Omega$ and $R(\Omega)$ are homeomorphic via $\ev_R$): then
\begin{equazione}\label{xxiii.3}
{({R\oplus S})(\Omega) \cap \bigl({\Lambda' \times_M \GL({\fifu S})}\bigr)} = ({R\oplus S})(\Omega')\text,
\end{equazione}
which proves our contention. From the remarks that precede Definition \refcnt{xxiii.2} we immediately conclude that the following property is satisfied by any smooth classical fibre functor \fifu: \textsl{when $(\Omega,R)$ is a representative chart of \tannakian\fifu, so must be $(\Omega,{R\oplus S})$ for each object $S \in \Ob(\Kt)$.}

\separazione
\noindent The converse holds:

\begin{proposizione}\label{xxiii.4}%
Let \fifu\ be a classical fibre functor. Then \fifu\ is smooth if and only if the following two conditions are satisfied:%
\begin{elenco}%
\item the domains of representative charts cover the space of arrows of the Tannakian groupoid \tannakian\fifu, ie for each $\lambda \in \tannakian\fifu$ there exists a representative chart $(\Omega,R)$ with $\lambda \in \Omega$;%
\item if $(\Omega,R)$ is a representative chart of \tannakian{\fifu} then the same is true of $(\Omega,{R\oplus S})$ for every object $S \in \Ob(\Kt)$.%
\end{elenco}%
\end{proposizione}%
\begin{proof}%
We have already proved that a smooth classical fibre functor satisfies conditions \textsl{i)} and \textsl{ii)}. Vice versa, suppose these conditions are satisfied: the crucial point now is to show that any representative chart $(\Omega,R)$ establishes an isomorphism of functionally structured spaces between $(\Omega,\Rf^\infty_\Omega)$ and the submanifold $X \bydef R(\Omega) \subset \GL({\fifu R})$ (endowed with the structure \smooth[X]).%

Since $\ev_R: \tannakian{} \to \GL({\fifu R})$ is a morphism of functionally structured spaces, it is clear that $f \in \C^\infty(X)$ implies ${f\circ \ev_R} \in \Rf^\infty(\Omega)$ (cfr.\ the proof of Proposition~\refcnt[N.20]{O.prp10}). The converse implication is less obvious: we will make use of the special properties of tame submanifolds we derived in the preceding section. Suppose $r = r_{S,\psi,\eta,\eta'} \in \Rf^\infty(\Omega)$ and let $f$ be the function on $X$ such that ${f\circ \ev_R} = r$; we must show that $f \in \C^\infty(X)$. Since $f = {q_{\psi,\eta,\eta'} \circ \ev_S \circ {\ev_R}^{-1}}$ where $q_{\psi,\eta,\eta'}$ is the smooth function on $\GL({\fifu S})$ given by $\nu \mapsto \smash{\bigsca{\nu\cdot\eta(\s[\nu])}{\eta'(\t[\nu])}}_\psi$ and ${\ev_R}^{-1}: X \xto\iso \Omega$ is the inverse map, it will be enough to show that ${\ev_S \circ {\ev_R}^{-1}}$ is a smooth mapping of $X$ into $\GL({\fifu S})$. Put $E = \fifu(R)$, $F = \fifu(S)$. Recall that ${\GL(E) \times_M \GL(F)}$ is the product of $\GL(E)$ and $\GL(F)$ in the category of Lie groupoids over $M$ (see Note \refcnt{O.sub6} below) and that therefore it comes equipped with two projections $\pr_E,\pr_F$ that are morphisms of Lie groupoids over $M$. One can build the following commutative diagram
\begin{equazione}\label{xxiii.5}%
\begin{split}%
\xymatrix@C=37pt@M=3.5pt{& ({R\oplus S})(\Omega)\ar[d]_\iso^{\text{homeo}}\ar@{^{(}->}[r]^-{e_{R,S}} & {\GL(E)\times_M \GL(F)}\ar[d]^{\pr_E} \\ \Omega\ar[ur]^(.45){\ev_{R\oplus S}}\ar[r]^-{\ev_R} & X = R(\Omega)\ar@{^{(}->}[r]^-{\text{submanifold}} & \GL(E)\text,\!\!}%
\end{split}%
\end{equazione}%
where $e_{R,S}$ is the smooth embedding whose composition with%
\begin{equazione}\label{xxiii.6}%
{\GL(E) \times_M \GL(F)} \into \GL({E\oplus F}) = \GL\bigl(\fifu({R\oplus S})\bigr)\text, \quad (\mu,\nu) \mapsto {\mu\oplus\nu}%
\end{equazione}%
equals the inclusion of $({R\oplus S})(\Omega)$ into $\GL(\fifu{R\oplus S})$. Now, $(\Omega,{R\oplus S})$ is a representative chart of \tannakian{\fifu} and hence $({R\oplus S})(\Omega)$ is a tame submanifold of $\GL(\fifu{R\oplus S})$, so we can apply Proposition~\refcnt[O.2.1]{O.prp101} to conclude that the transition homeomorphism in \refequ{xxiii.5} is in fact a diffeomorphism. This immediately implies the desired smoothness of the transition mapping ${\ev_S \circ {\ev_R}^{-1}}: X \to \GL(F)$, because of the commutativity of the following diagram:%
\begin{equazione}\label{xxiii.7}%
\begin{split}%
\xymatrix@C=41pt@M=4pt{& ({R\oplus S})(\Omega)\ar@{^{(}->}[r]^-{e_{R,S}}_(.42){\text{smooth}} & {\GL(E)\times_M \GL(F)}\ar[d]^{\pr_F} \\ X\ar[ur]^(.45){\text{trans.~diffeo~~~~}}_(.53){\!\iso}\ar[r]^(.57){{\ev_R}^{-1}} & \Omega\ar@{^{(}->}[r]^-{\ev_S}\ar[u]_{\ev_{R\oplus S}} & \GL(F)\text.\!\!}%
\end{split}%
\end{equazione}%

From condition~\textsl{i)} and what we have just proved, we see that $(\tannakian{},\Rf^\infty)$ is a smooth manifold and that each representative chart $(\Omega,R)$ induces a diffeomorphism $\ev_R|_\Omega$ of $\Omega$ onto $R(\Omega)$. Moreover, since on the domain of any representative chart $(\Omega,R)$ the source map of \tannakian\fifu\ is the composition of $\ev_R|_\Omega$ with the restriction to $R(\Omega)$ of the source map of $\GL({\fifu R})$, we also see that the source map of \tannakian\fifu\ is a submersion\nobreakdash---because such remains the source map of $\GL({\fifu R})$ when restricted to the tame submanifold $R(\Omega) \subset \GL({\fifu R})$. Proposition~\refcnt[N.21]{xxi3} allows us to finish the proof.%
\end{proof}%

There is yet one useful remark concerning Condition \textsl{ii)}: under the hypothesis that $(\Omega,R)$ is a representative chart, the evaluation map $\ev_{R\oplus S}$ establishes, as in \refequ{xxiii.3}, a homeomorphism between $\Omega$ and the subset $({R\oplus S})(\Omega)$ of the manifold $\GL(\fifu{R\oplus S})$, wherefore the pair $(\Omega,{R\oplus S})$ is a representative chart if and only if $({R\oplus S})(\Omega)$ is a tame submanifold of $\GL(\fifu{R\oplus S})$.

The usefulness of the last proposition will become evident in the study of weak equivalences of classical fibre functors (cfr Section \ref{O.2.4}) and in the study of classical fibre functors associated with proper Lie groupoids (Chapter \refcpt{6}).

\begin{corollario}\label{xxiii.8}
Let $\fifu: \Kt[C] \to \V[\infty]M$ be a classical fibre functor satisfying conditions \textsl{i)} and \textsl{ii)} of the preceding proposition.

Then there exists a unique manifold structure on the space of arrows of the groupoid \tannakian\fifu, that renders \tannakian\fifu\ a Lie groupoid and
$$
\ev_R: \tannakian\fifu \longto \GL({\fifu R})
$$
a smooth representation for each object $R$. Equivalently, the same manifold structure can be characterized as the unique manifold structure for which an arbitrary mapping $f: X \to \tannakian{}$ is smooth if and only if so is ${\ev_R\circ f}$ for all $R$. The correspondence $R \mapsto \bigl(\fifu(R),\ev_R\bigr)$, $a \mapsto \fifu(a)$ determines a faithful tensor functor $\ev$ of \Kt[C] into \R[\infty]{\tannakian\fifu}, which makes
\begin{equazione}\label{xxiii.9}%
\begin{split}%
\xymatrix@R=19pt{\Kt\ar[dr]_(.45)\fifu\ar[rr]^-\ev & & \R[\infty]{\tannakian\fifu}\ar[dl]\\ & \V[\infty]M &}%
\end{split}%
\end{equazione}%
commute as a diagram of tensor functors (where the unlabelled arrow is the standard forgetful functor of \refsez{N.13}).%
\end{corollario}
% - Da qui in poi tutto da rivedere - 22/02/2008
\begin{proof}
We only need to check the assertions concerning the uniqueness of the smooth structure. Thus, suppose $\ev_R$ smooth ${\forall R}$. For convenience, let $\tanngpd\fibfunc^*$ denote the ``unknown'' manifold structure on the set $\tanngpd\fibfunc$. Since the topology of $\tanngpd\fibfunc^*$ is necessarily finer than that of $\tanngpd\fibfunc$, an open subset of $\tanngpd\fibfunc$ must be in particular a tame submanifold of $\tanngpd\fibfunc^*$. Therefore if $(\Omega,R)$ is a representative chart, the homomorphism of Lie groupoids $\ev_R: \tanngpd\fibfunc^* \to \GL(\omega R)$ restricts to a smooth isomorphism of the open subset $\Omega \subset \tanngpd\fibfunc^*$ onto the (tame) submanifold $R(\Omega)$ of $\GL(\omega R)$. Thus, we see that the identity map is, locally in the domains of representative charts, a diffeomorphism between $\tanngpd\fibfunc$ and $\tanngpd\fibfunc^*$; since representative charts cover $\tanngpd\fibfunc$, we get $\tanngpd\fibfunc^* = \tanngpd\fibfunc$, as was to be proved.
\end{proof}

For the sake of completeness, we also record the following refinement of Lemma \refcnt[N.20]{O.lem1}, which may be regarded as a statement about the existence of representative charts of a special type:

\begin{corollario}\label{O.cor4}
Let \G\ be a proper Lie groupoid over a manifold $M$. Assume that $(E,\varrho)$ is a classical representation of \G, mapping a subset $\G(x,x')$ injectively into $\Lis(E_x,E_{x'})$.

Then there exist open balls $B$ and $B'$ in $M$, centred at $x$ and $x'$ respectively, such that the restriction
$$
\varrho: \G(B,B') \to \GL(E)
$$
is an embedding of manifolds.
\end{corollario}
\begin{proof}
To begin with, observe that for any given arrow $g \in \G(x,x')$ and open neighbourhood $\Gamma$ of $g$ in \G\ there is an open ball $P$ inside $\GL(E)$, centred at $\varrho(g)$, such that $\varrho^{-1}(P) \subset \Gamma$. To see this, we fix a sequence
$$
\cdots \subset P_{i+1} \subset P_i \subset \cdots \subset P_1
$$
of open balls inside $\GL(E)$, centred at $\varrho(g)$ and with $\lim_i{\mathrm{radius}(P_i)} = 0$, and then we argue as in the proof of Theorem \refcnt[N.20]{O.thm1}.

By Lemma \refcnt[N.20]{O.lem1}, every $g \in \G(x,x')$ admits an open neighbourhood $\Gamma_g$ in $\G$ such that $\varrho$ induces a smooth isomorphism between $\Gamma_g$ and a submanifold of $\GL(E)$. As observed above, one can then choose an open ball $P_g \subset \GL(E)$ at $\varrho(g)$ such that $\varrho^{-1}(P_g) \subset \Gamma_g$. Now, let $\Gamma = {\bigcup \varrho^{-1}(P_g)}$. We claim that $\varrho$ induces a smooth isomorphism between $\Gamma$ and a submanifold of $\GL(E)$. By construction, $\varrho$ restricts to an immersion of $\Gamma$ into $\GL(E)$. If $g \in \G(x,x')$ then
$$
{\varrho(\Gamma) \cap P_g} = \varrho\bigl(\varrho^{-1}(P_g)\bigr)
$$
is an open subset of the submanifold $\varrho(\Gamma_g) \subset \GL(E)$, because $\varrho$ is a smooth isomorphism of $\Gamma_g$ onto $\varrho(\Gamma_g)$. Since the open balls $P_g$ cover $\varrho(\Gamma)$ as $g$ ranges over $\G(x,x')$, $\varrho(\Gamma)$ is a submanifold of $\GL(E)$. Moreover, since $\varrho$ is a local smooth isomorphism of $\Gamma$ onto $\varrho(\Gamma)$, it will be also a global diffeomorphism provided it is globally injective over $\Gamma$: now, if $\varrho(\gamma') = \varrho(\gamma)$ then $\gamma', \gamma \in \varrho^{-1}(P_g) \subset \Gamma_g$ for some $g$ and therefore $\gamma' = \gamma$ because $\varrho$ is injective over $\Gamma_g$.

Finally, one further application of the usual properness argument will yield open balls $B,B' \subset M$ at $x,x'$ such that $\G(B,B')$ is contained in $\Gamma$ (this is an open neighbourhood of $\G(x,x')$ in \G).
\end{proof}

Note that the preceding corollary entails in particular that the image $\varrho(\G)$ is a submanifold of $\GL(E)$ for every proper Lie groupoid \G\ and faithful classical representation $(E,\varrho)$ of \G.

\sottosezione{Technical notes}

\begin{nota}\label{O.sub5}% - xxiii.14, Riveduta 22/02/2008
Suppose one is given an isomorphism $E \iso F$ of vector bundles over a manifold $M$. Then there is an induced isomorphism of Lie groupoids over $M$ (ie one that restricts to the identity mapping on $M$)
\begin{equazione}\label{xxiii.15}
\GL(E) \xto\iso \GL(F)\text,
\end{equazione}
given, for each $(x,x') \in {M\times M}$, by the bijection that makes the linear isomorphisms $\alpha$ and $\beta$ correspond to each other when they fit in the diagram%
\begin{equazione}\label{xxiii.16}%
\begin{split}%
\xymatrix@R=19pt{E_x\ar[d]^{\iso_x}\ar[r]^-\alpha & E_{x'}\ar[d]^{\iso_{x'}} \\ F_x\ar[r]^-\beta & F_{x'}\text.\!\!}%
\end{split}%
\end{equazione}%

In particular, if two objects $R, S \in \Ob(\Kt)$ are isomorphic, any induced isomorphism $\fifu(\iso): \fifu(R) \iso \fifu(S)$ will in turn yield an isomorphism of the corresponding linear groupoids $\GL({\fifu R}) \iso \GL({\fifu S})$ (identical on $M$), such that for each $\lambda \in \tannakian\fifu$ the linear mappings $\lambda(R)$ and $\lambda(S)$ correspond to one another\nobreakdash---because of naturality of $\lambda$:%
\begin{equazione}\label{xxiii.17}%
\begin{split}%
\xymatrix@C=33pt{({\fifu R})_x\ar@{=}[r]\ar[d]^{\fifu(\iso)_x} & \fifu_x(R)\ar[d]^{\fifu_x(\iso)}\ar[r]^-{\lambda(R)} & \fifu_{x'}(R)\ar[d]^{\fifu_{x'}(\iso)}\ar@{=}[r] & ({\fifu R})_{x'}\ar[d]^{\fifu(\iso)_{x'}} \\ ({\fifu S})_x\ar@{=}[r] & \fifu_x(S)\ar[r]^-{\lambda(S)} & \fifu_{x'}(S)\ar@{=}[r] & ({\fifu S})_{x'}\text.\!\!}%
\end{split}%
\end{equazione}%
Thus, the latter isomorphism will transform $\ev_R$ into $\ev_S$:%
\begin{equazione}\label{xxiii.18}%
\begin{split}%
\xymatrix@C=65pt@R=3pt{& \GL({\fifu R})\ar@{<->}[dd]^\iso \\ \tannakian\fifu\ar[dr]_{\ev_S}\ar[ur]^{\ev_R} & \\ & \GL({\fifu S})\text.\!\!}%
\end{split}%
\end{equazione}%
It follows that if $\Omega \subset \tannakian{}$ is any open subset then $R(\Omega)$ is a tame submanifold of $\GL({\omega R})$ if and only if $S(\Omega)$ is a tame submanifold of $\GL({\omega S})$ (see, for instance, Note~\refcnt[O.2.1]{O.sub3}) and that $R(\Omega)$ and $S(\Omega)$ are homeomorphic subsets; hence $\ev_R$ will induce a homeomorphism between $\Omega$ and $R(\Omega)$ if and only if $\ev_S$ induces one between $\Omega$ and $S(\Omega)$.%
\end{nota}%

% - Da qui in poi tutto da rivedere - 22/02/2008
\begin{nota}\label{O.sub6}%
Let $\G$ and $\H$ be two Lie groupoids over the manifold $M$. We want to construct, provided this is possible, their product in the category of Lie groupoids over $M$. It ought to be a Lie groupoid over $M$ endowed with canonical projections, satisfying the usual universal property%
\begin{equazione}\label{O.equ35}%
\begin{split}%
\xymatrix@C=55pt{ & & \G \\ \K\ar@/_1pc/[drr]^-\psi\ar@/^1pc/[urr]^-\varphi\ar@{-->}[r]^(.45){(\varphi,\psi)} & {\G\times_M\H}\ar[ur]_-{\pr_1}\ar[dr]^-{\pr_2} & \\ & & \H\text.\!\!}%
\end{split}%
\end{equazione}%
It must be kept in mind that all the arrows in this diagram are morphisms of Lie groupoids over $M$, ie they all induce the identity map $\id: M \to M$ at the base level.%

The construction of the product over $M$ can be obtained as a special case of the so-called ``strong fibred product construction'' for Lie groupoids, cfr.\ for example \textit{Moerdijk and Mr\v cun (2003),} \cite{Moerdijk&Mrcun'03}~p.~123.%

Namely, we regard the maps%
$$%
\xymatrix@C=30pt{\ar@{..>}[r]\ar@{..>}[d] & \G\ar[d]^{(\s,\t)} & \text{(viz.} & \G\ar[d]^{(\s,\t)}\ar[r]^-{(\s,\t)} & {M \times M}\ar[d]^{(\pr_1,\pr_2) \, = \, \id} & \\ \H\ar[r]^-{(\s,\t)} & {M \times M} && {M \times M}\ar[r]^{\id \times \id} & {M \times M} & \text{etc.)}}%
$$%
as morphisms of lie groupoids over $M$, where ${M \times M}$ is the pair groupoid, and apply the strong fibred product construction to them:%
$$%
\left\{\begin{aligned}%
\text{set of arrows} &= \bigl\{(g,h) \in {\G\times \H}: (\s,\t)(g) = (\s,\t)(h)\bigr\}\text,%
\\%
\text{set of objects} &= \bigl\{(m,m') \in {M\times M}: m=m'\bigr\} \can M\text.
\end{aligned}\right.%
$$%
Transversality criteria imply that this defines a Lie groupoid ${\G\times_M\H}$ over $\Delta(M) \can M$ whenever, for instance, one of the two maps is a submersion. (Terminology: we say that a Lie groupoid $\G \rightrightarrows M$ is \index{groupoid!locally transitive}\index{locally transitive groupoid}{locally transitive} if the map $(\s,\t): \G \to {M\times M}$ is a submersion. This appears to be reasonable, since $\G$ is said to be \index{groupoid!transitive}\index{transitive groupoid}{transitive} if that map is a surjective submersion.) Moreover, if the trasversality condition is satisfied, this construction gives a fibred product with the familiar universal property.%

Suppose that ${\G \times_M \H}$ makes sense, ie that the transversality condition is satisfied. We remark that the universal property (\ref{O.equ35}) is a consequence of the universal property of the pullback. Indeed, first of all, the two projections of the fibred product to its own factors are morphisms over $M$, as one sees directly at once. Secondly, if $\varphi: \K \to \G$ and $\psi: \K \to \H$ are morphisms over $M$, then the following diagram commutes (precisely by definition of morphism over $M$)%
$$%
\xymatrix@C=40pt{\K\ar[d]_{\psi}\ar@{-->}[dr]^{(\s,\t)}\ar[r]^{\varphi} & \G\ar[d]^{(\s,\t)} \\ \H\ar[r]^-{(\s,\t)} & {M \times M}}%
$$%
and therefore there exists a unique morphism of Lie groupoids $(\varphi,\psi): \K \to {\G \times_M \H}$ such that diagram (\ref{O.equ35}) commutes, so we need only verify that $(\varphi,\psi)$ is in fact a morphism over $M$. This follows at once from the commutativity of the diagram%
$$%
\xymatrix{\K\ar[dd]^{(\s,\t)}\ar[rr]^-{(\varphi,\psi)}\ar[dr]^{\varphi} & & {\G \times_M \H}\ar[dd]^{(\s,\t)}\ar[dl]_{\pr_1} \\ & \G\ar[dr]^{(\s,\t)} & \\ {M \times M}\ar[rr]^{\id \times \id} & & \,{M \times M}\text{.}}%
$$%

Observation. By construction, the manifold of arrows of ${\G \times_M \H}$ is a submanifold of the Cartesian product ${\G \times \H}$; it follows that the subsets of the form ${\Gamma \times \Lambda}$, for $\Gamma \subset \G$ and $\Lambda \subset \H$ open, form a basis for the topology of ${\G \times_M \H}$. (Of course, we write ${\Gamma \times \Lambda}$ but we mean ${(\Gamma \times \Lambda) \cap (\G \times_M \H)}$.) Thus, one sees immediately that, when the differentiable structure is discarded, the same construction yields the product in the category of topological groupoids over $M$.%

Now, we apply this general construction to the locally transitive Lie groupoids $\GL(E)$ associated to vector bundles $E \in {\Ob\,\V[\infty]M}$. (These are locally transitive since if $E_U \iso {U \times \mathbf{E}}$ and $E_V \iso {V \times \mathbf{F}}$ are local trivializations of $E$, then up to diffeomorphism the map $(\s,\t)$ coincides locally with a projection%
$$%
\GL(E)(U,V) \iso {U \times V \times \Lis(\mathbf{E},\mathbf{F})} \xto{\pr} {U \times V}%
$$%
and is in particular a submersion; note that this makes sense even when $\Lis(\mathbf{E},\mathbf{F}) = \varnothing$.)%
\end{nota}%

\sezione{Morphisms of Fibre Functors}\label{O.2.3}
A \index{morphism of fibre functors|emph}morphism of fibre functors, let us say one $(\Kt,\fifu) \to (\Kt',\fifu')$, consists of a smooth map $f: M \to M'$ of the respective base manifolds together with a linear tensor functor $\Phi^*: \Kt' \longto \Kt$ and a tensor preserving isomorphism $\alpha$
\begin{equazione}\label{O.equ14}
\begin{split}
\xymatrix{\Kt'\ar[d]^{\fifu'}\ar[r]^-{\Phi^*} & \Kt\ar[d]^\fifu_{\,}="b" \\ \V[\infty]{M'}\ar[r]^-{f^*}="a"\ar@{=>}@/^0.5pc/"a";"b"^\alpha & \V[\infty]M\text,\!\!}%
\end{split}%
\end{equazione}%
where $f^*$ = pullback along $f$. In place of the correct $(f,\Phi^*,\alpha)$, our preferred notation for morphisms of fibre functors will be the incorrect $(f^*,\Phi^*)$, in order to emphasize the algebraic symmetry.%

Composition of morphisms is defined as%
\begin{equazione}\label{xxiv.1}%
{(g^*,\Psi^*) \cdot (f^*,\Phi^*)} = \bigl(({g\circ f})^*,{\Phi^*\circ\Psi^*}\bigr)\text.%
\end{equazione}%
Note that if in our definition we required \refequ{O.equ14} to commute in the strict sense we would get into trouble because $({g\circ f})^* \can {f^*\circ g^*}$ are canonically isomorphic but not really identical tensor functors.%

Lemmas~\refcnt{O.lem102} and \refcnt{O.lem103} below apply directly to \refequ{O.equ14} to yield maps%
\begin{multiriga}\label{xxiv.2}%
\Hom^\otimes(\fifu_x,\fifu_y) \xto{\text{Lem.~\ref{O.lem102}}} \Hom^\otimes\bigl({\fifu_x\circ\Phi^*},{\fifu_y\circ\Phi^*}\bigr)%
\\%
= \Hom^\otimes\bigl({x^*\circ\fifu\circ\Phi^*},{y^*\circ\fifu\circ\Phi^*}\bigr)%
\\%
\xto{\iso\text{~(\ref{O.equ14})~+~Lem.~\ref{O.lem103}}} \Hom^\otimes\bigl({x^*\circ f^*\circ\fifu'},{y^*\circ f^*\circ\fifu'}\bigr)%
\\%
\xto{\can\text{~Lem.~\ref{O.lem103}}} \Hom^\otimes\bigl({{f(x)}^* \circ \fifu'},{{f(y)}^* \circ \fifu'}\bigr)%
\\%
= \Hom^\otimes\bigl(\fifu'_{f(x)},\fifu'_{f(y)}\bigr)\text.%
\end{multiriga}%
Moreover, since ${({\lambda\circ\mu}) \cdot \Phi^*} = {({\lambda\cdot\Phi^*}) \circ ({\mu \cdot \Phi^*})}$ and ${\id \cdot \Phi^*} = \id$, these can be pieced together in a functorial way, so that they form a homomorphism of groupoids%
\begin{equazione}\label{O.equ15}%
\begin{split}%
\xymatrix@C=35pt{\tannakian{\fifu}\ar[d]\ar[r]^-\Phi & \tannakian{\fifu'}\ar[d] \\ {M \times M}\ar[r]^-{f \times f} & {M' \times M'}\text,\!\!}%
\end{split}%
\end{equazione}%
which can be characterized as the unique map making%
\begin{equazione}\label{O.equ16}%
\begin{split}%
\xymatrix@C=35pt{\tannakian{\fifu}\ar[d]^{\ev_{\Phi^* R'}}\ar[r]^-\Phi & \tannakian{\fifu'}\ar[d]^{\ev_{R'}} \\ \GL(\fifu {\Phi^* R'})\ar[r]^-{\gamma \circ \alpha_*^{-1}} & \GL(\fifu' R')}%
\end{split}%
\end{equazione}%
commute for all $R' \in \Ob(\Kt')$, where the morphism $\gamma$ is the ``projection'' $\GL({f^*(\fifu'R')}) \can ({f\times f})^*(\GL(\fifu'R')) \to \GL(\fifu'R')$ and the isomorphism%
\begin{equazione}\label{xxiv.3}%
\alpha_*: \GL({f^*\fifu'R'}) \isoto \GL({\fifu\Phi^*R'})%
\end{equazione}%
comes from $\alpha_{R'}: {f^*\fifu'R'} \isoto {\fifu\Phi^*R'}$ according to Note~\refcnt[O.2.2+N.23]{O.sub5}. It is also immediate from \refequ{O.equ16} that such a solution $\Phi$ is necessarily a morphism of $\C^\infty$\nobreakdash-functionally structured spaces, so \refequ{O.equ15} proves to be a homomorphism of $\C^\infty$\nobreakdash-functionally structured groupoids.%

We shall refer to $\Phi$ as the {\em realization} of the morphism $(f^*,\Phi^*)$. This construction is functorial with respect to composition of morphisms of fibre functors, and therefore defines a functor into the category of $\C^\infty$\nobreakdash-structured groupoids, called the {\em realization functor}.%

\begin{proposizione}\label{O.prp102}%
Let $(\Kt,\fifu), (\Kt',\fifu')$ be smooth classical fibre functors and$$(f^*,\Phi^*): (\Kt,\fifu) \to (\Kt',\fifu')$$a morphism of fibre functors. Then the corresponding realization is a \index{homomorphism of groupoids}homomorphism of Lie groupoids.
\end{proposizione}%
\begin{proof}%
It follows from \refequ{O.equ16} that the composite ${\ev_{R'} \circ \Phi}$ is smooth for every object $R'$ of $\Kt'$. The map $\Phi$ is then smooth by the ``universal property'' of the Lie groupoid $\tannakian{\fifu'}$.%
\end{proof}%

\sottosezione{Notes}%

\begin{nota}\label{O.sub9}%
In this note we recall a couple of elementary properties of tensor functors and tensor preserving natural transformations.%
\begin{lemma}\label{O.lem102}
Let $F$, $G$, $S$, $T$ be tensor functors relating suitable tensor categories. Then%
\begin{elenco}%
\item[1.]the rule $\lambda \mapsto {\lambda\cdot S}$ maps $\Hom^\otimes(F,G)$ into $\Hom^\otimes({F\circ S},{G\circ S})$;%
\item[2.]the rule $\lambda \mapsto {T \cdot\lambda}$ maps $\Hom^\otimes(F,G)$ into $\Hom^\otimes({T\circ F},{T\circ G})$.%
\end{elenco}%
\end{lemma}%
\begin{proof}%
(1) The natural transformation $(\lambda \cdot S)(X) = \lambda(S X)$ is a morphism of tensor functors if such is $\lambda$, because%
$$%
\xymatrix{{{F {S X}} \otimes {F {S Y}}}\ar[d]^{\can}\ar[rr]^{\lambda(S X) \otimes \lambda(S Y)} & & {{G {S X}} \otimes {G {S Y}}}\ar[d]^{\can} & & \TU\ar[d]^{\can}\ar[r]^{\id} & \TU\ar[d]^{\can} \\ F({S X} \otimes {S Y})\ar[d]^{F \, \can}\ar[rr]^{\lambda({S X} \otimes {S Y})} & & G({S X} \otimes {S Y})\ar[d]^{G \, \can} & & {F \TU}\ar[d]^{F \, \can}\ar[r]^{\lambda(\TU)} & {G \TU}\ar[d]^{G \, \can} \\ {F S(X \otimes Y)}\ar[rr]^{\lambda(S(X \otimes Y))} & & {G {S (X \otimes Y)}} & & {F {S \TU}}\ar[r]^-{\lambda(S \TU)} & \,{G {S \TU}}\text{.}}%
$$%

\noindent(2) The same can be said of $(T \cdot \lambda)(X) = T(\lambda(X))$, since%
$$%
\xymatrix{{{T {F X}} \otimes {T {F Y}}}\ar[d]^{\can}\ar[rr]^{{T \lambda(X)} \otimes {T \lambda(Y)}} & & {{T {G X}} \otimes {T {G Y}}}\ar[d]^{\can} & & \TU\ar[d]^{\can}\ar[r]^{\id} & \TU\ar[d]^{\can} \\ T({F X} \otimes {F Y})\ar[d]^{T \, \can}\ar[rr]^{T(\lambda(X) \otimes \lambda(Y))} & & T({G X} \otimes {G Y})\ar[d]^{T \, \can} & & {T \TU}\ar[d]^{T \, \can}\ar[r]^{T(\id)} & {T \TU}\ar[d]^{T \, \can} \\ {T {F (X \otimes Y)}}\ar[rr]^{T \lambda(X \otimes Y)} & & {T {G (X \otimes Y)}} & & {T {F \TU}}\ar[r]^{T \lambda(\TU)} & \,{T {G \TU}}\text{.}}%
$$%
\end{proof}%

Let $(\mathcal{C},\otimes)$ and $(\mathcal{V},\otimes)$ be tensor categories. Suppose that$$F, F', G, G': \mathcal{C} \longto \mathcal{V}$$are tensor functors, and that $F \iso F', G \iso G'$ are tensor preserving natural isomorphisms. For every $X \in \Ob(\mathcal{C})$, there is an obvious bijective map $a \mapsto a'$ determined by the commutativity of%
\begin{equazione}\label{O.equ33}%
\begin{split}%
\xymatrix{{F X}\ar[d]^{\iso}\ar[r]^{a} & {G X}\ar[d]^{\iso} \\ {F' X}\ar[r]^{a'} & \,{G' X}\text{.}}%
\end{split}%
\end{equazione}%
Given a natural transformation $\lambda \in \Hom(F,G)$, we put $\lambda'(X) = {\lambda(X)}'$.%
\begin{lemma}\label{O.lem103}%
The rule which to $\lambda$ associates $\lambda'$ determines a bijective correspondence\begin{equazione}\label{O.equ34}\Hom^\otimes(F,G) \isoto \Hom^\otimes(F',G')\text{.}\end{equazione}%
\end{lemma}%
\begin{proof}%
Obvious.%
\end{proof}%
\end{nota}%

\sezione{Weak Equivalences}\label{O.2.4}
\begin{definizione}\label{xxv1}
A \index{weak equivalence|emph}\textit{weak equivalence}\footnote{Note on terminology: We shall reserve the term `weak equivalence' for the context of fibre functors. When dealing with Lie groupoids, we prefer to use the term `Morita equivalence'.} of fibre functors, symbolically $(\Kt,\fifu) \xto{\iso} (\Kt',\fifu')$, is a morphism of fibre functors
$$%
(f^*,\Phi^*): (\Kt,\fifu) \to (\Kt',\fifu')%
$$%
satisfying the following two conditions%
\begin{elenco}
\item[1.]the base mapping $f: M \to M'$ is a surjective submersion;
\item[2.]the functor $\Phi^*$ is a \index{tensor equivalence}tensor equivalence, ie there exist a tensor functor $\Phi_*: \Kt \longto \Kt'$ and tensor preserving natural isomorphisms
$$%
\left\{\begin{array}{l}%
{\Phi^*\circ\Phi_*} \iso \Id_{\Kt}%
\\%
{\Phi_*\circ\Phi^*} \iso \Id_{\Kt'}%
\text.%
\end{array}\right.%
$$%
\end{elenco}%
\end{definizione}
In order to conclude that $\Phi^*$ is a tensor equivalence, it suffices to know it to be an ordinary categorical equivalence. Every quasi-inverse equivalence $\Phi_*$ is then necessarily a linear functor. (Details may be found in Note~\refcnt{O.sub10}.) Weak equivalences of fibre functors are stable under composition of morphisms of fibre functors, as defined in Section \ref{O.2.3}.%

\begin{proposizione}\label{O.prp7}%
Let%
$$%
(f^*,\Phi^*): (\Kt,\fifu) \xto{\iso} (\Kt',\fifu')%
$$%
be a weak equivalence of fibre functors. Then its realization diagram%
\begin{equazione}\label{O.equ36}%
\begin{split}%
\xymatrix@C=35pt{\tannakian{\fifu}\ar[d]\ar[r]^-\Phi & \tannakian{\fifu'}\ar[d] \\ {M\times M}\ar[r]^-{f\times f} & {M'\times M'}}%
\end{split}%
\end{equazione}%
is a topological pullback, ie a pullback in the category of topological spaces, and $\Phi: \tannakian{\fifu} \epito \tannakian{\fifu'}$ is a surjective open mapping.%
\end{proposizione}%
\begin{proof}%
Let $T$ be a topological space, and suppose given a problem%
\begin{equazione}\label{O.equ22}%
\begin{split}%
\xymatrix{T\ar@/_1.5pc/[dddr]\ar@{..>}[dr]^a\ar@/^1pc/[drr] & & \\ & \tannakian{\fifu}\ar[d]^{\ev_{\Phi^*R'}}\ar[r]^-\Phi & \tannakian{\fifu'}\ar[d]^{\ev_{R'}} \\ & \GL(\fifu {\Phi^* R'})\ar[d]\ar[r]^-{\gamma \circ \alpha_*^{-1}} & \GL(\fifu' R')\ar[d] \\ & {M \times M}\ar[r]^-{f \times f} & {M' \times M'}}%
\end{split}%
\end{equazione}%
stated in the category of topological spaces and continuous mappings. There exists a unique set-theoretic solution $a$, because \refequ{O.equ36} is already known to be a set-theoretic pullback (by Note~\refcnt{O.sub10} again). Thus, we must check that $a$ is continuous. Note that ${\forall R}$ in $\Kt$, ${\ev_R \circ a}$ is continuous if and only if ${\ev_{\Phi^* \Phi_* R} \circ a}$ is continuous, because of the isomorphism ${\Phi^* \Phi_* R} \iso R$, see also the comments in Note~\refcnt[O.2.2+N.23]{O.sub5}. Therefore, if we put $R' = {\Phi_* R}$ in \refequ{O.equ22}, we conclude at once that ${\ev_{\Phi^* R'} \circ a}$ is continuous from the fact that the lower square of \refequ{O.equ22} is, by definition, a topological pullback.%

Next, observe that if one has a topological pullback%
\begin{equazione}\label{O.equ20}%
\begin{split}%
\xymatrix{X\ar[d]^p\ar[r]^f & Y\ar[d]^q \\ M\ar[r]^g & N}%
\end{split}%
\end{equazione}%
along a submersive morphism $g$ of smooth manifolds, there is the following local decomposition up to diffeomorphism%
\begin{equazione}\label{O.equ21}%
\begin{split}%
\xymatrix{X_U\ar[d]^p\ar[r]^f & Y_V\ar[d]^q \\ U\ar[r]^g & V}\qquad\xymatrix{X_U\ar[d]^p\ar@{-}[r]^-\iso & {Y_V \times P}\ar[d]^{q \times \id}\ar[r]^-\pr & Y_V\ar[d]^q \\ U\ar@{-}[r]^-\iso & {V\times P}\ar[r]^-{\pr} & V\text,\!\!}%
\end{split}%
\end{equazione}%
where $U \subset M$ is open and so small that, up to diffeomorphism, $g|_U$ is a projection ${V\times P} \to V = g(U)$ for some open ball $P$; of course, $X_U = p^{-1}(U)$ etc. (Note that in \refequ{O.equ21}, $U \iso {V\times P}$ is a diffeomorphism whereas $X_U \iso {Y_V\times P}$ is a homeomorphism.) It follows that $f$ is a `topological submersion', in particular an open mapping; in addition, if $g$ is surjective then it is clear that $f$ must be also surjective. This shows that the statement that $\Phi$ is an open mapping follows from the statement that \refequ[O.2.4]{O.equ36} is a topological pullback.%
\end{proof}%

Suppose a topological pullback \refequ{O.equ20} along a smooth submersion is given, and let $U \subset M$ be an open subset such that $g|_U$ is, up to diffeomorphism, a projection $U \iso {V\times P} \xto\pr V$ onto an open subset $V \subset N$. Let $A \subset X$ be an open subset, and put $B = f(A)$; $B \subset Y$ is open because $f$ is an open mapping. We shall be interested in the subspaces $p(A) \subset M$ and $q(B) \subset N$; note that $g$ restricts to a continuous mapping of $p(A)$ onto $q(B)$. Assume that $A$ has the following property: {\em the commutative square%
\begin{equazione}\label{O.equ26}%
\begin{split}%
\xymatrix{{A\cap p^{-1}(U)}\ar[d]^p\ar[r]^-f & {B\cap q^{-1}(V)}\ar[d]^q \\ U\ar[r]^-{g} & V}%
\end{split}%
\end{equazione}%
is a topological pullback}. Then there is a trivialization, analogous to \refequ{O.equ21}, which shows that the smooth iso $U \iso {V\times P}$ induces a correspondence between$${p(A)\cap U} = p\bigl(A\cap p^{-1}(U)\bigr)$$and$${\bigl(q(B)\cap V\bigr)\times P} = {q\bigl(B\cap q^{-1}(V)\bigr)\times P}\text{.}$$Thus, ${\forall u} \in U$ one has ${u \in p(A)} \aeq {g(u) \in q(B)}$. Note also that $p$ restricts to a homeomorphism of ${A\cap p^{-1}(U)}$ onto ${p(A)\cap U}$ if and only if $q$ restricts to a homeomorphism of ${B\cap q^{-1}(V)}$ onto ${q(B)\cap V}$. The two relevant cases for the present discussion occur, in the first place, when $A = f^{-1}(f(A))$, and secondly, when $A \subset p^{-1}(U)$ coincides with ${B\times P}$ in the trivialization \refequ{O.equ21}.%

\separazione%

Fix an object $R' \in \Ob(\Kt')$. Then the outer rectangle of \refequ{O.equ22} is a topological pullback\textemdash note that it coincides with \refequ{O.equ36}; the lower square enjoys the same property. Consequently, the upper square, viz \refequ[O.2.3]{O.equ16}, must be a topological pullback as well; moreover, since the smooth mapping ${\gamma \circ \alpha_*^{-1}}: \GL(\fifu{\Phi^*R'}) \to \GL(\fifu'R')$ is a (surjective) submersion, it is a pullback of the form \refequ{O.equ21}. Hence the preceding remarks apply, and we get:%

\vskip 5pt%
{\em 1. If $(\Omega',R')$ is a representative chart of $(\Kt',\fifu')$ then $(\Phi^{-1}(\Omega'),{\Phi^*R'})$ is a representative chart of $(\Kt,\fifu)$.} Since diagram \refequ[O.2.3]{O.equ16} is a topological pullback,%
$$%
\xymatrix{\Phi^{-1}(\Omega')\ar[d]^{\ev_{\Phi^*R'}}\ar[r]^-{\Phi} & \Omega'\ar[d]^{\ev_{R'}} \\ {\Phi^*R'}(\Phi^{-1}(\Omega'))\ar[r]^-{\gamma \circ \alpha_*^{-1}} & R'(\Omega')}%
$$%
is also a topological pullback and therefore $\ev_{\Phi^*R'}$ induces a homeomorphism between $\Phi^{-1}(\Omega')$ and its image ${\Phi^*R'}(\Phi^{-1}(\Omega'))$, because $\ev_{R'}$, on the right, does the same. Proposition~\refcnt[O.2.1]{O.prp1} implies that ${\Phi^*R'}(\Phi^{-1}(\Omega'))$ is a tame submanifold of $\GL(\fifu{\Phi^*R'})$ if and only if $R'(\Omega')$ is a tame submanifold of $\GL(\fifu'R')$, because ${\gamma \circ \alpha_*^{-1}}$ is a Morita equivalence and $\Omega' = \Phi(\Phi^{-1}(\Omega'))$.%

\vskip 5pt%
{\em 2. Let $\Omega \subset \tannakian{\fifu}$ be an open subset and $\lambda_0 \in \Omega$. For any given object $R' \in \Ob(\Kt')$, there is a smaller open neighbourhood $\lambda_0 \in \Omega_0 \subset \Omega$ such that $(\Omega_0,{\Phi^*R'})$ is a representative chart of $(\Kt,\fifu)$ if and only if $(\Phi(\Omega_0),R')$ is a representative chart of $(\Kt',\fifu')$.} Let $\Lambda$ be an open neighbourhood of $\lambda_0(\Phi^* R')$ in $\GL(\fifu\Phi^*R')$ such that ${\gamma \circ \alpha_*^{-1}}|_\Lambda$ is, up to diffeomorphism, a projection ${\Lambda' \times P} \to \Lambda' = {\gamma \circ \alpha_*^{-1}}(\Lambda)$. Making the open ball $P$, and thus $\Lambda$, smaller if necessary, we find an open neighbourhood $\Omega_0 \subset {\ev_{\Phi^* R'}^{-1}(\Lambda) \cap \Omega}$ of $\lambda_0$ such that the homeomorphism $\ev_{\Phi^* R'}^{-1}(\Lambda) \iso {\ev_{R'}^{-1}(\Lambda') \times P}$ of \refequ{O.equ21} produces a decomposition%
\begin{equazione}\label{O.equ25}%
\begin{split}%
\xymatrix{\Omega_0\ar[d]\ar[r]^-\Phi & \Phi(\Omega_0)\ar[d] \\ \Lambda\ar[r]^-{\gamma \circ \alpha_*^{-1}} & \Lambda'}\qquad\quad\xymatrix{\Omega_0\ar[d]\ar@{-}[r]^-\iso & {{\Phi(\Omega_0)}\times P}\ar[d]^{\times\id}\ar[r]^-\pr & \Phi(\Omega_0)\ar[d] \\ \Lambda\ar@{-}[r]^-\iso & {\Lambda'\times P}\ar[r]^-\pr & \Lambda'}%
\end{split}%
\end{equazione}%
Therefore, if we put $\Sigma = {\Phi^*R'}(\Omega_0) \subset \Lambda$ and $\Sigma' = R'(\Phi(\Omega_0)) \subset \Lambda'$ we have $\lambda \in \Sigma \aeq {\gamma \circ \alpha_*^{-1}\lambda} \in \Sigma'$ for all $\lambda \in \Lambda$, and Proposition~\refcnt[O.2.1]{O.prp1} implies that $\Sigma$ is a tame submanifold of $\GL(\fifu{\Phi^*R'})$ if and only if $\Sigma'$ is a tame submanifold of $\GL(\fifu'R')$, since ${\gamma \circ \alpha_*^{-1}}$ is a Morita equivalence.%

\vskip 5pt%
Clearly, these statements imply that whenever a weak equivalence of fibre functors $(\Kt,\fifu) \xto{\iso} (\Kt',\fifu')$ is given, Condition~\textsl{i)} of Proposition~\refcnt[O.2.2+N.23]{xxiii.4} holds for $(\Kt,\fifu)$ if and only if it holds for $(\Kt',\fifu')$. (As a consequence of the fact that $\Phi$ is surjective and open: Fix $\lambda'_0 = \Phi(\lambda_0)$. If $(\Omega,R)$ is a chart at $\lambda_0$, then $(\Phi(\Omega_0),{\Phi_*R})$ is a chart at $\lambda'_0$ for some open $\lambda_0 \in \Omega_0 \subset \Omega$; conversely, if $(\Omega',R')$ is a chart at $\lambda'_0$ then $(\Phi^{-1}(\Omega'),{\Phi^*R'})$ is a chart at $\lambda_0$.)%

On the other hand, they also imply invariance of Condition~\textsl{ii)} of the same proposition, p.~\pageref{xxiii.4}, as follows. Assume the condition holds for $(\Kt',\fifu')$: Let $(\Omega,R)$ be a chart of $(\Kt,\fifu)$ and $S \in \Ob(\Kt)$ an object. Choose a point $\lambda_0 \in \Omega$. There exists a neighbourhood $\Omega_0 \subset \Omega$ of $\lambda_0$ such that $(\Phi(\Omega_0),{\Phi_*R})$, and consequently $(\Phi(\Omega_0),{{\Phi_*R}\oplus{\Phi_*S}})$, is a chart of $(\Kt',\fifu')$. Since $\Omega_0 \subset {\Phi^{-1}{\Phi(\Omega_0)}}$ and $\Phi^*({\Phi_*R}\oplus{\Phi_*S}) \iso {R\oplus S}$, it follows that $(\Omega_0,{R\oplus S})$ is a chart of $(\Kt,\fifu)$. Since $\lambda_0$ was arbitrary, we conclude that $\Omega$ can be covered with open subsets $\Omega_0$ such that $(\Omega_0,{R\oplus S})$ is a chart, and therefore that $(\Omega,{R\oplus S})$ is a chart as well. Conversely, assume Condition 2 holds for $(\Kt,\fifu)$: Let $(\Omega',R')$ be a chart of $(\Kt',\fifu')$ and $S' \in \Ob(\Kt')$ an object. Fix a point $\lambda'_0 \in \Omega'$; since $\Phi$ is surjective, ${\exists \lambda_0}$ with $\lambda'_0 = \Phi(\lambda_0)$. Since $(\Phi^{-1}(\Omega'),{\Phi^*R'})$ is a chart, $(\Phi^{-1}(\Omega'),{{\Phi^*R'}\oplus{\Phi^*S'}})$ and, consequently, $(\Phi^{-1}(\Omega'),\Phi^*(R'\oplus S'))$ are charts of $(\Kt,\fifu)$ as well. Hence there exists a neighbourhood $\Omega_0 \subset \Phi^{-1}(\Omega')$ of $\lambda_0$ such that $(\Phi(\Omega_0),{R'\oplus S'})$ is a chart of $(\Kt',\fifu')$. As before, since $\lambda'_0$ was arbitrary it follows that $(\Omega',{R'\oplus S'})$ is a chart of $(\Kt',\fifu')$.

\separazione

\noindent We can collect our conclusions in the following
\begin{proposizione}\label{O.prp8}
Let%
$$%
(f^*,\Phi^*): (\Kt,\fifu) \longto (\Kt',\fifu')%
$$%
be a weak equivalence of fibre functors. Then $(\Kt,\fifu)$ is a \index{fibre functor!smooth}\index{smooth fibre functor}smooth classical fibre functor if and only if so is $(\Kt',\fifu')$. In this case,%
$$%
(f,\Phi): \tannakian{\fifu} \longto \tannakian{\fifu'}
$$
is a \index{Morita equivalence}Morita equivalence of Lie groupoids.
\end{proposizione}
\begin{proof}
That \refequ[O.2.3]{O.equ15} is a pullback in the category of manifolds of class $\C^\infty$ follows by the same argument used in the proof of Proposition~\refcnt{O.prp7}, because of the universal property of the Tannakian groupoid.%
\end{proof}%

\sottosezione{Notes}%

\begin{nota}\label{O.sub10}%
List of elementary facts.%

1. \textsl{Any quasi-inverse equivalence $\Phi_*$ is automatically a linear functor.} Indeed, the map%
$$%
\label{O.esp1}\Hom_{\Kt}(R,S) \to \Hom_{\Kt}({\Phi^* {\Phi_* R}},{\Phi^* {\Phi_* S}})\text{,} \quad a \mapsto {\Phi^* {\Phi_* a}}%
$$%
is a linear bijection, as it is clear from the commutativity of$$\label{O.esp2}\xymatrix{{\Phi^* {\Phi_* R}}\ar[d]_{\Phi^* {\Phi_* a}}\ar[r]^-{\iso_R} & R\ar[d]^{a} \\ {\Phi^* {\Phi_* S}}\ar[r]^-{\iso_S} & \,S\text{,}}$$and the functor $\Phi^*$ is linear and, being a categorical equivalence, faithful, hence the equality ${\Phi^* \Phi_*({\alpha a} + {\beta b})} = {{\alpha \, {\Phi^* {\Phi_* a}}} + {\beta \, {\Phi^* {\Phi_* b}}}} = \Phi^*({\alpha \, {\Phi_* a}} + {\beta \, {\Phi_* b}})$ implies the desired linearity $\Phi_*({\alpha a} + {\beta b}) = {{\alpha \, {\Phi_* a}} + {\beta \, {\Phi_* b}}}$.%

2. \textsl{The realization $\Phi: \tannakian\fifu \longto \tannakian{\fifu'}$ of a weak equivalence is a fully faithful morphism of groupoids, in other words \refequ{O.equ36} is a set-theoretic pullback.} This can be seen as follows.%

The tensor preserving isomorphism ${\Phi^* \circ \Phi_*} \iso \Id_{\Kt}$ gives, according to Lemma~\refcnt[O.2.3]{O.lem102} p.~\pageref{O.lem102}, a tensor preserving isomorphism%
\begin{equazione}\label{O.esp4}%
\fifu_x \iso {\fifu_x \circ {\Phi^* \circ \Phi_*}} \iso {\fifu'_{f(x)} \circ \Phi_*}\text{;}\end{equazione}%
similarly, ${\Phi_* \circ \Phi^*} \iso \Id_{\Kt'}$ yields another such isomorphism\begin{equazione}\label{O.esp5}%
\fifu'_{f(x)} \iso {\fifu'_{f(x)} \circ {\Phi_* \circ \Phi^*}}\text{.}\end{equazione}%
If now we apply Lemma~\refcnt[O.2.3]{O.lem103} p.~\pageref{O.lem103} to these, we conclude at once from the commutativity of the diagram%
$$%
\xymatrix{\Hom^\otimes(\fifu_x,\fifu_y)\ar[d]^{\text{(\ref{O.esp4})}}_{\iso}\ar[r]^-{\Phi_{x,y}} & \Hom^\otimes\bigl(\fifu'_{f(x)},\fifu'_{f(y)}\bigr)\ar[d]^{\text{(\ref{O.esp5})}}_{\iso}\ar[dl] \\ \Hom^\otimes\bigl({\fifu'_{f(x)} \Phi_*},{\fifu'_{f(y)} \Phi_*}\bigr)\ar[r] & \Hom^\otimes\bigl({\fifu'_{f(x)}{\Phi_*\Phi^*}},{\fifu'_{f(y)}{\Phi_*\Phi^*}}\bigr)}%
$$%
that the diagonal arrow is a surjective and injective map,  and hence that $\Phi_{x,y}$ is bijective. (The commutativity of the two triangles follows from the commutativity of the two squares%
$$%
\xymatrix@C=30pt{\fifu_x(\Phi^*{\Phi_*R})\ar[r]^{\lambda_{\Phi^*{\Phi_*R}}} & \fifu_y(\Phi^*{\Phi_*R}) & \fifu'_{f(x)}(R')\ar[d]^{\fifu'_{f(x)}\,\iso}\ar[r]^{\lambda_{R'}} & \fifu'_{f(y)}(R')\ar[d]^{\fifu'_{f(y)}\,\iso} \\ \fifu_x(R)\ar[u]_{\fifu_x\,\iso}\ar[r]^{\lambda_R} & \fifu_y(R)\ar[u]_{\fifu_y\,\iso} & \fifu'_{f(x)}(\Phi_*{\Phi^*R'})\ar[r]^{\lambda_{\Phi_*{\Phi^*R'}}} & \fifu'_{f(y)}(\Phi_*{\Phi^*R'})}%
$$%
expressing naturality of $\lambda, \lambda'$ respectively.)%
\end{nota}%

\begin{nota}\label{O.sub11}
Let $X$ and $Y$ be topological spaces, and let $M$ and $N$ be smooth manifolds. Suppose%
\begin{equazione}\label{O.esp8}%
\begin{split}%
\xymatrix{X\ar[d]^p\ar[r]^f & Y\ar[d]^q \\ M\ar[r]^g & N}%
\end{split}%
\end{equazione}%
is a pullback diagram in the category of topological spaces, where $g$ is a smooth mapping.%

1. Given an open subset $B \subset Y$, put $A = f^{-1}(B)$. Then the continuous maps in (\ref{O.esp8}) restrict to a commutative diagram of topological spaces%
\begin{equazione}\label{O.esp9}%
\begin{split}%
\xymatrix{A\ar[d]^{p}\ar[r]^{f} & B\ar[d]^{q} \\ p(A)\ar[r]^{g} & \,q(B)\text{,}}\end{split}\end{equazione}%
which is again a topological pullback. Observe that if the restriction $q|_B$ induces a homeomorphism of $B$ onto $q(B)$, then $p|_A$ induces one between $A$ and $p(A)$. (This is a general property of pullbacks. Indeed, from$$\label{O.esp11}\xymatrix{C\ar[r]^{g}\ar@{..>}[dr]_{p'}\ar@/_1pc/[ddr]_{\id} & D\ar@/^0.5pc/[dr]^{q^{-1}} & \\ & A\ar[d]^{p}\ar[r]^{f} & B\ar[d]^{q} \\ & C\ar[r]^{g} & D}$$and from the equalities ${f \, {p' p}} = f$ and ${p \, {p' p}} = p$, it follows that ${p' p} = \id$, thus $p$ is invertible.)

2. Given an open subset $U \subset M$ such that $V = g(U)$ is open,%
\begin{equazione}\label{O.esp10}%
\begin{split}%
\xymatrix{p^{-1}(U)\ar[r]^{f}\ar[d]^{p} & q^{-1}(V)\ar[d]^{q} \\ U\ar[r]^{g} & V}\end{split}\end{equazione}%
makes sense and is clearly also a topological pullback.%
\end{nota}%

\capitolo[Classical Tannaka Theory]{Study of Classical Tannaka Theory of Lie Groupoids}\label{6}

In this conclusive chapter we are ideally going back to the point where we started from, namely the theory of classical representations of Lie groupoids expounded in \refsez{O.3.1+.4.1}. We will try to see what can be said about such theory by the light of the general results of Chapters~\refcpt[4]{5}. In particular, we will study in detail the standard classical fibre functor associated with a Lie groupoid. Recall that in \refsez{O.3.1+.4.1} we introduced the category \index{R(T;k)@\R[\infty]{\mathcal T;k} (category of smooth representations on vector bundles)}\R[\infty]{\G} of \index{C infinity representation@\ensuremath{C^\infty}-representation}\index{representation!C infinity or smooth@\ensuremath{C^\infty}- or smooth}classical representations $R = (E,\varrho)$ of a Lie groupoid \G, along with the standard classical fibre functor \index{omega(G)@\forget[T]\G, \forget[\infty]{\G} (forgetful functor)|emph}\index{standard fibre functor!classical@(classical) \forget[\infty]\G}\forget[\infty]{\G} defined as the forgetful functor $(E,\varrho) \mapsto E$ of \R[\infty]{\G} into the category \V[\infty]M of smooth vector bundles of locally finite rank over the base $M$ of \G. Let us give a brief review of the items we will be interested in, so as to fix the tacit notational conventions to be followed throughout the chapter.

Let \index{T(G)@\tannakian[T]\G, \tannakian[\infty]{\G} (Tannakian groupoid associated with a Lie groupoid)|emph}\index{Tannakian groupoid!T(G) classical@(classical) \tannakian[\infty]\G|emph}\tannakian[\infty]{\G} denote the Tannakian groupoid \xtannakian{\forget[\infty]{\G};\nR} associated with the fibre functor \forget[\infty]\G. Note that it does not make any difference whether we use real or complex coefficients in our theory, because eventually the groupoid \tannakian[\infty]{\G} and the other related items discussed below will be exactly the same; in fact, all what we are going to say holds for real as well as for complex coefficients: for simplicity, we assume real coefficients whenever we need to write them down explicitly. Recall from \refsez{N.21} that the Tannakian construction defines an operation
$$%
\G \mapsto \tannakian[\infty]\G\text, \quad \bigl\{\text{Lie~groupoids}\bigr\} \longto \bigl\{\text{$\C^\infty$-func.~structured~groupoids}\bigr\}\text;%
$$%
also note that the source and target map of \tannakian[\infty]{\G} are {\em submersions,} in the sense that they admit local sections which are morphisms of functionally structured spaces: this follows from the existence of such sections for \G\ and the fact that the envelope homomorphism \envelope[\infty]{} (see below) is a morphism of functionally structured spaces.%

Next, observe that for each Lie groupoid homomorphism $\varphi: \G \to \H$ the constructions of \refsez{O.2.3} may be applied to the equation ${\forget[\infty]\G \circ \varphi^*} = {f^* \circ \forget[\infty]\H}$ (identity of tensor functors), so as to yield a \index{homomorphism of groupoids}homomorphism of $\C^\infty$\nobreakdash-functionally structured groupoids
$$
\tannakian[\infty]\varphi: \tannakian[\infty]\G \to \tannakian[\infty]\H\text.%
$$%
In spite of the lack of functoriality of the operation $\varphi \mapsto \varphi^*$, in other words in spite of $({\psi\circ\varphi})^* \can {\varphi^* \circ \psi^*}$ being canonically isomorphic but not equal, the correspondence $\varphi \mapsto \tannakian[\infty]\varphi$ actually turns out to be a functor, i.e.\ the identities $\tannakian[\infty]{\psi\circ\varphi} = {\tannakian[\infty]\psi \circ \tannakian[\infty]\varphi}$ and $\tannakian[\infty]\id = \id$ hold.%

We let \index{envelope homomorphism@envelope homomorphism \envelope[T]\G, \envelope[\infty]\G|emph}\index{pi(G)@\envelope[T]\G, \envelope[\infty]{\G} (envelope homomorphism)|emph}\envelope[\infty]{\G} or, when there is no ambiguity, \envelope[\infty]{} denote the envelope homomorphism $\G \to \tannakian[\infty]\G$ defined by ${\envelope[\infty]g}(E,\varrho) = \varrho(g)$. The results of \refsez{N.20} concerning envelope homomorphisms can be applied. In particular, \envelope[\infty]{\G} will be a morphism of $\C^\infty$\nobreakdash-functionally structured groupoids. The correspondence $\G \mapsto \envelope[\infty]\G$ determines, in fact, a natural transformation $\envelope[\infty]{\text-}: (\text-) \mapsto \tannakian[\infty]{\text-}$, that is to say the diagram below commutes for each Lie groupoid homomorphism $\varphi: \G \to \H$
$$
\xymatrix@C=45pt{\G\ar[d]^\varphi\ar[r]^-{\envelope[\infty]\G} & \tannakian[\infty]\G\ar[d]^{\tannakian[\infty]\varphi} \\ \H\ar[r]^-{\envelope[\infty]\H} & \tannakian[\infty]\H\text.\!\!}%
$$%

The main result of the present chapter, to be proved in \refsez{O.3.6}, is: \textsl{for \G\ proper and regular, the standard classical fibre functor \forget[\infty]{\G} is smooth; in fact, \tannakian[\infty]{\G} is a proper regular Lie groupoid although, in general, not one equivalent to \G.} Furthermore, in \refsez{O.3.5} we prove some partial results about the smoothness of the standard classical fibre functor, that are valid for arbitrary proper Lie groupoids; we also remark that the evaluation functor%
$$%
\ev: \R[\infty]\G \longto \xR[\infty]{\tannakian[\infty]\G}\text, \quad R = (E,\varrho) \mapsto (E,\ev_R)%
$$%
is an isomorphism of tensor categories for each proper \G\ (recall the definition of the category \xR[\infty]{\tannakian[\infty]\G} in \refsez{N.21}). Finally, in \refsez{O.4.3} we give a few examples of classically reflexive (proper) Lie groupoids.%

\sezione[The Classical Envelope of a Proper Groupoid]{On the Classical Envelope of a Proper Lie Groupoid}\label{O.3.5}
Let \G\ be a Lie groupoid. Recall from \refsez{N.21} that to each classical representation $R = (E,\varrho)$ of \G\ one can associate a \index{evaluation representation@evaluation representation \ensuremath{\ev_R}}\index{ev R@\ensuremath{\ev_R} (evaluation representation)}representation $\ev_R: \tannakian[\infty]\G \to \GL(E)$, given by evaluation at the object $R \in {\Ob\,\R[\infty]\G}$:
\begin{equazione}\label{xxvi.1}%
{\tannakian[\infty]\G}(x,x') \ni \lambda \mapsto \lambda(R) \in \Lis(E_x,E_{x'})\text,%
\end{equazione}%
which makes the following triangle commute%
\begin{equazione}\label{O.equ8}%
\begin{split}%
\xymatrix@C=40pt@R=23pt{\G\ar[dr]_\varrho\ar[rr]^-{\envelope[\infty]\G} & & \tannakian[\infty]\G\ar[dl]^{\ev_R} \\ & \GL(E)\text,\!\! &}%
\end{split}%
\end{equazione}%
where \envelope[\infty]\G\ denotes the envelope homomorphism ${\envelope[\infty]g}(E,\varrho) = \varrho(g)$.%

Throughout the present section we shall be interested mainly in proper Lie groupoids. Therefore, from now on we assume that \G\ is a proper Lie groupoid and we regard this assumption as made once and for all. As ever, $M$ will denote the base manifold of \G. When we want to state a result that is true under less restrictive assumptions on \G, we shall explicitly point it out. We are going to apply the general theory of representative charts (\refsez{O.2.2+N.23}) to the standard classical fibre functor \forget[\infty]\G.%

\begin{lemma}\label{O.lem7}%
Let $(E,\varrho)$ be a classical representation of a (not necessarily proper) Lie groupoid \G. Suppose we are given an open subset $\Gamma$ of the manifold of arrows of \G, such that the image $\Sigma = \varrho(\Gamma)$ is a submanifold of $\GL(E)$ and such that $\varrho$ restricts to an open mapping of $\Gamma$ onto $\Sigma$.%

Then $\Sigma$ is a \index{tame submanifold}tame submanifold of $\GL(E)$, and the restriction of $\varrho$ to $\Gamma$ is a submersion of $\Gamma$ onto $\Sigma$.

Moreover, when \G\ is proper then the assumption that $\varrho$ should restrict to an open mapping of $\Gamma$ onto $\Sigma$ is superfluous.%
\end{lemma}%
\begin{proof}%
We prove the statement in the proper case first, so without making the assumption that $\varrho$ is an open map of $\Gamma$ onto $\Sigma$.%

We start by observing that for each $x_0 \in M$ the image $\varrho\bigl(\G(x_0,\text-)\bigr)$ is a principal submanifold of $\GL(E)$ and the mapping%
\begin{equazione}\label{O.equ9}%
\G(x_0,\text-) \xto{\;\varrho\;} \varrho\bigl(\G(x_0,\text-)\bigr)%
\end{equazione}%
is a submersion. In particular, the latter will be an open mapping and this forces the open subset%
\begin{equazione}\label{xxvi.3}%
\Sigma(x_0,\text-) = \varrho\bigl({\G(x_0,\text-) \cap \Gamma}\bigr) \subset {\varrho\,\G(x_0,\text-)}%
\end{equazione}%
to be a principal submanifold of $\GL(E)$ as well.

Our argument is as follows. Fix $g_0$ in $\G(x_0,\text-)$ and let $\lambda_0 = \varrho(g_0)$. Choose an open subset $V \subset M$ containing $x_0' = \t(g_0)$, small enough to ensure that the principal bundle $\G(x_0,\text-)$ is trivial over $Z = {{\G x_0}\cap V}$, ie that a local equivariant chart $\G(x_0,Z) \iso {Z\times G_0}$ can be found, where $G_0$ denotes the isotropy group at $x_0$; it is no loss of generality to assume $g_0 \iso (x_0',e)$ in such a chart which we now use, along with the representation $\varrho$, to obtain a smooth section $z \mapsto (z,e) \iso g \mapsto \varrho(g)$ to the target map of $\GL(E)$ over $Z$. Next, the isotropy homomorphism $G_0 \to \GL(E)_0$ determined by $\varrho$ at $x_0$ canonically factors through the quotient Lie group obtained by dividing out the kernel, thus yielding a closed Lie subgroup $H \into \GL(E)_0$. As usual, this Lie subgroup and the target section above can be combined into an embedding of manifolds of type \refequ[O.2.1]{O.equ18}, which fits in the following square%
\begin{equazione}\label{xxvi.4}%
\begin{split}%
\xymatrix@M=4.5pt@C=40pt{{Z\times G_0}\ar[d]^{\id\times\pr}\ar[r]^-\iso & \G(x_0,Z)\ar[d]^\varrho \\ {Z\times H}\ar@{^{(}->}[r]^-{\refequ[O.2.1]{O.equ18}} & \GL(E)}%
\end{split}%
\end{equazione}%
and hence simultaneously displays ${\varrho\,\G(x_0,Z)}$ as a principal submanifold of $\GL(E)$ and, according to the initial remarks of Section \ref{O.2.1}, the mapping $\varrho: G(x_0,Z) \to {\varrho\,\G(x_0,Z)}$ as a submersion; since the subset%
\begin{equazione}\label{xxvi.5}%
{\varrho\,\G(x_0,Z)} = {{\varrho\,\G(x_0,\text-)} \cap \t^{-1}(V)} \subset {\varrho\,\G(x_0,\text-)}%
\end{equazione}%
is an open neighborhood of $\lambda_0$ in ${\varrho\,\G(x_0,\text-)}$, we can conclude.

At this point, in order to prove that $\Sigma$ is a tame submanifold of $\GL(E)$ we need only verify that the restriction $\Sigma \to M$ of the source map of $\GL(E)$ is a submersion. So, fix $\sigma_0 \in \Sigma$, say $\sigma_0 = \varrho(g_0)$ with $g_0 \in \Gamma$. There exists a local smooth source section $\gamma: U \to \Gamma$ through $g_0 = \gamma(\s g_0)$, hence we can also find a local smooth source section $\sigma = {\varrho\circ\gamma}: U \to \Sigma$ through $\sigma_0$.

Finally, we come to the statement that $\varrho: \Gamma \to \Sigma$ is a submersion. Fix $g_0 \in \Gamma$ and let $\sigma_0 = \varrho(g_0)$. Since both $\Gamma$ and $\Sigma$ are tame submanifolds, there exist local trivializations of the respective source maps around the points $g_0 \iso (x_0,0)$ and $\sigma_0 \iso (x_0,0)$, which yield a local expression for $\varrho|_{\Gamma_0}$%
\begin{equazione}\label{xxvi.6}%
\begin{split}%
\xymatrix@R=19pt{\Gamma_0\ar[d]^\iso\ar[r]^\varrho & \Sigma_0\ar[d]^\iso \\ {U\times B}\ar@{-->}[r] & {V\times C}}%
\end{split}%
\end{equazione}%
of the form $(u,b) \mapsto (u,c(u,b))$, where $U \subset V$ are open subsets of $M$ and $B, C$ are Euclidean balls. The partial map $b \mapsto c(x_0,b)$ is submersive at the origin because it is the local expression of \refequ{O.equ9}.

Now we turn to the general case where \G\ is not necessarily proper. Thus, assume that $\varrho$ restricts to an open mapping of $\Gamma$ onto $\Sigma$.%

As explained above, for any given $g_0 \in \G(x_0,\text-)$ there is a submanifold $Z \subset M$ contained in ${\G x_0}$\inciso{although, in general, this is no longer of the form $Z = {{\G x_0}\cap V}$}such that the subset $\G(x_0,Z) \subset \G(x_0,\text-)$ is open, the image ${\varrho\,\G(x_0,Z)}$ is a principal submanifold of $\GL(E)$ and the induced mapping $\varrho: \G(x_0,Z) \to {\varrho\,\G(x_0,Z)}$ is submersive. On the other hand, from the assumption that $\varrho: \Gamma \to \Sigma$ is open it follows that the restriction $\varrho: \Gamma(x_0,\text-) \to \Sigma(x_0,\text-)$ must be open as well, because one has%
\begin{equazione}\label{xxvi.7}%
\varrho\bigl(\mathrm X(x_0,\text-)\bigr) = {\varrho(\mathrm X)}(x_0,\text-)%
\end{equazione}%
for any subset $\mathrm X \subset \mca[1]\G$. Then, since ${\Gamma\cap \G(x_0,Z)} = {\Gamma(x_0,\text-) \cap \G(x_0,Z)}$ is an open subset of $\Gamma(x_0,\text-)$, it is evident that%
\begin{equazione}\label{xxvi.8}%
\Sigma(x_0,Z) = \varrho\bigl({\Gamma\cap \G(x_0,Z)}\bigr) \subset {\varrho\,\G(x_0,Z)}%
\end{equazione}%
is both an open neighbourhood of $\lambda_0$ in $\Sigma(x_0,\text-)$ and an open subset of the principal submanifold ${\varrho\,\G(x_0,Z)}$ of $\GL(E)$. This means that $\Sigma(x_0,\text-)$ is a principal submanifold of $\GL(E)$. Moreover, from what we said it is evident that $\varrho$ induces a submersion of $\Gamma(x_0,\text-)$ onto $\Sigma(x_0,\text-)$.%

The rest of the proof holds without modifications.%
\end{proof}%

Note that the preceding lemma holds for real as well as for complex coefficients\nobreakdash---that is, for $(E,\varrho)$ in \R[\infty]{\G,\nR} or in \R[\infty]{\G,\nC}.%

\separazione%

Our main goal in the present section is to show that the standard classical fibre functor \forget[\infty]\G\ associated with a proper Lie groupoid \G\ always satisfies condition \textsl{ii)} of Proposition~\refcnt[O.2.2+N.23]{xxiii.4}.%

First of all, note that in order that $(\Omega,R)$ may be a \index{representative chart}representative chart of \tannakian[\infty]\G, where $\Omega$ is an open subset of the space of arrows of \tannakian[\infty]\G\ and $R = (E,\varrho) \in {\Ob\,\R[\infty]\G}$, it is sufficient that $\ev_R$ establishes a one-to-one correspondence between $\Omega$ and a submanifold of $\GL(E)$. For if we set $\Gamma = (\envelope[\infty]{})^{-1}(\Omega)$, we have $\varrho(\Gamma) = R(\Omega)$ because of \refequ{O.equ8} and the surjectivity of \envelope[\infty]{}; then Lemma~\refcnt{O.lem7} implies that $R(\Omega)$ is a tame submanifold of $\GL(E)$ and that $\varrho: \Gamma \to R(\Omega)$ is a submersion\nobreakdash---so, in particular, that the map $\ev_R: \Omega \to R(\Omega)$ is open and hence a homeomorphism.%

Our claim about the condition \textsl{ii)} of Proposition~\refcnt[O.2.2+N.23]{xxiii.4} essentially follows from a simple general remark about submersions. Namely, suppose that a commutative triangle of the form%
\begin{equazione}\label{O.esp101}%
\begin{split}%
\xymatrix@C=70pt@R=5pt{& X\ar@{-->}[dd]^g \\ Y\ar[dr]_{f'}\ar[ur]^f & \\ & X'}%
\end{split}%
\end{equazione}%
is given, where $X$, $X'$ and $Y$ are smooth manifolds, $f$ is a submersion onto $X$, $f'$ is a smooth mapping and all we know about $g$ is that it is a set-theoretic solution which fits in the triangle. Then the map $g$ is necessarily smooth; in particular, in case $f'$ is also a surjective submersion, $g$ is a diffeomorphism if and only if it is a set-theoretic bijection.%

To see how this may be used to prove compatibility of charts, suppose we are given an arbitrary representative chart $(\Omega,R)$ of \tannakian[\infty]\G\ to start with, where let us say $R = (E,\varrho)$, and an arbitrary classical representation $S = (F,\sigma)$. Let $\Gamma = (\envelope[\infty]{})^{-1}(\Omega)$, so that $\Gamma$ is an open submanifold of \G. We have already observed that $\varrho$ induces a submersion of $\Gamma$ onto the submanifold $R(\Omega)$ of $\GL(E)$; also, the homomorphism of Lie groupoids%
\begin{equazione}\label{xxvi.9}%
(\varrho,\sigma): \G \longto {\GL(E) \times_M \GL(F)}%
\end{equazione}%
can be restricted to $\Gamma$ to yield a smooth mapping into ${\GL(E) \times_M \GL(F)}$. We get an instance of \refequ{O.esp101} by introducing the following map
\begin{equazione}\label{xxvi.10}%
s = {(\ev_R,\ev_S) \circ {\ev_R}^{-1}}: R(\Omega) \to {\GL(E) \times_M \GL(F)}
\end{equazione}%
(note that $\ev_R: \Omega \to R(\Omega)$ is invertible because we assume $(\Omega,R)$ to be a representative chart), which is then a {\em smooth} section to the projection%
\begin{equazione}\label{xxvi.11}%
{\GL(E) \times_M \GL(F)} \to \GL(E)
\end{equazione}%
and thus, in particular, an immersion. Now, if $s$ is indeed the embedding of a submanifold\inciso{ie if it is an open map onto its image}then we are done, since in that case $(R,S)(\Omega) = s(R(\Omega))$ is a submanifold of ${\GL(E) \times_M \GL(F)}$ and $(\ev_R,\ev_S)$ a bijective map onto it; equivalently, $({R\oplus S})(\Omega)$ is a submanifold of $\GL({E\oplus F})$ and $\ev_{R\oplus S}$ is a bijection of $\Omega$ onto it. (Cf.\ Section \refsez{O.2.2+N.23}. As observed above, this is enough to conclude that $(\Omega,{R\oplus S})$ is a representative chart.) For each open subset $\Lambda$ of $\GL(E)$,%
\begin{equazione}\label{xxvi.12}%
s\bigl({R(\Omega) \cap \Lambda}\bigr) = {s(R(\Omega)) \cap \bigl({\Lambda \times \GL(F)}\bigr)}%
\end{equazione}%
is in fact an open subset of the subspace $s(R(\Omega))$.%

\separazione%

\noindent We can summarize what we have concluded so far as follows:%
\begin{proposizione}\label{O.prp11}%
Let \G\ be a proper Lie groupoid.%

Then the standard classical fibre functor \index{smooth fibre functor}\index{fibre functor!smooth}\forget[\infty]{\G} is smooth if and only if the space of arrows of the classical Tannakian groupoid \tannakian[\infty]{\G} can be covered with open subsets $\Omega$ such that for each of them one can find some $R = (E,\varrho) \in {\Ob\,\R[\infty]\G}$ with the property that $\ev_R$ establishes a bijection between $\Omega$ and a submanifold $R(\Omega)$ of $\GL(E)$.

Moreover, in case the latter condition is satisfied then the envelope homomorphism $\envelope[\infty]\G: \G \longto \tannakian[\infty]\G$ will be a surjective submersion of Lie groupoids.%
\end{proposizione}%
\begin{proof}%
The first assertion is already proven.%

The second assertion follows from the (previously noticed) fact that for each representative chart $(\Omega,R)$ the mapping $\varrho: \Gamma \to R(\Omega)$ is a submersion, where as usual $R = (E,\varrho)$ and we put $\Gamma = (\envelope[\infty]{})^{-1}(\Omega)$. (Remember from the proof of Prop.~\refcnt[O.2.2+N.23]{xxiii.4} that $\ev_R$ establishes a diffeomorphism between $\Omega$ and the submanifold $R(\Omega)$ of $\GL(E)$.)%
\end{proof}%

Note that, for any proper Lie groupoid \G\ whose associated standard classical fibre functor \forget[\infty]{\G} is smooth, the preceding proposition allows us to characterize the familiar Lie groupoid structure on the Tannakian groupoid \tannakian[\infty]{\G} as the unique such structure for which the envelope homomorphism \envelope[\infty]{\G} becomes a submersion. Indeed, assume that an unknown Lie groupoid structure, making \envelope[\infty]{\G} a submersion, is assigned on the Tannakian groupoid of \G. Let \tannakian[*]{\G} indicate the Tannakian groupoid of \G\ endowed with the unknown smooth structure. Now, the identity homomorphism of the Tannakian groupoid into itself fits in the following triangle%
\begin{equazione}\label{xxvi.17}%
\begin{split}%
\xymatrix@C=90pt@R=7pt{ & \tannakian[\infty]\G\ar@{-->}[dd]^\id \\ \G\ar[dr]_{\envelope[*]{}}\ar[ur]^{\envelope[\infty]{}} & \\ & \tannakian[*]\G}%
\end{split}%
\end{equazione}%
where $\envelope[\infty]{} = \envelope[\infty]\G = \envelope[*]{}$ are surjective submersions. It follows that the identity $\id: \tannakian[\infty]\G = \tannakian[*]\G$ is a diffeomorphism.%

\separazione%

Under the assumption of properness, we can also say something useful about condition~\textsl{i)} of Proposition~\refcnt[O.2.2+N.23]{xxiii.4}:%
\begin{nota}\label{xxvi.21}%
Let \G\ be a proper Lie groupoid. Suppose that for each identity arrow $x_0$ of the Tannakian groupoid \tannakian[\infty]{\G} one can find a representative chart for \tannakian[\infty]{\G} about $x_0$. Then we contend that the condition~\textsl{i)} of Proposition~\refcnt[O.2.2+N.23]{xxiii.4} is satisfied by the classical fibre functor \forget[\infty]\G.%

Let an arbitrary arrow $\lambda_0: x_0 \to x_0'$ of \tannakian[\infty]{\G} be given. Because of properness, we have $\lambda_0 = \envelope[\infty]{g_0}$ for some arrow $g_0: x_0 \to x_0'$ of \G. Select any smooth local bisection $\sigma: U \to \mca[1]\G$, defined over a neighbourhood $U$ of $x_0$ and with $\sigma(x_0) = g_0$, and let $U' = \t(\sigma(U))$. Now, let $(\Omega,R)$ be a representative chart about $x_0$, let us say with $\Omega \subset \tannakian[\infty]\G|_U$ and $R = (E,\varrho)$. Notice that one has the following commutative square%
\begin{equazione}\label{xxvi.22}%
\begin{split}%
\xymatrix{\G|_U\ar[d]_\iso^{\sigma\text-}\ar[r]^-\varrho & \GL(E)|_U\ar[d]_\iso^{({\varrho\circ\sigma})\text-} \\ \G(U,U')\ar[r]^-\varrho & \GL(E)(U,U')\text,\!\!}%
\end{split}%
\end{equazione}%
where ${\sigma\text-}$ denotes the left translation diffeomorphism $g \mapsto {\sigma(\t(g)) \cdot g}$ and, similarly, ${({\varrho\circ\sigma})\text-}$ denotes the diffeomorphism $\mu \mapsto {\varrho(\sigma(\t[\mu])) \cdot \mu}$. Let $\Gamma = (\envelope[\infty]{})^{-1}(\Omega)$, so $\Gamma \subset \G|_U$ is an open subset. Then $\Gamma^\sigma = {\sigma\text-}(\Gamma)$ is an open neighbourhood of $g_0$, $\Omega^\sigma = {({\envelope[\infty]{}\circ\sigma})\text-}(\Omega)$ is an open neighbourhood of $\lambda_0$ and $\Gamma^\sigma = (\envelope[\infty]{})^{-1}(\Omega^\sigma)$. It follows that the subset%
\begin{equazione}\label{xxvi.23}%
R(\Omega^\sigma) = \varrho(\Gamma^\sigma) = {({\varrho\circ\sigma})\text-}(\varrho(\Gamma)) = {({\varrho\circ\sigma})\text-}(R(\Omega))%
\end{equazione}%
is a submanifold of $\GL(E)(U,U')$. Similarly, one sees that $\Omega^\sigma$ bijects onto $R(\Omega^\sigma)$ via $\ev_R$. So $(\Omega^\sigma,R)$ is a representative chart about $\lambda_0$.%
\end{nota}%

\separazione%

The next, conclusive result provides, in the special case under exam, a positive answer to the question raised in \refsez{N.21} about the evaluation functor being an equivalence of categories.%
\begin{proposizione}\label{O.prp12}%
Let \G\ be any proper Lie groupoid.%

Then the \index{evaluation functor@evaluation functor \ensuremath{\ev}}\index{ev@\ensuremath{\ev} (evaluation functor)}evaluation functor
$$%
\ev: \R[\infty]\G \longto \R[\infty]{\tannakian[\infty]\G}\text, \quad R = (E,\varrho) \mapsto (E,\ev_R)%
$$%
is an isomorphism of categories, having the pullback along the envelope homomorphism of \G\ as inverse.%
\end{proposizione}%
\begin{proof}%
This can be verified directly, since the envelope homomorphism of a proper Lie groupoid is already known to be surjective.%
\end{proof}%

\sezione[Proper Regular Groupoids]{Smoothness of the Classical Envelope of a Proper Regular Groupoid}\label{O.3.6}
We start by recalling a few basic definitions and properties. For additional information, see {\em Moerdijk (2003)}~\cite{Moerdijk'03}.%

Recall that a Lie groupoid \G\ over a manifold $M$ is said to be \index{regular groupoid}\index{groupoid!regular}{\em regular} when the rank of the differentiable map $\t_x: \G(x,\text-) \to M$ locally keeps constant as the variable $x$ ranges over $M$; an equivalent condition is that the anchor map of the Lie \index{Lie algebroid}\index{algebroid}algebroid of \G, let us call it $\rho: \mathfrak g \to \T M$, should have locally constant rank (as a morphism of vector bundles over $M$). If \G\ is regular then the image of the anchor map $\rho$ is a subbundle $F$ of the tangent bundle \T M; in fact, $F$ turns out to be an integrable subbundle of \T M and hence determines a foliation $\mathcal F$ of the base manifold $M$, called the \index{orbit foliation}{\em orbit foliation} associated with the regular groupoid \G.

Recall that a \index{leaf}{\em leaf} of a foliation $\mathcal F$ associated with an integrable subbundle $F$ of \T M is a maximal connected immersed submanifold $L$ of $M$ with the property of being everywhere tangent to $F$. The codimension of $L$ in $M$ coincides with the codimension of $F$ in \T M. Also recall that a \index{transversal (complete)}{\em transversal for $\mathcal F$} is a submanifold $T$ of $M$, everywhere transversal to $F$ and of dimension equal to the codimension of $F$. There always exist \index{complete transversal}{\em complete} transversals, i.e.\ transversals that meet every leaf of the foliation.

{\em Bundles of Lie groups,} that is to say Lie groupoids whose source and target map coincide, form a very special class of regular Lie groupoids. Proper bundles of Lie groups are also called {\em bundles of compact Lie groups.}%

\begin{lemma}\label{O.lem6}%
Let \G\ be a bundle of compact Lie groups over a manifold $M$. Let $R = (E,\varrho)$ be a classical representation of \G.%

Then the image $\varrho(\G)$ is a submanifold of $\GL(E)$.%
\end{lemma}%
\begin{proof}%
By a result of Weinstein~\cite{Weinstein'02}, every bundle of compact Lie groups is {\em locally trivial.} This means that for each $x \in M$ one can find an open neighborhood $U$ of $x$ in $M$ and a compact Lie group $G$ such that there exists an isomorphism of Lie groupoids over $U$ (viz.\ a local trivialization)%
\begin{equazione}\label{xxvii.2}%
\G|_U \iso {U\times G}\text.%
\end{equazione}%
At the expense of replacing $U$ with a smaller open neighborhood, one can also assume that there is a local trivialization $E|_U \iso {U\times \boldsymbol V}$, where $\boldsymbol V$ is some vector space of finite dimension; as explained in Note~\refcnt[O.2.2+N.23]{O.sub5}, such a trivialization will determine an isomorphism $\GL(E|_U) \iso {U\times \GL(\boldsymbol V)}$ of Lie groupoids over $U$. Then one can take the following composite homomorphism%
\begin{equazione}\label{xxvii.3}%
\begin{split}%
\xymatrix{{U\times G}\ar[r]^-\iso\ar[d] & \G|_U\ar[r]^-{\varrho|_U}\ar[d] & \GL(E|_U)\ar[r]^-\iso\ar[d] & {U\times \GL(\boldsymbol{V})}\ar[r]^-\pr\ar[d] & \GL(\boldsymbol{V})\ar[d] \\ {U\times U}\ar[r]^\id & {U\times U}\ar[r]^\id & {U\times U}\ar[r]^\id & {U\times U}\ar[r] & {\pt\times \pt}\text.\!\!}%
\end{split}%
\end{equazione}%
This yields a smooth family of representations of the compact Lie group $G$ on the vector space $\boldsymbol V$, parametrized by the connected open set $U$. We will denote such family by $\varrho_U: {U\times G} \to \GL(\boldsymbol V)$.%

Now, it follows from the so-called `homotopy property of representations of compact Lie groups' (Note~\refcnt[O.3.1+.4.1]{ii1}) that all the representations of the smooth family $\varrho_U$ are equivalent to each other; in particular, they all have the same kernel $K \subset G$. Hence there exists a unique map $\widetilde{\varrho_U}$ making the diagram%
\begin{equazione}\label{xxvii.4}%
\begin{split}%
\xymatrix{{U\times G}\ar[d]_{\id\times\pr}\ar[r]^-{\id\times\varrho_U} & {U\times \GL(\boldsymbol{V})} \\ {U\times ({G/K})}\ar@{-->}[ur]_{\widetilde{\varrho_U}} & }%
\end{split}%
\end{equazione}%
commute. Note that the map $\widetilde{\varrho_U}$ must be smooth, because ${\id\times\pr}$ is a surjective submersion; of course, the same map is also a faithful representation of the bundle of compact Lie groups ${U\times ({G/K})}$ on the trivial vector bundle ${U\times \boldsymbol V}$. Then Corollary~\refcnt[O.2.2+N.23]{O.cor4} implies that the image of $\widetilde{\varrho_U}$ is a submanifold of ${U\times \GL(\boldsymbol{V})}$. The latter submanifold coincides, via the diffeomorphism $\GL(E)|_{\Delta U} \iso {U\times \GL(\boldsymbol V)}$, with the intersection ${\varrho(\G)\cap \GL(E)|_U}$.%
\end{proof}%

It is evident from the above proof that the kernel of the envelope homomorphism $\envelope[\infty]{}: \G \to \tannakian[\infty]\G$ must be a (locally trivial) bundle of compact Lie groups \K, embedded into \G. Thus, if $U$ is a connected open subset of $M$ and $R = (E,\varrho)$ is a classical representation such that $\kernel{\varrho_u} = \K|_u$ at some point $u \in U$, it follows from the aforesaid homotopy property that $\kernel{\varrho|_U} = \K|_U$ and therefore\inciso{from the commutativity of \refequ[O.3.5]{xxvi.1}}that the evaluation representation $\ev_R$ is faithful on $\tannakian[\infty]\G|_U$.%

From the latter remark, the discussion about smoothness in the preceding section and Lemma~\refcnt{O.lem6} it follows immediately that \textsl{the standard classical fibre functor \forget[\infty]\G\ associated with a bundle of compact Lie groups \G\ is smooth.} Indeed, let an arbitrary arrow $\lambda_0 \in \tannakian[\infty]\G$ be fixed, let us say $\lambda_0 \in \tannakian[\infty]\G|_{x_0}$ with $x_0 \in M$. Take an object $R \in {\Ob\,\R[\infty]\G}$ with the property that the restriction of the evaluation representation $\ev_R$ to $\tannakian[\infty]\G|_{x_0}$ is faithful (this exists by Prop.~\refcnt[O.1.4+.3.3]{O.lem8}) and then choose any connected open neighbourhood $U$ of $x_0$ in $M$. Then the pair $\bigl(\tannakian[\infty]\G|_U,R\bigr)$ constitutes a representative chart for \forget[\infty]\G\ about $\lambda_0$.%

More generally, let \G\ be a proper Lie groupoid with the property that for each $x_0 \in M$ there exists an open neighbourhood $U$ of $x_0$ in $M$ such that $\G|_U$ is a bundle of compact Lie groups. By adapting the above recipe for the construction of representative charts about the arrows belonging to the isotropy of \tannakian[\infty]\G\ and by taking into account Note~\refcnt[O.3.5]{xxvi.21}, we see that \forget[\infty]\G\ is smooth also in the present case.%

\separazione%

We are going to generalize the latter remark to arbitrary proper regular Lie groupoids. The shortest way to do this is to apply the theory of weak equivalences of \refsez{O.2.4}.%

\begin{proposizione}\label{O.prp9}
Let \G\ be a proper regular Lie groupoid.

Then the standard classical fibre functor \forget[\infty]\G\ associated with \G\ is \index{smooth fibre functor}\index{fibre functor!smooth}smooth.

Recall that in view of Proposition~\refcnt[O.3.5]{O.prp11} this can also be expressed by saying that there exists a (necessarily unique) Lie groupoid structure on the Tannakian groupoid \tannakian[\infty]{\G} such that the \index{envelope homomorphism@envelope homomorphism \envelope[T]\G, \envelope[\infty]\G}envelope homomorphism \envelope[\infty]{\G} becomes a smooth submersion.
\end{proposizione}
\begin{proof}
Let $M$ be the base of \G. Select a complete transversal $T$ for the foliation of the manifold $M$ determined by the orbits of \G. Note that $T$ is in particular a {\em slice,} so the restriction $\G|_T$ is a proper Lie groupoid embedded into \G\ (by Note~\refcnt[N.4]{iv.4}). If $i: T \into M$ denotes the inclusion map then, by the remarks at the end of \refsez{N.4}, the embedding of Lie groupoids%
\begin{equazione}\label{xxvii.6}%
\begin{split}%
\xymatrix@C=29pt@M=5pt@R=17pt{\G|_T\ar[d]\ar@{^{(}->}[r]^-{\text{inclusion}} & \G\ar[d] \\ {T\times T}\ar@{^{(}->}[r]^-{i\times i} & {M\times M}}%
\end{split}%
\end{equazione}%
is a Morita equivalence. One may therefore find another (proper) Lie groupoid \K, along with Morita equivalences $\G|_T \xfrom{\:\text{M.e.}\:} \K \xto{\:\text{M.e.}\:} \G$ inducing surjective submersions at the level of base manifolds. The corresponding morphisms of standard classical fibre functors%
\begin{equazione}\label{xxvii.7}%
\bigl(\R[\infty]{\G|_T},\forget[\infty]{\G|_T}\bigr) \xfrom{\:\text{w.e.}\:} \bigl(\R[\infty]\K,\forget[\infty]\K\bigr) \xto{\:\text{w.e.}\:} \bigl(\R[\infty]\G,\forget[\infty]\G\bigr)%
\end{equazione}
are \index{weak equivalence}weak equivalences. Hence, by Proposition~\refcnt[O.2.4]{O.prp8}, one is reduced to showing that \forget[\infty]{\G|_T} is a smooth fibre functor.

Now, $\G|_T$ is a proper Lie groupoid over $T$ with the above-mentioned property of being, locally, just a bundle of compact Lie groups.%
\end{proof}%

Let $\mathfrak{ProReg}$ denote the category of proper regular Lie groupoids. One may summarize the conclusions of the present section as follows:

\sloppy
\begin{theorem}\label{O.thm4}
The classical Tannakian correspondence $\G \mapsto \tannakian[\infty]\G$ induces an idempotent functor
\begin{equazione}\label{xxvii.9}
\tannakian[\infty]{\text-}: \mathfrak{ProReg} \longto \mathfrak{ProReg}\text;
\end{equazione}
moreover, envelope homomorphisms fit together into a natural transformation
\begin{equazione}\label{xxvii.10}
\envelope[\infty]{\text-}: \Id \longto \tannakian[\infty]{\text-}\text.
\end{equazione}
\end{theorem}
\textit{Open Question.} It is natural to ask whether this result can be generalized to the whole category of proper Lie groupoids.

\fussy

\sezione[Classical Reflexivity: Examples]{A few Examples of Classically Reflexive Lie Groupoids}\label{O.4.3}
Recall that a Lie groupoid $\G \rightrightarrows X$ is said to be {\em \'etale} if the source and target maps $\s, \t: \G \to X$ are \'etale maps, that is to say local isomorphisms of smooth manifolds. An open subset $\Gamma \subset \G$ will be said to be {\em flat} if the source and target map restrict to open embeddings of $\Gamma$ into $X$. A Lie groupoid \G\ will be said to be {\em source-proper} or, for short, {\em\s\nobreakdash-proper} when the source map of \G\ is a proper map.%

\begin{proposizione}\label{xxviii.1}
Let \G\ be a source-proper \'etale Lie groupoid.

Then \G\ admits globally faithful classical \index{reflexive}\index{groupoid!reflexive}representations.
\end{proposizione}
\begin{proof}
The {\em regular representation $(R,\varrho)$} of \G\ exists and has locally finite rank. A couple of remarks before starting. Let $X$ be the base of \G.%

For every point $x$ of $X$, the \s-fiber $\s^{-1}(x)$ is a finite set. Indeed, it is discrete, because if $g \in \s^{-1}(x)$ then since \s\ is \'etale there exists a flat open neighborhood $\Gamma \subset \G$ and therefore $\{g\} = {\Gamma \cap \s^{-1}(x)}$ is a neighborhood of $g$ in the \s-fiber. It is also compact, because of \s-properness.%

Put $\ell(x) = \norma{\s^{-1}(x)}$, the cardinality of this finite set. Then the fiber $R_x$ of the vector bundle $R \to X$ is by definition the vector space%
\begin{equazione}\label{xxviii.2}%
\C^0(\s^{-1}(x);\nR) \can \nR^{\ell(x)}%
\end{equazione}%
of \nR-valued maps. This makes sense because%
\begin{lemma}%
The assignment $x \mapsto \ell(x)$ defines a locally constant function on $X$, with values into positive integers.%
\end{lemma}%
\noindent\textsl{Proof of the lemma.} Fix $x \in X$, and say $\s^{-1}(x) = \{g_1,\ldots,g_\ell\}$. For every $i = 1, \ldots, \ell$, there exists a flat open neighborhood $\Gamma_i \subset \G$ of $g_i$. Choosing an open ball $B \subset {\bigcap \s(\Gamma_i)}$ at $x$, we can assume $\s: \Gamma_i \isoto B$ to be an isomorphism ${\forall i}$. Moreover, it is no loss of generality to assume the open subsets $\Gamma_1, \ldots, \Gamma_\ell$ to be pairwise disjoint. (As a consequence of the fact that a finite union of open balls in any manifold---not necessarily Hausdorff---is a Hausdorff open submanifold.) Then, ${\forall i} = 1, \ldots, \ell$ and ${\forall z} \in B$, the intersection ${\s^{-1}(z) \cap \Gamma_i}$ consists of a single point $g_i(z)$, and these points $g_1(z), \ldots, g_\ell(z) \in \G$ are pairwise distinct, because the $\Gamma_i$ are pairwise disjoint. This shows $\ell(z) \geqq \ell(x)$ ${\forall z} \in B$. To prove the converse inequality, it will suffice to prove that ${\exists N} \subset B$, a smaller ball at $x$, such that $\s^{-1}(N) \subset \Gamma = {\Gamma_1 \cup \cdots \cup \Gamma_\ell}$. Consider a decreasing sequence of closed balls $C_{n+1} \subset C_n \subset B$ shrinking to $x$, and the corresponding decreasing sequence $\Sigma_n = {\s^{-1}(C_n) - \Gamma}$ of closed subsets of the compact subspace $\s^{-1}(C_1) \subset \G$; there $\exists n$ such that $\Sigma_n = \varnothing$, in other words $\s^{-1}(C_n) \subset \Gamma$. This concludes the proof of the lemma.%
\vskip 5pt%
Thus, it makes sense to regard $R \to X$ as the set-theoretic support of a \nR-linear vector bundle of locally finite rank. The proof of the lemma contains also a recipe for the construction of local trivializations. Namely, let $x \in X$ be fixed, and choose an ordering $\s^{-1}(x) = \{g_1,\ldots,g_\ell\}$ of the corresponding fiber; there exist an open ball $B \subset X$ centered at $x$ and disjoint flat open neighborhoods $\Gamma_1, \ldots, \Gamma_\ell \subset \G$ of $g_1, \ldots, g_\ell$ such that $\s^{-1}(B) = {\Gamma_1 \cup \cdots \cup \Gamma_\ell}$. Then one gets a bijection $R|_B \iso {B \times \nR^\ell}$ by setting, for $z \in B$ and $f \in \C^0(\s^{-1}(z);\nR)$,%
$$%
(z,f) \mapsto \bigl(z,f(g_1(z)),\ldots,f(g_\ell(z))\bigr)\text.%
$$%
(Cf.\ the notation used in the proof of the lemma.) The transition mappings are smooth, because locally they are given by constant permutations $(a_1,\ldots,a_\ell) \mapsto \bigl(a_{\tau(1)},\ldots,a_{\tau(\ell)}\bigr)$.%

The \nR-linear isomorphism%
$$%
\varrho(g) \in \Lis(R_x,R_y)\text,%
$$%
associated with $g \in \G(x,y)$, is defined by `translation'%
$$%
f \mapsto {\varrho(g)}(f) \equiv f({\text-\,g})\text.%
$$%
The resulting functorial map $\varrho: \G \longto \GL(R)$ is clearly faithful; it is also smooth, because in any trivializing local charts it looks like a locally constant permutation.%
\end{proof}%

\separazione%

\noindent If \G\ is any \'etale Lie groupoid with base manifold $X$, there is a morphism of Lie groupoids $\mathrm{Ef}: \G \longto {\Gamma X}$, where ${\Gamma X}$ is the \'etale Lie groupoid (with base $X$) of germs of smooth isomorphisms $U \isoto V$ between open subsets of $X$. It sends $g \in \G$ to the germ of the local smooth isomorphism associated with a flat open neighborhood of $g$. An {\em effective} Lie groupoid is an \'etale Lie groupoid such that $\mathrm{Ef}$ is faithful, in other words such that every $g \in \G$ is uniquely determined by its `local action' on the base manifold $X$. (Some of the simplest \'etale groupoids, such as for instance the trivial ones ${X\times K}$, $K$ a discrete group, are not effective at all!)%

The class of effective Lie groupoids is stable under weak equivalence among \'etale Lie groupoids. (Cf.\ \textit{Moerdijk and Mr\v cun (2003),} \cite{Moerdijk&Mrcun'03} p.~137.)%

The following conditions on a Lie groupoid \G\ are equivalent:%
\begin{elenco}%
\item[1.]\G\ is weakly equivalent to a proper effective groupoid;%
\item[2.]\G\ is weakly equivalent to the Lie groupoid associated with an orbifold.%
\end{elenco}%
(Cf.\ \textit{ibid.}\ p.~143.) The relevance of this theorem in the present context is that it tells that if one wants to study orbifolds through their associated Lie groupoid and Tannakian duality, it is sufficient to prove the duality result for proper effective groupoids.%

Any \'etale Lie groupoid $\G \rightrightarrows X$ has a canonical representation on the tangent bundle $\T{X} \to X$, which associates to $g \in \G(x,y)$ the invertible \nR-linear map $\T[x]{X} \to \T[y]{X}$ of tangent spaces given by the tangent map at $x$ of the germ of local smooth isomorphisms $\mathrm{Ef}(g)$. In general, this representation need not be faithful. However%
\begin{proposizione}%
If \G\ is a proper effective Lie groupoid with base $X$, the canonical representation on the tangent bundle \T{X} is faithful.\footnote{This was pointed out to me by I.~Moerdijk.}%
\end{proposizione}%
\begin{proof}%
If $\G \rightrightarrows X$ is a proper \'etale Lie groupoid and $x \in X$, there exist a neighborhood $U \subset X$ of $x$ and a smooth action of the isotropy group $G_x = \G|_x$ on $U$, such that the Lie groupoid $\G|_U \rightrightarrows U$ is isomorphic to the action groupoid ${G_x\ltimes U}$. I need to recall part of the proof. (Cf.\ \textit{Moerdijk and Mr\v cun (2003),} \cite{Moerdijk&Mrcun'03} p.~142.) Let $G_x = \{1,\ldots,\ell\}$. There are a connected open neighborhood $W \subset X$ of $x$ and \s-sections $\sigma_1, \ldots, \sigma_\ell: W \to \G$ with $\sigma_i(x) = i \in G_x$ ${\forall i}$, such that the maps $f_i = {\t\circ \sigma_i}$ send $W$ diffeomorphically onto itself and satisfy ${f_i\circ f_j} = f_{ij}$ for all $i, j \in G_x$.%

Since \G\ is also effective, the group homomorphism $i \mapsto f_i$, of $G_x$ into the group $\Aut(W;x)$ of smooth automorphisms of $W$ that fix the point $x$, is injective. Now, if $M$ is a connected manifold and $H \subset \Aut(M)$ is a finite group of smooth automorphisms of $M$, the group homomorphism which maps $f \in H_x = \{f \in H| f(x) = x\}$ to the tangent map $\T[x]{f} \in \Aut(\T[x]M)$ is injective ${\forall x} \in M$. (\textit{Ibid.}\ p.~36.) In the case $M = W$ and $H = \{f_i| i \in G_x\} = H_x$, this says precisely that the canonical representation of \G\ on the tangent bundle \T{X} restricts to a faithful representation $G_x \into \Aut(\T[x]X)$.%
\end{proof}%

\separazione%

\noindent Another simple example is offered by \index{action groupoid@action groupoid \ensuremath{G\ltimes M}}\index{G ltimes M@\ensuremath{G\ltimes M} (action groupoid)}action groupoids associated with compact Lie group actions.

Precisely, let $K$ be a compact Lie group acting smoothly on a manifold $X$, say from the left. We denote by ${K\ltimes X}$ the Lie groupoid over $X$ whose manifold of arrows is the Cartesian product ${K\times X}$, with the second projection $(k,x) \mapsto x$ as source map, the action ${K\times X} \to X$ as range map and%
$$%
{(k',{k\cdot x})\cdot(k,x)} = (k'k,x)%
$$%
as composition of arrows.%

If $\boldsymbol{V}$ is a faithful $K$-module (in other words a faithful representation $\varrho$ of the compact Lie group $K$ on a vector space $\boldsymbol{V}$), then we get a faithful representation of the groupoid ${K\ltimes X}$ on the trivial vector bundle ${X\times \boldsymbol{V}}$, defined by%
$$%
(k,x) \mapsto \bigl(x,{k\cdot x},\varrho(k)\bigr)\text.%
$$%

\pagestyle{headings}\sloppy\hbadness=1000\vbadness=1000

\refstepcounter{chapter}\addcontentsline{toc}{chapter}{Bibliography}
\bibliography{gtrentinaglia-2008}

\pagestyle{headings}\sloppy\hbadness=10000\vbadness=1000

\newpage
\refstepcounter{chapter}\addcontentsline{toc}{chapter}{Index}
\printindex

\end{document}